\date{} %
\title{Contour lines of the \\ two-dimensional discrete Gaussian free field}
\author{Oded Schramm \and Scott Sheffield}
\documentclass[12pt,naturalnames]{article}
\usepackage{amsmath}
\usepackage{amssymb}
\usepackage{amsthm}
\usepackage{amsfonts}
\usepackage{graphicx}
\input labelfig.tex
\input epsf.tex

\newif\iffigures\figurestrue
\figuresfalse

\newif\ifhyper\IfFileExists{hyperref.sty}{\hypertrue}{\hyperfalse}

\ifhyper\usepackage{hyperref}
\def\hitem#1#2{\item[\hypertarget{#1}{#2}]\expandafter\gdef\csname LBL#1ITM\endcsname{#2}}
\def\iref#1{\hyperlink{#1}{\csname LBL#1ITM\endcsname}}
\else
\def\hitem#1#2{\item[{#2}]\expandafter\gdef\csname LBL#1ITM\endcsname{#2}}
\def\iref#1{{\csname LBL#1ITM\endcsname}}
\fi

\newif\ifdraft
\drafttrue
\long\def\note#1/{\ifdraft {\bf [#1]}\fi}
\long\def\comment#1{}
\long\def\old#1{}

\numberwithin{equation}{section}
\numberwithin{figure}{section}

\newtheorem{theorem}{Theorem}
\numberwithin{theorem}{section}
\newtheorem{corollary}[theorem]{Corollary}
\newtheorem{lemma}[theorem]{Lemma}
\newtheorem{proposition}[theorem]{Proposition}

\theoremstyle{remark}\newtheorem{definition}[theorem]{Definition}
\theoremstyle{remark}\newtheorem{remark}[theorem]{Remark}
\def\eref#1{(\ref{#1})}
\let\qqed=\qed
\def\QED{\qqed\medskip}
\let\qed=\QED

\newcommand{\R}{\mathbb{R}}

\newcommand{\Z}{\mathbb{Z}}
\newcommand{\N}{\mathbb{N}}
\def\H{\mathbb{H}}
\def\U{\mathbb{U}}
\def\diam{\mathop{\mathrm{diam}}}
\def\rad{\mathop{\mathrm{rad}}}

\def\dist{\mathop{\mathrm{dist}}}

\def\Im{{\rm Im}\,}
\def\Re{{\rm Re}\,}
\def\SLEkk#1/{$\mathrm{SLE}(#1)$}
\def\SLEr#1/{$\mathrm{SLE(\kappa;#1)}$}
\def\SLEkr#1;#2/{$\mathrm{SLE(#1;#2)}$}
\def\SLEk/{\SLEkk{\kappa}/}
\def\SLEtwo/{\SLEkk2/}
\def\SLE/{$\mathrm{SLE}$}
\def\SLEab/{\SLEkr 4; {a/\hco-1}, {b/\hco-1}/}

\def\Ito/{It\^o}
\def \eps {\epsilon}
\def \P {{\bf P}}
\def\md{\mid}
\def\Bb#1#2{{\def\md{\bigm| }#1\bigl[#2\bigr]}}
\def\BB#1#2{{\def\md{\Bigm| }#1\Bigl[#2\Bigr]}}
\def\Bs#1#2{{\def\md{\mid}#1[#2]}}
\def\Pb{\Bb\P}
\def\Eb{\Bb\E}
\def\PB{\BB\P}
\def\EB{\BB\E}
\def\Ps{\Bs\P}
\def\Es{\Bs\E}

\def \p {{\partial}}
\def \E {{\bf E}}

\def\closure{\overline}
\def\ev#1{{\mathcal{#1}}}
\def \proof {{ \medbreak \noindent {\bf Proof.} }}
\def\proofof#1{{ \medbreak \noindent {\bf Proof of #1.} }}

\def\setminus{\smallsetminus}
\def\backslash{\smallsetminus}
\def\bl{\bigl}\def\br{\bigr}\def\Bl{\Bigl}\def\Br{\Bigr}

\def\VD{V_D}

\def\hatVD{\hat V_D}

\def\capacity{\mathrm{cap}_\infty}

\def\Qp{U_Q}
\def\oconf#1{\Phi_{{#1}}}
\def\oconftR{\oconf {3R}}
\def\oconftRp{\oconf {R'}}
\def\oconfr{\oconf r}
\def\conf#1{\Theta_{{#1}}}
\def\confR{\conf R}

\def\ball{\mathfrak{B}}
\def\bal#1{\ball_{#1}}

\def\Qual{Q}

\def\inr#1{\mathrm{rad}_{#1}}

\def\WW{W^*}
\def\WW_#1{{\tilde W}_{s_#1}}

\def\errprb{\eta}
\def\tlambda{\lambda_{\TG}}
\def\TG{TG}
\def\GG{\mathfrak G}
\def\DDp{{\mathfrak {D}}^+}
\def\DDm{{\mathfrak {D}}^-}

\def\dhaus{{d}_{\mathrm H}}
\def\dstrong{{d}_{\mathcal U}}
\def\dCKC{{d}_{\mathrm {CKC}}}

\def\hco{\lambda}
\def\lowerhco{\Lambda_0}
\def\upperhco{\bar \Lambda}

\def\bbeta{\xi}
\def\rr{{r_D}}

\def\DD{D}

\def\gSLE{\gamma_{\mathrm{SLE}}}
\def\cH{\closure{\H}}

\def\achi{\chi}
\def\tchi{{\tilde\achi}}

\def\bS{{\breve S}}
\def\Va{W^1}
\def\Vb{W_7}
\def\Vp{W_1^7}
\def\Vpp{\Va\cup\Vb}
\def\Vs{W_2^6}%
\def\Vss{W_3^5}
\def\Vtwo{W^4}
\def\Vtwop{W_4}

\def\bah{{\bar h}}
\def\rrho{{r_0}}
\def\vv{v_0}
\def\exter#1{\mathrm{ext}_{#1}}
\def\inter#1{\mathrm{int}_{#1}}
\def\sig{\sigma}
\def\zz{{S_\tau}}
\def\zzp{S'_{\tau'}}
\def\zzT{z_T}
\def\evZ{\evZs{\sig}}
\def\evZp{{\evZ}'}
\def\evZs#1{\ev Z_0^{#1}}
\def\evZa{\ev Z_0}

\def\th{{\tilde h}}
\def\hth{{h}}
\def\expo#1{{\zeta_{#1}}}
\def\beg{\hat\gamma_g}
\def\OC{O_{\upperhco}}
\def\patha{\alpha_a}
\def\pathb{\alpha_b}
\def\pathq{\alpha_q}

\def\density{{\mathfrak p}}

\def\SmirnovPerc{MR1851632}
\def\CamiaNewmanSLE{MR2249794}

\def\LSWi{MR2002m:60159a}
\def\LSWiii{MR1899232}

\def\LSWlesl{MR2044671}
\def\LSWrestriction{MR1992830}

\def\SchSLE{MR1776084}

\def\Ahlfors{MR50:10211}

\def\RSsle{MR2153402}
\def\Lbook{MR92f:60122}

\def\PommeBDRY{MR95b:30008}

\def\RevuzYor{MR2000h:60050}

\def\KaratsasShreve{MR89c:60096}

\def\WernerStFlour{MR2079672}
\def\LawlerSLbook{MR2129588}
\def\CardySLESurvey{MR2148644}
\def\LawlerSLEintro{MR2087784}
\def\LawlerStrict{MR1645225}
\def\SchrammSheffieldHE{MR2184093}
\def\LawlerConformalBook{MR2129588}
\def\DoyleSnell{MR920811}

\def\ConiglioPotts{1989PhRvL..62.3054C}
\def\DupSalWinding{1988PhRvL..60.2343D}
\def\DupSalDenseSAW{1987NuPhB.290..291D}
\def\DupSalPercolationhull{1987PhRvL..58.2325S}
\def\Kadanoff{1978JPhA...11.1399K}

\def\KondevHenleyninetyfive{1995PhRvL..74.4580K}
\def\KondevHenleytwothousand{2000PhRvE..61..104K}
\def\HuberKondevGeometryofloops{2001APS..DCM.Q2008H}
\def\HuberDuKondevGeometryofloops{2002APS..MAR.U4003K}
\def\Nienhuiseightytwo{1982PhRvL..49.1062N}

\def\KostThouThree{1973JPhC....6.1181K}
\def\KostRoughening{1977JPhC...10.3753K}
\def\Foltin{2001JPhA...34.5327F}

\def\FroSpenAbelianSpin{MR634447}
\def\NaddafSpencer{MR1461951}

\def\KenyonDominoGFF{MR1872739}

\def\Mandelbrot{1967Sci...156..636M}

\def\IsichenkoSurvey{1992RvMP...64..961I}
\def\SpencerSurvey{MR1460292}
\def\BriMelFroSurvey{MR833220}
\def\GiacominSurvey{MR1919512}
\def\DuplantierSurvey{1989PhR...184..229D}
\def\NienhuisKagerSurvey{MR2065722}
\def\NienhuisSurvey{MR751711}

\def\GFFSurvey{math.PR/0312099}
\def\GFFSurvey{MR2322706}

\def\GlimmJaffe{MR887102}
\def\GawcedzkiCFTStrings{MR1701610}
\def\diFrancescoetalCFT{MR1424041}
\def\BPZ{MR757857}
\def\denNijs{MR690541}
\def\Beffara{math.PR/0211322}

\def\noopsort#1{}

\begin{document}
\maketitle

\begin{abstract} We prove that the chordal contour lines of the discrete
Gaussian free field converge to forms of SLE$(4)$. Specifically, there
is a constant $\lambda > 0$ such that when $h$ is an interpolation
of the discrete Gaussian free field on a Jordan domain --- with boundary
values $-\lambda$ on one boundary arc and $\lambda$ on the complementary arc
--- the zero level line of $h$ joining the endpoints of these arcs
converges to SLE$(4)$ as the domain grows larger. If instead the boundary
values are $-a < 0$ on the first arc and $b > 0$ on the
complementary arc, then the convergence is to
SLE$(4;a/\lambda-1,b/\lambda-1)$, a variant of SLE$(4)$.
 \end{abstract}

\newpage
\tableofcontents
\newpage

\section{Introduction} \label{introsection}

\subsection{Main result}
The two-dimensional massless Gaussian free field (GFF) is a
two-dimensional-time analog of Brownian motion.  Just as Brownian
motion is a scaling limit of simple random walks and various other
one-dimensional systems, the GFF is a scaling limit of several
discrete models for random surfaces. Among these is the discrete
Gaussian free field (DGFF), also called the harmonic crystal. We
presently discuss the basic definitions and describe the main
results of the current work, postponing an overview of the history
and general context to Section \ref{historysection}.

Let $G=(V,E)$ be a finite graph and let $V_\p\subset V$ be some nonempty
set of vertices.
Let $\Omega$ be the set of functions $h:V\to\R$ that are zero on $V_\p$.
Clearly, $\Omega$ may be identified with $\R^{V\setminus V_\p}$.
The DGFF on $G$ with zero boundary values on $V_\p$
is the probability measure on $\Omega$ whose density with respect to
Lebesgue measure on $\R^{V\setminus V_\p}$ is proportional to
\begin{equation}\label{e.density}
\exp\Bigl(\sum_{\{u,v\}\in E} -\frac 12\,\bl(h(v)-h(u)\br)^2\Bigr).
\end{equation}
Note that under the DGFF measure, $h$ is a multi-dimensional Gaussian random
variable. Moreover, the DGFF is a rather natural discrete model for a random
field: the term $-(h(v)-h(u))^2/2$ corresponding to each edge $\{u,v\}$ penalizes
functions $h$ which have a large gradient along the edge.

Now fix some function $h_\p:V_\p\to\R$, and let $\Omega_{h_\p}$
denote the set of functions $h:V\to\R$ that agree with $h_\p$ on $V_\p$.
The probability measure on $\Omega_{h_\p}$ whose density with respect to
Lebesgue measure on
$\R^{V\setminus V_\p}$ is proportional to~\eref{e.density} is the DGFF
with boundary values given by $h_\p$.

Let $\TG$ be the usual triangular grid in the complex plane, i.e.,
the graph whose vertex set is the integer span of $1$ and
$e^{\pi\,i/3}=(1+\sqrt{3}\,i)/2$, with straight edges joining $v$ and
$w$ whenever $|v-w|=1$.  A {\bf $\TG$-domain} $D \subset \mathbb R^2
\cong \mathbb C$ is a domain whose boundary is a simple closed curve
comprised of edges and vertices in $\TG$.  Let $V=V_{\closure D}$ be the set of
$\TG$-vertices in the closure of $D$, let $G=G_D$ be the induced
subgraph of $\TG$ with vertex set $V_{\closure D}$, and write $V_\p =
\partial D \cap V_{\closure D}$.
While introducing our main results, we will
focus on graphs $G_D$ and boundary sets $V_\p$ of this form (though analogous
results hold if we replace $\TG$ with another doubly periodic planar
graph; see Section \ref{precisestatementsection}).

We may assume that any function $f:V \rightarrow \mathbb R$ is
interpolated to a continuous function on the closure of $D$ which is
affine on each triangle of $\TG$.  We often interpret $f$ as a
surface embedded in three dimensions and refer to $f(v)$ as the {\bf
height} of the surface at $v$.

\begin{figure}
\centerline{\epsfxsize=1.2\hsize%
\epsfbox{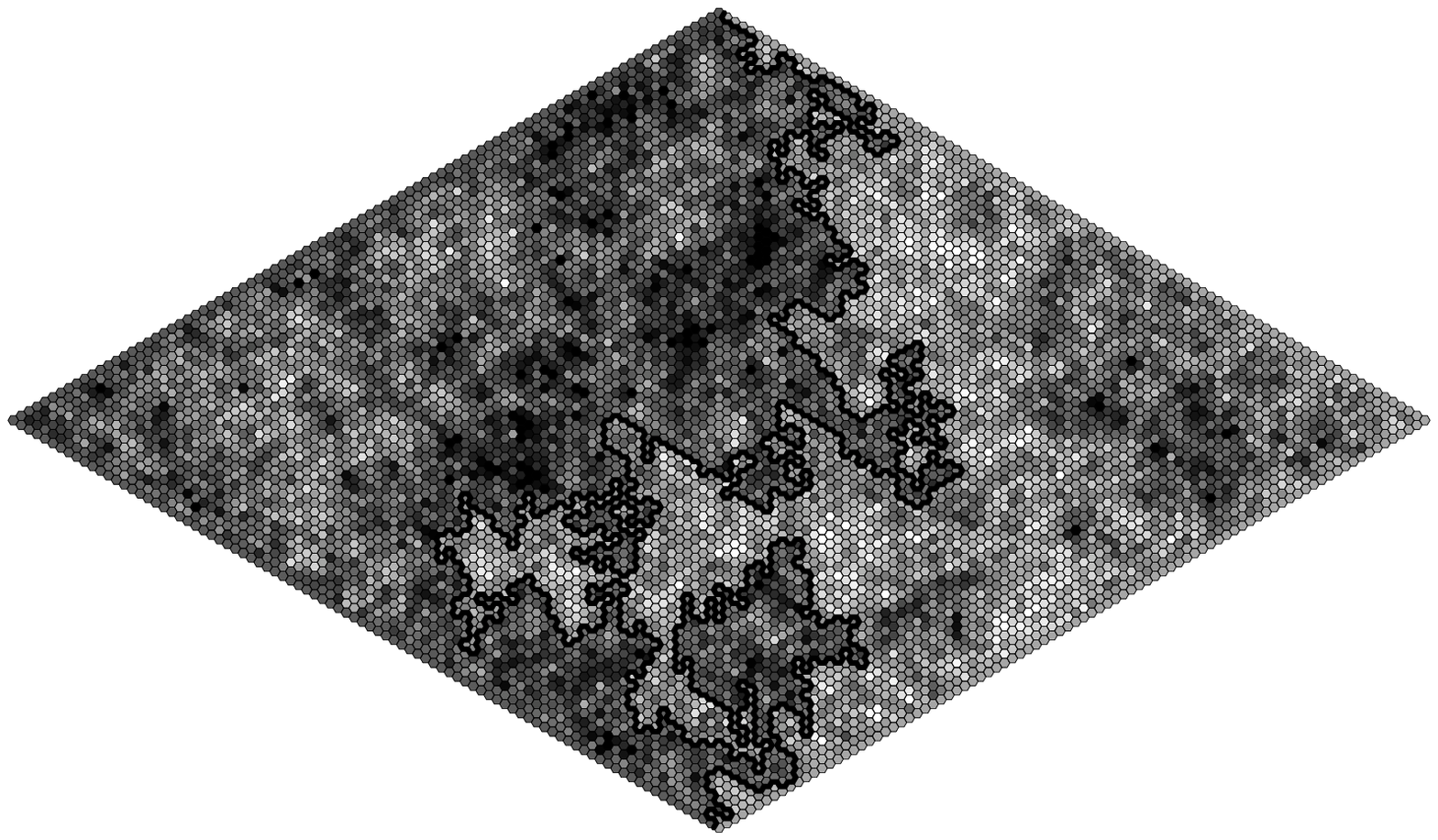}}
\vskip 0.5in
\centerline{\epsfxsize=1.2\hsize%
\epsfbox{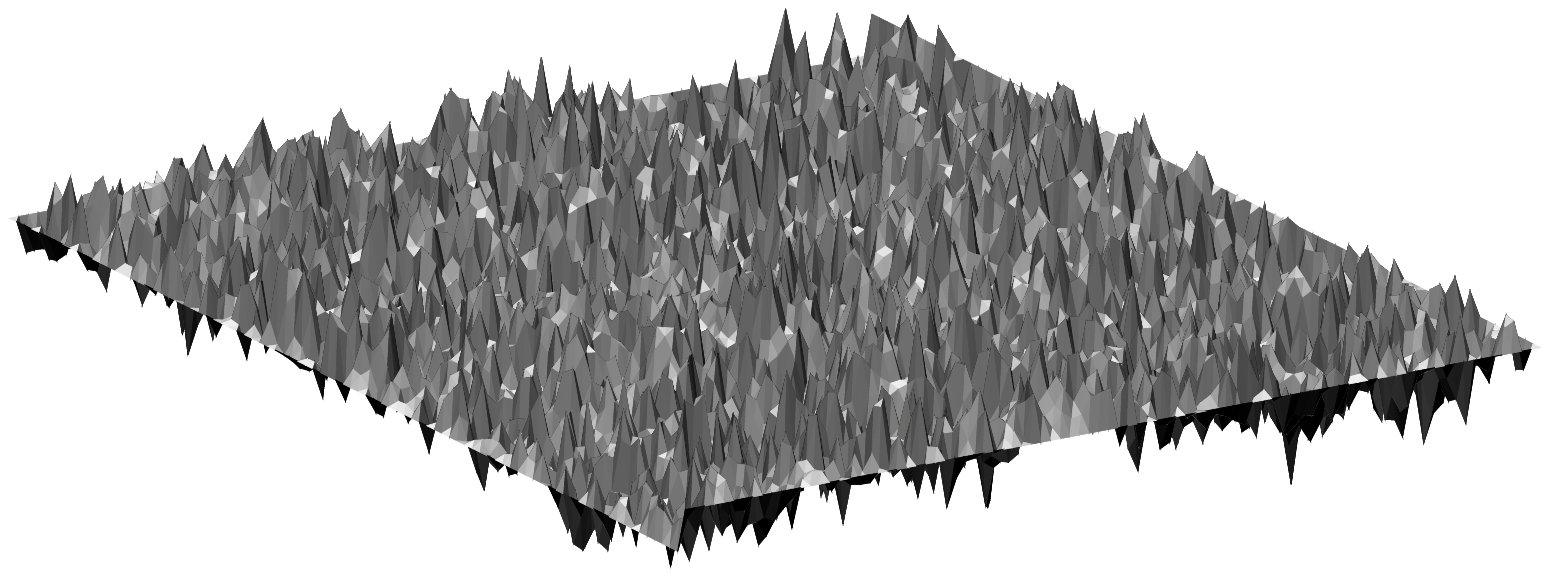}}
\vskip 0.5in
\old{
\caption{\label{fullDGFF}(a) DGFF on $90 \times 90$ hexagon array
with boundary values $-\tlambda$ on left, $\tlambda$ on right, where
$\tlambda = 3^{-1/4} \sqrt{\pi/8}$, faces shaded by height. (b)
Surface plot of DGFF.}
}
\caption{\label{fullDGFF}(a) DGFF on $90 \times 90$ hexagon array
with boundary values $\lambda$ on the right and
$-\lambda$ on the left; faces shaded by height. (b)
Surface plot of DGFF.}
\end{figure}

Let $\p D=\p_+\cup\p_-$ be a partition
of the boundary of a $\TG$-domain $D$ into two
disjoint arcs whose endpoints are midpoints of two distinct $\TG$-edges
in $\p D$.  Fix two constants $a,b>0$.
Let $h$ be an instance of the DGFF on $(G_D,V_\p)$, with boundary function
$h_\p$ equal to $-a$ on the vertices
in $\p_-$ and equal to $b$ on the vertices in
$\p_+$.
Then $h$ (linearly interpolated on triangles) almost surely assumes the
value zero on a unique piecewise linear path $\gamma_h$ connecting the two
boundary edges containing endpoints of $\p_+$.

\old{
Let $h$ be an instance of the DGFF on $(G_D,V_\p)$, with boundary function
$h_\p$ given by the constant $\lambda =
\tlambda:=3^{-1/4}\sqrt{\pi/8}$ (whose significance will be
explained later) on the vertices on one arc $\p_+\subset\partial D$
(whose endpoints are midpoints of edges of $\TG$) and $-\lambda$ on
vertices in the complementary arc $\p_-:=\p D\setminus\p_+$. Then
$h$ (linearly interpolated on triangles) almost surely assumes the
value zero on a unique piecewise linear path connecting the two
endpoints of $\p_+$.
}

In Section \ref{ss.definingsle}, we will briefly review the
definition of \SLEkk 4/ (a particular type of random chordal path
connecting a pair of boundary points of $D$ whose randomness comes
from a one-dimensional Brownian motion), along with the variants
of \SLEkk4/ denoted \SLEkr 4; \rho_1, \rho_2/.
Our main result, roughly stated, is the following:

\begin{theorem} \label{SLEconvergence}
Let $D$ be a $\TG$-domain, $\p D=\p_+\cup \p_-$
and let $h$ and $\gamma_h$ be as above.
There is a constant $\lambda>0$
 such that if $a=b=\lambda$, then as the triangular mesh gets finer,
the random path $\gamma_h$ converges in distribution to \SLEkk 4/.
If $a,b\ge \lambda$ are not assumed to equal $\lambda$,
 then the convergence is to \SLEab/.
\end{theorem}

See Section \ref{precisestatementsection} for a more precise version,
which describes the topology under which convergence is attained.
As explained there, we can also prove convergence in a
weaker form when the conditions $a,b\ge \lambda$ are relaxed.

We will elaborate on the role of the constant
$\lambda$ in \S\ref{heightgapsection}. This constant depends only
on the lattice used.
Although we do not prove it in this paper,
for the triangular grid the value of $\lambda$
is $\tlambda:=3^{-1/4}\sqrt{\pi/8}$ (see Section \ref{ss.sequel}).

\medskip
Figure \ref{fullDGFF} illustrates a dual perspective on an instance
of $\gamma_h$.  Here, each vertex in the closure of a
rhombus-shaped $\TG$-domain $D$ is replaced with a hexagon in the
honeycomb lattice.  Call hexagons {\bf positive} or {\bf negative}
according to the sign of $h$.  Then there is a cluster of positive
hexagons that includes the positive boundary hexagons, a similar
cluster of negative hexagons, and a path $\gamma$ forming the
boundary between these two clusters.  Figure \ref{fullDGFF} depicts
a computer generated instance of the DGFF---with $\pm \lambda$
boundary conditions---and the corresponding $\gamma$. Followed from
bottom to top, the interface $\gamma$ turns right when it hits a negative
hexagon, left when it hits a positive hexagon.  It closely tracks
the boundary-hitting zero contour line $\gamma_h$ in the following
sense: the edges in $\gamma$ are the duals of the edges of $\TG$
that are crossed by $\gamma_h$. This is because $h$ is almost
surely non-zero at each vertex in $V$, so whenever a zero contour
line contains a point on an edge of $\TG$, $h$ must be positive on
one endpoint of that edge and negative on the other; hence the dual
of that edge separates a positive hexagon from a negative hexagon.

In the fine mesh limit, there will be no difference between $\gamma$
and $\gamma_h$.  Thus (by Theorem \ref{SLEconvergence}) the path in
Figures \ref{fullDGFF}.a and~\ref{conditionedDGFF}.a approximates
\SLEkk 4/, while the path of Figure \ref{otherboundaryDGFF}.a
approximates \SLEkr 4; 2, 2/.  We will state and prove most of our
results in terms of the dual perspective displayed in the figures.

\subsection{Conditional expectation and the height gap} \label{heightgapsection}

We derive the following well known facts as a warm-up in
Section~\ref{ss.generalremarks}
(see also, e.g., \cite{\GiacominSurvey}):

\vspace{.1 in} {\bf
\begin{noindent}Boundary influence:\end{noindent}} The law of the DGFF with boundary
conditions $h_\p : V_\p \rightarrow \R$ is the same as that of the
DGFF with boundary conditions $0$ {\em plus} a deterministic
function $\tilde h_\p: V \rightarrow \R$ which is the unique
discrete harmonic interpolation of $h_\p$ to $V$. (By {\bf
discrete harmonic} we mean that for each $v \in V \backslash V_\p$,
the value $h(v)$ is equal to the average value of $h(w)$ over $w$
adjacent to $v$.) In particular, the expected value of $h(v)$ is
discrete harmonic in $V \backslash V_\p$.

\vspace{.1 in} {\bf \begin{noindent}Markov property:\end{noindent}}
Let $h:V \rightarrow \R$ be a random function whose law is the DGFF
on $G$ with some boundary values $h_\p$ on $V_\p$. Then
given the values of $h$ on a superset $V_0\supset V_\p$, the
conditional law of $h$ is that of a DGFF on $G$ with boundary set
$V_0$ and with boundary values equal to the given values. \vspace{.1
in}

\begin{figure}
\centerline{\epsfxsize=1.2\hsize%
\epsfbox{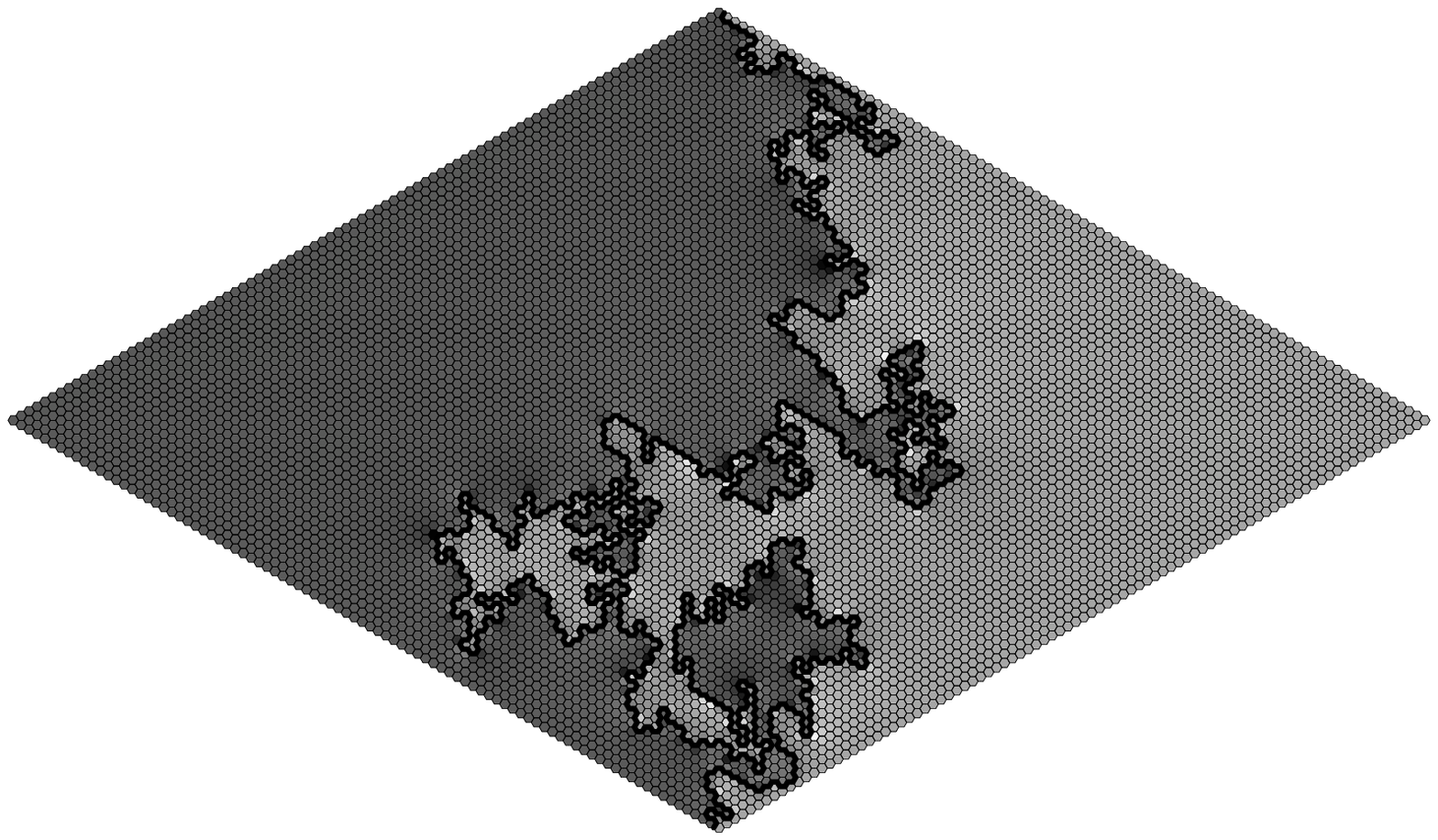}}
\vskip 0.5in
\centerline{\epsfxsize=1.2\hsize%
\epsfbox{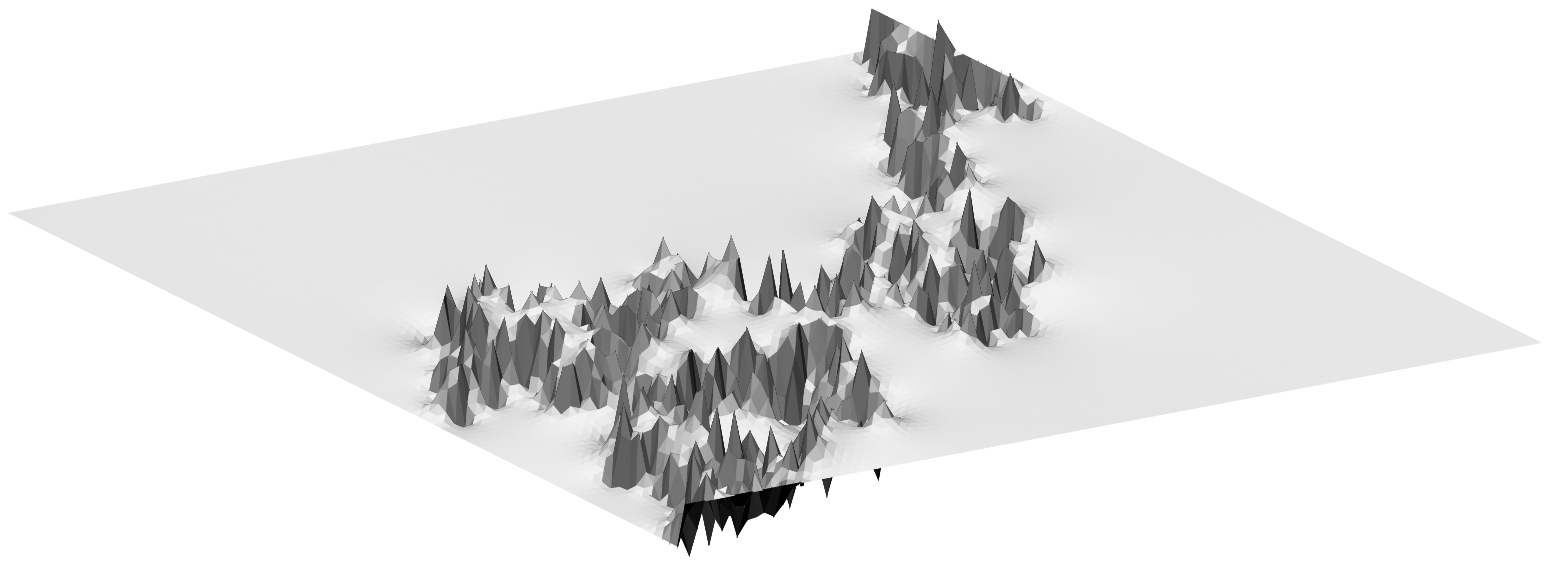}}
\vskip 0.5in
\caption{\label{conditionedDGFF}(a) Expectation of DGFF
with boundary values $\pm\lambda$ given its
values at hexagons bordering interface. (b) Surface plot of the
above.}
\end{figure}

\begin{figure}
\centerline{\epsfxsize=1.2\hsize%
\epsfbox{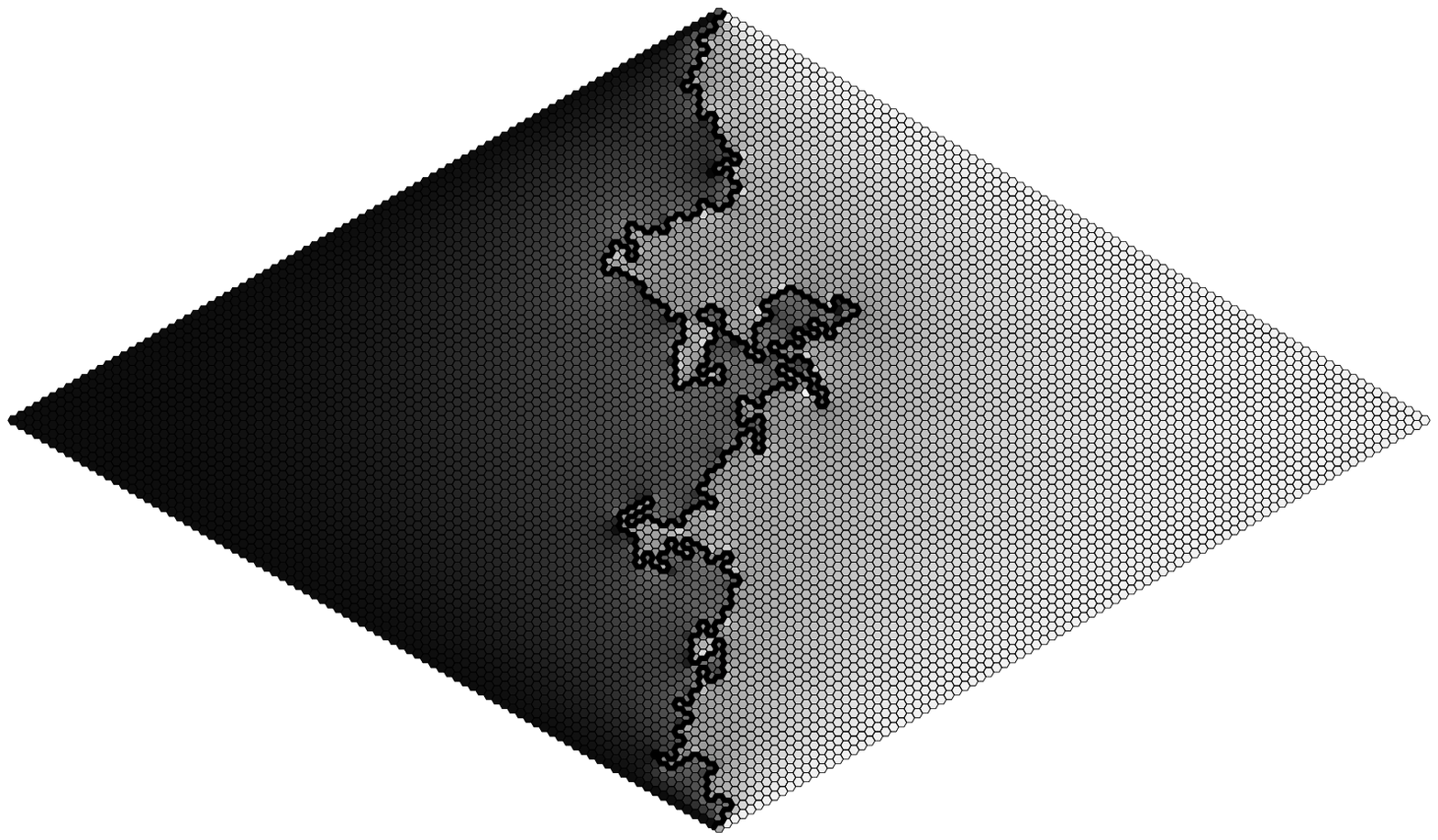}}
\vskip 0.5in
\centerline{\epsfxsize=1.2\hsize%
\epsfbox{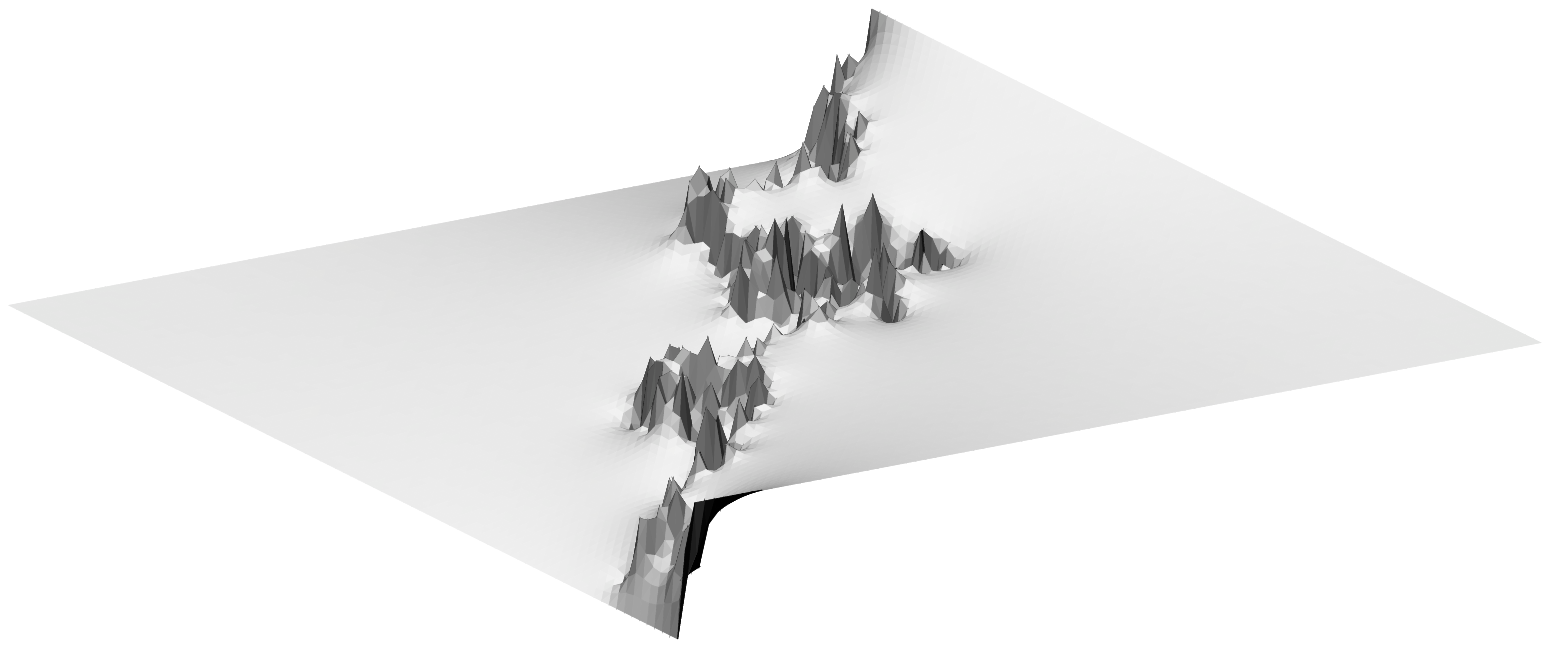}}
\vskip 0.5in
\caption{\label{otherboundaryDGFF} (a) Expectation of DGFF given its
values at hexagons bordering interface; exterior boundary values are
$-3\lambda$ on left, $3\lambda$ on right. (b) Surface plot of the
above.}
\end{figure}

From these facts it follows that conditioned on the path $\gamma$
described in the previous section and on the values of $h$ on the
hexagons adjacent to $\gamma$, the expected value of $h$ is discrete
harmonic in the remainder of $G_D$.  Figures \ref{conditionedDGFF}
and \ref{otherboundaryDGFF} illustrate the expected value of $h$
conditioned on the values of $h$ on the hexagons adjacent to
$\gamma$.

The reader may observe in Figure~\ref{conditionedDGFF} that although the expected
value of $h$ given the values along $\gamma$ varies a great deal
among hexagons close to $\gamma$, the expected value at five or ten
lattice spacings away from $\gamma$ appears to be roughly constant
along either side of $\gamma$.
On the other hand in Figure~\ref{otherboundaryDGFF},
away from $\gamma$, the expected height appears to be a smooth but
non-constant function.
 In a sense we make precise in Section
\ref{tailtrivialitysection} (see Theorem \ref{t.gap}),
the values $-\lambda$ and $\lambda$ describe the expected value of
$h$, conditioned on $\gamma$, at the vertices near (but not
microscopically near) the left and right sides of $\gamma$; in the
fine mesh limit there is thus an ``expected height gap'' of $2\,
\lambda$ between the two sides of $\gamma$.
In Figure~\ref{conditionedDGFF} the height expectation
appears constant away from $\gamma$, because
the boundary values of $\pm\lambda$ are the same as the
expected values near (but not microscopically near) $\gamma$.

Once we have established the height gap result, the proof of Theorem
\ref{SLEconvergence} (at least for the simplest case that the boundary
conditions are $-\lambda$ and $\lambda$) is similar to the
proof that the harmonic explorer converges to \SLEkk 4/, as given by
the present authors in \cite{\SchrammSheffieldHE},
which, in turn, follows the same strategy as the
proof of convergence of the loop-erased random walk to \SLEkk 2/
and the uniform spanning tree Peano curve to \SLEkk 8/
in~\cite{\LSWlesl}.

We will now briefly describe some of the key ideas in the proof of
the height gap result. The main step is to show that if one samples
a vertex $z$ on $\gamma$ according to discrete harmonic measure
viewed from a typical point far away from $\gamma$, then the
absolute value of $h(z)$ is close to independent of the values of
$h$ (and geometry of $\gamma$) at points that are {\em not}
microscopically close to $z$. In other words, if we start a random
walk $S$ at a typical point in the interior of $D$ and stop the first
time it hits a vertex $z$ which either belongs to $V_\p$ or
corresponds to a hexagon incident to $\gamma$, then (conditioned on
$z \not \in V_\p$), the random variable $|h(z)|$ (and in particular its
conditional expectation) is close to independent of the behavior of $\gamma$ and
$h$ at vertices far away from $z$.

To prove this, we will actually prove something stronger, namely
that (up to multiplication by $-1$) the collection of all of the
values of $h$ (and the geometry of $\gamma$) in a microscopic
neighborhood of $z$ is essentially independent of the values of
$h$ (and the geometry of $\gamma$) at points that are not
microscopically close to $z$. One consequence of our analysis is
Theorem \ref{t.limit}, which states that if one takes $z$ to be
the origin of a new coordinate system and conditions on the behavior
of $\gamma$ and $S$ outside of a ball of radius $R$ centered at
$z$ and $S$ starts outside that ball,
then as $R$ tends to infinity the conditional law of the
interface $\gamma$ has a weak limit (which is independent of the
sequence of boundary conditions chosen), which is the law of a
random infinite path $\gamma$ on
the honeycomb grid $\TG^*$ (almost surely containing
an edge adjacent to the hexagon centered at the origin $z = 0$).
We will define a function of such infinite paths $\gamma$ which (in
a certain precise sense) describes the expected value of $|h(z)|$
conditioned on $\gamma$; the value $\lambda$ is the expectation of
this function when $\gamma$ is chosen according to the limiting
measure described above.

We remark that many important problems in statistical physics
involve classifying the measures that can arise as weak limits of
Gibbs measures on finite systems.  In such problems, showing the
uniqueness of the limiting measure often involves proving that
properties of a random system near the origin are approximately
independent of the properties of the system far away from the
origin.  In our case, we need to prove that in some sense the
behavior of the triple $(h,\gamma,S)$ near the the origin (i.e., the
first point $S$ hits $\gamma$) is close to independent of the
behavior of $(h,\gamma,S)$ far from the origin.

Very roughly speaking, our strategy will be to describe the joint
law of $(h,\gamma,S)$ near the origin and $(h,\gamma,S)$ far from
the origin by considering a
different measure in which the two are independent
and weighting it by the probability that the inside and outside
configurations properly ``hook up'' with one another.  To get a
handle on these ``hook up'' probabilities, we will need to develop
various techniques to control the probabilities (conditioned on the
values of $h$ on certain sets) that certain zero height level lines
hook up with one another, as well as the probabilities that these
level lines avoid certain regions. We will also need bounds on the
probability that there exist clusters of positive or negative
hexagons crossing certain regions; these are roughly in the spirit
of the Russo-Seymour-Welsh theorems for percolation, but the proofs
are entirely different.  All of the height gap related results are proved in
Section \ref{tailtrivialitysection}.

\subsection{GFF definition and background} \label{historysection}

To help put our DGFF theorems in context and provide further
intuition, we now briefly recall the definition of the (continuum)
GFF and mention some basic facts described, e.g., in
\cite{\GFFSurvey}. Let $H_s(D)$ be the set of smooth functions
supported on compact subsets of a planar domain $D$, and let $H(D)$
be its Hilbert space completion under the {\bf Dirichlet inner
product} $(f, g)_{\nabla} = \int_{D}\nabla f \cdot \nabla g\;dx$,
where $dx$ refers to area measure. We define an instance of the
Gaussian free field to be the formal sum $$h = \sum_{i = 1}^{\infty}
\alpha_i f_i,$$ where $\alpha_i$ are i.i.d.\ one-dimensional
standard (unit variance, zero mean) Gaussians and the $f_i$ are an
orthonormal basis for $H(D)$. Although the sum does not converge
pointwise or in $H(D)$, it does converge in the space of
distributions \cite{\GFFSurvey}. In particular, the sum
$$(h, g)_\nabla := \sum_{i=1}^\infty\alpha_i\,\bigl( f_i, g\bigr)_\nabla$$
is almost surely convergent for every $g \in H_s(D)$.

It is worthwhile to take a moment to compare with the situation
where $D$ is one dimensional. If $D$ were a bounded open interval in
$\R$, then the partial sums of $h=\sum_{i=1}^\infty\alpha_i f_i$
would almost surely converge uniformly to a limit, whose law is that
of the Brownian bridge, having the value zero at the interval's
endpoints.  If $D$ were the interval $(0, \infty)$, then the partial
sums would converge (uniformly on compact sets) to a function whose
law is that of ordinary Brownian motion $B_t$, indexed by $t \in [0,
\infty)$, with $B_0=0$ \cite{\GFFSurvey}.

Let $g$ be a {\bf conformal} (i.e., bijective analytic) map from $D$
to another planar domain $D'$. When $g$ is a rotation, dilation, or
translation, it is obvious that
$$\int_{D'} \nabla (f_1 \circ g^{-1} ) \cdot \nabla (f_2 \circ g^{-1} )\,dx = \int_{D} (\nabla f_1 \cdot \nabla f_2)\,dx.$$
for any $f_1, f_2 \in H_s(D)$, and an elementary change of variables
calculation gives this equality for any conformal $g$. Taking the
completion to $H(D)$, we see that the Dirichlet inner product---and
hence the two-dimensional GFF---is invariant under conformal
transformations of $D$.

Up to a constant, the DGFF on a $\TG$-domain $D$ can be realized as
a projection of the GFF on $D$ onto the subspace of $H(D)$
consisting of functions which are continuous and are affine on each
triangle of $D$ \cite{\GFFSurvey}.  Note that if $f$ is such a
function, then
$$(f,f)_{\nabla} = \frac{\sqrt 3}{6} \sum\bigl( |f(j) - f(i)|^2 +
|f(k)-f(i)|^2 + |f(k) - f(j)|^2\bigr),$$ where the sum is over all
triangles $(i,j,k)$ in $V_{\closure D}$. This is because the area of each
triangle is $\sqrt 3/4$ and the norm of the gradient squared in the
triangle is
$$\frac{2}{3} \left( |f(j) - f(i)|^2 + |f(k)-f(i)|^2 + |f(k) -
f(j)|^2 \right).$$ Since each interior edge of $D$ is contained in
two triangles, for such $f$,
\begin{equation}
\label{e.affineenergy} \|f\|^2_\nabla =3^{-1/2} \sum_{\{i,j\}\in
E_I} |f(j) - f(i)|^2 + \frac{3^{-1/2}}{2} \sum_{\{i,j\}\in E_\p}
|f(j) - f(i)|^2,
\end{equation}
where
$E_I$ and $E_\p$ are the interior and boundary (undirected) edges of
$\TG$ in $\closure{D}$.  We will refer to the sum $\sum_{E_I} |f(j)
- f(i)|^2$ as the {\bf {\em discrete} Dirichlet energy} of $f$. It
is equivalent---up to the constant factor $3^{-1/2}$ and an additive
term depending only on boundary values of $f$---to the {\bf
Dirichlet energy} $(f,f)_\nabla$ of the piecewise affine
interpolation of $f$ to $D$.

\medskip

The above analysis suggests a natural coupling between the GFF and a sequence of DGFF approximations
to the GFF (obtained by taking finer mesh approximations of the same domain).
The GFF can also be obtained as a scaling limit of other discrete random
surface models (e.g., solid-on-solid, dimer-height-function, and $\nabla
\phi$-interface models)
\cite{\KenyonDominoGFF,\NaddafSpencer,\SpencerSurvey}. Its Laplacian
is a scaling limit of some Coulomb gas models, which describe random
electrostatic charge densities in two-dimensional domains
\cite{\Foltin,\FroSpenAbelianSpin,\KostRoughening,\KostThouThree,\SpencerSurvey}.
Physicists often use heuristic connections to the GFF to predict
properties of two-dimensional statistical physics models that are
not obviously random surfaces or Coulomb gases (e.g., Ising and
Potts models, $O(n)$ loop models)
\cite{\denNijs,\diFrancescoetalCFT,\DuplantierSurvey,\Kadanoff,
\NienhuisKagerSurvey,\Nienhuiseightytwo,\NienhuisSurvey}. As a model
for the field theory of non-interacting massless bosons,
the GFF is a starting point for many constructions in quantum field
theory, conformal field theory, and string theory
\cite{\BPZ,\diFrancescoetalCFT,\GawcedzkiCFTStrings,\GlimmJaffe}.

Because of the conformal invariance of the GFF, physicists and
mathematicians have hypothesized that discrete random surface models
that are believed or known to converge to the Gaussian free field
(e.g., the discrete Gaussian free field, the height function of the
oriented $O(n)$ loop model with $n=2$, height functions for domino
and lozenge tilings) have level sets with conformally invariant
scaling limits \cite{\ConiglioPotts,\DupSalDenseSAW,
\DupSalWinding,\DupSalPercolationhull,\HuberDuKondevGeometryofloops,\HuberKondevGeometryofloops,\KenyonDominoGFF,
\KondevHenleyninetyfive,\KondevHenleytwothousand,\Nienhuiseightytwo}.
Our results confirm this hypothesis for the discrete Gaussian free
field.

Various properties of the DGFF contour lines (such as
winding exponents and the fact that the fractal dimension is $3/2$) have
been predicted correctly
in the physics literature
\cite{\ConiglioPotts,\DupSalWinding,\DupSalDenseSAW,
\DupSalPercolationhull,\HuberKondevGeometryofloops,\HuberDuKondevGeometryofloops,\KondevHenleyninetyfive,\KondevHenleytwothousand,
\Nienhuiseightytwo}.  The techniques used to make these predictions
are also described in detail in the survey papers
\cite{\DuplantierSurvey,\NienhuisKagerSurvey,\NienhuisSurvey}.
Analogous results about winding exponents and fractal
dimension have now been proved rigorously for \SLE/ \cite{\SchSLE,
\RSsle, \Beffara}.

The study of level lines of the DGFF and related random surfaces is
also related to the study of equipotential lines of random charge
distributions in statistical physics.  The so-called {\bf
two-dimensional Coulomb gas} is a model for electrostatics in which
the force between charged particles is inversely proportional to the
distance between them. In this model, a continuous function $f \in
H_s(D)$ is the {\bf Coulomb gas electrostatic potential function}
(``grounded'' at the boundary of $D$) of $-\Delta f$, when $\Delta
f$ is interpreted as a charge density function.  The value
$(f,f)_{\nabla}$ is then the total potential energy---also called
the {\bf energy of assembly} of the charge distribution $-\Delta f$.
In the Coulomb gas model, this is the amount of energy required to
move from a configuration in which the charge density is zero
throughout $D$ to a configuration in which the charge density is
given by $-\Delta f$.

In statistical physics, it is often natural to consider a
probability distribution on configurations in which the probability
of a configuration with potential energy $H$ is proportional to
$e^{-H}$. If $\rho$ is a smooth charge distribution, then its energy
of assembly is given by $(-\Delta^{-1} \rho, -\Delta^{-1}
\rho)_\nabla = (\rho, -\Delta^{-1} \rho)$; if we define $\rho$ to be
the standard Gaussian in $\Delta H(D)$ determined by this quadratic
form, then $\rho$ is the Laplacian of the Gaussian free field
(which, like the GFF itself, is well defined as a random
distribution but not as a function). In other words, the Laplacian
of a Gaussian free field is a random distribution that we may
interpret as a model for random charge density in a statistical
physical Coulomb gas.

However, we stress that when physicists refer to the Coulomb gas method
for $O(n)$ model computations, they typically have in mind a more
complicated Coulomb gas model in which the charges are required to be
discrete (i.e., $\rho$ is required to be a sum of unit positive and
negative point masses) and hard core constraints may be enforced.

The surveys \cite{\BriMelFroSurvey,\GiacominSurvey,\SpencerSurvey}
contain additional references on lattice spin models that have the
GFF as a scaling limit and Coulomb gas models that have its
Laplacian as a scaling limit---for example, the {\bf harmonic
crystal} (a.k.a.\ the {\bf discrete Gaussian free field}) with
quadratic nearest-neighbor potential, the more general {\bf
anharmonic crystal}, the {\bf discrete-height Gaussian} (where $h$
is a function on a lattice, with values restricted to integers), the
{\bf Villain gas} (where $h$ is a function on a lattice and the
values of its discrete Laplacian $\density = - \Delta h$ are
restricted to integers), and the {\bf hard core Coulomb gas} (where
$h$ is a function on a lattice and its discrete Laplacian $\density
= - \Delta h$ is $\pm 1$ valued).

The physics literature on applications of the GFF to field theory
and statistical physics is large, and the authors themselves are
only familiar with parts of it.  Outside of these areas, there is a
body of experimental and computational research on contour lines of
random topographical surfaces, such as the surface of the earth.
Mandelbrot's famous {\em How Long Is the Coast of Britain?...}
\cite{\Mandelbrot}, which prefigured the notion of ``fractal''
introduced by Mandelbrot years later, is an early example.  The
results about contour lines in these studies (including fractal
dimension computations) are less detailed than the ones provided
here and are not all mathematically rigorous. However, some of the
models are similar in spirit to the GFF, involving functions whose
Fourier coefficients are independent Gaussians. An eclectic overview
of this literature appears in \cite{\IsichenkoSurvey}.

\subsection{\SLE/ background and prior convergence results}\label{ss.definingsle}

We now give a brief definition of (chordal) \SLEk/ for $\kappa > 0$.
See also the surveys~\cite{\WernerStFlour, \NienhuisKagerSurvey,
\LawlerSLbook, \CardySLESurvey}
or~\cite{\LawlerSLEintro}.  The discussion below along with further discussion
of the special properties of \SLEkk4/ appears in another paper by the current authors~\cite{\SchrammSheffieldHE}.
That paper shows that \SLEkk4/ is the scaling limit of a random interface called the harmonic explorer
(designed in part to be a toy model for the DGFF contour line addressed here).

Let $T>0$.
Suppose that $\gamma:[0,T]\to\overline\H$ is a continuous simple path in the closed
upper half plane $\overline\H$ which satisfies $\gamma[0,T]\cap\R=\{\gamma(0)\}=\{0\}$.
For every $t\in[0,T]$, there is a unique conformal homeomorphism
$g_t:\H\setminus\gamma[0,t]$ which satisfies the so-called {\bf hydrodynamic}
normalization at infinity
$$
\lim_{z\to\infty} g_t(z)-z=0\,.
$$
The limit
$$
\capacity(\gamma[0,t]):=\lim_{z\to\infty} z(g_t(z)-z)/2
$$
is real and monotone increasing in $t$.
It is called the (half plane) {\bf capacity} of $\gamma[0,t]$ from $\infty$, or
just capacity, for short.
Since $\capacity(\gamma[0,t])$ is also continuous in $t$,
it is natural to reparameterize $\gamma$ so that
$\capacity(\gamma[0,t])=t$.
Loewner's theorem states that in this case the maps $g_t$ satisfy
his differential equation
\begin{equation}\label{e.chordal} \p_t g_t(z) =
\frac {2}{g_t(z)-W_t}\,,\qquad g_0(z)=z\,, \end{equation}
where $W_t=g_t(\gamma(t))$.
(Since $\gamma(t)$ is not in the domain of definition of $g_t$, the expression
$g_t(\gamma(t))$ should be interpreted as a limit of
$g_t(z)$ as $z\to\gamma(t)$ inside $\H\setminus\gamma[0,t]$.  This limit
does exist.)   The function $t \to W_t$ is
continuous in $t$, and is called the {\bf driving parameter} for $\gamma$.

One may also try to reverse the above procedure.
Consider the Loewner evolution defined by the ODE~\eref{e.chordal},
where $W_t$ is a continuous, real-valued function.
For a fixed $z$, the evolution
defines $g_t(z)$ as long as $|g_t(z)-W_t|$ is bounded away from zero.
For $z\in\closure{\H}$
let $\tau_z$ be the first time $t\ge0$ in which $g_t(z)$ and $W_t$ collide,
or set $\tau_z=\infty$ if they never collide.
Then $g_t(z)$ is well defined on $\{z\in\closure{\H}: \tau_z \ge t\}$.
The set $K_t:=\{z\in\closure\H: \tau_z\le t\}$ is sometimes
called the {\bf evolving hull} of the evolution.
In the case discussed above where the evolution is generated
by a simple path $\gamma$ parameterized by capacity
and satisfying $\gamma(t)\in\H$ for $t>0$, we have $K_t=\gamma[0,t]$.

The path of the evolution is defined as $\gamma(t)=\lim_{z\to W_t}g_t^{-1}(z)$,
where $z$ tends to $W_t$ from within the upper half plane $\H$, provided that the limit exists
and is continuous. However, this is not always the case.
The process (chordal)
\SLEkk\kappa/ in the upper half plane, beginning at $0$ and ending at $\infty$, is the path
$\gamma(t)$ when $W_t$ is $\sqrt\kappa\,B_t$, where $B_t=B(t)$ is a standard one-dimensional Brownian
motion.
(\lq\lq Standard\rq\rq\ means $B(0)=0$ and $\E[B(t)^2]=t$, $t\ge 0$.
Since $(\sqrt k\,B_t:t\ge 0)$ has the same distribution as
$(B_{\kappa\,t}:t\ge 0)$,
taking $W_t=B_{\kappa\,t}$ is equivalent.)
In this case a.s.\ $\gamma(t)$ does exist and is a continuous path. See~\cite{\RSsle} ($\kappa\ne 8$)
 and~\cite{\LSWlesl} ($\kappa=8$).

We now define the processes \SLEr \rho_1, \rho_2/.  Given a
Loewner evolution defined by a continuous $W_t$, we will let $x_t$
and $y_t$ be defined by $x_t:=\sup\{g_t(x):x<0, x\notin K_t\}$ and
$y_t:=\inf\{g_t(x):x>0, x\notin K_t\}$. When the Loewner evolution
is generated by a simple path $\gamma(t)$ satisfying
$\gamma(t)\in\H$ for $t>0$, these points $x_t$ and $y_t$ can be
thought of as the two images of $0$ under $g_t$. Note that
by~\eref{e.chordal}
\begin{equation}
\label{e.xy}
  \p_t x_t=2/(x_t-W_t),\qquad \p_t y_t=2/(y_t-W_t)\,,
\end{equation}
for all $t$ such that $x_t < W_t < y_t$.  Beginning from an
initial time $r$ for which $x_r < W_r < y_r$, we define \SLEr
\rho_1, \rho_2/ to be the evolution that makes $(x_t,W_t, y_t)$ a
solution to the SDE system
\begin{equation}\label{e.sdesystem}
dW_t
=\sqrt\kappa\,dB_t+\frac{\rho_1\,dt}{W_t-x_t}+\frac{\rho_2\,dt}{W_t-y_t},
\,\,\, dx_t=\frac{2\,dt}{x_t-W_t}, \,\,\,
dy_t\,=\,\frac{2\,dt}{y_t-W_t},
\end{equation}
noting that existence and uniqueness of solutions to this SDE (at
least from the initial time $r$ until the first $s>r$ for which
either $x_s=W_s$ or $W_s=y_s$) follow easily from standard results
in~\cite{\RevuzYor}.  (The $x_t$ and $y_t$ are called {\bf force
points} because they apply a ``force'' affecting the drift of the
process $W_t$ by an amount inversely proportional to their
distance from $W_t$.)

Some subtlety is involved in extending the definition of \SLEr
\rho_1, \rho_2/ beyond times when $W_t$ hits the force points, and
in starting the process from the natural initial values
$x_0=W_0=y_0 = 0$. This is closely related to the issues which come
up when defining the Bessel processes of dimension less than
two and will be discussed in more detail in Section \ref{s.drive}.

\bigskip

Although many random self-avoiding lattice paths from the statistical physics literature
are conjectured to have forms of \SLE/ as scaling limits, rigorous proofs have thus far appeared only for
a few cases: site
percolation cluster boundaries on
the hexagonal lattice (\SLEkk6/,~\cite{\SmirnovPerc};
see also~\cite{\CamiaNewmanSLE}), branches (loop-erased random walk) and outer boundaries
(random Peano curves) of uniform spanning
tree (forms of \SLEkk2/ and \SLEkk8/ respectively,~\cite{\LSWlesl}), the harmonic explorer
(\SLEkk4/,~\cite{\SchrammSheffieldHE}), and boundaries of
simple random walks (forms of \SLEkk{8/3}/,~\cite{\LSWrestriction}).

In the latter case, conformal invariance properties follow almost immediately from the conformal invariance of
two dimensional Brownian motion.  In each of the other cases listed above, the initial step of the proof is to show
that a certain function of the partially generated paths $\gamma([0,t])$, which is a martingale in $t$ when
$\gamma$ is \SLEk/ for the appropriate $\kappa$, has a discrete analog which is (approximately or exactly)
a martingale for the discrete paths and is approximately equivalent to the continuous version in the fine mesh limit.
For loop-erased random
walk, harmonic explorer, and uniform spanning tree Peano curves, this initial step is the easy part of the
argument; it follows almost immediately from the fact that simple random walk converges to Brownian motion.
The analogous step for site percolation on the hexagonal lattice, as
given by~\cite{\SmirnovPerc}, is an ingenious but nonetheless short and simple argument.

By contrast, the analogous step in this paper (which requires the proof of
the height gap lemma,
as given in Section \ref{tailtrivialitysection}) is quite involved; it is
the most technically
challenging part of the current work and includes many new techniques and lemmas about the
geometry of DGFF contours that we hope are interesting for their own sake.

Another way in which the DGFF differs from percolation, the harmonic
explorer, and the uniform spanning tree
is that it has a natural continuum analog (the GFF) which can be easily rigorously constructed without any
reference to \SLE/, and which is itself (like Brownian motion) an object of great significance.
It becomes natural to ask
whether the DGFF results enable us to define the ``contour lines'' of the continuum GFF in a canonical way;
we plan to answer this question (affirmatively) in a subsequent work (see Section \ref{ss.sequel}).

A final difference is that, for the DGFF, there is a continuum of choices for left and right boundary
conditions ($a$ and $b$) which are equally natural a priori, so we are led to consider a family of paths \SLEab/
instead of simply \SLEkk4/.  (The case $a=b=0$ is particularly natural; see Figure \ref{boundaryzeroset}.)
In these processes, the driving parameters $W_t$ are generally no longer
Brownian motions (rather, they are continuous semimartingales with constant quadratic variation and a
drift term that can become singular on a fractal set).  Proving driving parameter convergence
to these processes requires
some rather general convergence infrastructure (Section \ref{s.approximatediffusions}), which
we hope will be useful in other settings as well.

\subsection{Precise statement of main result} \label{precisestatementsection}

Let $\H$ be the upper half plane.  Let $D$ be any $\TG$-domain and let
$\p D=\p_+\cup\p_-$ be a partition of the boundary of $D$
into two disjoint arcs whose endpoints are midpoints of two
$\TG$-edges contained in $\p D$.
As before, let $V$ denote the vertices of $\TG$ in $\closure D$.
Let $h_\p=-a$ on $\p_-\cap V$ and $h_\p=b$ on $\p_+\cap V$,
where $a$ and $b$ are positive constants.

Let $h:V\to \R$ be an instance of the
DGFF with boundary conditions $h_\partial$.  Let $\phi_D$ be any
conformal map from $D$ to $\H$ that maps $\p_+$ bijectively onto the
positive real ray $(0, \infty)$. (Note that $\phi_D$ is unique up to
positive scaling.)

There is almost surely a unique interface $\gamma\subset\closure D$
between hexagons in the dual grid $\TG^*$
containing $\TG$-vertices where $h$ is positive and
such hexagons where $h$ is negative such that the endpoints
of $\gamma$ are on $\p D$. In fact, the endpoints of $\gamma$
are the same as the endpoints of $\p_+$.
(As mentioned above, this interface $\gamma$ stays within
a bounded distance from the
zero height contour line $\gamma_D$ of the affine interpolation of $h$.)
Now,
$\phi_D \circ \gamma$ is a random path on $\H$ connecting $0$ to
$\infty$. We will show that this path converges to a form of
\SLEkk4/. Rather than considering a fixed domain $\hat D$ and a
sequence of discrete domains $D_n$ approximating $\hat D$, with the mesh
tending to $0$, we will employ a setup that is more general in which
the mesh is fixed (the triangular lattice will not be rescaled), and
we consider domains $D$ that become \lq\lq larger\rq\rq. The correct
sense of \lq\lq large\rq\rq\ is measured by
$$\rr=r_{D,\phi} := \inr{\phi_D^{-1}(i)}(D)\,,$$ where $\inr{x}(D)$ denotes
the radius of $D$ viewed from $x$, i.e., $\inf_{y \not \in D}|x-y|$.
Of course, if $\phi_D^{-1}(i)$ is at bounded distance from $\p D$,
then the image of the triangular grid under $\phi_D$ is not fine
near $i$, and there is no hope for approximating \SLE/ by
$\phi_D\circ\gamma$.

We have chosen to use $\H$ as our canonical domain (mapping all
other paths into $\H$), because it is the most convenient domain in
which to define chordal \SLE/.  However, to make
the completion of $\H$ a compact
metric space, we will endow $\H$ with the metric it inherits from
its conformal map onto the unit disc $\U$. Namely, we let $d_*(\cdot,\cdot)$
be the metric on $\closure\H\cup\{\infty\}$ given by
$d_*(z,w)=|\Psi(z)-\Psi(w)|$, where $\Psi(z):=(z-i)/(z+i)$ maps
$\closure\H\cup\{\infty\}$ onto $\closure\U$. If $z\in\closure\H$,
then $d_*(z_n,z)\to 0$ is equivalent to $|z_n-z|\to 0$, and
$d_*(z_n,\infty)\to 0$ is equivalent to $|z_n|\to\infty$.

If $\gamma_1$ and $\gamma_2$ are distinct unparametrized simple
paths in $\closure \H$, then we
define $\dstrong(\gamma_1, \gamma_2)$ to be the infimum over all
pairs $(\eta_1, \eta_2)$ of parameterizations of $\gamma_1$ and
$\gamma_2$ in $[0,1]$ (i.e., $\eta_j:[0,1]\to \closure \H$
is a simple path satisfying $\eta_j([0,1])=\gamma_j$ for $j=1,2$)
 of the uniform distance
$\sup\bl\{d_*\bl(\eta_1(t),\eta_2(t)\br): t\in[0,1]\br\}$
with respect to the metric $d_*$.

Our strongest result is in the case where $a,b\ge \lambda$.
We prove the following:

\begin{theorem}\label{preciseSLEconvergence}
There is a constant $\lambda>0$ such that
if $a,b\ge\lambda$,
then as $\rr\to\infty$, the random paths
$\phi_D\circ\gamma$ described above converge in distribution to
\SLEab/ with respect to
the metric $\dstrong$.

In other words, for every $\eps>0$ there is some $R=R(\eps)$ such
that if $\rr>R$, then there is a coupling of $\phi_D \circ \gamma$
and a path $\gSLE$ whose law is that of
\SLEab/ such that
$$
\PB{\dstrong(\phi_D \circ \gamma, \gSLE)>\eps}<\eps\,.
$$
\end{theorem}

We first comment that it follows that when $r_D$ is large,
$\gamma$ is \lq\lq close\rq\rq\ to $\phi_D^{-1}\circ \gSLE$.
For example, if $r_D\to\infty$ and $r_D^{-1}\,D$ tends to a bounded domain
$\hat D$ whose boundary is a simple closed path in such a
way that the boundaries of the domains may be parameterized to
give uniform convergence of parameterized paths
and if $r_D^{-1}\p_+$ converges, then $r_D^{-1}\,\gamma$ converges
in law to the corresponding SLE in $\hat D$.
To prove this from Theorem~\ref{preciseSLEconvergence}, we
only need to note that
in this case the maps $r_D^{-1}\phi_D^{-1}$
converge uniformly in $\closure\H\cup\{\infty\}$
(see, e.g.,~\cite[Proposition 2.3]{\PommeBDRY}).

\medskip

When we relax the assumption $a,b\ge\lambda$ to $a,b>0$, we still prove
some sort of convergence to \SLEab/, but with respect to
a weaker topology.
In fact, we can allow $a$ and $b$ to be zero or even slightly negative,
but in this case we need to
appropriately adjust the above definition of the interface $\gamma$.
Say that a hexagon in the hexagonal grid $\TG^*$ dual to $\TG$ is positive
if either the center $v$ of the hexagon is in $D$ and $h(v)>0$,
or $v\in \p_+$.
Likewise, say that the hexagon is negative if $v\in D$ and $h(v)<0$,
or $v\in\p_-$.
Let $\gamma$ be the unique oriented path in $\TG^*$ that joins
the two endpoints of $\p_+$, has only positive hexagons
adjacent to its right hand side
and only negative hexagons adjacent to its left hand side.
(If $a,b>0$, this definition clearly agrees with the previous definition
of $\gamma$.)
We prove:

\begin{theorem}\label{absmall}
For every constant $\upperhco>0$ there is a
constant $\lowerhco=\lowerhco(\upperhco)>0$ such that if
$a,b\in[-\lowerhco,\upperhco]$
and $\gamma$ is the DGFF interface defined above,
then as $r_D\to\infty$
the Loewner driving term $\hat W_t$ of
$\phi_D\circ\gamma$ parameterized by capacity from $\infty$
converges in law to the driving term $W_t$ of
\SLEab/ with respect to the topology of locally uniform
convergence.
That is,
for every $T,\eps>0$ there is some $R>0$ such that if $r_D>R$,
then $\gamma$ and \SLEab/ may be coupled so that with probability
at least $1-\eps$
$$
\sup\bl\{|\hat W_t-W_t|:t\in[0,T]\br\}<\eps\,.
$$
\end{theorem}

Some (essentially well-known) geometric consequences of this
kind of convergence are proved in \S\ref{CKC}.

\begin{figure} \leavevmode \put(-65,-10){ \resizebox{18 cm}{!}
{\includegraphics{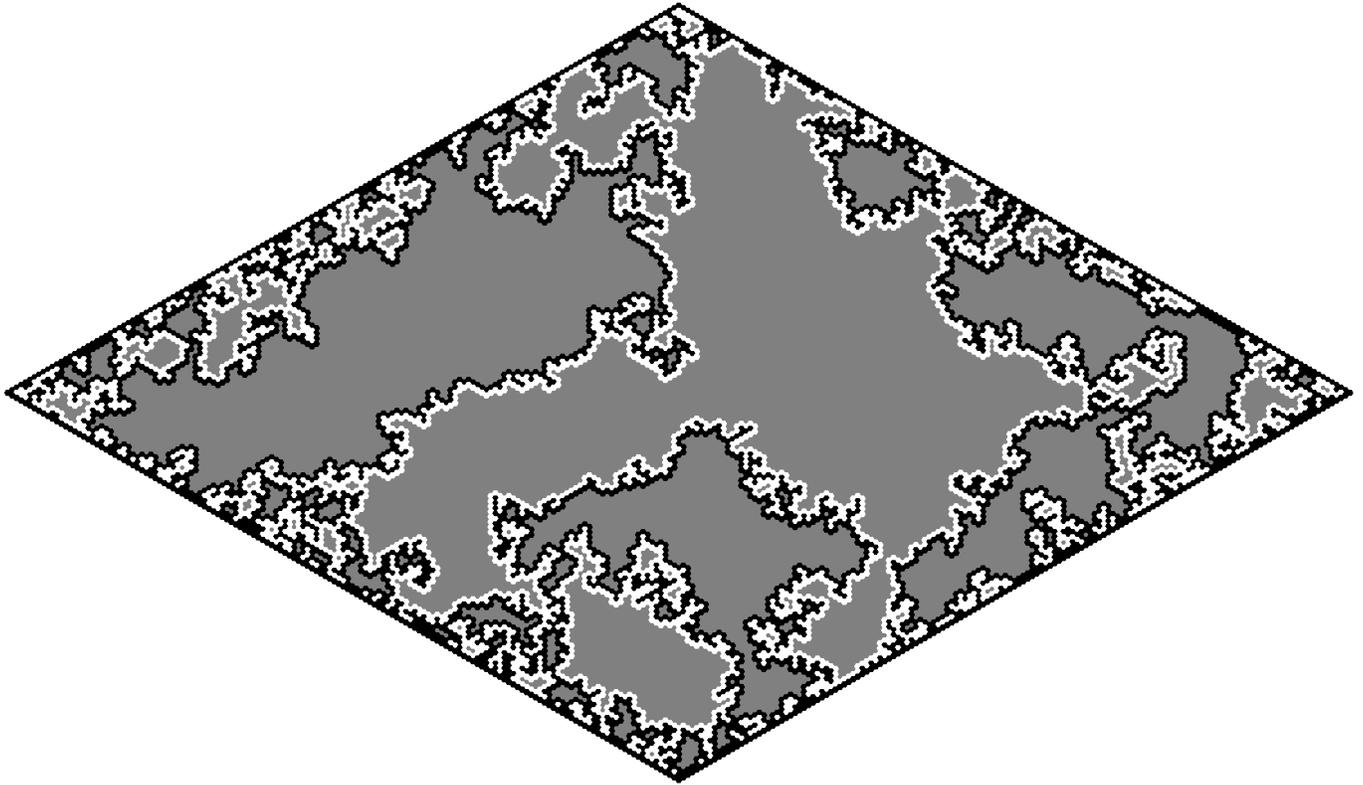}}}
\medskip
\caption{\label{boundaryzeroset} Zero height interfaces starting and ending
on the boundary, shown for the discrete GFF on a $150 \times 150$ hexagonal array
with zero boundary.   Interior white hexagons have height greater
than $0$; interior black hexagons have height less than $0$;
boundary hexagons (of height zero) are black.  Hexagons that are not incident
to a zero height interface that reaches the boundary are gray.}
\end{figure}

A particularly interesting case of Theorem \ref{absmall} is the case $a=0$ and $b=0$, corresponding to the DGFF with zero boundary values.  In this case,
when $h$ is interpolated linearly to triangles, its zero level set will almost surely include a
finite number of piecewise linear arcs in $D$ whose endpoints on $\partial D$ are vertices of $\TG$.
A dual representation of this set of arcs is shown in Figure \ref{boundaryzeroset}.
For any fixed choice of endpoints on the boundary, the interface connecting those endpoints will
converge to \SLEkr 4; -1, -1/.  The limit of the complete set of arcs in Figure \ref{boundaryzeroset}
is in some sense a coupling of \SLEkr 4; -1, -1/ processes, one for each pair of boundary points.

\medskip
Finally, we discuss the generalizations: replace $\TG$ with an arbitrary
weighted doubly periodic planar lattice---i.e., a
connected planar graph $\GG\subset\R^2$ invariant under
two linearly independent translations, $T_1,T_2$,
 such that every compact subset of $\R^2$
meets only finitely many vertices and edges, together with
a map $w$ from the
edges of $\GG$ to the positive reals, which is invariant under
$T_1$ and $T_2$.

A {\bf $\GG$-domain} $D \subset \mathbb R^2$ is a domain whose boundary is a simple closed curve
comprised of edges and vertices in $\GG$.  Let $V=V_{\closure D}$ be the set of
$\GG$-vertices in the closure of $D$, let $G=G_D=(V,E)$ be the induced
subgraph of $\GG$ with vertex set $V_{\closure D}$, and write $V_\p =
\partial D \cap V_{\closure D}$.  Given a boundary value function $h_\p:V_\p\to\R$, the edge weighted DGFF on $G$
has a density with respect to Lebesgue measure on $\R^{V\setminus V_\p}$ which is proportional to
\begin{equation*}
\exp\Bigl(\sum_{\{u,v\}\in E} -\frac 12\,w\bl(\{u,v\}\br)\bl(h(v)-h(u)\br)^2\Bigr).
\end{equation*}
If every face of $\GG$ has three edges, then every vertex in the dual graph is an endpoint
of exactly three edges, and the boundary between positive and negative faces can be defined
as a simple path in this dual lattice, similar to the one shown in Figures \ref{fullDGFF}.a.
If not every face of $\GG$ has three edges, then we may ``triangulate'' $\GG$ by
adding additional edges to $\GG$, while maintaining the
invariance under $T_1$ and $T_2$,
to make this the case (and set $w$ to zero on these edges so that their presence does not affect
the law of the DGFF).

We define the weighted random walk on $\GG$ to be
the Markov chain with transition probability
$w(\{u,v\})/\sum_{v'} w(\{u,v'\})$ from $u$ to $v$,
where we take $w(\{u,v\})=0$ unless $u$ and $v$ are neighbors in $\GG$.
It is well known and easy to prove that such a walk
on the re-scaled lattice $\epsilon \GG$ converges to a linear transformation of
time-scaled two-dimensional Brownian motion
when $\epsilon$ tends to zero
(but since we could not find a reference, we very briefly explain this
in \S\ref{s.finalremarks}).
It is convenient to replace the embedding of $\GG$ into $\R^2$
described above with a linear transformation of that embedding that causes this limit to be
standard Brownian motion.

\begin{theorem} \label{t.otherlattices}
Both Theorem \ref{absmall} and Theorem \ref{preciseSLEconvergence} continue to hold if $\TG$
is replaced by a general weighted doubly periodic planar lattice
 $\GG$, as described above, provided
that $\GG$ is embedded in $\R^2$ in such a way that the weighted random walk converges to Brownian motion.
\end{theorem}

If $\GG$ is the grid $\mathbb Z^2$, then one natural way to triangulate $\GG$ is to add
all the edges of the form $\{(x,y), (x+1, y+1) \}$.  Another would be to add edges of the form $\{(x,y),
(x+1, y-1)\}$.  The above theorem implies, perhaps surprisingly, that the limiting law of
the zero height interface is the same in either case, with no need for a linear change
of coordinates.

\subsection{Outline}

Section \ref{s.preliminaries} introduces the basic notation and assumptions that
are necessary for the height gap results proved in Section \ref{tailtrivialitysection}.
Sections \ref{ss.estimates}, \ref{ss.nearindependence},
\ref{ss.narrows}, and \ref{ss.barriers} develop bounds and estimates
related to the geometry of zero height interfaces.  The random walk
$S$ comes into the picture in Section \ref{ss.enterRW}, where we
develop results about the near-independence of the triple $(h,
\gamma, S)$ on microscopic and macroscopic scales.  Section
\ref{ss.coupling} applies these results to prove uniqueness of the
limiting measure, namely Theorem \ref{t.limit}, and Section
\ref{ss.boundaryvalues} applies this to prove our
main height gap result, Theorem \ref{t.gap}.

Based on Theorem~\ref{t.gap}, the convergence in the case where
the boundary values are $\pm\lambda$ is not too hard, using the
method from~\cite{\LSWlesl}. However, to prove the convergence to
\SLEab/ in Section~\ref{s.drive}, we need to contend with a few
other issues which stem from the fact that the driving parameter of
\SLEab/ is the solution to an SDE with a drift term that blows up
to infinity on a fractal set of times.  To overcome these
difficulties, we change coordinates to a coordinate system in
which the drift terms stay bounded. In
\S\ref{s.approximatediffusions} we define and study {\bf
approximate diffusions}. These are random processes that are not
necessarily Markov, but satisfy an approximate discrete version of
an SDE. The Loewner driving term of the DGFF interface (before
going to the scaling limit) is the approximate diffusion we are
interested in. The main point is that an approximate diffusion of
an SDE is shown to be close to the corresponding true diffusion
satisfying the same SDE, under appropriate regularity conditions.
This is how the convergence of the driving term of the DGFF
interface to the driving term of the corresponding \SLEkr
\kappa;\rho_1,\rho_2/ is established. In sections \S\ref{CKC} and
\S\ref{ss.improve} more geometric convergence results are deduced
from the convergence of the Loewner driving term of the interface.

Finally, the rather brief Section~\ref{s.finalremarks} then
describes the (very minor) modifications required for the
generalization to other lattices, Theorem \ref{t.otherlattices}.

\subsection{Sequel} \label{ss.sequel}

This paper is actually the first of two papers the current authors
are writing about this subject.  In the second paper we will make sense of
the ``contour lines'' of the continuum GFF.  An instance $h$ of the
continuum GFF is a random distribution, not a random function;
however, given an instance of the GFF on a domain $D$, we can
project $h$ onto the space of functions which are piecewise linear on
a triangulation of $D$ to yield an instance of the DGFF which is, in
some sense, a piecewise linear approximation to $h$.  We can then
define the level lines of the GFF to be the limits of
the level lines of its piecewise
linear approximations (after proving that these limits exist).  We
will also characterize these random paths directly---without
reference to discrete approximations---by showing that they are the
unique path-valued functions of $h$ which satisfy a simple Markov
property.  Similar techniques allow us to describe the contour lines
of $h$ that form loops (instead of starting and ending at points on
the boundary of $D$).

The determination of the value of $\lambda$ for a given lattice
is not too hard, but fits better with the general spirit of our
next paper on the subject, in which we will prove, in particular,
$\tlambda = 3^{-1/4} \sqrt{\pi/8}$.
If the DGFF is scaled so that its fine mesh limit is the
ordinary GFF, we have $\lambda = \sqrt{\pi/8}$.

\medbreak {\noindent\bf Acknowledgments.} We wish to thank John Cardy,
Julien Dub\'edat, Bertrand Duplantier, Richard Kenyon, Jan\'{e} Kondev,
and David Wilson for inspiring and useful conversations,
and thank Alan Hammond and an
anonymous referee for numerous remarks on an earlier
version of this paper.

\section{Preliminaries} \label{s.preliminaries}

\subsection{A few general properties of the DGFF}
\label{ss.generalremarks}

In this subsection, we recall a few well-known properties of the
DGFF that are valid on any finite graph.
Let $(V,E)$ be a finite graph, and let $V_\p\subset V$ be a nonempty
set of vertices.
Let $h$ denote a sample from the DGFF with boundary values
given by some function $h_\p:V_\p\to\R$.

When $f$ is a function on $V$, the discrete gradient $\nabla f$ is the
function
on the set of ordered pairs $(v,u)$ such that $\{v,u\}\in E$ defined by
$\nabla f\bl((u,v)\br) = f(v)-f(u)$.
When
defining the norm of the gradient we sum over undirected edges,
i.e., we write
\begin{equation} \| \nabla f(v)\|^2 =
\sum_{\{u,v\}\in E} \bl(f(v)-f(u)\br)^2.
\end{equation}
Thus,  the probability density of $h$
is proportional to $e^{-\|\nabla
h(v)\|^2/2}$.  Therefore, when $h_\p=0$ on $V_\p$, $h$ is a
standard Gaussian with respect to the norm $\|\nabla h \|$ on $\Omega$.
  The (discrete) Dirichlet inner product that defines this norm can be
written
\begin{equation} \label{e.nablaprod} ( f, g )_\nabla =
\sum_{\{u,v\}\in E} \bl(f(v)-f(u)\br)\bl(g(v)-g(u)\br).
\end{equation}

Now write $\Delta f(v)=\sum_{u\sim v} \bigl(f(u)-f(v)\bigr)$, where the sum
is over all neighbors $u$ of $v$.
By expanding and rearranging the summands in~\eref{e.nablaprod},
we find
\begin{equation}\label{e.nablaDelta}
(f,g)_\nabla = -(\Delta f,g)\,.
\end{equation}

\medskip

Let $V_0\subset \VD$.
We claim that the vector space of functions $f:V\to \R$
that are zero on $V\setminus V_0$, and the vector space
of functions $f:V\to\R$ that are discrete-harmonic on $V_0$ (i.e.,
$\Delta f=0$ on $V_0$)
are orthogonal to each other
with respect to the inner product $({\cdot},{\cdot})_\nabla$,
and together they span $\R^V$.
This basic observation will be used frequently below.
Indeed, that they are orthogonal follows immediately
from~\eref{e.nablaDelta},
and a dimension count now shows that the two
spaces together span $\R^V$.

The following consequence of this orthogonality property
will be used below. Let $V_0\subset V$
satisfy $V_0\supset V_\p$ and let
$h_0$ denote the function that is discrete harmonic in
$V\setminus V_0$ and equal to $h$ in $V_0$.
Then $h-h_0$ and $h_0$ are independent random variables,
because $h=h_0+(h-h_0)$ is the corresponding
orthogonal decomposition of $h$.
It also follows that
\begin{equation}\label{e.bdmark}
h-h_0 \text{ is the DGFF in }
V\setminus V_0 \text{ with zero boundary values on }
V_0\,.
\end{equation}
Observe that the Markov property and the effect of boundary conditions
that were mentioned in \S\ref{heightgapsection} 
are immediate consequences of~\eref{e.bdmark}.
An additional useful consequence is that the law of $h_0$ is
proportional to $e^{- (h_0,h_0)_\nabla^2/2}$
times Lebesgue measure on $\R^{V_0}$.

\medskip

We now derive a useful well-known
expression for the expectation of $h(v)\,h(u)$:
\begin{equation} \label{e.greensrandomwalk}
\Eb{h(v)\,h(u)}-\Eb{h(v)}\,\Eb{h(u)}=G(u,v)/\operatorname{deg}(v)\,,
\end{equation}
where $G(u,v)$ is the expected number of
visits to $v$ by a simple random walk started at $u$ before it hits
$V_\p$ and $\operatorname{deg}(v)$ is the degree of $v$; that is, the number
of edges incident with it.
(The function $G$ is known as the Green function.)
As we have noted in the introduction,
$h$ is the sum of the
discrete harmonic extension of $h_\p$ and a DGFF with zero boundary values.
It therefore suffices to prove~\eref{e.greensrandomwalk}
in the case where $h_\p=0$ on $V_\p$.
In this case, $\Eb{h(v)}=0=\Eb{h(u)}$.
Setting $G_v(u)=G(u,v)$, we observe (or recall) that
$\Delta G_v(u) = -\operatorname{deg}(v)\,1_v(u)$.
Thus,
$$
h(v)=(h,1_v)= -
\frac{(h,\Delta G_v)}{\operatorname{deg}(v)}
\overset {\eref{e.nablaDelta}}=
\frac{(h,G_v)_\nabla}{\operatorname{deg}(v)}
\,.
$$
If $X$ is a standard Gaussian in $\R^n$, and
$x,y\in\R^n$, then
$\Eb{(X\cdot x)\,(X\cdot y)} = x\cdot y$.
Consequently, we have when $h_\p=0$
\begin{multline*}
\Eb{h(v)\,h(u)}=
\frac{\Eb{(h,G_v)_\nabla\,(h,G_u)_\nabla}}
{\operatorname{deg}(v) \operatorname{deg}(u)}
=
\frac{(G_v,G_u)_\nabla}
{\operatorname{deg}(v) \operatorname{deg}(u)}
\\
=
\frac{-(G_v,\Delta G_u)}
{\operatorname{deg}(v) \operatorname{deg}(u)}
=
\frac{(G_v,1_u)}
{\operatorname{deg}(v)}
=
\frac{G(u,v)}
{\operatorname{deg}(v)}\,.
\end{multline*}
This proves~\eref{e.greensrandomwalk}.

\subsection{Some assumptions and notations}

We will make frequent use of the following notations and
assumptions:
\begin{enumerate}
\hitem{i.h}{(h)}
  A bounded domain (nonempty, open, connected set) $D\subset\R^2$
  whose boundary $\p D$ is a subgraph of $\TG$ is fixed.
  The set of vertices of $\TG$ in $D$ is denoted by $\VD$ and $V_\p$ denotes
  the set of vertices in $\p D$.
  A constant $\upperhco>0$ is fixed, as well as a function $h_\p:V_\p\to\R$
  satisfying $\|h_\p\|_\infty\le\upperhco$.
  The DGFF on $D$ with boundary values given by $h_\p$ is denoted by $h$.
  Also set $V=V_{\closure D}=\VD\cup V_\p$.
\end{enumerate}

We denote by $\TG^*$ the hexagonal grid which is dual to the
triangular lattice $\TG$---so that each hexagonal face of $\TG^*$ is
centered at a vertex of $\TG$.
Generally, a $\TG^*$-hexagon will mean
a closed hexagonal face of $\TG^*$.
  Denote by $\ball_R$ the union of all
$\TG^*$-hexagons that intersect the
ball $B(0,R)$.

Sometimes, in addition to~\iref{i.h} we will need to assume:
\begin{enumerate}
\hitem{i.D}{(D)} The domain $D$ is simply connected, and (to avoid minor but annoying trivialities)
$\p D$ is a simple closed curve.
We fix two distinct midpoints of $\TG$-edges $x_\p$ and $y_\p$ on
$\p D$. Let the counterclockwise [respectively, clockwise] arc of
$\p D$ from $x_\p$ to $y_\p$ be denoted by $\p_+$ [respectively,
$\p_-$].
\end{enumerate}

If $H\subset D$ is a $\TG^*$-hexagon, we write $h(H)$ as
a shorthand for the value of $h$ on the center of $H$ (which is a vertex
of $\TG$).
Assuming~\iref{i.h} and~\iref{i.D},
let $\DDp$ denote the union of all $\TG^*$-hexagons contained in $D$ where $h$ is positive
together with the intersection of $\closure D$ with $\TG^*$-hexagons centered at vertices in $\p_+$.
Let $\DDm$ be the closure of $\closure D\setminus\DDp$ (which a.s.\ consists of
$\TG^*$-hexagons in $D$ where $h<0$ and
the intersection of $\closure D$ with $\TG^*$-hexagons whose center is in $\p_-$).
Then $\p \DDm\cap\p\DDp$ necessarily consists of the interface we previously called $\gamma$,
and a collection of disjoint simple closed paths.
We use the term {\bf interface} (or {\bf zero height interface}) to
describe a simple (or simple closed) path in
$\p\DDm\cap\p\DDp$ oriented so that $\DDp$ is on its right.
(That is, oriented clockwise around $\DDp$.)

Throughout, the notation $O(s)$ represents any quantity $f$ such
that $|f|\le C\,s$ for some absolute constant $C$. We use the notation
$\OC(s)$ if the constant also depends on $\upperhco$.  When
introducing a constant $c$, we often write $c=c(a,b)$ as shorthand
to indicate that $c$ may depend on $a$ and $b$.

\subsection{Simple random walk background}

We need to recall a very useful property of the harmonic measure of simple random walk.

\begin{lemma}[Hit near]\label{l.hitnear}
Let $v$ be a vertex of the grid $\TG$, and let $H$ be a connected
subgraph of $\TG$. Set $d=\diam(H)$.  The probability that a simple
random walk on $\TG$ started from $v$ exits the ball $B(v,d)$ before
hitting $H$ is at most $c\,(\dist(v,H)/d)^\expo1$, where $c$ and
$\expo1\in(0,1)$ are absolute constants.

Likewise, the same bound applies to the probability that a simple
random walk started at some vertex outside $B(v,d)$ will hit
$B\bl(v,\dist(v,H)\br)$ before $H$.
\end{lemma}

In fact, we may take $\expo1=1/2$. The continuous version of this
statement is known as the Beurling projection theorem (the extremal
case is when $H$ is a line segment). The above statement can
probably be deduced from the discrete Beurling theorem as given
in~\cite[2.5.2]{\Lbook}, though the setting there is slightly
different. In any case, since we do not require any particular value
for $\expo1$, the lemma is rather easily proved directly
(see~\cite[Lemma 2.1]{\SchSLE}).

\medskip
We will also use the (well-known) discrete Harnack principle, in the following
form.

\begin{lemma}[Harnack principle]\label{l.harnack}
Let $D,V,V_\p$ and $\VD$ be as in~\iref{i.h},
let $v,u\in\VD$, and
let $f: V\to\R$ be discrete-harmonic in $\VD$.
\begin{enumerate}
\item
If $v$ and $u$ are neighbors, then
$$
|f(v)-f(u)|\le O(1)\,\|f\|_\infty/\dist(v,\p D)\,.
$$
\item If we assume that $f\ge 0$ and that there is a
 path of length $\ell$ from $v$ to
 $u$ whose minimal distance to $\p D$ is $\rho$,  then
$$
f(u)\ge f(v) \,\exp\bl(-O(\ell+\rho)/\rho\br).
$$
\end{enumerate}
\end{lemma}

The proof of Statement 1 can be obtained by noting that
$f(v)$ is the expected value of $f$ evaluated
at the first hitting point on $V_\p$ of a random
walk started at $v$, and observing that a random
walk started at $v$ and a random walk started at $u$
may be coupled so that they meet before reaching a
distance of $r:=\dist(v,\p D)$ with probability $1-O(1/r)$
and walk together after they meet.
(See also the more general~\cite[Lemma 6.2]{\LSWlesl}, for example.)

To prove Statement 2, let $W:=\{w\in V:f(w)\ge f(v)\}$,
and observe that the maximum principle implies
that $W$ contains a path from $v$ to $\p D$.
If we assume that $\ell\le \rho/2$, then a random walk
started at $u$ has probability bounded away from zero to
hit $W$ before $\p D$, which by the Optional Sampling Theorem
implies 2 in this case. The case
$\ell>\rho/2$ now follows by induction on
$\lceil 2\,\ell/\rho\rceil$.

We also need the following well-known estimate on the Green function $G(u,v)$.

\begin{lemma}\label{l.green}
Let $D,\VD$ and $V$ be as in~\iref{i.h},
let $u,v\in \VD$, and let $\eps>0$.
Set $r:=\dist(u,\p D)$, and
suppose that within distance $r/\eps$
from $u$ there is a connected component of $\p D$
of diameter at least $\eps\, r$.
If $|u-v|\le 2\,r/3$, say, then
$G_D(v,u)$ (the expected number of visits
to $u$ by a simple random walk starting at $v$ before hitting $\p D$)
satisfies
$$
 G_D(v,u) = \exp\bl(O_\eps(1)\br) \log \frac {r+1}{|u-v| +1}\,.
$$
\end{lemma}

The probability $H_D(v,u)$ that random walk started from $v$ hits $u$ before $\p D$
can be expressed as $G_D(v,u)/G_D(u,u)$ and hence by the lemma
\begin{equation}\label{e.hitprob}
H_D(v,u) = \exp\bl(O_\eps(1)\br)\,\Bigl(1-\frac{\log(|u-v|+1)}{\log(r+1)}\Bigr)\,.
\end{equation}

Let $p_R$ denote the probability that simple random
walk started at $0$ does not return to  $0$ before
exiting the ball $B(0,R)$.
We now show that the lemma follows from the
well-known estimate
\begin{equation}\label{e.escape}
p_R = \exp(O(1))/\log (R+2)\,.
\end{equation}
In the setting of the lemma, consider some vertex $w\in \VD$
such that $|w-u|>r$.
It is easy to see that with probability bounded away from zero
(by a function of $\eps$), a simple random walk started from $w$
will hit $\p D$ before $u$.
Hence, $q:= \min\{1- H_D(w,u):w\in\VD\setminus B(u,r)\}$
is bounded away from $0$ by a positive function of $\eps$.
Clearly, 
$$
(q\,p_r)^{-1} \ge G_D(u,u)\ge (p_r)^{-1}\,.
$$
This proves the case $u=v$ of the lemma from~\eref{e.escape}.
Now start a simple random walk at a vertex
$v\ne u$ satisfying $|v-u|\le 2\,r/3$, and let
the walk stop when it hits $\p D$.
It is easy to see that the probability that the
walk visits $v$ after exiting the ball $B(v,|v-u|/4)$, say,
is within a bounded factor from $H_D(v,u)$.
Therefore, the expected number of visits to $v$ after
exiting $B(v,|v-u|/4)$ is withing a bounded factor
of $G_D(v,u)$. But the former is the same as
$G_D(v,v)- G_{B(v,|v-u|/4)}(v,v)$.
The lemma follows.

We have not found a reference proving~\eref{e.escape} 
in a way that generalizes to the setting of Theorem~\ref{t.otherlattices},
though the result is well-known. In fact, it easily follows from
Rayleigh's method, as explained in~\cite[\S2.2]{\DoyleSnell}.
In that book, the goal is to show that $p_R\to 0$
as $R\to\infty$ for lattices in the plane but not in $\R^3$.
However, the method easily yields the quantitative bounds~\eref{e.escape}.

\section{The height gap in the discrete setting} \label{tailtrivialitysection}
\subsection{A priori estimates} \label{ss.estimates}

This subsection contains some technical (and uninspiring) estimates
that are necessary to carry out the technical (but hopefully
interesting) coupling argument of the later parts of the section.

Suppose that $\beta$ is some path in the hexagonal grid $\TG^*$,
which has a positive probability to be a subset of a zero height
interface $h$. We will need to understand rather well the behavior
of $h$ conditioned on $\beta$ being a subset of a contour line. This
conditioning amounts to conditioning $h$ to be positive on vertices
adjacent to $\beta$ on one side and negative on vertices adjacent to
$\beta$ on the other side. (Here and in the following, a vertex $v$
of $\TG$ is {\bf adjacent} to $\beta$ if $\beta$ intersects the
interior of one of the six boundary edges of the
$\TG^*$-hexagon centered at $v$.) Thus, the following lemma will be rather
useful.

\begin{lemma}[Expectation bounds]\label{l.ibd}
There is a finite $c=c(\upperhco)>0$ such that the following holds.
Assume~\iref{i.h}.
Let $V_+$ and $V_-$ be nonempty disjoint subsets of $\VD$,
and set $U:=V_\p\cup V_+\cup V_-$.
Suppose that every vertex in $V_+$ has a neighbor in $V_-\cup V_\p$
and every vertex in $V_-$ has a neighbor in $V_+\cup V_\p$.
Let $\ev K$ be the event that $h>0$ on $V_+$ and $h<0$ on $V_-$.
Then for every $v\in V_+\cup V_-$,
\begin{equation}\label{e.hbexp}
\Eb{e^{ |h(v)|} \md\ev K}<c\,,
\end{equation}
and
\begin{equation}\label{e.atv}
c^{-1}< \Eb{|h(v)|\md \ev K}<c\,.
\end{equation}
Moreover
\begin{equation}\label{e.btv}
\Eb{|h(v)|^{-1/2}\md \ev K}<c\,.
\end{equation}
Let $B\subset D$ be a disk whose radius $r$ is smaller than its distance
to $U$. Assume that $B\cap \VD\ne\emptyset$. Then
\begin{equation}\label{e.maxbh}
\EB{\max\{|\bar h(v)| : v\in \VD\cap B\} \md \ev K} < c\,,
\end{equation}
where $\bar h$ denotes the discrete harmonic extension of the restriction of $h$ to $U$
(which is also the conditional expectation of $h$ given its restriction to $U$).
Moreover, if $\eps>0$ and $U$ has a connected component whose distance from $B$ is $R$ and
whose diameter is at least $\eps\,R$, then
\begin{equation}\label{e.av}
\EB{\Bigl(|\VD\cap B|^{-1}\sum_{v\in \VD\cap B} h(v)\Bigr)^2\md\ev K}
\le c+c'\,\log \frac{R}{r}\,,
\end{equation}
where $c'=c'(\eps)$.
\end{lemma}

\proof
The proof of~\eref{e.hbexp} is a maximum principle type argument.
Suppose that $v_1$ maximizes
$\Eb{e^{|h(v)|}\md \ev K}$ among $v\in U$, and let $M:= \Eb{e^{|h(v_1)|}\md \ev K}$.
Clearly, $M<\infty$.
Assume first that $v_1\in V_+$.
Set $\tilde V:= U\setminus\{v_1\}$, and for $v\in\tilde V$ let $p_v$ be the
probability that simple random walk starting at $v_1$ first hits $\tilde V$ in $v$. We may write
$$
h(v_1)= X + Y\,,
$$
where
\begin{equation}\label{e.Yexpansion}
Y:=\sum_{v\in \tilde V} p_v\,h(v)
\end{equation}
and $X$ is a centered Gaussian independent from $(h(v):v\in\tilde
V)$.
Moreover, by~\eref{e.greensrandomwalk}, $\Eb{X^2}$ is $1/6$
times the expected number of
visits to $v_1$ by a simple random walk starting from $v_1$ until it
first hits $\tilde V$.
Set $a:=\Es{X^2}$. Since $v_1$ has a neighbor in $\tilde V$, it
follows that $a=O(1)$. Also, clearly, $a$ is bounded away from zero,
since this random walk has at least one visit to $v_1$.

For every $y\in\R$ we have
\begin{multline*}
\Eb{e^{h(v_1)} \md Y=y,\ev K}
=
\Eb{e^{X+y}\md Y=y,\ev K}
\\
= e^y\, \Eb{e^X\md X+y>0}
=
e^y\,\frac
{\int_{-y}^\infty\exp(x-x^2/(2a)) \,dx}
{\int_{-y}^\infty \exp({-x^2/(2\,a)})\,dx}
\,.
\end{multline*}
If $y>0$, the right hand side is bounded by $O(e^y)$.
If $y\le 0$, then both integrals on the right are comparable
to their value when the upper integration limit is reduced to $-y+1$.
But then the ratio between the corresponding integrands is bounded by
$e^{1-y}$, which implies the same for the ratio of the integrals.
Consequently, the right hand side is $O(1)$ when $y\le 0$.
We take expectations conditioned on $\ev K$ and get
\begin{equation}\label{e.Mbd}
M = O(1)\,\Eb{e^Y\md\ev K}+O(1)\,.
\end{equation}
Repeated use of H\"older's inequality shows that for
non-negative random variables $x_1,\dots,x_n$
and non-negative constants $p_1,\dots,p_n$ such that
$\sum_1^np_j\le 1$ we have $\Eb{\prod_1^n x_j^{p_j}}\le
\prod_1^n \Es{x_j}^{p_j}$.
Thus,~\eref{e.Yexpansion} and~\eref{e.Mbd} give
$$
M\le O(1)\,\prod_{v\in \tilde V} \Eb{e^{h(v)}\md \ev K}^{p_v}+O(1)\,.
$$
Clearly,
$$
\Eb{e^{h(v)}\md \ev K}\le \begin{cases}
1 & v\in V_-\,,\\
M & v\in V_+\,,\\
e^{\upperhco} & v \in V_\p\,.
\end{cases}
$$
Setting
$p_+:=\sum_{v\in\tilde V\cap V_+} p_v$, we therefore obtain
\begin{equation}\label{e.Mbd2}
M \le O(1) \, M^{p_+}\,e^{(1-p_+)\upperhco}+O(1)\,.
\end{equation}
Since $v_1$ has a neighbor in $U\setminus V_+$, we have $p_+\le 5/6$
and $M\le O(1)\,e^{\upperhco}$ follows.
A symmetric argument applies if $v_1\in V_-$. If $v_1\in V_\p$,
then obviously $M\le e^{\upperhco}$.
Thus~\eref{e.hbexp} holds with $c=O(e^{\upperhco})$.
The right hand inequality in~\eref{e.atv} is an immediate
consequence of~\eref{e.hbexp} (possibly with a different $c$).

We now prove~\eref{e.maxbh}.
Given $v\in \VD$ and $u\in U$ let $H(v,u)$ denote the
probability that simple random walk started at $v$ first hits
$U$ at $u$.  Suppose $v,v'\in B\cap\VD$ and $u\in U$.
The discrete Harnack principle (Lemma~\ref{l.harnack})
then gives
$H(v',u)\le H(v,u) \,O(1)$.
Thus,
$$
|\bar h(v')|=\Bigl|\sum_{u\in U} H(v',u)\, h(u)\Bigr|
\le \sum_{u\in U} O(1)\,H(v,u)\, |h(u)|\,.
$$
The right hand side therefore bounds
$\max\{ |\bar h(v')|:v'\in \VD\cap B\}$.
Now,~\eref{e.maxbh} follows from the right hand inequality of~\eref{e.atv}
and $\sum_{u\in U}H(v,u)=1$.

We now prove \eref{e.btv}.  Consider some $v\in V_+$.  We first
claim that if $w$ neighbors with $v$ and $w\notin U$, then
$\Eb{h(w)^2\md\ev K}=\OC(1)$.  To see this, first note that
$\overline h(w)$ (where as above $\overline h$ is the discrete
harmonic extension of the values of $h$ on $U$) is a linear
combination of the values $\overline h(v)$ for $v \in U$
(which are the same as the values of $h$ there); since each
of these values has mean and variance which are $\OC(1)$ (by
\eref{e.hbexp}), the mean and variance of $\overline h(w)$ are also
$\OC(1)$.
By~\eref{e.greensrandomwalk}, conditioned on $\overline h$, the value $h(w)$ is
Gaussian with variance $1/6$ times the expected number of
times a random walk started at $w$ visits $w$ before hitting $U$,
which is $O(1)$ because $U$ contains a neighbor of $w$.

Let $Z$ denote the average of $h$ on the neighbors of $v$, and let
$Z':=h(v)-Z$. The above implies that $\Eb{Z^2\md \ev K}=\OC(1)$.
Since $Z'$ is a centered gaussian with variance $1/6$,
for every $z\in\R$,
\begin{multline*}
\Eb{h(v)^{-1/2}\md \ev K,\,Z=z}=
\Eb{(Z'+z)^{-1/2}\md Z'+z>0}
\\
= \frac{\int_{0}^\infty x^{-1/2}\,e^{-6(x-z)^2/2}\,dx}
{\int_{0}^\infty e^{-6(x-z)^2/2}\,dx}\,.
\end{multline*}
If $z\ge -2$, then this is clearly bounded. Assume therefore that
$z< -2$. It is easy to verify that the integrals in the numerator
and denominator are comparable to the same integrals restricted to
the range $x\in\bl[0,|z|^{-1}\br]$. But in this range, the maximum
value of $\exp(-6(x-z)^2/2)$ is comparable to the minimum value of
the same quantity. Consequently, when $z< -2$,
\begin{multline*}
\Eb{h(v)^{-1/2}\md \ev K,\,Z=z}
\\
= O(1)\, \frac{\int_{0}^{|z|^{-1}} x^{-1/2}\,e^{-6(x-z)^2/2}\,dx}
{\int_{0}^{|z|^{-1}} e^{-6(x-z)^2/2}\,dx} = O(1)\,
\frac{\int_{0}^{|z|^{-1}} x^{-1/2}\,dx} {\int_{0}^{|z|^{-1}}\,dx}
=O(1)\,|z|^{1/2},
\end{multline*}
which gives
$$ \Eb{h(v)^{-1/2}\md \ev K}
= O(1)+O(1)\,\Eb{|Z|^{1/2}\md\ev K}.
$$
Since $\Eb{Z^2\md\ev K}=\OC(1)$, we certainly have
$\Eb{|Z|^{1/2}\md\ev K}=\OC(1)$.
This proves~\eref{e.btv}. Now the left hand inequality
in~\eref{e.atv} is an immediate consequence.

We now prove \eref{e.av}.  Consider two vertices $v,u\in B\cap \VD$.
Assume that within distance $R$ from $B$ there is a connected
component of $U$ of diameter at least $\eps\,R$. Since
$\Eb{h(v)\md\ev K,\bar h}=\bar h(v)$, we have
\begin{equation}\label{e.ehh}
\Eb{h(v)h(u)\md\ev K}
= \Eb{(h(v)-\bar h(v))(h(u)-\bar h(u))\md\ev K}
+
\Eb{\bar h(v)\bar h(u)\md\ev K}
\,.
\end{equation}
Now, $\bar h(v)$ is just a weighted average of $h(w)$ with $w\in U$
(according to harmonic measure from $v$).
Consequently,
$$
\Eb{\bar h(v)^2\md\ev K}\le \max_{w\in\U} \Eb{h(w)^2\md \ev K} \le
\max_{w\in U} 2\,\Eb{e^{|h(w)|}\md\ev K}\le 2\,c
$$
and Cauchy-Schwarz implies that the last summand in~\eref{e.ehh} is also bounded by $2\,c$.
Since $h-\bar h$ is the Gaussian free field with zero boundary values on $U$
and is independent from the restriction of $h$ to $U$,
by~\eref{e.greensrandomwalk}
$$
\Eb{(h(v)-\bar h(v))(h(u)-\bar h(u))\md\ev K}
$$
is $1/6$ times the expected number of visits to $u$ by a random walker
started at $v$ which stops when it hits $U$.
Thus
$$
\sum_{u\in B\cap \VD}
\Eb{(h(v)-\bar h(v))(h(u)-\bar h(u))\md\ev K}
$$
is $1/6$ times the expected number of steps that the walker spends
in $B$, which is
$$
O_\eps(1)\,\left|B\cap \VD\right|\,\log(R/r)\,,
$$
by Lemma~\ref{l.green}.
The estimate~\eref{e.av} is now obtained by averaging~\eref{e.ehh} over all $v,u\in B\cap \VD$.
\QED

The following lemma provides a variant of the right hand bound
in~\eref{e.atv} in the case where instead of looking for a zero
height interface of $h$, we consider instead a zero height interface
of $h-g$, for some fixed function $g$.

\begin{lemma}[Further expectation bounds]\label{l.ibd2}
In the setting of Lemma~\ref{l.ibd},
suppose that $g: \closure D\to\R$ is zero on $V_\p$
and Lipschitz in $\closure D$.
Let $\ev K_g$ be the event that $h>g$ on $V_+$ and $h<g$ on $V_-$.
Then for every $v\in V_+\cup V_-$,
\begin{equation}\label{e.atv2}
\Eb{|h(v)-g(v)|\md \ev K_g}\le
c\,(1+\upperhco)+c\,\frac{\|g\|_\infty}{\log (q+2)}\,,
\end{equation}
where $c>0$ is a universal constant
and $q=\|g\|_\infty/\|\nabla g\|_\infty$.
\end{lemma}

\proof
The proof is similar to the proof of Lemma~\ref{l.ibd}.
Suppose that $v_1$ maximizes
$\E[h(v)-g(v)\md \ev K_g]$ among $v\in V_+$,
and let $M:= \E[h(v_1)-g(v_1)\md \ev K_g]$. By symmetry, it is enough
to get a bound on $M$.
We define $\tilde V, X, Y$ and $p_v$ as in the proof of Lemma~\ref{l.ibd}.

For every $y\in\R$ we have
\begin{multline*}
\Eb{h(v_1)-g(v_1)\md Y=y,\ev K_g} = \Eb{X+y-g(v_1)\md X+y-g(v_1)>0}
\\
\le O(1)+(y-g(v_1))_+ \,,
\end{multline*}
where $(x)_+:=\max\{0,x\}$.
Consequently,
\begin{align*}
M
&
\le O(1)+ \Eb{(Y-g(v_1))_+\md\ev K_g}
\\&
\le
O(1)+ \sum_{v\in\tilde V} p_v\, \Eb{(h(v)-g(v_1))_+\md\ev K_g}
\\&
\le
O(1)+ \sum_{v\in\tilde V} p_v\, \Eb{(h(v)-g(v))_+\md\ev K_g}
+ \sum_{v\in\tilde V} p_v \bigl|g(v)-g(v_1)\bigr|
\\&
\le
O(1)+\upperhco+ M\, \sum_{v\in V_+} p_v
+ \sum_{v\in\tilde V} p_v \bigl|g(v)-g(v_1)\bigr|
\,.
\end{align*}
Consider a simple random walk started at $v_1$ and stopped when it first
hits $\tilde V$. Denote the vertex where it first hits $\tilde V$ by
$w$. (Then, $\P[v=w]=p_v$.)
Set $r:=q/\log (q+2)$.
Since $v_1$ has a neighbor in $\tilde V$, we have
$\Pb{|v-v_1|>r}\le O(1)/\log (r+2)$, by standard random walk estimates.
When $|v-v_1|\le r$, we have $|g(v)-g(v_1)|\le O(1)\, r\,\|\nabla g\|_\infty$.
Therefore,
$$
\sum_{v\in\tilde V} p_v \bigl|g(v)-g(v_1)\bigr|
\le
O(1)\, r\,\|\nabla g\|_\infty+O(1)\,\frac{\|g\|_\infty}{\log (r+2)}
\le
O(1)\,\frac{\|g\|_\infty}{\log (q+2)}\,.
$$
This gives
$$
M\,\sum_{v\in\tilde V\setminus V_+} p_v \le
O(1)+\upperhco + O(1)\,\frac{\|g\|_\infty}{\log (q+2)}\,.
$$
Because
$\sum_{v\in\tilde V\setminus V_+} p_v$ is bounded away from zero, the proof
is now complete.
\QED

Our next result establishes the continuity of the conditional distribution
of $h$ in the specified data. More precisely,

\begin{proposition}[Heights interface continuity]\label{p.hcont}
For every $\eps>0$ there is some $R=R(\eps,\upperhco)>1/\eps$ such
that the following holds. Let $D$, $\VD$, $ V_\p$, $h_\p$, $V_+$,
$V_-$, $U$ and $\ev K$ be as in Lemma~\ref{l.ibd}, and let $\hat D$,
$\hatVD$, $\hat V_\p$, $\hat h_\p$, $\hat V_+$, $\hat V_-$, $\hat U$
and $\hat{ \ev K}$ be another such system, which is also assumed to
satisfy $\|\hat h_\p\|_\infty\le\upperhco$. Let $h_{\ev K}$ be a
DGFF in $D$ with boundary values given by $h_\p$ conditioned on $\ev
K$, and let $h_{\hat {\ev K}}$ be a DGFF in $\hat D$ with boundary
values given by $\hat h_\p$ conditioned on $\hat{\ev K}$. Suppose
that within $\ball_R$, the two systems are the same; that is, $D\cap
\ball_R=\hat D\cap \ball_R$, $h_\p|_{\ball_R}=\hat h_\p|_{\ball_R}$,
$V_+\cap\ball_R=\hat V_+\cap\ball_R$ and $V_-\cap\ball_R=\hat
V_-\cap\ball_R$. Further suppose that $0\in U$. Then there is a
coupling of $h_{\ev K}$ and of $h_{\hat{\ev K}}$ such that for every
vertex $v\in\bal {1/\eps}$ we have $\Eb{|h_{{\ev K}}(v)-h_{\hat{\ev
K} }(v)|}<\eps$.
\end{proposition}

The following lemma will be needed in the proof.

\begin{lemma}\label{l.hcont}
Let $X$ be a one dimensional Gaussian of zero mean and unit
variance. Let $x,\hat x\in \R$, let $Z$ be a random variable whose
distribution is the same as that of $X+x$ conditioned on $X+x>0$,
and let $\hat Z$ be a random variable whose distribution is the same
as that of $X+\hat x$ conditioned on $X+\hat x>0$. Then there is a
coupling of $Z$ and $\hat Z$ such that $|Z-\hat Z|<|x-\hat x|$
almost surely if $x\ne\hat x$. Moreover, there is a continuous
function $\delta(x,\hat x)$ satisfying $\delta(x,\hat x)<1$ such
that $\Eb{|Z-\hat Z|}\le \delta(x,\hat x)\,|x-\hat x|$ under this
coupling.
\end{lemma}

The coupling that we use is what is known as the quantile coupling of $Z$ and $\hat Z$.

\proof
Let $F(s)=\Ps{X<s}$, and let $G=F^{-1}$.
Set $t:=F(-x),\,\hat t:=F(-\hat x)$.
Let $p$ be a random variable uniformly distributed in $[0,1]$.
Then $Z(t):=x+G\bl(t+p\,(1-t)\br)=G(t+p-t\,p)-G(t)$ has
the same distribution as $Z$.
Therefore, $\bl(Z(t),Z(\hat t)\br)$ is a coupling of
$Z$ and $\hat Z$.  Consequently, to verify the first claim it is
sufficient to show that $|\p_t Z(t)|< \p_t G(t)$.
In fact, we will prove the stronger statement,
$ -\p_t G(t)<\p_t Z(t) < 0$ for all $p\in(0,1)$,
which is equivalent to
$$
0< \p_tG(t+p-t\,p)< \p_t G(t)\,.
$$
The left hand inequality is immediate, because $G'>0$ on $(0,1)$.
The right hand inequality translates to
$(1-p)\, G'(t+p-t\,p)< G'(t)$, which we rewrite as
$\bl(1-(t+p-t\,p)\br) \,G'(t+p-t\,p)<(1-t)\,G'(t)$.
This is equivalent to $(1-t_p)/F'(G(t_p))<(1-t)/F'(G(t))$, where $t_p:=t+p-t\,p>t$.
Now, note that
$$
\frac{1-t}{F'(G(t))}
=\frac
{\int_{-x}^\infty \exp(-s^2/2)\,ds}
{\exp(-x^2/2)}
=
{\int_{-x}^\infty e^{(x^2-s^2)/2}\,ds}
=
{\int_{0}^\infty e^{xs-s^2/2}\,ds}
$$
is strictly decreasing in $t$, because $x$ is strictly decreasing in $t$.
This proves the first claim.
The second claim follows with
$$
\delta(x,\hat x) :=
\int_0^1
\frac{Z(t)-Z(\hat t)}
{x-\hat x}\,dp
\qquad
\text{ when } x\ne \hat x
$$
and
$$
\delta(x,x) :=
\int_0^1
\frac{Z'(t)}
{-G'(t)}\,dp \,.
$$
\QED

\proofof{Proposition~\ref{p.hcont}} Fix a coupling of $h_{\ev K} $
and $h_{\hat{\ev K}}$ that minimizes $\sum_{v\in \VD\cap\hatVD}
\Eb{|h_{\ev K}(v)-h_{\hat{\ev K}}(v)|}$. Standard continuity and
compactness arguments show that there is such a coupling. Set
$f(v):=\Eb{|h_{\ev K}(v)-h_{\hat{\ev K}}(v)|}$ for vertices $v\in
(\VD\cup V_\p)\cap(\hatVD\cup\hat V_\p)$.

First, we claim that $f$ is discrete subharmonic on vertices in
$\bal R\setminus U$. Indeed, fix a vertex $w\in \bal R\setminus U$.
The conditional distribution of $h_{\ev K}(w)$ given the value of
$h_{\ev K}$ at every vertex but $w$ is that of $x+ A\,X$, where $x$ is
the average of $h_{\ev K}$ on the neighbors of $w$, $X$ is a
standard Gaussian, and $A$ is the lattice-dependent constant
$1/\sqrt{6}$ (since each vertex has six neighbors in $TG$). Similarly,
$h_{\hat{\ev K}}(w)=\hat x+A\hat X$. By the choice of coupling, when
we fix the values of $h_{\ev K}$ and $h_{\hat{\ev K}}$ off of $w$,
the corresponding conditioned coupling of $h_{\ev K}(w)$ and
$h_{\hat{\ev K}}(w)$ minimizes the conditioned expectation of
$|h_{\ev K}(w)-h_{\hat{\ev K}}(w)|$. But one such conditioned
coupling is obtained by taking $X=\hat X$. Thus, $f(w)\le
\Eb{|x-\hat x|}$, which implies that $f$ is subharmonic at $w$,
since
\begin{equation}\label{e.ftb}
\sum_{u\sim w} |h_{\ev K}(u)-h_{\hat{\ev K}}(u)|
\ge
\Bl|\sum_{u\sim w} h_{\ev K}(u)-\sum_{u\sim w} h_{\hat{\ev K}}(u)\Br|,
\end{equation}
where the sums are over the neighbors of $w$.

Next, consider any vertex $v\in V_+\cap \bal R$. As before, we may
write $h_{\ev K}(v)=x+A\,X$, where $x$ is the average of $h_{\ev K}$
on the neighbors of $v$ and $X$ is a random variable whose conditional
law given $x$ is that of a standard Gaussian 
conditional on $x+A\,X>0$.  Lemma~\ref{l.hcont} applied with $x/A$ instead of $x$ and
$\hat x/A$ instead of $\hat x$ implies that $f$ is also subharmonic
at $v$. We claim that there is a constant $b=b(\upperhco,\eps)>0$
such that
\begin{equation}\label{e.strictsub}
\Delta f(v) \ge b\,,\qquad\text{ if }f(v)\ge \eps/2\,,
\end{equation}
where
$\Delta$ denotes the discrete Laplacian on $\TG$.
Indeed, the optimality of the coupling gives
$$
\Eb{|h_{\ev K}(v)-h_{\hat{\ev K}}(v)|\md x,\hat x}\le
\delta(x/A,\hat x/A)\,|x-\hat x|\,,
$$
where $\delta(\cdot,\cdot)<1$ is as in the lemma.
Thus,
$$
f(v) \le \Eb{ \delta(x/A,\hat x/A)\,|x-\hat x|}.
$$
By~\eref{e.ftb} with $v$ in place of $w$ we have
\begin{equation}\label{e.Delta}
\Delta f(v) \ge \Eb{|x-\hat x|}-f(v) \ge \EB{\bl(1- \delta(x/A,\hat
x/A)\br)\,|x-\hat x|}.
\end{equation}
For every neighbor $u$ of $v$ we have by~\eref{e.av} $\Eb{h_{\ev
K}(u)^2}< \OC(1)$. It easily follows that $\Eb{x^2}<\OC(1)$. Hence,
Cauchy-Schwarz implies that there is a constant
$b_0(\eps,\upperhco)>0$ such that $\Eb{|x|\,1_{\ev A}}<\eps/8$ for
every event $\ev A$ satisfying $\Pb{\ev A}<b_0$. There is a constant
$b_1=b_1(\eps,\upperhco)$ such that $\Pb{|x|>b_1}<b_0/2$. The same
inequalities will hold with $\hat x$ and $h_{\hat{\ev K}}$ in place
of $x$ and $h_{\ev K}$. Therefore,
$$
\begin{aligned}
\Eb{|x-\hat x|}
&
=
\Eb{|x-\hat x|\,1_{|x|\vee|\hat x|\le b_1}}
+
\Eb{|x-\hat x|\,1_{|x|\vee|\hat x|> b_1}}
\\&
\le
\Eb{|x-\hat x|\,1_{|x|\vee|\hat x|\le b_1}}
+
\Eb{(|x|+|\hat x|)\,1_{|x|\vee|\hat x|> b_1}}
\\&
\le
\Eb{|x-\hat x|\,1_{|x|\vee|\hat x|\le b_1}}
+
\eps/4\,.
\end{aligned}
$$
Thus, if we assume that $f(v)\ge \eps/2$, then also $\Eb{|x-\hat x|}\ge \eps/2$
and therefore
$ \Eb{|x-\hat x|\,1_{|x|\vee|\hat x|\le b_1}}\ge \eps/4$.
Therefore,~\eref{e.Delta} gives~\eref{e.strictsub} with
$$
b:=(\eps/4)\,\min\{1-\delta(x/A,\hat x/A):x,\hat x\in [-b_1,b_1]\}.
$$
Clearly, $f$ is also discrete-subharmonic on
$V_-\cap\bal R$ and~\eref{e.strictsub} also holds for $v\in V_-\cap\bal R$ and
(trivially) for $v\in V_\p\cap\bal R$.

Next, we prove that for all vertices $w\in \bal R$
we have
\begin{equation}\label{e.upbd}
f(w)\le \OC(1)\,\sqrt{\log R}\,.
\end{equation}
Fix such a $w$, and assume that $w\notin U$. We may decompose
$h_{\ev K}(w)$ as a sum $h_{\ev K}(w)=y+Y$, where $y$ is the
value at $w$ of the
discrete harmonic extension of the restriction of $h_{\ev K}$ to
$U$, and $Y$ is a centered Gaussian whose variance is
$1/6$ times the expected
number of visits to $w$ by a simple random walk started at $w$ that
is stopped when it hits $U$. A simple random walk on $\TG$ started
at $w$ has probability at least a positive constant times $1/\log R$
to reach distance $R$ from $w$ before returning to $w$, and once it
does reach this distance it has probability bounded away from zero
to hit $0$ before returning to $w$. Since $0\in U$, it follows that
$\Eb{Y^2}=O(\log R)$. Thus, $\Eb{|Y|}=O(\sqrt{\log R})$. Since $y$
is the average of the value of $h_{\ev K}$ on $U$ with respect to
harmonic measure from $w$, it follows from~\eref{e.atv} that
$\Eb{|y|}=\OC(1)$. Thus, we have $\Eb{|h_{\ev K}(w)|}\le
\OC(\sqrt{\log R})$. This certainly also holds if $w\in U$, and a
similar estimate holds for $h_{\hat{\ev K}}(w)$. Now~\eref{e.upbd}
follows, since $f(w)\le \Eb{|h_{\ev K}(w)|}+ \Eb{|h_{\hat{\ev
K}}(w)|}$.

We now show that the established properties of $f$ imply that
$f\le\eps$ on $\bal{1/\eps}$ if $R$ is sufficiently large. Fix some
vertex $w\in\bal{1/\eps}$, and let $S_t$ be a simple random walk on
$\TG$ started at $w$.
Let $t_1$ be the first time $t$ such that $|S_t|>R/2$ or $S_t=0$.
Since $f$ is discrete-subharmonic on $\VD\cap\bal R$,
$t\mapsto f(S_{t\wedge t_1})$ is a submartingale.
The optional sampling theorem implies that $f(w)\le \Eb{f(S_{t_1})}$.
By standard random walk estimates, $\Pb{|S_{t_1}|>R/2}\le O(1)\,|\log\eps|/\log R$.
(We assume, with no loss of generality, that $\eps<1/2$, say.)
Consequently,
\begin{multline*}
f(w)\le
\Eb{f(S_{t_1})}\le f(0)+ O(1)\,\frac{|\log\eps|}{\log R} \,\max\{f(u):u\in \VD\cap\bal R\}
\\
\le
f(0)+ \OC(1)\,\frac{|\log\eps|}{\sqrt{\log R}}\,.
\end{multline*}
This proves that $f(w)\le \eps$ if $f(0)\le\eps/2$ and $R$ is sufficiently large.

Now assume that $f(0)\ge \eps/2$. Let $\tilde S_t$ be a simple random walk
starting at $0$.
Let $t_*:=\min\{t:|\tilde S_t|>R/2\}$ and
let $n_s$ be the number of $t\in\{0,\dots,s-1\}$ such that $\tilde S_t=0$.
By~\eref{e.strictsub} and our assumption that $f(0)\ge\eps/2$,
$f(\tilde S_{t\wedge t_*})-b\,n_{t\wedge t_*}$ is a submartingale.
Thus,
\begin{multline*}
0\le f(0)\le
\Eb{f(\tilde S_{t_*})}-b\,\Eb{n_{t_*}}
\le \max\{f(u):u\in\VD\cap\bal R\} - b\,\Eb{n_{t_*}}
\\
\le \OC(1)\,\sqrt{\log R} - b\,\Eb{n_{t_*}}.
\end{multline*}
Now note that as $R\to\infty$ while $\eps$ is fixed
$\Eb{n_{t_*}}$ grows at least as fast as a positive
constant times $\log R$, because the probability for $\tilde S_t$ not to return to
$0$ after any specific visit to $0$ is bounded by $O(1/\log R)$.
 Thus, the above rules out the possibility that $f(0)\ge \eps/2$ if
$R$ is sufficiently large. This completes the proof.
\QED

As a corollary of the proposition, we now show that
the correlation in the values of $h$ at two vertices
in $U$ decays as the distance between them tends to infinity.

\begin{corollary}[Correlation decay]\label{c.corr}
For every $\eps>0$ there is some $R=R(\eps,\upperhco)$ such that the following
holds.
Let $D$, $\VD$, $ V_\p$, $h_\p$, $V_+$, $V_-$, $U$ and $\ev K$ be as in Lemma~\ref{l.ibd},
and let $v_1,v_2\in U$ satisfy $|v_1-v_2|>R$.
Then
$$
\Bl|\Eb{h(v_1)\,h(v_2)\md\ev K}-
\Eb{h(v_1)\md\ev K}\,\Eb{h(v_2)\md\ev K}\Br|<\eps\,.
$$
\end{corollary}

\proof
Suppose, without loss of generality, that $v_2\in V_+$.
Fix some $a>0$ and let $X:=1_{\{0<h(v_2)\le a\}}$.
We may apply Proposition~\ref{p.hcont} to our present setup and to the setup
where the value of $h(v_2)$ is fixed at some constant $y\in(0,a]$ and $v_2\in \p D$.
Thus, the proposition would apply, provided that
$\upperhco$ is replaced by $\upperhco\vee a$.  Consequently, we find that there is an
$R'=R'(\eps,\upperhco,a)$ such that if $|v_1-v_2|>R'$, then
$$
\Bl|\Eb{h(v_1)\md h(v_2),\ev K}-
\Eb{h(v_1)\md \ev K}
\Br| X\le \eps/(2a)\,.
$$
Since $h(v_2)\,X\le a$, $X^2=X$ and $h(v_2)\,X$ is $h(v_2)$ measurable, this gives
\begin{multline*}
\eps/2\ge
\Bl|\Eb{h(v_1)\md h(v_2),\ev K}
\,
h(v_2)
\,
X
-
\Eb{h(v_1)\md \ev K}\,
h(v_2)
\,
X
\Br|
\\
=
\Bl|\Eb{h(v_1)\,h(v_2)\,X\md h(v_2),\ev K}
-
\Eb{h(v_1)\md \ev K}\,
h(v_2)
\,
X
\Br| .
\end{multline*}
Taking expectations conditioned on $\ev K$ now gives
\begin{equation}\label{e.restricted}
\Bl|
\Eb{h(v_1)\,h(v_2)\,X\md \ev K}
-
\Eb{h(v_1)\md \ev K}\,
\Eb{h(v_2) \, X\md\ev K}\Br|\le\eps/2\,.
\end{equation}
Since
$$
(1-X)\,h(v_2)^2\le (1-X)\, |h(v_2)|^3/a\le 6\, e^{|h(v_2)|}/a
$$
and $h(v_1)^2\le 2\,e^{|h(v_1)|}$,
if $c$ denotes the constant satisfying~\eref{e.hbexp}, then
Cauchy-Schwarz gives
\begin{multline*}
\Eb{h(v_1)\,h(v_2)\,(1-X)\md\ev K}^2
\le
\Eb{h(v_1)^2\md\ev K}\,
\Eb{h(v_2)^2\,(1-X)\md\ev K}
\\
\le (2\,c)\,(6\,c/a)\,.
\end{multline*}
Similarly,
$$
0\le\Eb{|h(v_2)| \, (1-X)\md\ev K}
\le \Eb{2\,e^{|h(v_2)|}/a\md\ev K}\le 2\,c/a\,.
$$
Consequently, if $a$ is chosen sufficiently large then
\begin{multline*}
\Bl|
\Eb{h(v_1)\,h(v_2)\,X\md \ev K}
-
\Eb{h(v_1)\,h(v_2)\md \ev K}\Br|
\\
=
\Bl|
\Eb{h(v_1)\,h(v_2)\,(1-X)\md \ev K}
\Br|
<\eps/4
\end{multline*}
and
$$
\Bl|
\Eb{h(v_1)\md \ev K}\,
\Eb{h(v_2) \, X\md\ev K}
-
\Eb{h(v_1)\md \ev K}\,
\Eb{h(v_2) \md\ev K}
\Br|<\eps/4\,.
$$
The corollary now follows from~\eref{e.restricted}.
\QED

Next, we provide a simple lemma which bounds the amount in which 
adding a function to $h$ affects its distribution.

\begin{lemma}[DGFF distortion]\label{l.distort}
Assume~\iref{i.h}.
Let $f:\VD\cup V_\p \to \R$ satisfy $f=0$ on $V_\p$.
Let $\mu$ be the law of $h$,
 and let $\mu_f$ be the law of
$\tilde h:=h+f$.  Then for every event $\ev A$,
$$
\mu_f[\ev A] \le \exp\bl(\|\nabla f\|^2/2\br)\, \mu[\ev A]^{1/2}.
$$
\end{lemma}

\proof Suppose that $X$ is a standard Gaussian in $\R^n$, $y\in\R^n$
is some fixed vector, and $A\subset\R^n$ is measurable. Then
\begin{multline*}
\Pb{X+y\in A}= c_n\int 1_A\, \exp(- \|x-y\|^2/2)
\\
=
c_n\int 1_A\,\exp(x\cdot y-\|y\|^2/2)\,\exp(-\|x\|^2/2)
\,,
\end{multline*}
where $ c_n^{-1} =\int \exp(-\|x\|^2/2) $ and the integrals are with
respect to Lebesgue measure in $\R^n$.
(This is the Cameron-Martin formula.)
We may think of the right
hand side as the inner product of $1_A$ and $\exp(x\cdot
y-\|y\|^2/2)$ with respect to the Gaussian measure. Hence,
Cauchy-Schwarz gives
\begin{multline*}
\Pb{X+y\in A}\le\Pb{X\in A}^{1/2}\,\Eb{\exp(2\,X\cdot y-\|y\|^2)}^{1/2}
\\
=
\Pb{X\in A}^{1/2}\,
\exp(\|y\|^2/2)\,.
\end{multline*}

Let $\bah$ denote the discrete harmonic extension of $h_\p$. Then
$h-\bah$ is the DGFF with zero boundary values,
and hence is a standard Gaussian on $\R^{\VD}$ with respect to the
norm $g\mapsto \|\nabla g\|_2$. The lemma follows. \QED

\subsection{Near independence} \label{ss.nearindependence}

In this subsection we build on the infrastructure developed above
to prove that under appropriate
assumptions the shape of an interface inside a ball does
not depend too strongly on the shape of an interface outside
a slightly larger ball. More precisely, we have:

\begin{proposition}[Near independence]\label{p.decoup}
Let $C>1$ and let $R>10^3\,C$.
  Assume~\iref{i.h} and $\ball_{5R}\subset D$.
Let $R_1,R_2,R_3\in [R,5R]$ satisfy $R_1+C^{-1}\,R<R_2$, $R_2+C^{-1}\,R<R_3$ and $R_3+C^{-1}\,R<5R$.
 Let $V_+^3$, $V_-^3$ be disjoint sets of vertices in $D\setminus \ball_{R_3}$
and let $V^1_+$, $V^1_-$ be disjoint sets of vertices in $\ball_{R_1}$.
Suppose that every vertex of $V^1_+$ neighbors with a vertex
in $V^1_-$, every vertex of $V^1_-$ neighbors with a vertex in
$V^1_+$, and similarly for $V^3_-$ and $V^3_+$.
Also suppose that a random walk started at $0$ has probability
at least $1/C$ to hit $V^3_-\cup V^3_+$ before exiting $\ball_{5R}$.
Let $\ev K_1$ be the event that $h>0$ on $V^1_+$ and
that $h<0$ on $V^1_-$, and let $\ev K_3$ be the corresponding
event for $V^3_-$ and $V^3_+$.
Let $a(V^1_+,V^1_-)$ be the probability of $\ev K_1$ for the DGFF
on $\ball_{R_2}$ with zero boundary values outside $\ball_{R_2}$.
Then there is a constant $c=c(\upperhco,C)>0$,
such that
$$
c^{-1}a(V^1_+,V^1_-) \le \Pb{\ev K_1\md \ev K_3} \le c\, a(V^1_+,V^1_-)\,.
$$
\end{proposition}

\proof
For $j=0,1,2,\dots,6,7,8$ set $r_j=R_1+j\,(R_2-R_1)/8$.
Then $r_8=R_2$ and $r_{j+1}\ge r_j+C^{-1}\,R/8$.
We will use the notations $W^j:=(\VD\cup V_\p)\cap \bal {r_j}=\VD\cap\bal {r_j}$,
$W_j:=(\VD\cup V_\p)\setminus W^j$ and $W_j^i:=W_j\cap W^i$.
Let $\th$ denote the harmonic extension of the restriction of $h$ to
$\Vpp$. We may identify $\th$ with a point in $\R^{\Vpp}$; namely,
its restriction to $\Vpp$.
As we have noted after~\eref{e.bdmark} 
 the probability density of $\th$
with respect to Lebesgue measure on $\R^{\Vpp}$ is proportional to
$\exp(-\|\nabla\th\|^2/2)$.
Hence,
\begin{equation}\label{e.fii}
\Pb{\ev K_1\md \ev K_3} = \frac {\int_{\ev K_1\cap\ev K_3}
\exp\bl(-\|\nabla\th\|^2/2\br)} {\int_{\ev K_3}{
\exp\bl(-\|\nabla\th\|^2/2\br)}} \,,
\end{equation}
where the integrals are with respect to Lebesgue measure on $\R^{\Vpp}$.
Let $h_1$ be the function that agrees with $h$ on $\Va$, is harmonic
on $\Vp$, and is zero in $\Vb$.  Let $h_3$ be the function that agrees
with $h$ on $\Vb$, is harmonic on $\Vp$, and is zero in $\Va$.
Clearly, $\th = h_1+h_3$.

We claim that
\begin{equation}\label{e.prod}
{\int_{\ev K_1\cap\ev K_3} \exp\Bl(-\frac{\|\nabla\th\|^2}{2}\Br)}
\asymp \int_{\ev K_1\cap\ev K_3} \exp\Bl(-\frac{\|\nabla
h_1\|^2}{2}\Br) \exp\Bl(-\frac{\|\nabla h_3\|^2}{2}\Br) \,,
\end{equation}
where $\asymp$ means equivalence up to multiplicative constants
depending on $\upperhco$ and $C$. Let $V_+:=V_+^1\cup V_+^3$ and
$V_-:=V_-^1\cup V_-^3$. Fix some $v\in\Vs$. Let $\bah$ denote the
harmonic extension of the restriction of $h$ to $V_+\cup V_-$. Then
$h(v)-\th(v)$ and $\th(v)-\bah(v)$ are independent Gaussian random
variables, and both are independent of $\bah$
(by the orthogonality property noted in
\S\ref{ss.generalremarks}).
By~\eref{e.greensrandomwalk},
the variance of $h(v)-\th(v)$ is $1/6$ times the
expected number of visits to $v$ by a random walk started at $v$,
which is stopped when it hits $\Vpp$, and the variance of
$h(v)-\bah(v)$ is $1/6$ times the expected number of visits to $v$
by the same random walk stopped when it hits $V_+\cup V_-$.
Consequently,
the variance of $\th(v)-\bah(v)$ is $1/6$ times the expected number of
visits to $v$ after the first hit of $\Vpp$ and before the first hit
of $V_+\cup V_-$. Note that the probability to hit $v$ by a random
walk started in $\Vpp$ before exiting $\bal {5R}$ is $O(1/\log R)$
and conditioned on hitting $v$ before exiting $\bal {5R}$ the number
of visits to $v$ prior to exiting $\bal{5R}$ is $O(\log R)$. Our
assumption on the probability to hit $V_-^3\cup V_+^3$ therefore
easily implies that $ \Eb{\bl(\th(v)-\bah(v)\br)^2} = O_{C}(1) $ and
hence $\Eb{|\th(v)-\bah(v)|} = O_C(1)$. Since $\ev K_1\cap\ev K_3$
is determined by $\bah$, it is independent of $\th(v)-\bah(v)$, and
consequently $\Eb{|\th(v)-\bah(v)|\md \ev K_1,\ev K_3} = O_C(1)$.
Now, by~\eref{e.maxbh}, we have $\Eb{|\bah(v)|\md \ev K_1,\ev
K_3}=\OC(1)$. Combining these estimates, we get $\Eb{|\th(v)|\md \ev
K_1,\ev K_3}=O_{C,\upperhco}(1)$.

We will now apply the argument used to prove~\eref{e.maxbh}
in order to establish
\begin{equation}\label{e.vss}
\EB{\max\bl\{|\th(v)|:v\in \Vss\br\}\md \ev K_1,\ev K_3}=O_{C,\upperhco}(1).
\end{equation}
Indeed, let $A$ denote the set of vertices in $\Vs$ neighboring
with some vertex outside $\Vs$, and let $H(v,u)$ denote the probability
that simple random walk started at $v$ first hits $A$ in $u$.
As in the proof of~\eref{e.maxbh}, $H(v,u)\le O_C(1)\,H(v',u)$ for 
$v,v'\in \Vss$.  Now~\eref{e.vss} follows as in the proof of~\eref{e.maxbh}.

Next, we want to show that~\eref{e.vss} holds with $h_1$ and $h_3$ replacing
$\th$; that is,
\begin{equation}\label{e.vssj}
\Es{M\md \ev K_1,\ev K_3}=O_{C,\upperhco}(1)\,,
\end{equation}
where $M:= \max\bl\{|h_j(v)|:v\in \Vss,\, j=1,3\br\}$.
Let $v_1$ be the vertex $v\in \Vss$ where $|h_1(v)|$ is maximized,
and let $v_3$ be the vertex $v\in \Vss$ where $|h_3(v)|$ is maximized.
Then $M=\max\{|h_1(v_1)|,|h_3(v_3)|\}$.  Assume, for now, that $M= |h_3(v_3)|$.
The maximum principle for discrete harmonic functions implies
that $v_3$ neighbors with a vertex outside $\ball_{r_5}$ and
$v_1$ neighbors with a vertex in $\ball_{r_3}$.
Let $p$ be the probability that simple random walk started at $v_3$
exits $\ball_{R_3}$ before hitting a vertex neighboring with
a vertex in $\ball_{r_3}$.
Then $p$ is bounded away from $0$ by a function of $C$.
Since $h_1$ composed with a simple random walk is a martingale while
the walk stays in $\Vp$,
we get $|h_1(v_3)|\le (1-p)\, |h_1(v_1)| \le (1-p)\, M$.
Since $\th=h_1+h_3$, we get,
$ |\th(v_3)|\ge |h_3(v_3)|-|h_1(v_3)| = M-|h_1(v_3)| \ge p\,M $.
The case $M=|h_1(v_1)|$ is similarly treated. 
Using~\eref{e.vss}, we then get~\eref{e.vssj}.

Next, we want to prove
\begin{equation}\label{e.ipm}
|\nabla h_1\cdot\nabla h_3| = O(M^2)\,.
\end{equation}
Since $h_1$ is harmonic in $\Vp$, if $v\in \Vp$,
$$
\sum_{u\sim v}\bl(h_1(v)-h_1(u)\br)\,h_3(v)=0,
$$
where the sum is over the neighbors of $v$.  This is also true for
$v\in \Va$, since $h_3$ is zero there.  Consequently,
\begin{equation}\label{e.v2}
\sum_{v\in \Vtwo}\sum_{u\sim v}\bl(h_1(v)-h_1(u)\br)\,h_3(v)=0.
\end{equation}
Similarly, we find,
\begin{equation}\label{e.v2p}
\sum_{u\in \Vtwop}\sum_{v\sim u}\bl(h_3(u)-h_3(v)\br)\,h_1(u)=0.
\end{equation}
Set $\p \Vtwo:= \{(v,u)\in\Vtwo\times\Vtwop:u\sim v\}$.
By considering the contribution of each edge $[v,u]$ to $\nabla h_1\cdot\nabla h_3$,
we compare the sum of the left hand sides of~\eref{e.v2} and~\eref{e.v2p} to
$\nabla h_1 \cdot\nabla h_3$ and conclude that
\begin{equation}\label{e.nablas}
\begin{aligned}
&
\nabla h_1 \cdot\nabla h_3
=
\sum_{(v,u)\in\p\Vtwo}
\bl(h_3(v)\,h_1(u)-h_1(v)\,h_3(u)\br)
\\&\;\;
=
\sum_{(v,u)\in\p\Vtwo}
\Bigr(\bl(h_1(u)-h_1(v)\br)h_3(u)+
\bl(h_3(v)-h_3(u)\br)h_1(u)\Bigr).
\end{aligned}
\end{equation}
The number of summands is clearly $O(R)$. Note that for every $v\in
\Vtwo$ neighboring with a vertex $u\in \Vtwop$, there is a disk of
radius proportional to $R/C$ such that all the vertices in that disk
are in $\Vss$.  Consequently, the Discrete Harnack Principle~\ref{l.harnack}
gives $\bl|h_1(u)-h_1(v)\br|=O(M/R)$ and
$\bl|h_3(u)-h_3(v)\br|=O(M/R)$. Hence,~\eref{e.nablas}
gives~\eref{e.ipm}.

Now,~\eref{e.vssj} implies that the expectation
of $ |\nabla h_1\cdot\nabla h_3|^{1/2}$ conditioned on $\ev K_1\cap
\ev K_3$ is bounded by a function of $C$ and $\upperhco$. In
particular, there is a constant $c_1=c_1(\upperhco,C)>0$ such that
$$
\Pb{|\nabla h_1\cdot\nabla h_3|<c_1\md \ev K_1,\,\ev
K_3}>(c_1)^{-1}\,.
$$ In terms of Lebesgue measure, this may be
written as
$$
{\int_{\ev K_1\cap\ev K_3}
\exp\Bl(\frac{-\|\nabla\th\|^2}{2}\Br)}<c_1 {\int_{\ev K_1\cap\ev
K_3} \exp\Bl(\frac{-\|\nabla\th\|^2}{2}\Br)}\, 1_{\{|\nabla
h_1\cdot\nabla h_3|<c_1\}}\,.
$$
Since $\th=h_1+h_3$, this implies
$$
{\int_{\ev K_1\cap\ev K_3} \exp\Bl(\frac{-\|\nabla\th\|^2}{2}\Br)} <
c_1\,e^{c_1}\, \int_{\ev K_1\cap\ev K_3}
\exp\Bl(-\frac{\|\nabla h_1\|^2 +\|\nabla h_3\|^2}{2}\Br) ,
$$
which gives one side of~\eref{e.prod}.

The other direction is proved in essentially the same way. Under the
probability measure weighted by $\exp\bl(-(\|\nabla h_1\|^2
+\|\nabla h_3\|^2)/2\br)$ (with respect to Lebesgue measure on
$\R^{\Vpp}$), $h_1$ restricted to $\Va$ has the law of the DGFF with
zero boundary values on $\Vb$ restricted to $\Va$. Similarly, with
this weighting, $h_3$ restricted to $\Vb$ has the law of the DGFF
with zero boundary values on $\Va$ and with boundary values given by
$h_\p$ on $\p D$, restricted to $\Vb$. Moreover, under this measure
$h_1$ and $h_3$ are clearly independent.  The above arguments show
that under this measure~\eref{e.vssj} holds (where $\ev K_1$ refers
to $h_1$ while $\ev K_3$ refers to $h_3$). Since~\eref{e.ipm} is
still valid, the opposite inequality in~\eref{e.prod} is then easily
established.

We may also apply~\eref{e.prod} in the case where $V_+^1=V_-^1=\emptyset$,
and hence $\ev K_1$ has full measure.
Since
\begin{multline*}
\int_{\ev K_1\cap\ev K_3} \exp\Bl(-\frac{\|\nabla h_1\|^2}{2}\Br)
\exp\Bl(-\frac{\|\nabla h_3\|^2}{2}\Br)
\\
= \int_{\R^{\Va}} 1_{\ev K_1} \exp\Bl(-\frac{\|\nabla
h_1\|^2}{2}\Br) \int_{\R^{\Vb}} 1_{\ev K_3} \exp\Bl(-\frac{\|\nabla
h_3\|^2}{2}\Br),
\end{multline*}
from~\eref{e.fii} we get,
$$
\Pb{\ev K_1\md \ev K_3} \asymp \frac {\int_{\ev K_1}
\exp\bl(-\|\nabla h_1\|^2/2\br)} {\int{ \exp\bl(-\|\nabla
h_1\|^2/2\br)}} \,,
$$
which completes the proof.
\QED

\subsection{Narrows and obstacles} \label{ss.narrows}

The present subsection and the next will use the infrastructure
developed in~\S\ref{ss.estimates} to prove some bounds on the
probabilities that contour lines cross certain regions in specified
ways. This is roughly in the spirit of the Russo-Seymour-Welsh
theorem for percolation, though the proofs are entirely different.

The following lemma is an estimate for having a crossing by $\TG^*$-hexagons
where $h$ is negative between two arcs on the boundary of some
subset of the domain, conditioned on some zero height interface
paths. The statement below is slightly complicated, because we need
to keep the geometric assumptions quite general. In percolation,
boundary values do not play a role, of course. But in our case we
need the crossing estimate in the case where one boundary arc of the
domain is conditioned to be an interface.

\begin{figure}
\SetLabels
\L(.51*.08)$\gamma_1$\\
\L(.13*.35)$A$\\
\L(.54*.5)$\alpha'$\\
\L(.76*.21)$\alpha_1$\\
\B(.5*.77)$\alpha_2$\\
\endSetLabels
\centerline{\epsfysize=2.5in%
\AffixLabels{%
\epsfbox{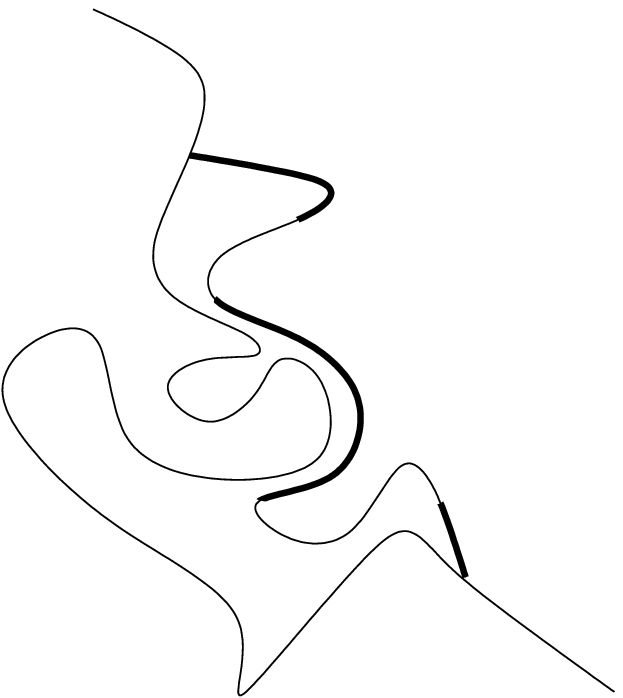}}
}
\begin{caption} {\label{f.narrows}Setup in the narrows lemma.}
\end{caption} \end{figure}

\begin{lemma}[Narrows]\label{l.narrows}
For every $\eps>0$ there is a $\delta=\delta(\upperhco,\eps)>0$ such
that the following crossing estimate holds. Assume~\iref{i.h} and~\iref{i.D}. Let
$\ev K$ be the event that a fixed collection
$\{\gamma_1,\dots,\gamma_k\}$ of oriented paths in $\TG^*$ are
contained in oriented zero height interfaces of $h$, and suppose
that $\P[\ev K]>0$. Let $\alpha\subset D\setminus
(\gamma_1\cup\cdots\cup\gamma_k)$ be a simple path
that has both its endpoints
on the right hand side (positive side) of $\gamma_1$.
Let $A$ be the domain bounded by $\alpha$ and a subarc of $\gamma_1$, and assume
that $A$ does not meet the left side of $\gamma_1$
and $\closure A\cap (\gamma_2\cup\gamma_3\cup\cdots\cup\gamma_k\cup\p D)=\emptyset$.
Let $\alpha_1,\alpha_2,\alpha'$ be three disjoint subarcs
of $\alpha$, where $\alpha_1$ contains one endpoint of $\alpha$
and $\alpha_2$ contains the other endpoint of $\alpha$.
(See Figure~\ref{f.narrows}.)
Suppose that each point in $\alpha'$ is contained in a $\TG^*$-hexagon
whose center is outside $A$.
Set $d_1:=\sup_{z\in\alpha'} \dist(z,\gamma_1)$.
Let $\ev C$ be the event
that there is a path crossing from
$\alpha_1$ to $\alpha_2$ in $A$ inside hexagons where
$h$ is negative.
Let $d^*$ be the infimum diameter of any path connecting
$\alpha'$ to $\gamma_1\cup\cdots\gamma_k\cup\p D$ which does
not contain a sub-path connecting $\alpha\setminus (\alpha_1\cup\alpha_2)$ to $\gamma_1$ in $A$.
If
\begin{equation}\label{e.assume}
d_1+1 \le \delta\,\min\bl\{d^*,
\dist(\alpha_1,\alpha'), \dist(\alpha_2,\alpha'),\diam \alpha' \br\},
\end{equation}
then
$$\Pb{\ev C\md\ev K}<\eps\,.
$$
\end{lemma}

The idea of the proof is to observe the effect that such a crossing would
have on certain averages of heights of vertices, and thereby conclude that
it is unlikely. The challenge in the implementation of this strategy is
to condition on a crossing in such a way that would make the expected
heights easy to estimate.

\proof
Set $N:=\lfloor \diam \alpha'/2\rfloor$.
Assume that~\eref{e.assume} holds.
Note that $\diam \gamma_1 \ge \diam \alpha'-2\,d_1$.
Choose some point $z_0\in\alpha'$.
For $j=1,2,\dots,N$ let $z_j$ be a point in $\alpha'$ at distance
$j$ from $z_0$. Now for each $z_j$ we let $s_j$ be some center of
a hexagon that contains $z_j$ satisfying $s_j\notin A$.
Then $|s_j-s_i|\ge |j-i|-O(1)$.
Let $U$ be the union of $V_\p$ and the set of vertices adjacent to any one of the
paths $\gamma_1,\dots,\gamma_k$.
(These are precisely the vertices $v$ where $h$ takes boundary values or the
sign of $h(v)$ is determined by $\ev K$.)
Set
$$
X:=N^{-1}\sum_{j=1}^N  h(s_j)\,,
$$
$b:=\Eb{X\md\ev K}$ and
fix some $\eps'>0$. We first claim that
\begin{equation}\label{e.Xcondvar}
\Eb{(X-b)^2\md\ev K}<\eps'
\end{equation}
if $\delta=\delta(\eps',\upperhco)>0$ is sufficiently small.
Let $h_U$ denote the harmonic extension of the restriction of $\hth$ to $U$
and set $X_U:=N^{-1}\sum_{j=1}^N h_U(s_j)$.
Note that $\Eb{X-X_U\md\ev K,h_U}=0$, and hence $\Eb{X_U\md\ev K}=b$.
For  each $u\in U$ and $j\in\{1,2,\dots,N\}$ let $p(j,u)$
denote the probability that simple random walk started from $s_j$
first hits $U$ at $u$. Also set $p(u):=N^{-1}\sum_{j=1}^N p(j,u)$.
Then $X_U=\sum_{u\in U} p(u)\,h(u)$.
Consequently,
\begin{multline*}
\Eb{(X_U-b)^2\md\ev K}
\\
=\sum_{u,u'\in U}
p(u)\,p(u')\,\Bl(\Eb{h(u)\,h(u')\md\ev K}-
\Eb{h(u)\md\ev K}\,\Eb{h(u')\md\ev K}\Br).
\end{multline*}
Let $Z(u,u')$ denote the term in parentheses corresponding to
the summand involving $u$ and $u'$.
Then $\Eb{(X_U-b)^2\md\ev K}$
is just the average of the
conditioned covariances $Z(u,u')$ weighted by $p(u)\,p(u')$.
We know from~\eref{e.hbexp} that $\Eb{h(u)^2\md\ev K}$
is bounded by a constant depending only on
$\upperhco$. It then follows by Cauchy-Schwarz that the same is
true for $Z(u,u')$.
Consequently, to prove that $\Eb{(X_U-b)^2\md\ev K}$ is
small, it suffices to show that
when $(u,u')$ is chosen with probability
$p(u)\,p(u')$ it is very likely that $|Z(u,u')|$ is small.
Suppose that we select $j$ from $\{1,2,\dots,N\}$ uniformly
at random and given $j$ select $u\in U$ with probability
$p(j,u)$. Independently, we also select $(j',u')$ with
the same distribution. It suffices to show that
$|Z(u,u')|$ is likely to be small, and
by Corollary~\ref{c.corr} it suffices to show that
the distance between $u$ and $u'$ is likely
to be large. Since $|s_j-s_{j'}|=|j-j'|+O(1)$ is unlikely
to be much smaller than $\diam \alpha'$,
which is larger than $\delta^{-1}\,(d_1+1)$,
it follows from Lemma Hit Near~\ref{l.hitnear}
that for any fixed $R$ the probability that $|u-u'|<R$ tends to zero as $\delta\to0$.
Consequently,
$$
\Eb{(X_U-b)^2\md\ev K}<\eps'/2\,,
$$
provided that $\delta$ is sufficiently small.

Set $X_j:=h(s_j)-h_U(s_j)$. Recall from \S\ref{ss.generalremarks} that
given the restriction of $h$ to $U$ the function $h-h_U$ is the DGFF
on $\VD\setminus U$ with zero boundary values on $U$.
By~\eref{e.greensrandomwalk}, therefore
$\Eb{X_iX_j\md \ev K,h_U}=G(s_i,s_j)/6$, where $G(v,u)$ is the
expected number of visits to $u$ by a random walker started at $v$
and stopped when it hits $U$.
From Lemma Hit Near~\ref{l.hitnear} and
Lemma~\ref{l.green}
$$
G(s_i,s_j)\le\begin{cases}O(1)\log(d_1/(|s_i-s_j|\vee 1)) & |s_i-s_j|<d_1/2,\\
O(1)\bl(d_1/( |s_i-s_j|\vee 1)\br)^{\expo1} & |s_i-s_j|\ge d_1/2.
\end{cases}
$$
Since $|s_j-s_i|\ge |j-i|-O(1)$ and
$\expo1\in(0,1)$, these estimates give
\begin{equation*}
\Eb{(X-X_U)^2\md \ev K}=N^{-2}\,\sum_{i,j=1}^N G(s_j,s_i)/6=O(1)\,(d_1/N)^\expo1.
\end{equation*}
Now~\eref{e.Xcondvar} follows for sufficiently small $\delta=\delta(\eps',\upperhco)>0$,
since $X-X_U$ is independent from $X_U$ and they are also independent given $\ev K$.

\medskip

We now claim that $b\ge 0$ if $\delta$ is sufficiently small.
Let $c$ be the constant given by Lemma~\ref{l.ibd}.
If $u$ is a fixed vertex adjacent to $\gamma_1$ on the right,
 then $\Eb{h(u)\md\ev K}>1/c$, by~\eref{e.atv}.
On the other hand,
$\Eb{|h(u)|\md\ev K}<c$ for every $u\in U$.
By~\eref{e.assume} and Lemma Hit Near~\ref{l.hitnear}
it follows that when $\delta$ is small with high probability a random walk
starting at any $s_j$ is likely to first hit $U$ at a vertex adjacent
to the right hand side of $\gamma_1$.
Thus, when $\delta$ is small, we have $b=\Eb{X_U\md\ev K}>0$.
Also, clearly, $b\le c$.

\medskip

Now set $a=b+1$. Let $Q$ denote the union of the closed hexagons in
$\TG^*$ for which $\hth(H)\in [0,a]$ and let $Q'$ denote the union
of the edges in $\TG^*$ that are  on the common boundary of two
hexagons $H_1$ and $H_2$ satisfying $\hth(H_1)<0,\, a<\hth(H_2)$.
Let $Q_0$ denote the connected component of $(Q\cup Q')\cap
\closure{A}$ that contains $\gamma_1\cap\p A$. Let $\ev Q$ denote
the event $Q_0\cap \alpha\setminus(\alpha_1\cup\alpha_2)=\emptyset$.
If $\ev C$ holds, then the corresponding crossing by hexagons where
$h$ is negative separates $\gamma_1\cap \p A$ from
$\alpha\setminus(\alpha_1\cup\alpha_2)$ in $A$, and hence $\ev Q$
holds as well. Thus, $\ev C\subset\ev Q$.

\begin{figure}
\SetLabels
(.1*.09)$1$\\
(.28*.09)$2$\\
(.37*.552)$Q_0$\\
(.46*.09)$3$\\
(.55*.24)$4,6$\\
(.46*.398)$5$\\
(.91*.552)$9$\\
(.55*.86)$12$\\
(.73*.244)$7$\\
(.93*.25)$A'$\\
(.82*.398)$8$\\
(.82*.704)$10$\\
(.64*.704)$11$\\
\endSetLabels
\centerline{\epsfysize=2.5in%
\AffixLabels{%
\epsfbox{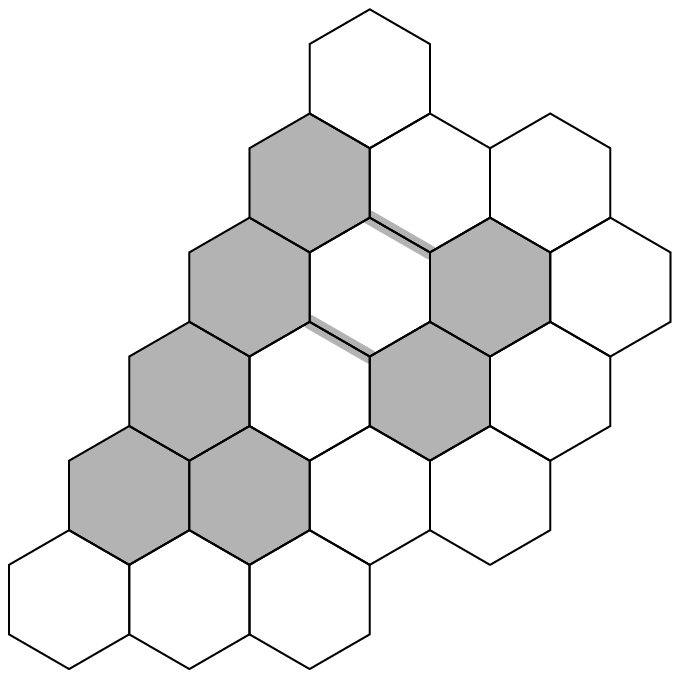}}
}
\begin{caption} {\label{f.hexseq}A portion of the sequence of
hexagons adjacent to $\p A'\setminus\p A$.}
\end{caption} \end{figure}

On the event $\ev Q$, let $A'$ be the connected component of $\closure A\setminus Q_0$
that contains $\alpha'$, and let $\Qp$ denote the set of centers of
hexagons $H$ such that $\closure{H\cap A'}\cap Q_0\ne\emptyset$.
Clearly, $\hth(v)<0$ or $\hth(v)>a$ for each $v\in \Qp$.
Since $A'$ and $Q_0$ are connected, it is immediate to verify (using
the Jordan planar curve theorem) that $A'$ is simply connected and
$\p A'\setminus\p A\subset Q_0$ is connected.
The closed hexagons of $\TG^*$
with centers in  $\Qp$
form a sequence (possibly with repetitions) with each
pair $H,H'$ of consecutive hexagons along the sequence satisfying
$H\cap H'\setminus Q_0\ne\emptyset$ and $H\cap H'\cap Q_0\ne\emptyset$.
(See Figure~\ref{f.hexseq}.)
If $v,u\in \Qp$ are centers of consecutive hexagons in this sequence,
then it is impossible that
$\hth(v)<0$ and $\hth(u)>a$ (otherwise, the boundary
between the hexagons would be in $Q'$).
Thus, either $\hth(\Qp)\subset (a,\infty)$,
or $\hth(\Qp)\subset(-\infty,0)$.
Let $\ev Q_+$ be the event that $\ev Q$ occurs and
$\min\hth(\Qp)>a$, and let $\ev Q_-$ be the event that
$\ev Q$ occurs and $\max\hth(\Qp)<0$.

We now want to estimate $\Eb{X\md \ev K,\ev Q_-}$ and
$\Eb{X\md \ev K,\ev Q_+}$.
Let $U'$ be the set of vertices that are either
in hexagons adjacent to $Q_0$ or in $U$.
Since $a=\OC(1)$,
it is clear that the proof of~\eref{e.atv} gives
$\Eb{|h(u)|\md \ev K,U',\ev Q_\pm}=\OC(1)$
for $u\in U'$. On the other hand, Lemma Hit Near~\ref{l.hitnear}
shows that at least $1-O(\delta^{\expo1})$ of the discrete harmonic
measure on $U'$ starting from every $s_j$
is in $\Qp$. If $\ev Q_-$ holds and $u\in\Qp$,
then $\Eb{h(u)\md\ev K,\ev Q_-,\Qp}$ is negative and bounded
away from zero, by the corresponding analog of the left
hand side of~\eref{e.atv}.
Thus, we find that
$\Eb{X\md\ev K,\ev Q_-}$ is negative and bounded away from
zero when $\delta=\delta(\eps',\upperhco)>0$ is small.
Since $\Eb{(X-b)^2\md \ev K}<\eps'$ and
$b\ge 0$, we conclude that by choosing $\eps'>0$ small
it can be guaranteed that $\Pb{\ev Q_-\md\ev K}<\eps/2$.

On $\ev Q_+$, we clearly have $h(u)>a\ge b+1$ on every
$u\in\Qp$. Thus, as above, it follows that when $\delta$ is
small $\Eb{X\md\ev K,\ev Q_+}> b+1/2$.
This again implies that $\Pb{\ev Q_+\md\ev K}$ can
be made smaller than $\eps/2$.
Since $\ev C\subset\ev Q_+\cup\ev Q_-$, this completes
the proof.
\QED

Next, we formulate an analogous lemma for crossings near the
boundary of the domain.

\begin{lemma}[Domain boundary narrows]\label{l.bdnarrows}
There is a constant $\lowerhco=\lowerhco(\upperhco)>0$ such that for
every $\eps>0$ there is a $\delta=\delta(\upperhco,\eps)>0$ such
that the following crossing estimate holds. Assume~\iref{i.h} and~\iref{i.D},
assume that $\p_+$ is a simple path contained in $\p D$ and that
$h_\p\ge -\lowerhco$ on $\p_+\cap V_\p$. Set $\p_-:=\p D\setminus
\p_+$. Let $\ev K$ be the event that a fixed collection
$\{\gamma_1,\dots,\gamma_k\}$ of oriented paths in $\TG^*$ are
contained in oriented zero height interfaces of $h$, and suppose
that $\P[\ev K]>0$. Let $\alpha\subset D\setminus
(\gamma_1\cup\cdots\cup\gamma_k)$ be a simple path that has both its
endpoints on $\p_+$. Let $A$ be the domain bounded by $\alpha$ and a
subarc of $\p_+$, and assume that $\closure A\cap
(\gamma_1\cup\gamma_2\cup\cdots\cup\gamma_k\cup\p_-)=\emptyset$. Let
$\alpha'\subset\alpha$ be a subarc. Suppose that each point in
$\alpha'$ is contained in a hexagon whose center is outside $A$. Set
$d_1:=\sup_{z\in\alpha'} \dist(z,\p_+)$. Let $\alpha_1$ be a subarc
of $\alpha$ that contains one of the endpoints of $\alpha$, and let
$\alpha_2$ be a subarc of $\alpha$ that contains the other endpoint
of $\alpha$. Let $\ev C$ be the event that there is a path crossing
from $\alpha_1$ to $\alpha_2$ in $A$ inside hexagons where $h$ is
negative. Let $d^*$ be the infimum diameter of any path connecting
$\alpha'$ to $\gamma_1\cup\cdots\gamma_k\cup\p D$ which does not
contain a sub-path connecting $\alpha\setminus
(\alpha_1\cup\alpha_2)$ to $\p_+$ in $A$. If
\begin{equation}\label{e.assumeD}
d_1+1 \le \delta\,\min\bl\{d^*,
\dist(\alpha_1,\alpha'), \dist(\alpha_2,\alpha'),\diam \alpha' \br\},
\end{equation}
then
$$\Pb{\ev C\md\ev K}<\eps\,.
$$
\end{lemma}

\proof
The proof is slightly simpler but essentially the same as
the proof of the Narrows Lemma~\ref{l.narrows}.
We use the same notations as in that lemma,
and only indicate the few differences in the proof.
In the present setting $b=\Eb{X\md\ev K}$ can be made
larger than $-2\,\lowerhco$ by taking $\delta>0$ small.
Here, we define $Q_0$ as the connected component of
$(Q\cup Q'\cup \p_+)\cap \closure{A}$
that contains $\p_+\cap\p A$.
Observe that $\Qp\cap\p D=\emptyset$ on $\ev Q_-\cap\ev K$.
It follows that $\Eb{X\md\ev K,\ev Q_-}$ is negative
and bounded away from zero (by a function of $\upperhco$)
when $\delta>0$ is small.
By taking $\lowerhco>0$ sufficiently small, we can make sure that
$\Eb{X\md\ev K,\ev Q_-}<-3\,\lowerhco$. But since
$\Eb{(X-b)^2\md\ev K}$ is arbitrarily small and $b\ge -2\,\lowerhco$,
this makes $\Pb{\ev Q_-\md\ev K}$ small.
The rest of the argument is essentially the same.
\QED

The previous lemmas will help us control the behavior of the
continuation of contours near existing contours or the boundary of
the domain. The next lemma will help us control the behavior in the
interior away from existing contours.

\begin{lemma}[Obstacle]\label{l.obstacle}
For every $\eps>0$ there is some constant $c=c(\eps,\upperhco)>0$
such that the following estimate holds. Assume~\iref{i.h} and~\iref{i.D}. Let $\ev
K$ be the event that a fixed collection
$\{\gamma_1,\dots,\gamma_k\}$ of oriented paths in $\TG^*$ are
contained in oriented zero height interfaces of $h$, and suppose
that $\P[\ev K]>0$. Let $U$ be the union of $V_\p$ and the vertices
of $\TG$ adjacent to
 $\bar\gamma:=\gamma_1\cup\gamma_2\cup\cdots\cup\gamma_k$.
Let $g$ be a function defined on the vertices of $\TG$ that is $0$ on $U$.
Let $\beg$ denote the union of the interfaces of $h+g$ that contain
any one of the paths $\gamma_1,\gamma_2,\dots,\gamma_k$.
Let $B(z_0,r)$ be a disk of radius $r$ 
that is centered at some vertex
$z_0$ satisfying $|g(z_0)|\ge \|g\|_\infty/2$.
Let $d>0$ and suppose that at distance at most $\eps ^{-1} \,d$ from $z_0$
there is a connected component of $\bar\gamma\cup\p D$ whose
diameter is at least $\eps\,d$.
Also assume $\|g\|_\infty/\|\nabla g\|_\infty>c\,r>c$.
Then
$$
\Pb{\beg\cap B(z_0,r)\ne\emptyset\md\ev K} \le
c\, \log(d/r)\,\|g\|_\infty^{-2}.
$$
\end{lemma}

\proof
With no loss of generality, we assume that $g(z_0)>0$.
Set $q:=\|g\|_\infty/\|\nabla g\|_\infty$ and $r_1:= q/10$.
Since between any two vertices $z,z'$ in $\TG$ there is a path in
$\TG$ whose length is at most $2\,|z-z'|$,
$$
\min\{g(z):z\in B(z_0,r_1)\}\ge g(z_0)- 2\, r_1\,\|\nabla g\|_\infty
= g(z_0) - \|g\|_\infty/5 \ge \|g\|_\infty/4\,.
$$
Since $g=0$ on $U$, it follows that $\eps^{-1}\,d\ge r_1$.
Since we are assuming that $ q>c\,r$, and we may assume that $c$ is a large constant
which may depend on $\eps$,
it follows that $d/r>100$, say.
Thus, we also assume, with no loss of generality, that $\|g\|_\infty\ge \sqrt c$,
since the required inequality is trivial otherwise.

Let $X$ denote the average value of $\hth$ on the vertices in $B(z_0,r)$.
The inequality~\eref{e.av} and $d/r>100$ give
\begin{equation}\label{e.condvar}
\Eb{X^2\md \ev K} \le O_{\eps,\upperhco}(1)\log (d/r)\,.
\end{equation}

If $\gamma_1$ is not a closed path,
we start exploring the interface of $h+g$ containing $\gamma_1$ starting from
one of the endpoints of $\gamma_1$
until that interface is completed or $B(z_0,r)$ is hit, whichever occurs
first. (This may entail going through several of the
interfaces $\gamma_j$, $j>1$.)
If that interface is completed before we hit $B(z_0,r)$, we continue and
explore the interface of $h+g$ containing $\gamma_2$, and so forth,
until finally either all of $\beg$ is explored or
$B(z_0,r)$ is hit.
Let $\ev Q$ denote the event that $B(z_0,r)$ is hit, and let
$\beta$ be the interfaces explored up to the time  when
the exploration terminates.

Let $U'$ be the union of $U$ with the vertices
adjacent to $\beta$.
Since we are assuming
$q>c\,r$, $r>1$ and $\|g\|_\infty\ge \sqrt c$, and since $c$ may be chosen
an arbitrarily large constant,
Lemma~\ref{l.ibd2} shows that for every
vertex $v\in U'$, we have
\begin{equation}\label{e.bdval}
\EB{|\hth(v)+g(v)|\md \ev K,\,\ev Q,\,\beta}
\le \|g\|_\infty/100\,.
\end{equation}
(Note that conditioning on $\ev Q$ and $\beta$ amounts to
conditioning that $\hth+g>0$ on vertices adjacent
to the right hand side of $\beta$ and
$\hth+g<0$ on vertices adjacent to the left hand side.
Consequently, the lemma applies.)

Now let $x$ be any vertex in $B(z_0,r)$. For each $u\in U'$,
let $p_u$ denote the probability that simple random
walk started at $x$ will first hit $U'$ in $u$; that is,
the discrete harmonic measure from $x$.
Then~\eref{e.bdval} gives
\begin{equation}\label{e.hkqb}
\EB{\hth(x)\md \ev K,\,\ev Q,\,\beta} \le
\frac{\|g\|_\infty}{100}
-
 \sum_{u\in U'} p_u\, g(u)\,.
\end{equation}
Since $\beta\cap B(z_0,r)\ne\emptyset$ and
$\beta$ intersects the complement of $B(z_0,r_1)$,
Lemma Hit Near~\ref{l.hitnear} gives
$$
\sum_{u\in U'\setminus B(z_0,r_1)}p_u \le O(1)\, (r/r_1)^\expo1.
$$
Since we are assuming that $q\ge c\,r$, we may assume that the right
hand side is less than $1/10$.
Recall that $g\ge \|g\|_\infty/4$ inside $B(z_0,r_1)$.
Outside $B(z_0,r_1)$, the trivial estimate $g\ge -\|g\|_\infty$ applies.
When these estimates are applied to~\eref{e.hkqb}, one gets,
\begin{multline*}
\EB{\hth(x)\md \ev K,\,\ev Q,\,\beta} -\frac{\|g\|_\infty}{100}
\le
-\sum_{u\in U'\cap B(z_0,r_1)}p_u\, g(u)
- \sum_{u\in U'\setminus B(z_0,r_1)}p_u\, g(u)
\\
\le
-(\|g\|_\infty/4)(1 - 1/10)
+\|g\|_\infty/10
=
-\|g\|_\infty/8
\,.
\end{multline*}
We may take expectation with respect to $\beta$ and average with respect to $x$ to conclude that
$ \Eb{X\md \ev K,\,\ev Q} \le -\|g\|_\infty/9$, which implies by Jensen's inequality
$$
 \Eb{X^2\md \ev K,\,\ev Q} \ge \|g\|^2_\infty/81\,.
$$
Since
$$
\Pb{\ev Q\md \ev K} \le \Eb{X^2\md \ev K}/\Eb{X^2\md\ev Q,\ev K},
$$
the lemma now follows from the above and~\eref{e.condvar}.
\QED

\subsection{Barriers} \label{ss.barriers}

In this subsection we apply Lemmas~\ref{l.narrows}, 
\ref{l.obstacle} and~\ref{l.distort} to get a flexible (though slightly complicated)
criterion giving lower bounds for the probability that contours
avoid certain sets.
The complications arise from the need to handle pre-existing contours that
are highly non-smooth on large scales.

The following relative notion of distance will
sometimes be used below:
\begin{equation}\label{e.distdef}
\dist(A,B;X):=
\inf\bl\{\diam \alpha:\alpha\text{ is a path in }\closure X\text{ connecting }A\text{ and }B\br\}\,.
\end{equation}
If $\bar\gamma$ is a collection of paths in $\closure D$, let
$\DD(\bar\gamma)$ denote the complement in $D$ of the union
of the closed triangles of $\TG$ meeting $\bar\gamma$.

We now define the notion of a barrier.
Assume~\iref{i.h}.
Let $\gamma_1,\gamma_2,\dots,\gamma_k$
be disjoint simple paths or simple closed paths in $\closure D$.
Set $\bar\gamma:=\gamma_1\cup\cdots\cup\gamma_k$.
Let $\Upsilon\subset\closure{ \DD(\bar\gamma)}$ be a path,
which is contained in $\DD(\bar\gamma)$,
except possibly for its two endpoints.
Fix some $R>0$ and $\eps>0$.
We call $\Upsilon$ an
$(\eps,R)$-{\bf barrier} for the configuration
$(D,\bar\gamma)$ if the following conditions hold:
\begin{enumerate}
\item \label{i.diam} $\eps\,R< \diam\Upsilon\le R$,
\item \label{i.dist} within distance $\eps ^{-1} \,R$ from $\Upsilon$ there is
a connected component of $\bar\gamma\cup\p D$ whose diameter
is at least $\eps R$,
\item \label{i.A}
if $z\in\p \DD(\bar\gamma)$ is
an endpoint of $\Upsilon$
and $\achi_z$ is the connected component of
$\p B(z,\eps\,R)\cap \DD(\bar\gamma)$
first encountered when traversing $\Upsilon$ from $z$
(which exists by~\ref{i.diam}),
then the connected component $A_z$
of $\DD(\bar\gamma)\setminus\chi_z$ that contains
points of $\Upsilon$ arbitrarily close to $z$
satisfies (a) $\p A_z$ consists of $\achi_z$ and a simple path
contained in $\p \DD(\bar\gamma)$ and
(b) $\Upsilon\cap A_z\cap \p B(z,r)$
consists of a single point for every $r\in (0,\eps\,R)$.
\item \label{i.away}
for every point $w\in\Upsilon$ such that
$\dist(w,\p\DD(\bar\gamma))\le \eps\,R/5$
there is an endpoint $z\in \Upsilon\cap\p\DD(\bar\gamma)$
such that $w\in A_z$ and $|w-z|<\eps\,R/2$
(roughly, $\Upsilon$ does not get close to
$\p \DD(\bar\gamma)$ except near its endpoints on $\p\DD(\bar\gamma)$).
\end{enumerate}

\begin{figure}
\SetLabels
(.4*.7)$A_z$\\
\L(.63*.15)$\achi_z$\\
\B(.8*.29)$\Upsilon$\\
\L(.63*.6)$\gamma_1$\\
\B(.48*.305)$z$\\
\R(.2*.7)$\gamma_2$\\
\endSetLabels
\centerline{\epsfysize=2.5in%
\AffixLabels{%
\epsfbox{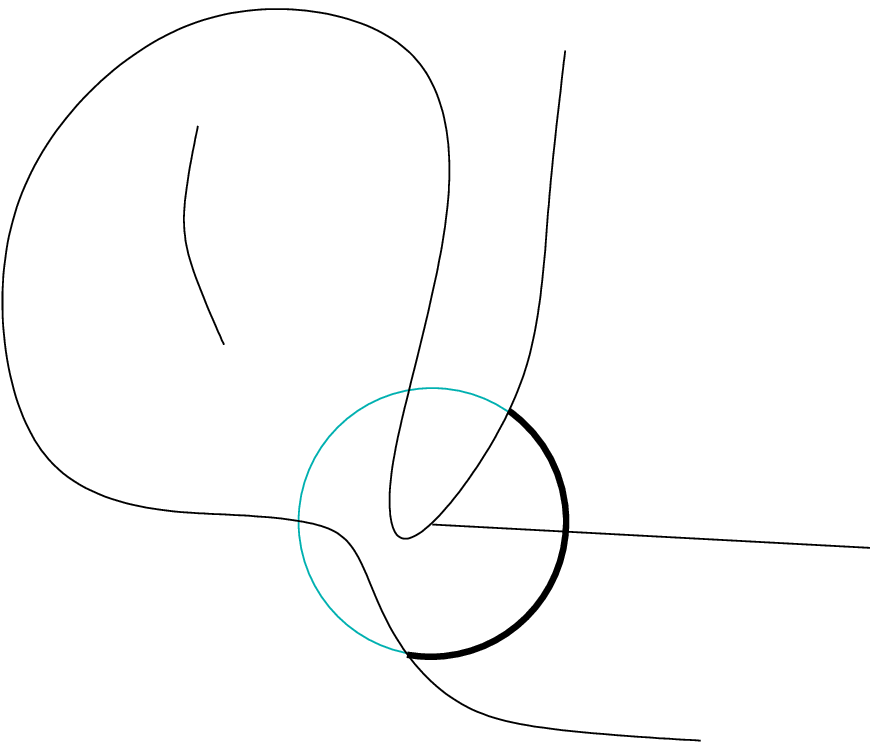}}
}
\begin{caption} {\label{f.A}The domain $A_z$ in a situation
where condition~\ref{i.A} fails.
(The figure does not show detail on the scale of the lattice;
that is, $\DD(\bar\gamma)$ appears as $D\setminus\bar\gamma$.)
}
\end{caption} \end{figure}

One example where condition~\ref{i.A} fails is given
in Figure~\ref{f.A}. If we remove the strand of $\gamma_2$
from that figure, we get an example that illustrates that
$A_z\subset\closure{B(z,\eps\,R)}$ does not follow from
the above conditions.
Note that it may happen that $\p A_z\setminus \achi_z$,
which is a simple path, by~\ref{i.A},
consists of an arc in $\p D$ together with one or two arcs
in $\bar\gamma$ that have endpoints in $\p D$.

We now use barriers to \lq\lq manipulate\rq\rq\ contours of the DGFF.

\begin{theorem}[Barriers]\label{t.barrier}
For every $\eps>0$, $m\in\N_+$ and $\upperhco>0$
there
is a $p=p(\eps,\upperhco,m)>0$
such that the following estimate holds.
Assume~\iref{i.h} and~\iref{i.D}.
Let $\ev K$ be the event that a fixed collection
$\{\gamma_1,\dots,\gamma_k\}$ of oriented paths in $\TG^*$ is
contained in the oriented zero height interfaces of $h$, and suppose
that $\P[\ev K]>0$. Set $\bar\gamma:=\bigcup_{j=1}^k\gamma_j$.
Let $V_+$ [respectively, $V_-$] be the set of vertices adjacent
to $\bar\gamma$
on the right [respectively, left] hand side.
Let $R>0$ and let $Y=Y_+\cup Y_-$
be a collection of $m$ $(\eps,R)$-barriers
for the configuration $(D,\bar \gamma)$.
Assume that the endpoints of these barriers are not on $\p D$
and
that for every $\Upsilon\in Y_+$ [respectively, $\Upsilon\in Y_-$] and
every endpoint $z\in \Upsilon\cap\p\DD(\bar\gamma)$, the vertices in
$\closure {A_z(\Upsilon)}\cap \p \DD(\bar\gamma)$ are in $V_+$
[respectively, $V_-$], where $A_z(\Upsilon)$ is as in
Condition~\ref{i.A}. Also assume that $\dist\bl(\cup Y_+,\cup
Y_-;\DD(\bar\gamma)\br)\ge 2\,\eps\,R$. In the situation where
$\eps\,R=O(1)$, we also need to assume that there is no hexagon in
$\TG^*$ meeting both $V_+\cup(\cup Y_+)$ and $V_-\cup(\cup Y_-)$ and
there is no hexagon meeting $\cup Y$ and $\p D$. Let $\hat\gamma$
denote the union of the zero height interfaces of $h$ which contain
any one of the arcs $\gamma_1,\dots,\gamma_k$. Then
$$
\Pb{(\hat\gamma\setminus\bar\gamma)\cap(\cup Y)=\emptyset\md\ev K}>p\,.
$$
\end{theorem}

The basic idea of the proof is as follows.
We define a function $g$ that is large (positive) near $\cup Y_+$
away from $\p \DD(\bar\gamma)$ and is negative and large in absolute
value near $\cup Y_-$ away from $\p \DD(\bar\gamma)$.
The Obstacle Lemma~\ref{l.obstacle} will
then imply that $\beg$, as defined there, is unlikely to hit $\Upsilon$, except near
endpoints of barriers. The Narrows Lemma~\ref{l.narrows} will be used to show that
$\beg$ is also unlikely to hit $\cup Y$ near endpoints.
Finally, the Distortion Lemma~\ref{l.distort} will be used to conclude that with probability
bounded away from zero, $\hat\gamma$ will not hit $\Upsilon$.

\proof
Let
$$
\hat A^1_+:=\bigcup\bl\{A_z(\Upsilon)\cap B(z,\eps\,R): \Upsilon\in Y_+,z\in\Upsilon\cap\bar\gamma\br\},
$$
$\hat A^2_+:=\bl\{z\in \DD(\bar\gamma):\dist\bl(z,\cup Y_+;\DD(\bar\gamma)\br)\le\eps\,R/10\br\}$,
and $\hat A_+:=\hat A^1_+\cup\hat A^2_+$.
Similarly define $\hat A_-$, with $Y_-$ replacing $Y_+$.
We fix constants $c_0>0$, large, and $\delta>0$ much smaller than
$\eps$, and set
$$
g(z):=
\begin{cases}
c_0\,\delta^{-1}\,R^{-1}\,\min\{\delta\,R,\dist(z,\R^2\setminus \hat A_+)\} &\text{ if } z\in\hat A_+,\\
-c_0\,\delta^{-1}\,R^{-1}\,\min\{\delta\,R,\dist(z,\R^2\setminus \hat A_-)\} &\text{ if } z\in\hat A_-,\\
0&\text{ otherwise.}
\end{cases}
$$
Note that
$\hat A_+\cap\hat A_-=\emptyset$ and $\|g\|_\infty=c_0$.

Let $Z:=\{z\in\R^2:|g(z)|= c_0\}$.
Let $c_1>1$ be some large constant and set
$r:=\delta\,R/c_1$.
Let $W$ be a maximal collection of vertices
in $\{v:|g(v)|\ge c_0/2\}$ such that the distance between
any two distinct vertices in $W$ is at least $r/3$,
and for $a\in W$ let $B_a$ denote the disk of radius
$r$ centered at $a$. Then the disks $B_a$, $a\in W$,
cover $Z$, assuming that $r>1$ and $c_1>100$, say.
Note that
$$
|W|=O(m)\,(R/r)^2=O(m)\,(c_1/\delta)^2.
$$

Fix some $a\in W$. We wish to invoke Lemma~\ref{l.obstacle}
to get a good upper bound on
$ \Pb{\beg\cap B_a\ne\emptyset\md\ev K} $.
We now verify the assumptions of the lemma.
We note that in our case $\|g\|_\infty=c_0$ and
$\|\nabla g\|_\infty = c_0/(\delta\,R)$.
Thus, $q:=\|g\|_\infty/\|\nabla g\|_\infty = \delta\,R$.
Consequently, we set $c_1$ to be the maximum of $100$ and twice
the constant $c$ in the lemma, and the assumption $q>c\,r$ is satisfied.
The assumption $r>1$ will hold once $R$ is large enough,
which we assume for now
(we promise that $\delta$ will be a constant depending only
on $\eps,m$ and $\upperhco$).
For the $d$ in the lemma we may take $d=R$.
Thus, $\log(d/r)=\log(c_1/\delta)$,
and the lemma gives
$$
 \Pb{\beg\cap B_a\ne\emptyset\md\ev K}
\le O_{\eps,\upperhco}(1)\,\log(c_1/\delta)\,c_0^{-2}.
$$
Since $\bigcup_{a\in W}B_a\supset Z$ and
$|W|= O(m)\,(c_1/\delta)^2$, we conclude that
$$
 \Pb{\beg\cap Z\ne\emptyset\md\ev K}
\le O_{\eps,\upperhco}(m)\,\log(c_1/\delta)\,(c_1/\delta)^2\,c_0^{-2}.
$$
Although we have not specified $\delta$ yet, we choose
$c_0=c_0(c_1,\delta,\eps,m,\upperhco)$ so that
\begin{equation}\label{e.say}
 \Pb{\beg\cap Z\ne\emptyset\md\ev K}
\le 1/10\,.
\end{equation}

\medskip

Now that we have established that it is unlikely that $\beg$
intersects $Z$, we need to worry about the case in which
$\beg$ circumvents $Z$ but hits $\cup Y\setminus Z$.
This can only happen near endpoints of barriers.
Let us fix some $\Upsilon\in Y_+$ that has an endpoint, say $z_1$,
on $\p \DD(\bar\gamma)$.
Let $\tchi_{z_1}$ be the arc $\p A_{z_1}(\Upsilon)\setminus\achi_{z_1}$.
(It is an arc, by Condition~\ref{i.A} of the definition of barrier.)
We will now prepare the geometric setup that will enable the use of Lemma~\ref{l.narrows}
to prove that $\Pb{\beg\cap A_{z_1}(\Upsilon)\cap \Upsilon\setminus Z\ne\emptyset\md\ev K}$ is small.

\begin{figure}
\SetLabels
\L(.11*.65)$\gamma_1$\\
\R(.38*.52)$z_1$\\
\B(.9*.63)$\Upsilon$\\
(.6*.7)$\tilde S$\\
\T(.2*.50)$\alpha_0$\\
\B\L(.68*.1)$\alpha_1$\\
\B(.65*.3)$S'$\\
\B(.55*.35)$\alpha_2$\\
(.1*.3)$A$\\
\endSetLabels
\centerline{\epsfysize=2.5in%
\AffixLabels{%
\epsfbox{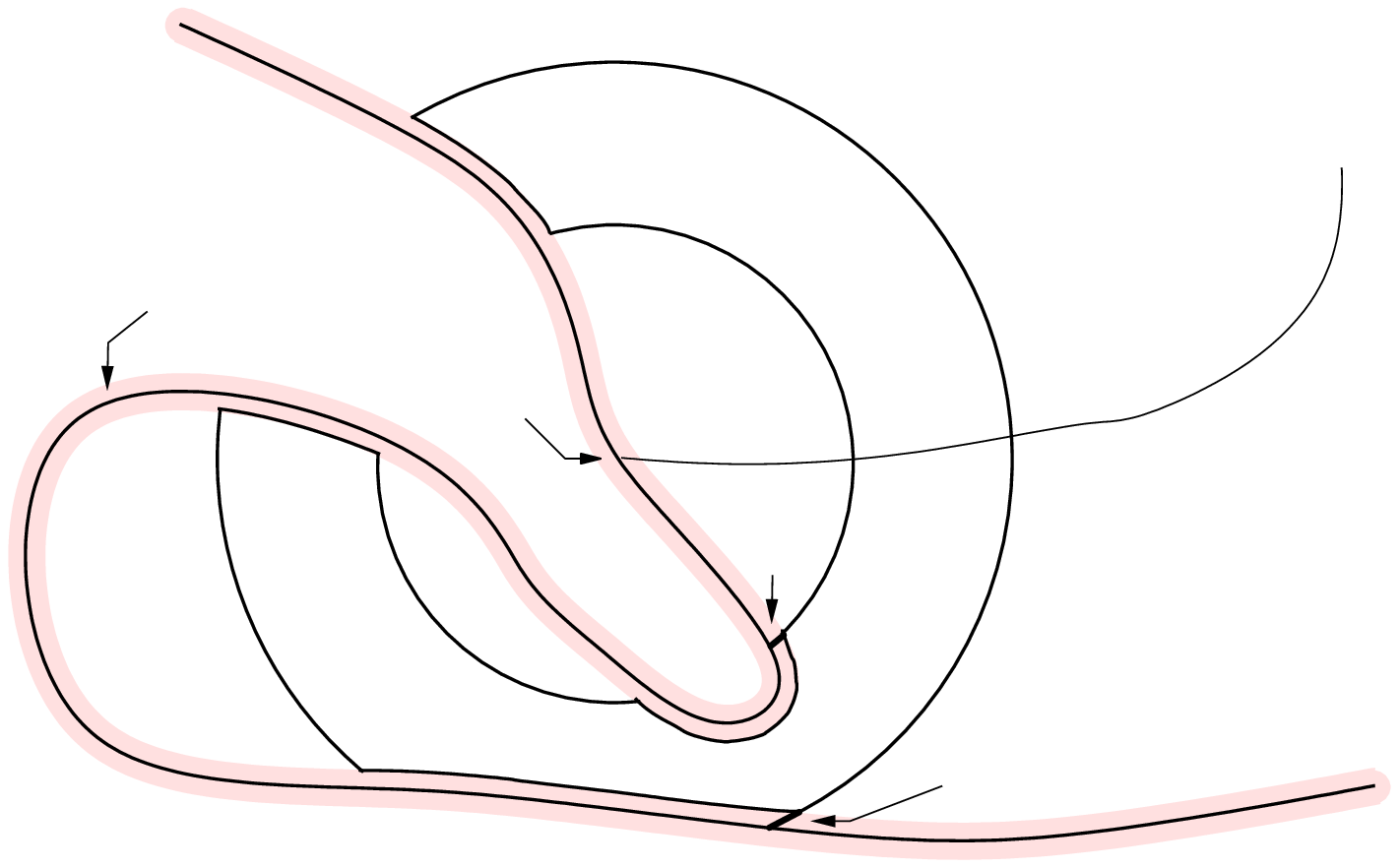}}
} 
\begin{caption} {\label{f.S}The set $S$, and the paths $\alpha_0$ and $\tilde\alpha_0$.
The shaded region is $Q$.}
\end{caption} \end{figure}

Let $Q$ be the set of all hexagons of the grid $\TG^*$ whose
distance from $\p \DD(\bar\gamma)$ is at most $2\,\delta\,R$.
Let $S$ be the connected component of
$$A_{z_1}\cap B(z_1,7\,\eps\,R/8)\setminus (Q\cup B(z_1,5\,\eps\,R/8))$$
that intersects $\Upsilon$.  (See Figure~\ref{f.S}.)
Condition~\ref{i.away} in
the definition of barriers and our assumption that $\delta\ll \eps$ guarantees that
there is a unique such component $S$.
We have $S\subset Z$, and so $\beg$ is unlikely to hit $S$.

Let $S'$ and $\tilde S$ be the two connected components
of $S\setminus \Upsilon$.
Consider the connected component $M_1$
of $\p B(z_1\,7\,\eps\,R/8)\setminus \tchi_{z_1}$
that intersects $\Upsilon$.
Let $\alpha_1$ be the arc in $M_1\setminus \closure S$
that has one endpoint in $\closure {S'}$ and the other in $\tchi_{z_1}$.
Likewise, let $M_2$ be the connected component of
$\p B(z_1\,5\,\eps\,R/8)\setminus \tchi_{z_1}$
that intersects $\Upsilon$, and let
$\alpha_2$ be the arc in $M_2\setminus \closure S$
that has one endpoint in $\closure {S'}$ and the other in $\tchi_{z_1}$.
Let $\alpha$ denote the union of $\alpha_1\cup\alpha_2$
with the arc of $\p S'\setminus\Upsilon$ that connects the endpoint of
$\alpha_1$ with the endpoint of $\alpha_2$.
Let $A\subset A_{z_1}(\Upsilon)$ be the domain whose boundary
consists of $\alpha_1\cup\alpha_2$, an arc of $\p S'$ connecting
$\alpha_1$ and $\alpha_2$ and an arc of $\tchi_{z_1}$ connecting
$\alpha_1$ and $\alpha_2$.
Let $\alpha_0$ be the unique arc of
$\p A\setminus \bl(\p B(z_1\,5\,\eps\,R/8)\cup \p B(z_1\,7\,\eps\,R/8)\br)$
connecting $\p B(z_1\,5\,\eps\,R/8)$ with $\p B(z_1\,7\,\eps\,R/8)$.
Finally, let $\alpha'$ be any sub-arc of $\alpha_0$
whose diameter is $\eps\,R/16$ that intersects
the circle $\p B(z_1,3\,\eps\,R/4)$.
We use $A,\alpha_1,\alpha_2,\alpha'$ in Lemma~\ref{l.narrows}.
In the present situation,
the value of $d_1$ of that lemma is $d_1=2\,\delta\,R+O(1)$.
If we have a path connecting $\alpha'$ to
$\p D\cup\bar\gamma$, it must either
exit $B(z_1,7\,\eps\,R/8)\setminus B(z_1,5\,\eps\,R/8)$,
hit $\Upsilon\setminus B(z_1,5\,\eps\,R/8)$ (whose distance from
$\bar\gamma$ is at least $\eps\,R/5$)
or connect $\alpha_0$ to $\tchi_{z_1}$ inside $A$.
Consequently, the minimum on the right hand side in~\eref{e.assume} is presently
at least $\eps\,R/16$.
The lemma now implies that if we choose our current
$\delta=\delta(\eps,m,\upperhco)>0$ sufficiently small
and make sure that $R>1/\delta$, then
$$
\Pb{\ev C\md\ev K} \le m^{-1}/10\,,
$$
where $\ev C$ is
the event that there is a crossing of hexagons satisfying
$h<0$ between $\alpha_1$ and $\alpha_2$ inside $A$.
If there is no such crossing,
then also $\beg$ does not make such a crossing,
because $g\ge 0$ in $A$.

Likewise, we may define
$\tilde\alpha_1,\tilde\alpha_2,\tilde\alpha', \tilde A$ and $\tilde {\ev C}$,
when we replace $S'$ by $\tilde S$
in the above paragraph. The same argument shows that
$ \Pb{\tilde{\ev C}\md\ev K} \le m^{-1}/10$.

\begin{figure}
\SetLabels
\R(.195*.94)$\tilde A$\\
(.2*.3)$A$\\
(.7*.4)$S$\\
\endSetLabels
\centerline{\epsfxsize=2.8in%
\AffixLabels{%
\epsfbox{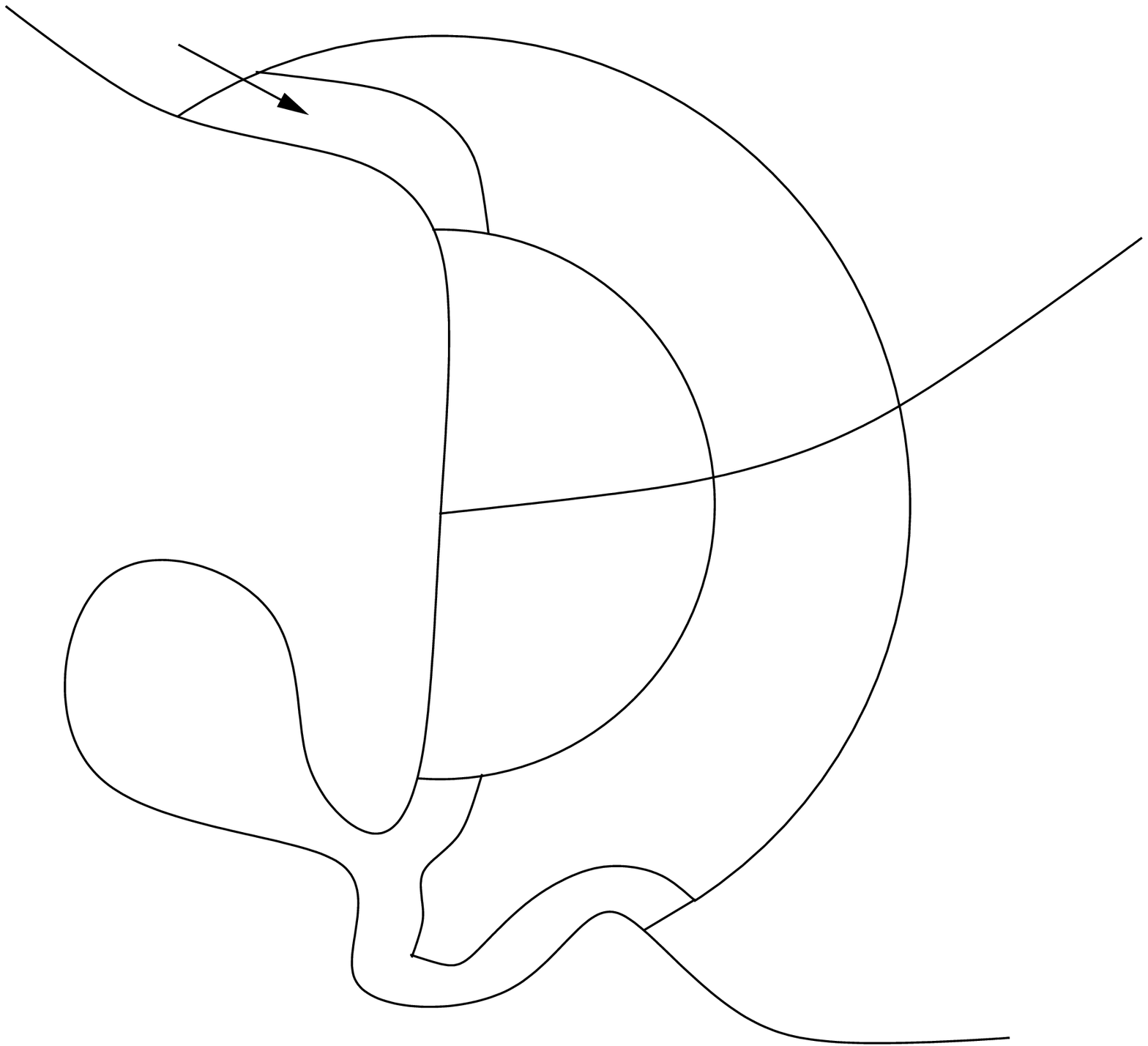}}
}
\begin{caption} {\label{f.noway}No way to penetrate.}
\end{caption} \end{figure}

Condition~\ref{i.A} of the definition of a barrier
implies that $\closure S\cup A\cup \tilde A\cup\chi_{z_1}$
separates $\Upsilon\cap A_{z_1}(\Upsilon)$ from all the endpoints of
the strands $\gamma_1,\dots,\gamma_k$.
Consequently, in order for $\beg$ to hit $\Upsilon\cap A_{z_1}(\Upsilon)$,
we must have $\beg\cap S\ne \emptyset$
or $\ev C\cup\tilde{\ev C}$.  See Figure~\ref{f.noway}
(also compare Figure~\ref{f.S}).
A similar argument applies for every other endpoint of a barrier.
Since $S\subset Z$, we conclude
\begin{equation}\label{e.Ups}
\Pb{\beg\cap(\cup Y)\ne \emptyset\md \ev K} \le 3/10\,.
\end{equation}

We are really more interested in $\hat\gamma$ than in $\beg$.
To do the translation, we will appeal to Lemma~\ref{l.distort}.
For this purpose, note that
$$
\|\nabla g\|_\infty=O(1)\,c_0\,\delta^{-1}\,R^{-1}=
O_{\eps,\upperhco,m}(1/R)
$$
and that $g$ is supported in a union of $m$ sets of diameter $O(R)$.
Consequently, $\|\nabla g\|_2 = O_{\eps,\upperhco,m}(1)$.
Let $U$ be the set of vertices in $\p\DD(\bar\gamma)$,
and let $h_U$ denote the restriction of $h$ to $U$.
By~\eref{e.Ups}, even if we further condition on
$h_U$, the left hand side stays bounded away from
one on an event whose probability is bounded away from zero, namely,
$$
\PB{
\Pb{\beg\cap(\cup Y) =\emptyset \md\ev K, h_U}>1/10\md\ev K}>1/10\,.
$$
Because $g=0$ on $U$,
we may apply Lemma~\ref{l.distort} and conclude that
for $h_U$ such that the inner inequality above holds,
$$
O_{\eps,\upperhco,m}(1)\,
\Pb{\hat\gamma\cap(\cup Y) =\emptyset \md\ev K, h_U} \ge 1\,.
$$
Since this set of $h_U$ has conditioned probability at least $1/10$,
the theorem follows.

It remains to remove the assumption that $R$ is larger than some fixed
constant $R_0=R_0(\eps,m,\upperhco)$.
Assume now that $R$ is bounded. It is not too hard to see that
the event $\ev H$ that $h(H)\in(0,1)$ for every hexagon meeting $\cup Y_+$
and $h\in(-1,0)$ for every hexagon meeting $\cup Y_-$ has
probability bounded below by (a rather small)
positive constant. This is proved by considering these $O(m\,R^2)$ hexagons
one by one.
On the event $\ev H$, we have $\hat\gamma\cap(\cup Y)=\emptyset$.
This completes the proof.
\QED

\begin{remark}\label{r.bdbarrier}
There is a corresponding analog of the Barriers Theorem~\ref{t.barrier}
in the case where the endpoints of the barriers are permitted to
land on $\p D$. In that case, it is necessary to assume that
$h_\p\ge -\lowerhco$ [respectively, $h_\p\le \lowerhco$]
on $\p A_z(\Upsilon)\cap V_\p$
if $\Upsilon\in Y_+$ [respectively, $Y_-$] and $z\in \p D\cap\Upsilon$,
where $\lowerhco=\lowerhco(\upperhco)>0$ is the constant given by Lemma~\ref{l.bdnarrows}.
We refrain from stating a complete formulation of this variant, though it will be useful.
The proof is the same, except that Lemma~\ref{l.bdnarrows} is used to deal with the narrows
near $\p D$, instead of the Narrows Lemma~\ref{l.narrows}.
\end{remark}

\subsection{Meeting of random walk and interface}\label{ss.enterRW}

We now need to further develop the basic setup and introduce some
notations. If $\alpha$ is a path in the hexagonal grid $\TG^*$, we
let $V(\alpha)$ denote the set of $\TG$-vertices adjacent to it. If
$\alpha$ is an arc of an oriented zero height interface of $h$, let
$V_+(\alpha)$ denote the vertices adjacent to it on its right hand
side, and let $V_-(\alpha)$ denote the vertices adjacent to it  on
its left hand side.

In addition to our previous assumptions~\iref{i.h}
and~\iref{i.D}, we will use the following setup:
\begin{enumerate}
\hitem{i.S}{(S)}
Let $\gamma$ denote the interface of $h$ from
$x_\p$ to $y_\p$.
 Let $\vv$ be some
 vertex of $\TG$ in $D$, and let $S$ be a simple random walk on
the vertices of $\TG$ started at $\vv$ that is independent from $h$.
Let $\tau$ be the first time $t$ such that $S_t\in \p \DD(\gamma)$.
\end{enumerate}

The point $\zz$ will play a special role. Essentially, we will be
interested in the configuration \lq\lq as viewed from $\zz$\rq\rq;
that is, in the coordinate system where $\zz$ is translated to $0$.
In order to eliminate too much additional notation, it will
be convenient to consider the event $\zz=0$ instead.
Let $\tau_0$ be the first $t$ such that $S_t=0$, and
let $\bS$ denote the reversed
walk $\bS_t:= S_{\tau_0-t}$, $t=0,1,2,\dots,\tau_0$.

For $\sig=0,1,\dots,5$,
let $e_\sig$ denote the edge $\bl[0,\exp(\pi\, i\,\sig/3)\br]$ of
the triangular grid $\TG$, and let $e^*_\sig$ denote the dual edge in $\TG^*$.
Let $\evZ$ denote the event
$\evZ:=\{\zz=0\}\cap\{e^*_\sig\subset\gamma\}$.

Fix some large $R$. Suppose that $\ball_R\subset D$ and
$\vv\notin\ball_R$. Let $\exter R\gamma$ denote the union of the
components of $\gamma\setminus \bal R$ containing $x_\p$ and $y_\p$.
(If $\gamma\cap\ball_R=\emptyset$, then $\exter R\gamma=\gamma$.) If
there is an interface of $\hth$ containing $e^*_\sig$, denote it by
$\hat \beta=\hat\beta^\sig$, and let $\beta=\beta_R^\sig$ be the
connected component of $\hat\beta\cap\ball_R$ that contains
$e^*_\sig$. Otherwise, set $\beta=\hat\beta=\emptyset$. Let $\inter
R\bS$ denote the part of $\bS$ up to the first exit of $\ball_R$,
and let $\exter R S$ denote the part of $S$ up to the first entry to
$\ball_R$.

Set $\oconf R:= (D,\p_+,h_\p,\vv,\exter {R}\gamma,\exter{R}S)$ and 
$\confR=\confR(\sig):= (\beta_R^\sig,\inter{R}\bS)$.
Our goal is to show that conditioned on $\evZ$, the distribution of $\beta$
does not depend strongly on $\oconf {4R}$.
(A precise version of this statement is given in Corollary~\ref{c.couple2} below.)
To this end we will use something like
\begin{equation}\label{e.resamp}
\PB{\confR=\vartheta\md \oconftR,\evZ}
=
\frac{\Pb{\confR=\vartheta\md \oconftR} \,\Pb{\evZ\md \confR=\vartheta,\oconftR}}
{\Pb{\evZ\md \oconftR}}\,.
\end{equation}
This equality is obtained by applying Bayes' formula to the measure
$\Pb{\cdot\md \oconftR}$.
The following lemma takes care of the first factor in the numerator
on the right hand side.

\begin{lemma}\label{l.aleq}
Assume~\iref{i.h}, \iref{i.D} and~\iref{i.S}.
There is a constant $c=c(\upperhco)>0$ and a function $p_R(\cdot)$
 such that if
$R>50$, $R'\in[5\,R/4,3\,R]$,
$D\supset \ball_{4R}$ and $\vv\notin\ball_{4R}$, then for all $\vartheta=(\tilde\beta,\tilde S)$ such
that $\tilde\beta\ne\emptyset$
$$
 c^{-1}\,p_R(\vartheta)\le
\Pb{\confR=\vartheta\md\gamma\cap \bal{4R}\ne\emptyset,\oconf{R'}}\le c\,p_R(\vartheta)\,.
$$
The function $p_R$ may depend on $R$ and $\vartheta$, but not on anything else
(in particular, not on $D$, $\vv$, $\oconf{R'}$ or $h_\p$).
\end{lemma}

\proof The corresponding statement with $\confR$ replaced by $\beta$, the first coordinate
of $\confR$, is an immediate consequence of Proposition~\ref{p.decoup}.

We assume $\vv\notin\ball_{4R}$.  The configuration $\oconf{R'}$
determines the first vertex, say $q$, inside $\ball_{R'}$ visited by
the random walk $S$. The continuation of the walk is just simple
random walk starting at $q$. Suppose that we had another such walk
starting at a vertex $q' \in\ball_{R'}$. It is easy to see that with
probability bounded away from zero the walk starting at $q$ visits
$q'$ before $0$. If that happens, we couple the continuation of the
walk to be the same as the walk which starts at $q'$ (otherwise, we
let them be independent). On the event that the walk started at $q$
hits $q'$ before $0$, the corresponding $\inter R\bS$ for both walks
will be the same. This proves the corresponding statement about the
second coordinate of $\confR$. Since the two coordinates are
independent given $\oconf{R'}$, the lemma follows. \QED

Proving an analogous result for the second factor in the
numerator of the right hand side of~\eref{e.resamp} will be
considerably more difficult. To this end,
we now define a measure of the {\bf quality} $\Qual=\Qual_R$ of
the configurations $\oconf R$ and $\conf R$.

If $\gamma\cap\bal R\ne\emptyset$,
let $x^R$ [respectively, $y^R$] denote the endpoint in
$\p\bal R$ of the component of $\exter R\gamma$ containing $x_\p$  [respectively, $y_\p$].
When $\vv\notin\bal R$, let $q^R$ denote the vertex in $\bal R$ first visited by $S$.
If $\gamma\cap\bal R=\emptyset$ or $\exter R S$
visits $\p \DD(\exter R\gamma)$, then set $\Qual(\oconf R)=0$.
Otherwise, define
$$
\Qual(\oconf R):=\frac{\dist(x,\exter R S)\wedge
\dist(y,\exter R S)\wedge \dist (q,\exter R\gamma)
\wedge |x-y|}{R}\wedge\frac 1 {100}\,,
$$
where $x=x^R$, $y=y^R$ and $q=q^R$.
This is a measure of the separation between the strands comprising $\oconf R$.
Similarly, define $\Qual (\conf R)$, as follows.
Suppose that $\vv\notin\bal R\subset D$ and
let $\hat q^R$ be the first vertex outside of $\bal R$ visited by $\bS$.
Fix an orientation of $e^*_\sig$.
If $\hat\beta\not\subset \bal R$, let $\hat x^R$ and $\hat y^R$ be the two
endpoints of the component of $\beta=\beta_R$ containing
$e^*_\sig$, chosen so that the orientation of the arc of $\beta$ from
$\hat x^R$ to $\hat y^R$ agrees with that of $e^*_\sig$.
If $\hat\beta\subset\bal R$ or
if $\bS$ visits any vertex in $\p\DD(\beta_R)\setminus\{0\}$,
then set $\Qual(\conf R)=0$.
Otherwise, set
$$
\Qual(\conf R):=\frac{\dist(\hat x,\inter R \bS)\wedge
\dist(\hat y,\inter R \bS)\wedge \dist (\hat q,\beta)
\wedge |\hat x-\hat y|}{R}\wedge\frac 1 {100}\,,
$$
where $\hat x=\hat x^R$, $\hat y=\hat y^R$ and $\hat q=\hat q^R$.

\begin{lemma}[Compatibility]\label{l.factor}
Assume~\iref{i.h}, \iref{i.D} and~\iref{i.S}.
For every $\eps>0$ there is a constant $c=c(\eps,\upperhco)>0$
such that
$$
c^{-1}\,1_{\{\Qual(\conf R)>\eps\}}\, 1_{\{\Qual(\oconf {R'})>\eps\}}
\le
\Pb{\evZ\md \confR,\oconf{R'}}\log R\le c
$$
holds whenever $R>c$, $5\,R>R'>9\,R/8$ and $\vv\notin\bal {6R}\subset D$.
\end{lemma}

\proof We start by proving the lower bound on $\Pb{\evZ\md
\confR,\oconf {R'}}$, which is the harder estimate. Assume that
$\Qual(\conf R)\wedge \Qual(\oconf {R'})>\eps$, $R>c>10^{10}/\eps$,
$5\,R>R'>9\,R/8$ and $\vv\notin\bal {6R}\subset D$. Let $\hat D_S$
denote the connected component of $\ball_R\setminus\beta_R$ that
intersects $\inter R\bS$ and let $D_S$ denote the connected
component of $D\setminus(\exter{R'}\gamma\cup\ball_{R'})$ that
contains $\vv$ and therefore $\exter{R'} S$. Note that the sign of
$h$ on vertices in $\hat D_S$ adjacent to $\beta$ is constant, as is
the sign of $h$ on vertices in $D_S$ adjacent to $\exter{R'}\gamma$.
Let $\ev D$ denote the event that these signs are the same, namely,
the sign of $h$ on vertices in $V(\exter {R'}\gamma)\cap D_S$ is the
same as on vertices in $V(\beta_R)\cap \hat D_S$. Using symmetry,
Proposition~\ref{p.decoup} immediately implies that $\Pb{\ev D\md
\oconf{R'},\confR}$ is bounded away from zero. (Although $\beta$ is
determined by $\confR$, its orientation as a subarc of an oriented
zero height interface of $h$ is not determined by $\confR$.)

\begin{figure}
\SetLabels
\R(.61*.5)$\bal R$\\
(.5*.87)$\bal {R'}$\\
\L(.91*.3)$\exter {R'}\gamma$\\
\R(.06*.93)$\exter {R'}\gamma$\\
\R(.335*.62)$\beta$\\
\L(.805*.04)$\exter {R'} S$\\
\R(.19*.45)$\inter {R} \bS$\\
\T(.74*.55)$A^b$\\
\T(.24*.74)$A^a$\\
\endSetLabels
\centerline{\epsfysize=2.9in%
\AffixLabels{%
\epsfbox{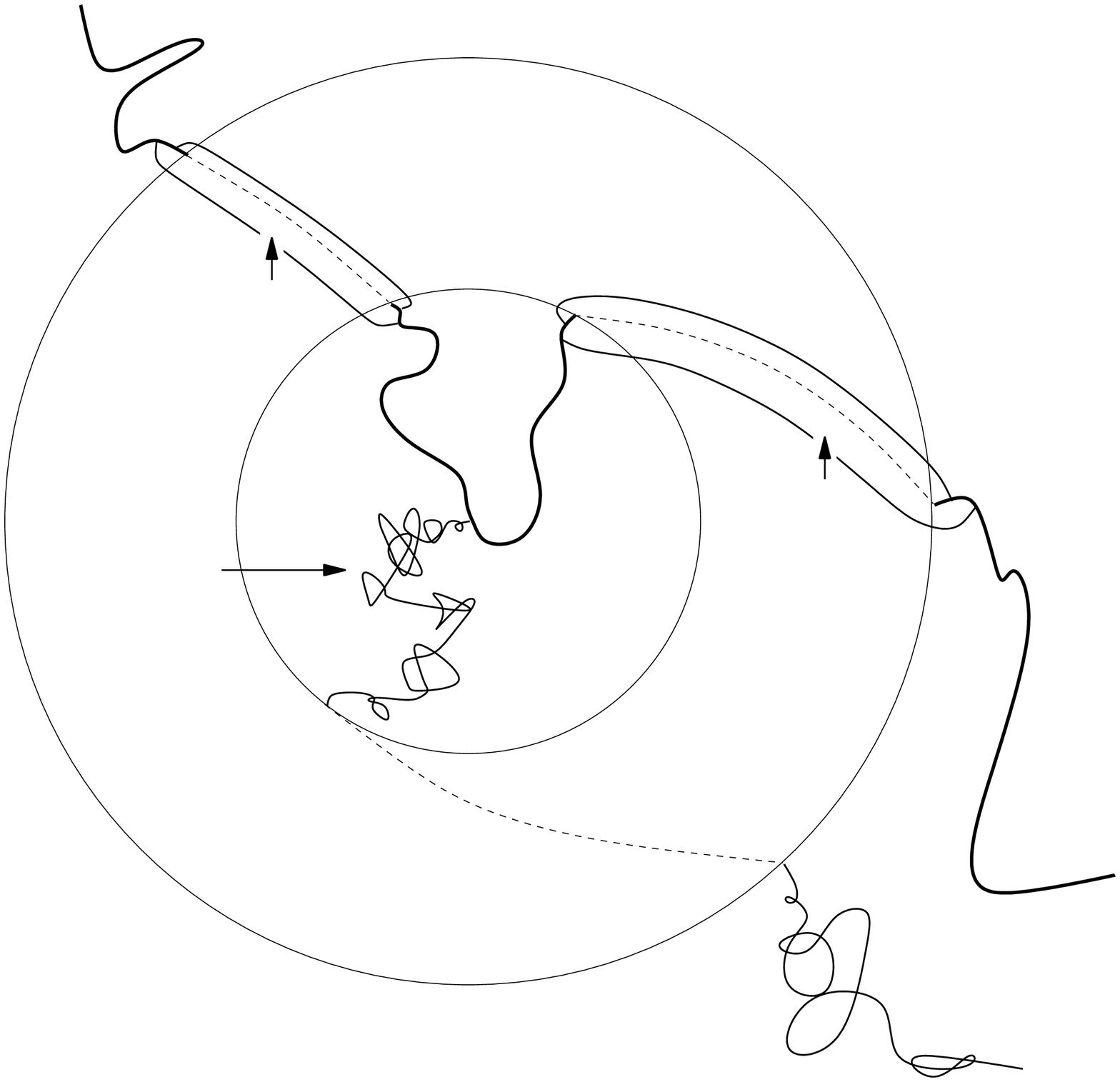}}
}
\begin{caption} {\label{f.paths}The construction of the barriers.}
\end{caption} \end{figure}

We now construct some barriers, as illustrated in Figure~\ref{f.paths}.
Let $a$ be the initial point  of the arc $\p\bal {R'} \cap \closure D_S$,
when the arc is oriented counterclockwise around $\bal {R'}$,
and let $b$ be the other endpoint of this arc.
Likewise,
let $\hat a$ be the initial point of the arc $\p\bal R\cap\closure {\hat D}_S$,
when the arc is oriented counterclockwise around $\bal R$,
and let $\hat b$ be the other endpoint of this arc.
Note that $\{\hat a,\hat b\}=\{\hat x^R,\hat y^R\}$ and
$\{a,b\}=\{x^{R'},y^{R'}\}$.

We now describe a  path $\patha$ connecting $a$ to $\hat a$ and
a path $\pathb$ connecting
$b$ to $\hat b$ and a path $\pathq$ connecting $\hat q^R$ to $q^{R'}$ such
that these paths do not come too close to each other.
For example, if the arguments $\arg a, \arg b, \arg q^{R'}, \arg \hat a,\arg \hat b$
and $\arg \hat q^R$
are chosen so that
$\arg \hat a<\arg \hat q^R<\arg \hat b<\arg \hat a+2\pi$,
$\arg a<\arg q^R<\arg b<\arg a+2\pi$,
and $|\arg a-\arg\hat a|\le \pi$,
then we may take
$\patha$ to be defined in polar coordinates by
$\theta = s+t\,r$, with $s$ and $t$ chosen
so that $a$ and $\hat a$ are on the path, and similarly
for $\pathb$ and $\pathq$.
It is easy to check that our assumptions guarantee
that the distance between any two of these paths is at least $c_1\,\eps\,R$ for
some constant $c_1>0$. Set $\eps' = c_1\,\eps/10$.

Let $D_+$ be the connected component of $D\setminus(\exter{R'}\gamma\cup\beta_R\cup\patha\cup\pathb)$
that contains $V_+(\exter{R'}\gamma)$,
and let $D_-$ be the other connected component.
Let $\alpha_a'$ be the connected component of $\{z\in D_+:\dist(z,\patha)=\eps'\,R\}$
that meets the circle $\p B\bl(0,(R+R')/2\br)$.
Note that $\alpha_a'$ is a simple path and
intersects the circle $\p B\bl(0,(R+R')/2\br)$ at one point.
Let $a_0$ denote that point.
We want to construct a pertubation of $\alpha_a'$, which will be some
$(\eps'',12\,R)$-barrier, with $\eps''$ not much smaller than $\eps'$.
Let $a_1$ be the closest point to $a_0$ along $\alpha_a'$ such that the distance
from $a_1$ to $\p\DD(\exter {R'}\gamma)$ is $\eps'\,R/10$,
and let $\hat a_1$ be the closest point to $a_0$ along $\alpha_a'$ such that
the distance from $\hat a_1$ to $\p\DD(\beta_R)$ is $\eps'\,R/10$.
Let $\alpha_+^a$ be the path which is the union of the arc of
$\alpha_a'$ connecting $a_1$ and $\hat a_1$ together with a
shortest line segment connecting $a_1$ to $\p\DD(\exter{R'}\gamma)$ and
a shortest line segment connecting $\hat a_1$ to $\p\DD(\beta_R)$.

We claim that $\alpha_+^a$ is an $(\eps'',12\,R)$-barrier
for the configuration $(D,\exter {R'}\gamma\cup\beta_R)$
with $\eps''=\eps'/1000$.
Indeed, Conditions~\ref{i.diam}, \ref{i.dist}, \ref{i.away} and~\ref{i.A}.(b)
in the definition of the barrier clearly
hold.
To verify Condition~\ref{i.A}.(a), let
$z_1$ be the endpoint of $\alpha_+^a$ on $\p\DD(\exter{R'}\gamma)$,
and let $z_1'$ and $z_1''$ be the two endpoints
of $\achi_{z_1}$ on $\p\DD(\exter{R'}\gamma)$.
Consider the simple arc $\tchi$ connecting
$z_1'$ to $z_1''$ in $\p \DD(\exter{R'}\gamma)$.
By the Jordan curve theorem, $\tchi\cup\achi_{z_1}$
separates the plane into two connected components.
Since $\alpha_+^a$ crosses $\achi_{z_1}$
it follows that the part of $\alpha_+^a$ inside $\ball_{R'}$
is outside of $A_{z_1}$, and thus the endpoints
$x^{R'},y^{R'}$ and also $\beta_R$ are all outside $A_{z_1}$.
It follows that $\p A_{z_1}=\achi_{z_1}\cup\tchi$,
as required. A similar argument applies
near the endpoint of $\alpha_+^a$ on $\p\DD(\beta_R)$.
Thus, $\alpha_+^a$ is indeed an $(\eps'',R')$-barrier.
Note also that the above easily implies that $\p A_{z_1}\subset D_+$.
This will be useful below when we apply the Barriers Theorem~\ref{t.barrier}.

We similarly construct a path $\alpha_-^a$ in $D_-$
close to $\patha$.
Likewise, we construct barriers $\alpha_+^b$ and
$\alpha_-^b$ near the path $\pathb$.
The construction is the same, except that we replace
$\patha$ by $\pathb$.

On the event $\ev D$, we may apply the Barriers Theorem~\ref{t.barrier}
with $Y_+=\{\alpha_+^a,\alpha_+^b\}$,
$Y_-=\{\alpha_-^a,\alpha_-^b\}$,
$Y=Y_+\cup Y_-$
and $\bar\gamma=\exter{R'}\gamma\cup\beta_R$.
(Here we use the assumption that $R>c$.)
Note that conditioning on $\confR$, $\oconftRp$ and $\ev D$,
amounts to conditioning on $\ev K$ in the theorem
and on the behavior of $\inter R\bS\cup\exter{R'}S$,
which is anyway independent from $h$.
Therefore, there is a $p=p(\eps,\upperhco)>0$ such that
$$
\Pb{\ev Y\md \confR,\oconftRp,\ev D}
\ge p\,,
$$
where $\ev Y$ denotes the event
$$
\ev Y:=\{(\gamma\setminus\bar\gamma)\cap
(\cup Y)=\emptyset\}\,.
$$
Let $A^a$ be the connected component of
$\DD(\exter{R'}\gamma\cup \beta_R)\setminus (\alpha_+^a\cup\alpha_-^a)$
that contains $\patha$, and let
$A^b$ be the connected component of
$\DD(\exter{R'}\gamma\cup \beta_R)\setminus (\alpha_+^b\cup\alpha_-^b)$
that contains $\pathb$.
Again, using the Jordan curve theorem,
it is easy to verify that $A^a\cap A^b=\emptyset$.
On the event $\ev Y$, there is
no other choice for the strand of $\gamma$ extending
$\exter{R'}\gamma$ at $a$, but to be confined to $A^a$
until it hooks up with $\beta$ at $\hat a$,
since every other exit from $A^a$ is blocked.
Consequently, on $\ev Y$,
we have $\gamma\supset \beta$.
A similar argument applies to $A^b$, and we get
$$
\beta\subset\gamma\subset \exter{R'}\gamma\cup\beta\cup A^a\cup A^b
\qquad\text{on }\ev Y
\,.
$$

\medskip
We now turn to the random walk $S$.
For $\evZ$ to hold, we must make sure that
$\{S_t:t<\tau_0\}$ does not meet any vertex neighboring
with $\gamma$.
First consider $\exter{R'}S$. Note that $\exter{R'}S$ does
not intersect $\p A^a$, because $\p A^a$ is
contained in $\ball_{R'}\cup B(a,2\,\eps'\,R)$,
and we are assuming that $\Qual(\oconf {R'})>\eps$.
(Recall, $\{a,b\}=\{x^{R'},y^{R'}\}$.)
Thus, $\exter{R'}S\cap A^a=\emptyset$,
and we may also conclude that $\exter{R'}S$ does
not visit any vertex adjacent to $A^a$ when $R$ is
large.
Similar arguments apply to $A^b$ and to  $\inter R\bS$.
Thus, $\inter R\bS\cup\exter{R'}S$ does not visit any vertex
adjacent to $\gamma$ on the event $\ev Y\cap \{
\Qual(\oconf {R'})>\eps,\,\Qual(\conf R)>\eps\}$.

Now let $S^*$ be the walk $S$ from the first time it visits
$q^{R'}$ until the first time it visits $\hat q^R$.
Then, conditioned on $\confR$ and $\oconf{R'}$,
$S^*$ is just simple random walk started at $q^{R'}$ conditioned
to hit $\hat q^R$ before hitting $0$.
Let $A^q$ denote the $\eps'\,R$-neighborhood of $\pathq$.
Clearly, $\dist\bl(A^q,A^a\cup A^b\br)\ge\eps'\,R$.
The probability that $S^*$ gets within distance $\eps'\,R/4$ of
$\hat q^R$ before exiting $A^q$ is at least some (perhaps small)
positive constant depending only on $\eps'$ (and hence on $\eps$).
Conditional on this event, the probability that
$S^*$ visits $\hat q^R$ before exiting $A^q$ is
within a constant multiple of $1/\log (\eps'\,R)$,
by~\eref{e.hitprob}.
Now let $S^{**}$ be the walk $S$ from the first visit of
$\hat q^R$ to the last visit of $\hat q^R$ before time $\tau_0$.
Note that $S^{**}$ and $\inter R \bS$ are independent given $\hat q^R$.
Thus, given $\oconf {R'},\confR,\ev D$ and $S^*$,
 we may sample $S^{**}$ by starting a random walk from $\hat q^R$,
stopping when it hits $0$, and then removing the part of that
walk after the last visit to $\hat q^R$.
When the latter walk first gets to distance
$\eps'\,R/2$ from $\hat q^R$, it has probability bounded away
from zero (by a constant depending only on $\eps'$)
to hit $0$ before $\hat q^R$.
(This follows, for example, from Lemma~\ref{l.harnack}
applied to the function giving for every vertex
the probability to hit $0$ before $\hat q^R$ for a random
walk started at that vertex.)
 Thus,
conditioned on $(S^*, \confR, \oconf {R'})$, with probability bounded
away from $0$, $S^{**}\subset A^q$.
Since
$$
\evZ\supset
\{S^*\cup S^{**}\subset A^q\}\cap\ev Y\,,
$$
we conclude that
$$
O_{\eps,\upperhco}(1)\,\Pb{\evZ\md \oconf {R'},\confR,\ev D}\ge 1/\log(R)\,.
$$
Above, we have argued that
$\Pb{\ev D\md \oconf {R'},\confR}$ is bounded away from zero,
and so we conclude that the lower bound estimate in the proposition
holds.

\medskip

It remains to prove the upper bound.
Conditional on $\gamma$, $\inter R\bS$ and $\exter{R'}S$, the probability
that $S^*$ (as defined in the proof of the lower bound)
hits $\hat q^R$ before hitting $\gamma$ is clearly
$O(1)/\log R$, since the conditional law of $S^*$ is that of random walk
started at $q^{R'}$ and conditioned to hit $\hat q^R$ before $0$,
and the probability that ordinary random walk started at $q^{R'}$
hits $\hat q^R$ before $0$ is bounded away from zero.
The upper bound now follows, and the proof is complete.
\QED

To make the previous lemma useful, we will need to argue that
configurations with quality bigger than $\eps$ are not too rare,
in an appropriate sense. This is achieved by the following lemma.

\begin{lemma}[Separation]\label{l.separation}
Assume~\iref{i.h}, \iref{i.D} and~\iref{i.S}.
Let $p<1$.
There is some constant $c=c(p,\upperhco)>0$ such that
if $R>1/c$ and $\vv\notin \bal {6R} \subset D$, then
\begin{equation}
\label{e.sepa}
\Pb{\Qual(\oconf {3R})\wedge \Qual(\conf{2R})>c\md\conf R, \oconf {4R},\evZ}>p,
\end{equation}
provided that $\Pb{\evZ\md\conf R, \oconf {4R}}>0$.
\end{lemma}

The proof of this lemma is modeled after Lawler's Separation Lemma
for Brownian motions from~\cite[Lemma 4.2]{\LawlerStrict}.

\begin{figure}
\SetLabels
(.64*.6)$\exter{r_j}\gamma$\\
\R(-0.01*.82)$\exter{r_j}S$\\
\endSetLabels
\centerline{\epsfysize=2.2in%
\AffixLabels{%
\epsfbox{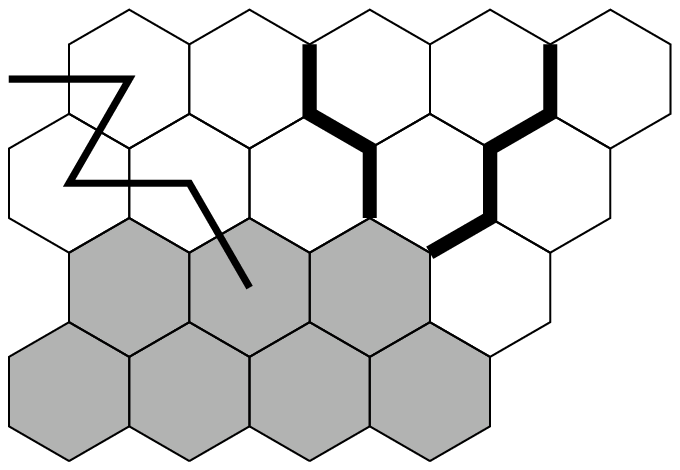}}
}
\begin{caption} {\label{f.discreteQ}An example for $\Qual_j=0\ne\Qual(\oconf{r_j})$.
The set $\bal {r_j}$ is shaded.}
\end{caption} \end{figure}

\proof
To keep the notations simple,
we start by proving a simpler version of the lemma, where we also assume
that $\Qual(\conf R)\ge 1/100$, say (and therefore
$\Qual(\conf R)=1/100$), and we prove
\begin{equation}
\label{e.qtar}
\Pb{\Qual(\oconf {3R})>c\md\conf R, \oconf {4R},\evZ}>p\,.
\end{equation}
We define inductively a random sequence $r_0,r_1,\dots$ as follows.
Set $r_0:=4\,R$.
Suppose that $r_j$ is defined.
Set
$$
\Qual_j:=\begin{cases}\Qual(\oconf {r_j}) &\text{if }
\Pb{\evZ \md \oconf {r_j},\conf R}>0,\\
0&\text{otherwise},
\end{cases}
$$
and $r_{j+1}:=(1-10\,Q_j)\,r_j $.
Note that $\Qual_j=0$ implies $ \Pb{\evZ \md \oconf {r_j},\conf R}=0$.
(An example showing that $\Qual_j=0\ne\Qual(\oconf{r_j})$ is possible is
given in Figure~\ref{f.discreteQ}.
Such a situation can only occur when $|x^{r_j}-y^{r_j}|=O(1)$.)
Also note that
$\Pb{\evZ \md \oconf {r_j},\conf R}>0$ if and only if
there are paths $\gamma_*$ and $S_*$ satisfying the following:
(1) $\gamma_*$ is a simple $\TG^*$-path in $\closure D$ containing $\beta^R$ and
$\exter {r_j}\gamma$,
(2) $S_*$ is a $\TG$-path in $D$ containing $\exter {r_j} S$ and the reversal of
$\inter R \bS$ and
(3) $S_*$ does not visit any vertex in $\p\DD(\gamma_*)$, except for $0$.

\begin{figure}
\SetLabels
\R(.14*.6)$\p\bal {r_j}$\\
(.6*.9)$\bal{r_{j+1}}$\\
\R(.1*.28)$\exter{r_j} S$\\
(.52*.14)$\exter{r_j}\gamma$\\
\endSetLabels
\centerline{\epsfysize=2.3in%
\AffixLabels{%
\epsfbox{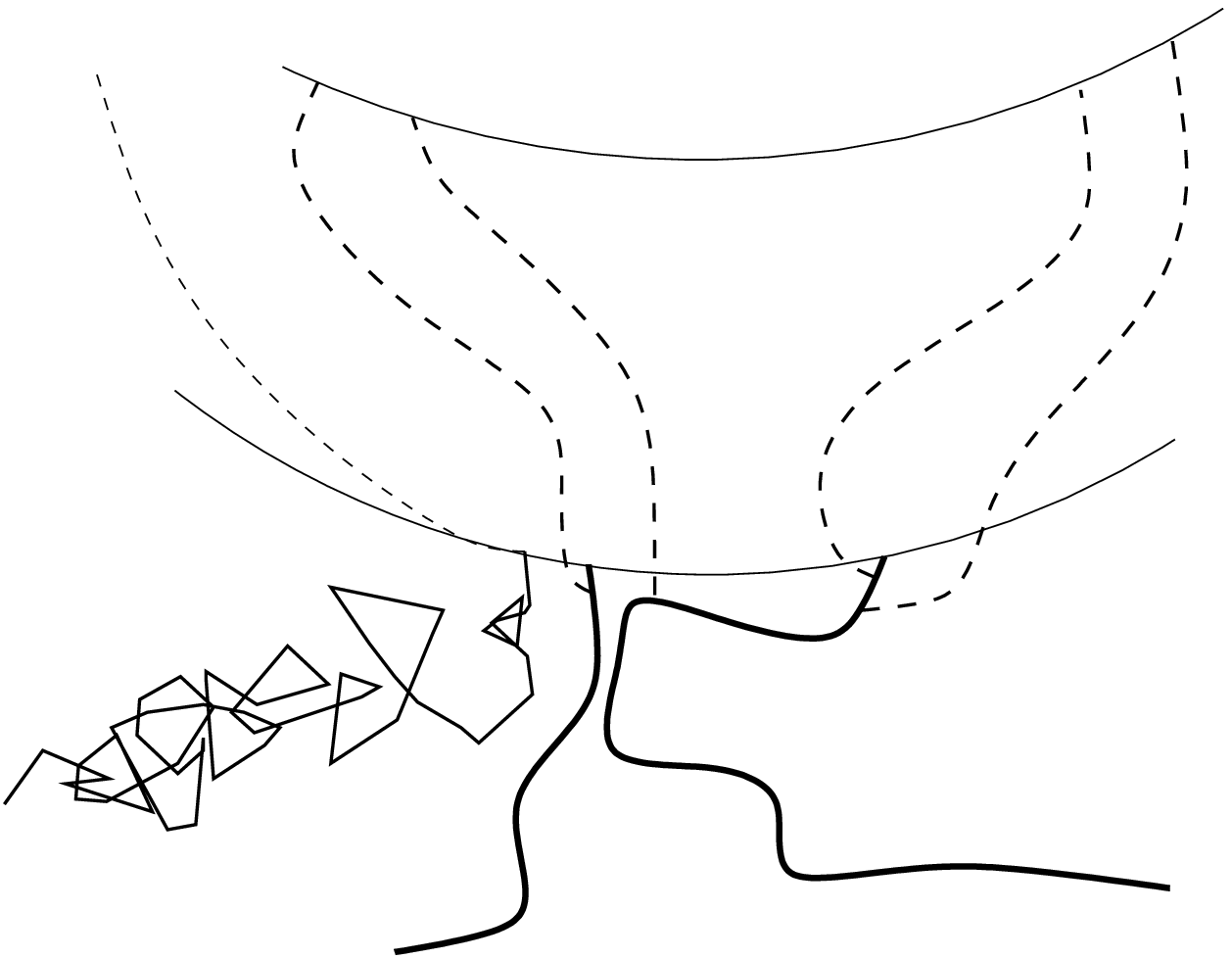}}
}
\begin{caption} {\label{f.improve}The construction of the barriers giving~\eref{e.improve}.}
\end{caption} \end{figure}

We claim that  for every $j\in\N$
\begin{equation}\label{e.improve}
\PB{\Qual_{j+1}\ge (2\,\Qual_j)\wedge (1/{100})\md \oconf {r_j},r_{j}>2\,R,\conf R}>c_0
\end{equation}
for some constant $c_0=c_0(\upperhco)>0$. Clearly it suffices to
prove this in the case $\Qual_j>0$. The gap between $\p \bal {r_j}$
and $\p\bal {r_{j+1}}$ is larger than but comparable to
$5\,r_j\,\Qual_j$. Note also that because $\bal {r_j}$ is a union of
$\TG^*$-hexagons, $r_j\,\Qual_j\ge 1/\sqrt 3$. If $r_j\,\Qual_j>20$,
say, then, we can easily use barriers as in the proof of the
Compatibility Lemma~\ref{l.factor} (see Figure~\ref{f.improve}) to
direct and separate the two strands of $\exter {r_{j+1}}
\gamma\setminus\exter {r_j}\gamma$ and the walk $\exter {r_{j+1}}
S\setminus\exter {r_j} S$ so as to obtain~\eref{e.improve}.

\begin{figure}
\SetLabels
(.70*.39)$\eta_a$\\
\L(.27*.33)$\eta_b$\\
\T(0.04*.28)$\p\ball{r_j}$\\
\B(0.03*.55)$\p\ball{r_{j+1}}$\\
\L(.51*.24)$S$\\
\endSetLabels
\centerline{\epsfysize=2.3in%
\AffixLabels{%
\epsfbox{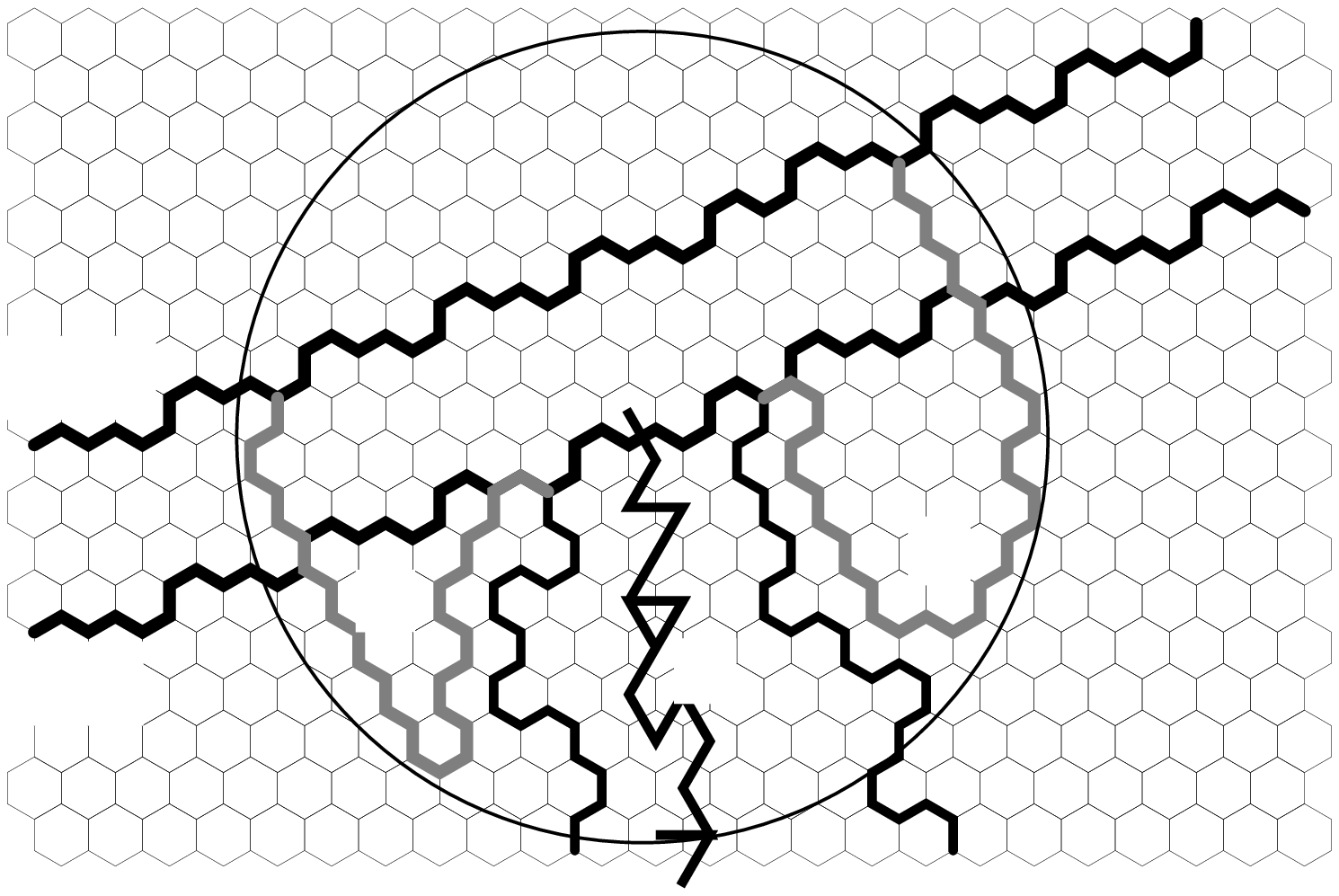}}
}
\begin{caption} {\label{f.spreadout} The interfaces spreading out.}
\end{caption} \end{figure}

Now assume that $r_j\,\Qual_j\le 20$. In this case the discrete
structure of the lattice is \lq\lq visible\rq\rq. Let $q'$ be the
point on $\p\bal{r_j}$ crossed by $\exter{r_j}S$ in its last step,
and let $\alpha$ be a longest arc among the three connected
components of $\p\bal{r_j}\setminus \{q',x^{r_j},y^{r_j}\}$. Suppose
first that $q'$ is not an endpoint of $\alpha$. Let $\alpha$ be
oriented counterclockwise around $\bal{r_j}$, and let $a$ and $b$ be
the initial and terminal points of $\alpha$, respectively.
Let $\eta_a$ [respectively $\eta_b$] denote the connected
component of $\exter {r_{j+1}}\gamma\setminus\exter {r_j}\gamma$
that has $a$ [respectively, $b$] as an endpoint.
Assume that $|q'-a|<500$ and $|q'-b|<500$.
Consider the event $\ev X$ that $\eta_a$ goes as far to the right as
possible subject to the conditions that it remains inside
$B(q',550)$ and avoids $\exter {r_j}\gamma$ and that $\eta_b$ goes
as far to the left as possible subject to the conditions that it
remains inside $B(q',550)$ and avoids $\exter{r_j}\gamma$. See
Figure~\ref{f.spreadout}. 
It is
easy to see that $\Pb{\evZ \md \oconf {r_j},\conf R}>0$ implies that
on $\ev X$ there is a simple $\TG$-path in $B(q',250)$ from $q'$ to
$\bal{r_{j+1}}$ that avoids $\p\DD(\exter{r_{j+1}}\gamma)$. Since
the number of edges traversed by these paths is bounded, it is easy
to see that the probability that $S$ follows the latter path and $\ev
X$ holds given $\oconf{r_j}$ and $\conf R$ satisfying the above
assumptions is bounded away from zero
(for $\ev X$, note that we can extend the interfaces one step at
a time and the probability for every specific step given the previous
ones is bounded away from zero). This gives~\eref{e.improve}
in this case. Similar, or simpler, arguments apply if one or more of
the assumptions $|q'-a|<500$ and $|q'-b|<500$ does not hold.

If $q'$ is an endpoint of $\alpha$ a similar argument may be used.
Suppose, for example, that $q'$ is the initial point of $\alpha$,
that $y^{r_j}$ is the other endpoint of $\alpha$ and that
$|q'-y^{r_j}|<300$. Then we may consider the possibility that the
connected component of $\exter{r_{j+1}}\gamma\setminus\exter {r_j}\gamma$
that has $y^{r_j}$ as an endpoint
goes as far to the left as possible subject to the
requirements that it stays inside $B(y^{r_j},5000)$ and avoids
hexagons containing vertices visited by $\exter{r_j}S$, and that the
connected component of $\exter{r_{j+1}}\gamma\setminus\exter {r_j}\gamma$
that has $x^{r_j}$ as an endpoint
goes as far to
the left as possible subject to the requirements that it stays
inside $B(x^{r_j},4000)$ and avoids the previous strand extending
$\exter{r_j}\gamma$ at $y^{r_j}$ and finally, the random walk avoids
$\p\DD(\exter{r_{j+1}}\gamma)$ and stays in $B(q',300)$ until it
hits $\bal {r_{j+1}}$. A similar argument applies if
$|q'-y^{r_j}|\ge 300$. This proves~\eref{e.improve}.

We now prove  for every $j\in\N$
\begin{equation}\label{e.die}
\Pb{\Qual_{j+2}=0\md
\oconf {r_j},r_{j}>2\,R,\conf R,\Qual_{j+1}<2\,\Qual_j}
\ge c_1\,,
\end{equation}
for some $c_1=c_1(\upperhco)>0$.
Let $x=x^{r_{j+1}}$, $y=y^{r_{j+1}}$ and $q=q^{r_{j+1}}$.
Let $S^\dagger $ be the
part of the walk $S$ from the first visit to $q^{r_{j}}$ up to the
first visit to $q$.
Note that $\dist(x,\exter {r_j} S)\ge 2\,r_j\,\Qual_j> 2\,r_{j+1}\,\Qual_j$,
if $\Qual_j>0$, and similarly for $y$. Thus, the event
$\Qual_{j+1}<2\,\Qual_j$ is the union of the following
five events
$$
\begin{aligned}
\ev M_0&:=\{\Qual_{j+1}=0\},\\
\ev M_1&:=\bl\{\dist(S^\dagger,\{x,y\})<(2\,\Qual_j\,r_{j+1})\wedge |x-y|\br\},\\
\ev M_2&:=\bl\{r_{j+1}\,\Qual_{j+1}=\dist(q,\exter {r_{j+1}}\gamma),\Qual_{j+1}<2\,\Qual_j\br\},\\
\ev M_3&:=\{0<|x-y|=r_{j+1}\,\Qual_{j+1},\Qual_{j+1}<2\,\Qual_j\},\\
\ev M_4&:=\{\Qual_{j+1}=1/100, \Qual_{j+1}<2\,\Qual_j\} .
\end{aligned}
$$
(In the definition of $\ev M_1$, $\dist (S^\dagger,\{x,y\})$ means the least distance
from a vertex visited by $S^\dagger$ to $x$ or $y$, of course.)
Clearly,
\begin{equation}\label{e.cases}
\OC(1)\,
\Pb{\Qual_{j+2}=0\md \oconf {r_j},\conf R,\ev M_i}
\ge 1
\end{equation}
holds for $i=0$.
The same is also true for $i=4$, because
$\Qual_{j+1}=1/100$ implies that the random walk started at $q$
has conditional probability bounded away from zero to hit $\exter {r_{j+1}}\gamma$
before $\p\bal{r_{j+2}}$.
A similar argument gives~\eref{e.cases} when $i=2$.

Now condition on $\ev M_1$, and let $v$ be the vertex first
visited by $S^\dagger$ that is at distance less than
$( 2\,\Qual_j\,r_{j+1})\wedge |x-y|$ from $\{x,y\}$.
Conditioned additionally on $\conf R$,
$\exter{r_{j+1}}\gamma$ and the walk $S^\dagger$ until
it hits $v$, there is clearly probability bounded away from zero
that $S^\dagger$ hits a vertex adjacent to $\exter{r_{j+1}}\gamma$ before
$\bal {r_{j+2}}$, and in this case we have $\Qual_{j+2}=0$.
Consequently,~\eref{e.cases}
also holds for $i=1$.

Now condition on $\ev M_3, \oconf{r_{j+1}}$ and $\conf R$.
Let $z$ be the midpoint of the segment $[x,y]$, and consider the
circle $\p B(z,2\,|x-y|)$.
We may build a barrier by using the connected component
of $\p B(z,2\,|x-y|)\setminus \p\DD(\exter {r_{j+1}}\gamma)$
that intersects $\bal {r_{j+1}}$ (and possibly perturbing it
slightly near its endpoints).
If $\gamma$ does not cross this barrier, then $\Qual_{j+2}=0$
holds. Thus, we get from the Barriers Theorem~\ref{t.barrier} that~\eref{e.cases} also holds
for $i=3$.
Since $\{\Qual_{j+1}<2\,\Qual_j\}=\bigcup_{i=0}^4 \ev M_i$,
and~\eref{e.cases} holds for $i=0,1,\dots,4$, it follows that~\eref{e.die}
holds as well.

Set $s_n:=2\,R\,\prod_{k=0}^{n-1} (1- 2^{-k}/10)^{-1}$.
It follows from the Compatibility Lemma~\ref{l.factor}
and our assumption that $\Qual(\conf R)\ge 1/100$ that
for any $j\in \N$
\begin{equation*}
\Pb{\evZ \md\oconf{r_j},\conf R}
\ge \frac{c_2} {\log R}\, 1_{\{\Qual_j\ge 1/100\}}\,1_{\{ r_j\ge s_0\}}\,,
\end{equation*}
for some $c_2=c_2(\upperhco)>0$.
An appeal to~\eref{e.improve} therefore implies
\begin{equation*}
\Pb{\evZ \md\oconf{r_j},\conf R}
\ge \frac{c_0\,c_2} {\log R}\, 1_{\{\Qual_j\ge 2^{-1}/100\}}\,1_{\{ r_j\ge s_1\}}\,.
\end{equation*}
Continuing inductively, we get for every $n\in\N$
\begin{equation*}
\Pb{\evZ \md\oconf{r_j},\conf R}
\ge \frac{c_0^n\,c_2} {\log R}\, 1_{\{\Qual_j\ge 2^{-n}/100\}}\,1_{\{ r_j\ge s_n\}}\,.
\end{equation*}
Since
$$
\sup_n s_n < { 2\,R}/\Bigl(1-\sum_{k=0}^\infty 2^{-k}/10\Bigr)=
\frac52\,R\,,
$$
we get for every $j\in\N$
\begin{equation}
\label{e.zlb}
\Pb{\evZ \md\oconf{r_j},\conf R}
\ge \frac{c_0^{-\log_2\Qual_j}\,c_2} {\log R}\, 1_{\{ r_j\ge \frac52\,R\}}
=
c_2\,\frac{\Qual_j^{-\log_2 c_0}} {\log R}\, 1_{\{ r_j\ge \frac52\,R\}}
\,.
\end{equation}

 From~\eref{e.die} we get that for every $j$ the conditional probability
that there is some $k>j$ such that $r_k>2\,R$,
$\Qual_k<2\,\Qual_j$ and $\Qual_{k+1}>0$ given $\oconf{r_j}$ and
$\conf R$ is at most $1-c_1$.
Let $m_n$ denote the number of $k\in\N$ such that
$r_k>3\,R$ and $\Qual_k\in (2^{-n},2^{1-n}]$.
Fix some $n\in\N$, and suppose that $\Pb{m_n>0\md\oconf {4R},\conf R}>0$.
On the event $m_n>0$, let $k_n$ be the first $k$ such that $\Qual_k\in(2^{-n},2^{1-n}]$.
By induction and~\eref{e.die} for every $m\in\N$
$$
\Pb{m_n> 2\,m\md m_n>0,\oconf {r_{k_n}},\conf R}\le (1-c_1)^m.
$$
An appeal to the upper bound in the Compatibility Lemma~\ref{l.factor} gives
$$
\Pb{\evZ,m_n> 2\,m\md m_n>0,\oconf {r_{k_n}},\conf R}\le c_3\, (1-c_1)^m/\log R.
$$
for some $c_3=c_3(\upperhco)$.
On the other hand,~\eref{e.zlb} gives
$$
\Pb{\evZ\md m_n>0,\oconf {r_{k_n}},\conf R}\ge
c_2\,c_0^n/\log R\,.
$$
Comparing the last two inequalities, we get
$$
\Pb{m_n> 2\,m\md\evZ,\oconf {r_{0}},\conf R}\le c_3\,c_0^{-n}\, (1-c_1)^m/c_2\,.
$$
In particular, there is an $n_0=n_0(\upperhco,p)\in\N$ and a
$c_4=c_4(\upperhco)$ such that
$$
\Pb{\exists n\ge n_0: m_n> c_4\,n \md\evZ,\oconf {r_{0}},\conf R}\le (1-p)/3\,.
$$
Let $n_1$ be the least integer larger than $3$ such that
$ \prod_{n=n_1}^{\infty} (1-10\cdot 2^{1-n})^{c_4n}> 7/8 $,
and let $n_2=n_1\vee n_0$.
On the event $\evZ\cap\bigcap_{n>n_2}\{m_n\le c_4\,n\}$
we must have some $j\in\N$ with $r_j>7\,R/2$ and $\Qual_j\ge 2^{-n_2}$
(because $r_{j+1}/r_j= 1-10\,\Qual_j$ and $n_2\ge n_1$).
Consequently,
\begin{equation}
\label{e.72}
\Pb{\exists j: \Qual_j> c_5,r_j>7\,R/2\md\evZ,\oconf {4R},\conf R}\ge 1-(1-p)/3
\end{equation}
holds for some $c_5=c_5(p,\upperhco)>0$.
Note that this almost achieves our goal of proving~\eref{e.qtar}. The difference
between~\eref{e.72} and~\eref{e.qtar} is that in the latter the radius $r$ at which
$\Qual(\oconf r)$ is bounded from below is variable.

Set $\ev A$ be the event that there is a $j\in\N$ with
$\Qual_j>c_5$ and $r_j>7\,R/2$, and on $\ev A$, let
$j_0$ denote the first such $j$.
We have from~\eref{e.zlb} that
\begin{equation}
\label{e.evZA}
\Pb{\evZ \md\ev A,\oconf{r_{j_0}},\conf R}
\ge
c_2\,c_5^{-\log_2 c_0}/ \log R
\,.
\end{equation}
Fix some $s>0$ small. We now argue that
\begin{equation}
\label{e.spow}
\Pb{\evZ \md\ev A,\oconf{r_{j_0}},\conf R,|x^{3R}-y^{3R}|<s\,R}
\le c_7\, s^{c_6}/\log R
\end{equation}
for some positive constants $c_7$ and $c_6$ depending only on $\upperhco$.
The argument is similar to the one given in the proof of the case $i=3$
in~\eref{e.cases}.
Let $z$ be the midpoint of the segment
$[x^{3R},y^{3R}]$. We construct a barrier as a perturbation of
the connected component of $\p B(z,2\,s\,R)\setminus\p\DD(\exter {3R}\gamma)$
that intersects $\bal {3R}$. If that barrier is hit by the extension
of $\exter {3R}\gamma$ (which happens with probability bounded away
from $1$), then we condition on the extension up to that barrier,
and construct another barrier at radius $4\,s\,R$, instead.
We continue in this manner, constructing  barriers at radii $2^n\,s\,R$
up to the least $n$ such that $2^n\,s>1/1000$, say.
Because the probability of avoiding the $n$'th barrier given
that the $(n-1)$'th barrier has been breached is bounded away from one,
we find that $\Pb{\gamma\supset\beta_R \md \ev A, \oconf{r_{j_0}},\conf R,|x^{3R}-y^{3R}|<s\,R}$
is bounded by a constant times some positive power of $s$.
The estimate~\eref{e.spow} follows by considering the behavior of $S$.

Suppose now that the random walk $S$ after its first hit to
$\bal {r_{j_0}}$ but before its first hit to
$\bal {3R}$ gets within distance $s\,R$ of $\exter{3R}\gamma$.
Then by Lemma Hit Near~\ref{l.hitnear}, conditional on $S$ up to the first
time this has happened and on $\ev A$, $\oconf{r_{j_0}}$, $\conf R$
and $\exter{3R}\gamma$, the conditional probability for
$S$ hitting $\bal {5\,R/2}$ before hitting $\p\DD(\exter {3R}\gamma)$
is at most $c_8\,s^\expo1$, for some universal constant $c_8$.
Thus, the conditional probability for $\evZ$ is at most
$c_8\,s^\expo1/\log R$.
Combining this with~\eref{e.spow}, one gets
$$
\Pb{\evZ \md\ev A,\oconf{r_{j_0}},\conf R,\Qual(\oconf {3R})<s}
\le (c_7\, s^{c_6}+c_8\,s^\expo1)/\log R\,.
$$
Comparison with~\eref{e.evZA} now gives
\begin{multline*}
\Pb{\Qual(\oconf{3R})<s\md \evZ, \ev A,\oconf{r_{j_0}},\conf R}
=
\frac
{\Pb{\Qual(\oconf{3R})<s,\evZ \md \ev A,\oconf{r_{j_0}},\conf R}}
{\Pb{\evZ \md \ev A,\oconf{r_{j_0}},\conf R}}
\\ 
\le
\frac
{\Pb{\evZ \md \ev A,\oconf{r_{j_0}},\conf R,\Qual(\oconf{3R})<s}}
{\Pb{\evZ \md \ev A,\oconf{r_{j_0}},\conf R}}
\le
\frac{c_7\, s^{c_6}+c_8\,s^\expo1} { c_2\,c_5^{-\log_2 c_0}}\,.
\end{multline*}
Thus, we obtain for all $s$ sufficiently small
$$
\Pb{\Qual(\oconf{3R})<s\md \evZ, \ev A,\oconf{4R},\conf R}
\le (1-p)/3\,.
$$
Taking~\eref{e.72} into account, this gives~\eref{e.qtar},
and completes the proof of the simplified case.

The argument in the general case proceeds as follows.
We define inductively two sequences $r_j$ and $\hat r_j$, starting
with $r_0=4\,R$ and $\hat r_0=R$.
At each step $j$,
we set
$$
\Qual_j:= \Qual(\oconf {r_j})\wedge \Qual(\conf {\hat r_j})
\wedge 1_{\{ \Ps{\evZ\md\oconf {r_j},\conf {\hat r_j}}>0\}}
$$
and take $r_{j+1}=(1-10\,\Qual_j)\,r_j$ and
$\hat r_{j+1}=(1+10\,\Qual_j)\,\hat r_j$.
The proof proceeds essentially as above.
The straightforward details are left to the reader.
\QED

\begin{corollary}\label{c.couple2}
There is a constant $c=c(\upperhco)>0$ such that the
following estimate holds.
Let $(D',\p'_+,\p'_-,h'_\p,h',\gamma',\vv')$ satisfy the same assumptions
as we have for $(D,\p_+,\p_-,h_\p,h,\gamma,\vv)$.
Let $R>c$, and assume that $\ball_{6R}\subset D'\cap D$ and
$\vv,\vv'\notin\ball_{6R}$.
Let $\confR',\oconf{4R}'$ and $\evZp$ be the objects corresponding to
$\confR,\oconf{4R}$ and $\evZ$ for the system in $D'$
(with the same $\sigma$, that is, $\sigma'=\sigma$).
Then
\begin{equation}
\label{e.couple2}
\Pb{\confR=\vartheta\md\evZ,\oconf{4R}}\le
c\,
\Pb{\confR'=\vartheta\md\evZp,\oconf{4R}'}
\end{equation}
holds for every $\vartheta$ and for every
$\oconf{4R}$ and $\oconf{4R}'$ satisfying
$\Pb{\evZ\md \oconf{4R}}>0$ and
$\Pb{\evZp\md \oconf{4R}'}>0$, respectively.
Consequently, under the same assumptions,
there exists a coupling of the conditional laws of $\confR$ and $\confR'$ such that
$$
\Pb{\confR=\confR'\md \evZ,\evZp,\oconf{4R},\oconf{4R}'}\ge 1/c\,.
$$
\end{corollary}

\proof
It is enough to prove the first claim, since the latter claim immediately follows.
Let $c'>0$ be the constant denoted as $c$ in the Separation Lemma~\ref{l.separation} with $p=1/2$.
Let $\ev Q$ denote the event $\Qual(\oconf {3R})\wedge\Qual(\conf{2R})\ge c'$.
Let $X$ be the collection of all $\theta$
such that $\conf {2R}=\theta$ is possible and $\Qual(\theta)\ge c'$,
and let $X_\vartheta$ be the collection of all $\theta\in X$
that are compatible with $\conf R=\vartheta$; that is,
such that $\{\conf {2R}=\theta\text{ and }\conf R=\vartheta\}$ is possible.
In the following, $f\approx g$ will mean that
$f/g$ is contained in $[1/c,c]$ for some constant $c=c(\upperhco)>0$.
By Lemma~\ref{l.separation} and the choice of $p$,
\begin{multline} \label{e.first}
\Pb{\confR=\vartheta\md\evZ,\oconf{4R}}
\approx
\Pb{\confR=\vartheta\md\ev Q,\evZ,\oconf{4R}}
\\
=
\sum_{\theta\in X_\vartheta}
\Pb{\conf{2R}=\theta\md\ev Q,\evZ,\oconf{4R}}
\end{multline}
Now if $\Qual(\oconf{3R})>c'$ and $\theta\in X$, then
\begin{align*}
&
\Pb{\conf{2R}=\theta\md\ev Q,\evZ,\oconf{3R}}
=
\frac{
\Pb{\conf{2R}=\theta,\ev Q,\evZ\md\oconf{3R}}
}
{
\Pb{\evZ,\ev Q\md\oconf{3R}}
}
\\&\qquad
=
\frac{
\Pb{\conf{2R}=\theta,\evZ\md\oconf{3R}}
}
{
\Pb{\evZ,\ev Q\md\oconf{3R}}
}
=
\frac{
\Pb{\conf{2R}=\theta\md\oconf{3R}}
\Pb{\evZ\md\conf{2R}=\theta,\oconf{3R}}
}
{
\Pb{\evZ,\ev Q\md\oconf{3R}}
}\,.
\end{align*}
We apply Lemma~\ref{l.aleq} to the first factor in the numerator
and Lemma~\ref{l.factor} to the second factor, and get
$$
\Pb{\conf{2R}=\theta\md\ev Q,\evZ,\oconf{3R}}
\approx
\frac{
p_{2R}(\theta)/\log R
}
{
\Pb{\evZ,\ev Q\md\oconf{3R}}
}
\,.
$$
The sum of the left hand side over all $\theta\in X$ is $1$.
Consequently,
$$
\Pb{\conf{2R}=\theta_0\md\ev Q,\evZ,\oconf{3R}}
\approx
\frac{p_{2R}(\theta_0)}
{\sum_{\theta\in X} p_{2R}(\theta)}\,.
$$
By taking expectation conditioned on $\ev Q,\evZ$ and $\oconf {4R}$, it
follows that the same relation holds when we replace
$\oconf{3R}$ by $\oconf{4R}$.
We now sum over $\theta_0\in X_\vartheta$ and invoke~\eref{e.first}, to obtain
$$
\Pb{\confR=\vartheta\md\evZ,\oconf{4R}}
\approx
\frac
{\sum_{\theta\in X_\vartheta} p_{2R}(\theta)}
{\sum_{\theta\in X} p_{2R}(\theta)}\,.
$$
This implies~\eref{e.couple2}, and completes the proof.
\QED

Our intermediate goal to show that the dependence between the local behavior near
$\zz$ and the global behavior far away is now accomplished.
Roughly, the next objective will be to show that it is unlikely that
$\gamma$ contains an arc with a very large diameter whose endpoints
are both relatively close to $\zz$.

For $R>r>0$ let $\ev J=\ev J(r,R)$ denote the event that
there are more than $2$ disjoint arcs of $\gamma$
connecting $\bal r$ and $\p\bal R$
or that $S$  exits $\bal R$ between the time it
first hits $\bal r$ and $\tau_0$.
Set $\evZa:=\bigcup_{\sig=0}^5\evZ$.
Our next objective is to show that conditioned on $\evZ$
or $\evZa$,
$\ev J(r,R)$ is unlikely if $R\gg r>0$.
More precisely, the claim is as follows.

\begin{lemma}\label{l.direct}
Assume~\iref{i.h}, \iref{i.D} and~\iref{i.S}.
For every $p>0$ there is some $a=a(p,\upperhco)>10$ such that if
$r>1$, $R>a\,r$, $\vv\notin\ball_{4R}\subset D$
and $\Pb{\evZ\md\oconf R,\conf r}>0$,
then
$$
\PB{\ev J(r,R)\md \oconf R,\conf r,\evZ}<p
$$
and also
$$
\PB{\ev J(r,R)\md \oconf R,\evZa}<p\,.
$$
\end{lemma}

One may first think that this
can be proved by repeating the argument in the proof of
the Compatibility Lemma~\ref{l.factor}.
The difficulty in carrying out this idea is that the sets
$A^a$ and $A^b$ described in the proof of Lemma~\ref{l.factor} may extend
beyond $\bal {sR'}$ for large $s$ 
(where $R'$ is as in that lemma) if there are more than two disjoint
arcs of $\exter {R'}\gamma$ connecting $\p\bal {R'+\eps'R}$ and
$\p\bal{sR'}$.

Since the proof of the lemma is a bit involved and somewhat indirect, we
take a few moments to give an overview of the strategy.
First, it is established that
under the conditioning the simple random walk $S$ is unlikely to
backtrack to $\p\bal R$ after hitting $\bal r$ and before $\tau_0$.
 Next, we identify a pair of
arcs $\alpha_1$ and $\alpha_2$ that are defined from
$\oconf{3r}$ each of which has one endpoint on $\exter {3r} S$
and the other on $\exter{3r}\gamma$.
A barriers argument is then used to show that with high
conditional probability $\gamma\setminus\exter{3r}\gamma$
does not hit $\alpha_1\cup\alpha_2$.
In this case, we see that $\alpha_1\cup\alpha_2$
have an alternative definition in terms of $\gamma$ and $S$,
which is in some sense more symmetric.
Next we define another pair of arcs $\tilde\alpha_1\cup\tilde\alpha_2$,
which have a similar definition as $\alpha_1\cup\alpha_2$,
except that they are defined from $\conf{R/3}$.
Again, the same barriers argument can be used
to show that with high conditional probability these
arcs are not visited by $\gamma\setminus\beta_{R/3}$.
In this case, these arcs have a more symmetric definition,
which leads us to conclude that with high conditional probability
$\alpha_1\cup\alpha_2=\tilde\alpha_1\cup\tilde\alpha_2$.
This is then used to establish that the endpoints
of these arcs on $\gamma$ belong to $\exter{3r}\gamma$
as well as $\beta_{R/3}$. Next, we prove that these
two endpoints belong to different connected components
of $\gamma\setminus e_\sig^*$. This then implies that
each of the two strands of $\exter{3r}\gamma$ merges with
$\beta_{R/3}$, which implies that there are no more than
two disjoint crossings between $\p \bal R$ and $\bal {3r}$
in $\gamma$.

\proof
The second claimed inequality with $p$ replaced by $6p$ follows from the first
inequality and taking conditional expectation,
since $\evZa=\bigcup_{\sig=0}^5\evZ$. Thus, we only need to prove the first.
By the Separation Lemma~\ref{l.separation}, there is a constant
$c_0=c_0(\upperhco,p)>0$ such that
\begin{equation}\label{e.goodq}
\Pb{\Qual(\oconf{3R/4})\wedge\Qual(\conf {R/2})\wedge \Qual(\oconf{3r})\wedge\Qual(\conf{2r})\le c_0\md
\oconf R,\conf r,\evZ}<p/10\,.
\end{equation}

Let $r'$ be in the range $\bl[\sqrt {r\,R},\sqrt{r\,R}+1\br]$, chosen so that
the circle $\p B(0,r')$ does not contain any $\TG$ vertices nor $\TG^*$ vertices.
Let $S_*$ denote the part of the walk $S$ from its first visit to
$q^{3r}$ until its last visit to $\hat q^r$ prior to $\tau_0$.
Conditional on $\oconf {3r}$ and $\conf r$,
the probability that $S_*$ exits $\bal {r'/3}$ without hitting
$\exter{3r}\gamma$ decays to zero as $a\to\infty$
(by Lemma~\ref{l.hitnear}).
On the event that this happens, let
$S_{**}$ be the initial segment of the walk $S_*$
until it exits $\bal {r'/3}$. By the proof of the upper bound in the Compatibility Lemma~\ref{l.factor},
$$
\Pb{\evZ\md \conf {2r},\oconf{3r},S_*\not\subset\bal{r'/3},S_{**}}
\le O(1)/\log r\,.
$$
By the lower bound in that lemma, on the event $\Qual(\oconf{3r})\wedge\Qual(\conf{2r})\ge c_0$
we have
$$
O_{p,\upperhco}(1)\,\Pb{\evZ\md \conf {2r},\oconf{3r}} \ge 1/\log r\,.
$$
Consequently, $a$ may be chosen sufficiently large so that
$$
1_{\{ \Qual(\oconf{3r})\wedge\Qual(\conf{2r})\ge c_0\}}\,
\Pb{S_*\not\subset\bal{r'/3}\md \conf{2r},\oconf{3r},\evZ}<p/10\,.
$$
Hence,~\eref{e.goodq} implies
\begin{equation}\label{e.Sdirect}
\Pb{S_*\not\subset\bal{r'/3}\md \conf{r},\oconf{R},\evZ}<p/5\,.
\end{equation}

Let $S'$ be the path traced by $\exter{3r}S$ from the last
time in which $\exter {3r} S$ was outside $\bal{R/3}$
until its terminal point $q^{3r}\in\bal{3r}$.
Let $\tilde S$ be the path traced by $\inter{R/3}\bS$ from the last
time in which $\inter {R/3}\bS$ was inside $\bal{3r}$
until its terminal point $\hat q^{R/3}$.
Observe that $S'$ is $\oconf{3r}$-measurable, $\tilde S$ is
$\conf{R/3}$-measurable, and when
$S_*\subset \bal{r'/3}$, we have $S'=\tilde S$ as unoriented paths.

Observe that there is a connected component
$\alpha$ of $\p B(0,r')\setminus S'$ such that
$S'\cup\alpha$ separates $\p D$ from $\bal{3r}$.
We fix such an $\alpha$, and if there is more than one possible
choice, we choose one in a way which depends only on $S'$.

On the event $\Qual(\oconf{3r})>0$ each strand of
$\exter{3r}\gamma$ connects $\p D$ with $\bal{3r}$,
and hence $\exter{3r}\gamma$ intersects $\alpha$.
Thus, there are precisely two connected components of $\alpha\cap\DD(\exter{3r}\gamma)$ which have
one endpoint in $S'$ and the other in $\p\DD(\exter{3r}\gamma)$.
Let $\alpha_1$ and $\alpha_2$ be these two arcs.

We now argue that
\begin{equation}\label{e.staynear}
1_{\{ \Qual(3r)\wedge\Qual(2r)\ge c_0\}}\,
\Pb{\gamma\cap \alpha_1\ne\emptyset\md \evZ,\oconf{3r},\conf{2r}}
< p/10
\end{equation}
if $a$ is sufficiently large.

The basic idea of the proof of~\eref{e.staynear} is to construct a sequence of barriers
separating $\bal {3r}$ from $\alpha_1$ in $\DD(\exter{3r}\gamma)$
such that if $\gamma$ hits a barrier
in the sequence, the conditional probability that it will hit the
next barrier is bounded away from $1$.

\begin{figure}
\SetLabels
(.5*.58)$D'$\\
(.63*.48)$\bal {3r}$\\
\L(.57*.33)$\bbeta_2$\\
\L(.61*.66)$\bbeta_1$\\
\T(.65*.22)$\p B(0,r')$\\
\R(.39*.22)$\alpha_1$\\
\endSetLabels
\centerline{\epsfysize=2.3in%
\AffixLabels{%
\epsfbox{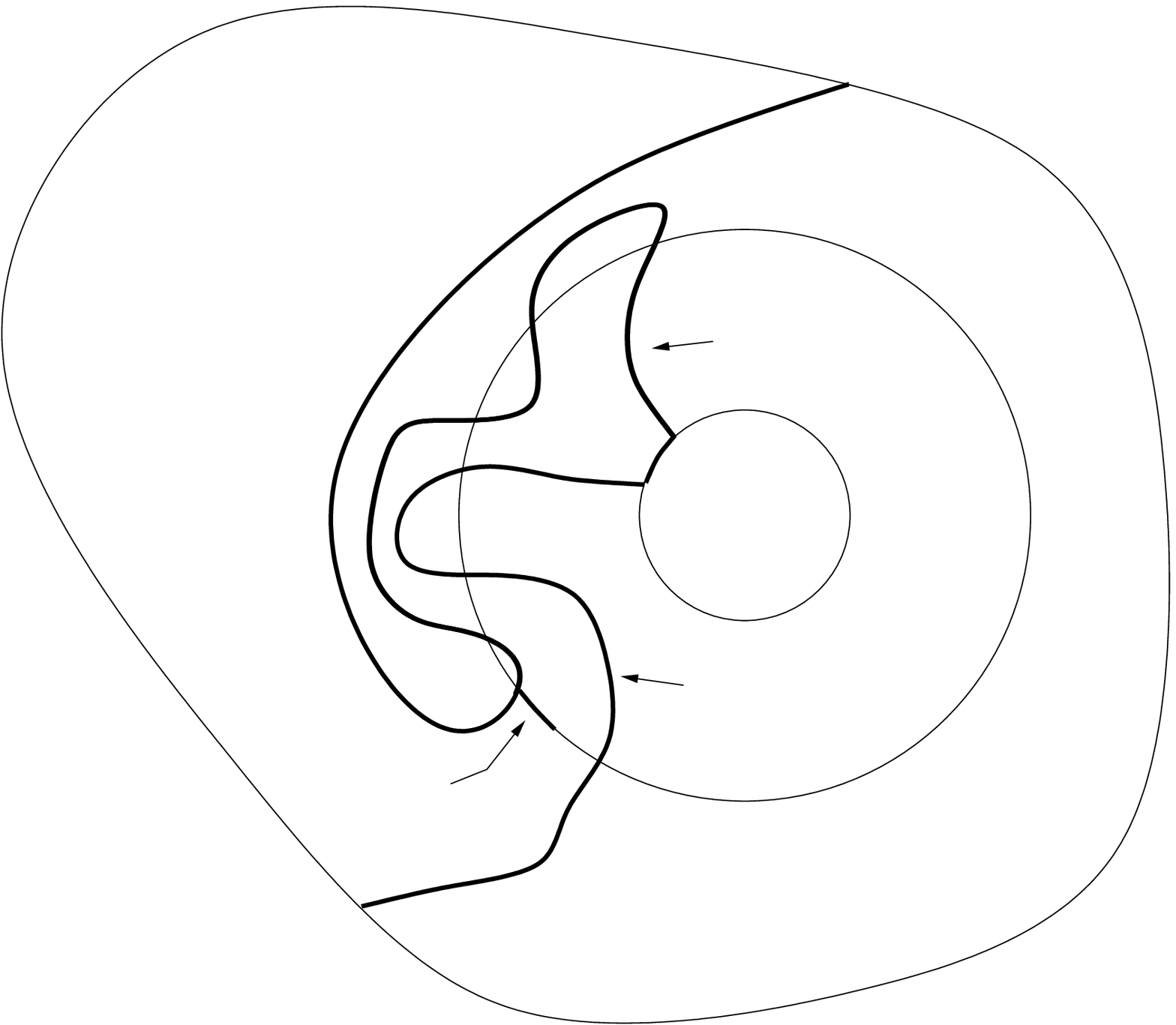}}
}
\begin{caption} {\label{f.Dp}The domain $D'$.}
\end{caption} \end{figure}

Let $D'$ be the connected
component of $\DD(\exter{3r}\gamma)\setminus \bal{3r}$ that contains
$\alpha_1$. Note that
$D'$ is a simply connected domain.
  See Figure~\ref{f.Dp}.
 Let $z_{\alpha_1}$ denote the endpoint of $\alpha_1$ on
$\p D'$ and let
$\bbeta_1$ and $\bbeta_2$ be the two connected components of
$\p D'\setminus\bl(\{z_{\alpha_1}\}\cup\closure{\bal {3r}}\br)$.
(Both have $z_{\alpha_1}$ as an endpoint and the
other endpoint in $\p \bal{3r}$.)
 Note that any path
in $D'\setminus \alpha_1$ connecting $\bbeta_1$ and $\bbeta_2$
separates $\alpha_1$ from $\p\bal{3r}$ in $\DD(\exter{3r}\gamma)$.
For each $\rho\in (3\,r+3,r'-3)$ let $A(\rho)$
denote the connected component of $D'\setminus \closure{B(0,\rho)}$ that
contains $\alpha_1$, and let $\alpha(\rho)$ denote the connected
component of $\p A(\rho)\cap \p B(0,\rho)$ that separates $\alpha_1$
from $\bal {3r}$ in $D'$.
Observe that $\alpha(\rho)$ has one endpoint on $\bbeta_1$ and the other
on $\bbeta_2$.
 If $3\,r+3<\rho<\rho'<r'-3$, then
$A(\rho)\supset A(\rho')$ and therefore $\alpha(\rho')$ separates
$\alpha(\rho)$ from $\alpha_1$ in $\DD(\exter{3r}\gamma)$ and $\alpha(\rho)$ separates
$\alpha(\rho')$ from $\bal {3r}$. When $3\,r+3<\rho<\rho'<r'-3$, let
$A(\rho,\rho')$ denote the connected component of
$D'\setminus(\alpha(\rho)\cup\alpha(\rho'))$ whose boundary contains
$\alpha(\rho)\cup\alpha(\rho')$.

For $n\in\N$, let $\rho_n:=2^{n}\,(3\,r+3)$,
and let $N$ be the largest $n$ such that $8\,\rho_n<r'$.
Fix some $n\in\{1,2,\dots,N\}$
and some small $\delta>0$ ($\delta=1/100$ should do).
Set $A_n^\delta:=A\bl((1-\delta)\,\rho_n,(1+\delta)\,\rho_n\br)$.
By continuity, $A^\delta_n$ contains points
$w$ such that $\dist(w,\bbeta_1)=\dist(w,\bbeta_2)$.
Set
$$
\eta(n):=\min\bl\{\dist(w,\bbeta_1):w\in \closure{A^\delta_n},\dist(w,\bbeta_1)=\dist(w,\bbeta_2)\br\}\,.
$$
First, assume that $\eta(n)>\delta\,\rho_n$.
In that case, a barrier $\Upsilon_n$ is defined as follows.
By continuity, there is a subarc $\Upsilon'$ of $\alpha(\rho_n)\subset A^\delta_n$
with endpoints $z_1$ and $z_2$
such that $\dist(z_j,\bbeta_j)=\delta\,\rho_n/10$, $j=1,2$, and $\dist(\Upsilon',\p D')=\delta\,\rho_n/10$.
Let $z_j'$ be a point in $\bbeta_j$ at distance $\delta\,\rho_n/10$ from $z_j$, $j=1,2$.
Then we take $\Upsilon_n$ as the union of $\Upsilon'$ with the two line segments
$[z_1,z_1']$ and $[z_2,z_2']$.
Recall the definition of $\dist(\cdot,\cdot;\cdot)$ from~\eref{e.distdef},
and
note that $\dist(z_j',\bbeta_{3-j};D')>\delta\,\rho_n/10$, $j=1,2$,
for otherwise, by continuity again, there would be a point $w\in\closure{D'}$
satisfying $\dist(w,z_j';D')\le\delta\,\rho_n/10$ (and therefore $w\in\closure{ A^\delta_n}$)
that is at equal distance from $\bbeta_1$ and from $\bbeta_2$, which would contradict
our assumption $\eta(n)>\delta\,\rho_n$.
It easily follows that in this case $\Upsilon_n$ is a $(\delta/20,2\,\rho_n)$-barrier.

We now assume that $\eta(n)\le \delta\,\rho_n$.
Let $w_1\in\closure{A^\delta_n}$ be a point satisfying
$\dist(w_1,\bbeta_1)=\dist(w_1,\bbeta_2)\le \delta\,\rho_n$.
For $j=1,2$, let $p_j\in\bbeta_j$ be a point satisfying
$|w_1-p_j|=\dist(w_1,\p D')$.
If $\dist(p_j,\bbeta_{3-j};D')\ge \delta^2\,\rho_n$ for $j=1,2$,
then we may take as our barrier the union $[p_1,w_1]\cup[w_1,p_2]$.
This will be a $(\delta/4,2\,\delta\,\rho_n)$-barrier.
Otherwise, fix a point $w_2$ satisfying $\dist(w_2,p_j;D')\le\delta^2\,\rho_n$
and $\dist(w_2,\bbeta_1)=\dist(w_2,\bbeta_2)\le \delta^2\,\rho_n$ and consider
the above construction with $w_2$ in place of $w_1$.
It may happen that the construction succeeds now, and we construct
a $(\delta/4,2\,\delta^2\,\rho_n)$-barrier.
Otherwise, we find a point $w_3\in D'$ satisfying $\dist(w_3,w_1;D')\le (\delta+\delta^2+\delta^3)\,\rho_n$
such that $\dist(w_3,\bbeta_1)=\dist(w_3,\bbeta_2)\le \delta^3\,\rho_n$.
We continue this procedure until some $(\delta/4,2\,\delta^m\,\rho_n)$ barrier is
obtained.  The procedure must terminate successfully at some finite $m$,
for otherwise the points $w_m$ would converge to some point in $\bbeta_1\cap\bbeta_2$
within distance $2\,\delta\,\rho_n$ from $\closure{A^\delta_n}$, which is clearly
impossible. Note that the barrier $\Upsilon_n$ thus constructed is contained in
$A^{2\,\delta}_n$.
Thus, when $1\le n'<n\le N$, $n,n'\in\N$, we have
$\dist(\Upsilon_n,\Upsilon_{n'};D')>\rho_n/4$
and $\Upsilon_n$ separates $\alpha_1$ from $\Upsilon_{n'}$ in $D'$.

Suppose $n\in\{1,2,\dots,N-1\}$.
Note that (contrary to what appears in Figure~\ref{f.Dp}, which does not show the scale of the
lattice) the endpoints of $\Upsilon_n$ are not on $\gamma$,
since $\bbeta_1\cup\bbeta_2$ are disjoint from $\gamma$, by construction.
On the event $\gamma\cap\Upsilon_n\ne\emptyset$,
let $\gamma_1$ and $\gamma_2$ be the two arcs
of $\gamma$ extending from the endpoints of
$\p_+$ to the first encounter with $\Upsilon_n$.
Now we apply the Barriers Theorem~\ref{t.barrier}
with $ Y=\{\Upsilon_{n+1}\}$.
Our careful construction above ensures that
$\Upsilon_{n+1}$ is an $(\eps,\diam\Upsilon_{n+1})$-barrier for some universal
constant $\eps>0$.
Note that $\bbeta_1$ and $\bbeta_2$ contain vertices
on which $h$ takes the same sign.  We conclude from the theorem that
$$
\Pb{\gamma\cap\Upsilon_{n+1}\ne\emptyset\md \gamma_1,\gamma_2,\oconf{3r},\conf{2r}}
<1-c_1
$$
for some $c_1=c_1(\upperhco)>0$.
The above implies
$$
 \Pb{\gamma\cap\Upsilon_{n+1}\ne\emptyset\md \gamma\cap\Upsilon_n\ne\emptyset,
\oconf{3r},\conf{2r}} <1-c_1 \,,
$$
which gives
\begin{equation}\label{e.pen}
\Pb{\gamma\cap\alpha_1\ne\emptyset\md\oconf{3r},\conf{2r}} \le (1-c_1)^{N-1}.
\end{equation}
Conditioned on $\gamma,\oconf{3r}$ and $\conf{2r}$, the probability of $\evZ$ is at most
$O(1)/\log r$,
by the proof of the upper bound in Lemma~\ref{l.factor}.
Thus,
$$
\Pb{\gamma\cap\alpha_1\ne\emptyset,\evZ\md\oconf{3r},\conf{2r}} \le O(1)\,(1-c_1)^{N-1}/\log r\,.
$$
On the other hand, the lower bound tells us that on the event
$ \{ \Qual(3r)\wedge\Qual(2r)\ge c_0\}$,
we have
$O_{\upperhco,p}(1)\,\Pb{\evZ\md\oconf{3r},\conf{2r}}\ge 1/\log r$.
Thus,~\eref{e.staynear} follows.
\medskip

Clearly,~\eref{e.staynear} also holds for $\alpha_2$.
On $\evZ$ let $\alpha_1^*$ and $\alpha_2^*$
be the two connected components of $\alpha\cap \DD(\gamma)$
that have one endpoint in $S'$ and the other in $\p \DD(\gamma)$.
Note that when $\evZ$ holds and $\gamma\cap(\alpha_1\cup\alpha_2)=\emptyset$,
we have $\alpha_1^*\cup\alpha_2^*=\alpha_1\cup\alpha_2$.
Thus,~\eref{e.staynear} for $\alpha_1$ and for $\alpha_2$ together with~\eref{e.goodq}
now gives
\begin{equation}\label{e.n2}
\Pb{ \alpha_1^*\cup\alpha_2^*\ne \alpha_1\cup\alpha_2
\md \oconf R,\conf r,\evZ}<3\,p/10\,.
\end{equation}

We now follow an analog of the above argument with the roles of inside
and outside switched.
Observe that there is a connected component
$\tilde \alpha$ of $\p B(0,r')\setminus \tilde S$ such that
$\tilde S\cup\tilde\alpha$ separates $\p D$ from $\bal{3r}$.
We fix such an $\tilde \alpha$, and if there is more than one possible
choice, we choose in the same way in which $\alpha$ was chosen from $S'$;
that is, we make sure that $\alpha=\tilde\alpha$ if $S'=\tilde S$ (as unoriented paths).
The point is that although $\alpha$ is $\oconf{3r}$-measurable
and $\tilde\alpha$ is $\conf{R/3}$-measurable, we have
$\alpha=\tilde\alpha$ if $S'=\tilde S$.

On the event $\Qual(\conf{R/3})>0$, let  $\tilde \alpha_1$ and $\tilde\alpha_2$ be the two
connected components of $\tilde\alpha\cap\DD(\beta_{R/3})$
that have one endpoint on $\tilde S$ and the other on
$\p \DD(\beta_{R/3})$.
On the event $\evZ$, let $\tilde\alpha_1^*$ and $\tilde\alpha_2^*$
be the two connected components of
$\tilde\alpha\cap\DD(\gamma)$
that have one endpoint on $\tilde S$ and the other on $\p\DD(\gamma)$.
Essentially the same proof which gave~\eref{e.n2} now gives
\begin{equation}\label{e.n3}
\Pb{ \tilde\alpha_1^*\cup\tilde\alpha_2^*\ne \tilde\alpha_1\cup\tilde\alpha_2
\md \oconf R,\conf r,\evZ}<3\,p/10\,.
\end{equation}
But observe that when $\evZ$ and $S'=\tilde S$ hold,
we clearly have
$\tilde\alpha_1^*\cup\tilde\alpha_2^*= \alpha_1^*\cup\alpha_2^*$.
Also recall that $S'=\tilde S$ when $S_*\subset\bal{r'/3}$.
Thus, we get from~\eref{e.Sdirect},~\eref{e.n2} and~\eref{e.n3}
$$
\Pb{ \tilde\alpha_1\cup\tilde\alpha_2\ne \alpha_1\cup\alpha_2
\text{ or }
S_*\not\subset \bal{r'/3}
\md \oconf R,\conf r,\evZ}<4\,p/5\,.
$$

Assume that
$\evZ$, $ \tilde\alpha_1\cup\tilde\alpha_2=\alpha_1^*\cup\alpha_2^*= \alpha_1\cup\alpha_2$
and $S_*\subset\bal{r'/3}$ hold.
It remains to show that in this case
the path $\gamma$ has no more than two disjoint arcs connecting
$\bal r$ and $\p\bal R$.
Recall that $\tilde\alpha_1$ has an endpoint on
$\p\DD(\beta_{R/3})$. This endpoint is on a $\TG$-triangle
containing a $\TG^*$-vertex $v_1\in\beta_{R/3}$.
Similarly, there is a $\TG^*$ vertex $v_2\in\beta_{R/3}$
for which the $\TG$-triangle containing it
has an endpoint of $\tilde\alpha_2$.
{}From $ \tilde\alpha_1\cup\tilde\alpha_2= \alpha_1\cup\alpha_2$
we conclude that $v_1,v_2\in \exter{3r}\gamma$ as well.

Shortly, we will prove that $v_1$ and $v_2$ are in separate connected
components of $\beta_{R/3}\setminus e^*_\sig$.
This implies that each connected component of
$\beta_{R/3}\setminus e^*_\sig$ intersects $\exter{3r}\gamma$.
Since $\beta_{R/3}\cup\exter{3r}\gamma\subset\gamma$,
and $\gamma$ is a simple path,
it easily follows that $\gamma=\beta_{R/3}\cup\exter{3r}\gamma$,
which implies that there are at most two disjoint crossings
in $\gamma$ between $\bal r$ and $\p\bal{R}$.

It remains to prove that $v_1$ and $v_2$ are in different connected
components of $\beta_{R/3}\setminus e_\sig^*$.
This will be established using planar topology arguments.
Let $\hat \alpha$ consist of
$\tilde\alpha_1\cup\tilde\alpha_2$, a simple path
$\tilde S_0\subset \tilde S$ connecting them, and short line segments
(contained in the $\TG$-triangles containing $v_1$ and $v_2$)
from the endpoints of $\tilde\alpha_1$ and
$\tilde\alpha_2$ on $\p \DD(\beta_{R/3})$
to $v_1$ and $v_2$. Then $\hat\alpha$ is a simple
path and only the endpoints of $\hat\alpha$ are on $\gamma$.
Let $\tilde\beta$ be the connected component of $\beta_{R/3}\setminus\{v_1,v_2\}$
with endpoints $v_1$ and $v_2$.
Then $\tilde\beta\cup\hat\alpha$ is a simple closed path
and it suffices to show that $e^*_\sig\subset \tilde\beta$.
If $\tilde\beta\cup\hat\alpha$ separates $0$ from $\p D$, then $\tilde\beta$
must contain $e_\sig^*$, because each connected component of
$\gamma\setminus\{e_\sig^*\}$ connects $\p D$ to
$\TG^*$ vertices adjacent to $0$ and is disjoint from
$\hat\alpha\setminus\{v_1,v_2\}$.

Suppose that $\tilde\beta\cup\hat\alpha$ does not separate $0$ from $\p D$.
Recall that $\tilde \alpha\cup\tilde S$ separates $\p D$ from $\bal {3r}$
and therefore from $0$. Since $\tilde S$ itself does not separate
$\p D$ from $0$, it follows that
the winding number of
$\tilde\alpha\cup \tilde S_0$ around $0$ is $\pm 1$
(depending on orientation).
Since
$\tilde\beta\cup\hat\alpha$ does not separate $0$ from $\p D$, its winding number
around $0$ is zero.
If we remove from the  union of the
two paths $\tilde\alpha\cup \tilde S_0$ and $\tilde\beta\cup\hat\alpha$
all the nontrivial arcs where they agree,
we get a closed curve $\chi$,
which consists of $\tilde\beta$, a segment of $\tilde\alpha$
and the two short connecting segments near $v_1$ and $v_2$,
and $\chi$ has odd winding number around $0$.
Consequently, it separates $0$ from $\p D$.
But observe that $ \tilde S$ is disjoint from $\chi$.
Moreover, since we are assuming $S_*\subset \bal{r'/3}$,
it follows that $\inter{R/3}\bS$ is also
disjoint from it. But this contradicts  the fact that
$\chi$ separates $\p D$ from $0$, since $\inter{R/3}\bS$
can be extended to a path disjoint from $\chi$ and connecting
$0$ and $\p D$. Thus, the proof is now complete.
\QED

\subsection{Coupling and limit} \label{ss.coupling}

In this subsection, we retain our previous assumptions~\iref{i.h}, \iref{i.D}
and~\iref{i.S} about the system $D,\p_+,\p_-,h_\p,\vv,h,\gamma$.
Moreover, we also consider another such system $D'$, $\p_+'$,
$\p_-'$, $h_\p'$, $\vv'$, $h'$, $\gamma'$, which is assumed to
satisfy the same assumptions. In particular,
$\|h_\p'\|_\infty\le\upperhco$. Generally, we will use $'$ to denote
objects related to the system in $D'$. For example, $\ev J'(r,R)$
will denote the event corresponding to $\ev J(r,R)$.

\begin{definition}
Fix $R>r>0$, and suppose that $\bal R\subset D\cap D'$.
Consider the intersection $\exter r\gamma\cap\bal R$
as a collection of oriented paths, oriented so
as to have vertices in $V_+(\gamma)$ on the right.
We say that $\oconf r$ and $\oconf r'$ {\bf match}
in $\bal R$ if
the set of vertices in $\bal R$ visited by
$\exter r S$ is the same as the corresponding
set for $S'$ and
$\exter r\gamma\cap\bal R=\exter r\gamma'\cap \bal R$
with all the orientations agreeing or with all the
orientations reversed.
\end{definition}

We now show that if the configurations match in a big annulus, then it is
likely that the interfaces agree in the inner disk; more precisely, we have:

\begin{lemma}\label{l.genseq}
For every $\delta>0$, $r\ge 10$ and $R^*>r+3$ there is an
$R=R(\delta,r,R^*,\upperhco)>R^*$ such that the following holds.
Suppose that $\vv',\vv\notin\ball_R\subset D\cap D'$, $\Pb{\neg\ev
J(r,R^*), \evZa\md\oconfr}>0$ and $\Pb{\neg\ev J'(r,R^*),
\evZa'\md\oconfr'}>0$. Assume that $\oconf r$ and $\oconf r'$ match
in $\bal R$.  In particular, the endpoint $q^r$ of $\exter{r} S$ in
$\ball_{r}$ is the same as that of $\exter{r} S'$. Let $\nu$ be the
law of $\gamma_*:=\gamma\setminus\exter r\gamma$ (as an unoriented
path) conditioned on $\evZa,\oconfr$ and $\neg\ev J(r,R^*)$, and let
$\nu'$ be the law of $\gamma'_*:=\gamma'\setminus\exter r\gamma'$
conditioned on $\evZa',\oconf r'$ and $\neg\ev J'(r,R^*)$. Then
$\|\nu-\nu'\|< \delta$.
\end{lemma}

Here, $\|\nu-\nu'\|$ denotes the total variation norm
$\sum_{\vartheta}\bl|\nu[\gamma_*=\vartheta]-\nu'[\gamma'_*=\vartheta]\br|$.

\proof
Assume that the orientation of
$\exter r\gamma\cap\bal R $ agrees with that of $\exter r\gamma'\cap\bal R$.
This involves no loss of generality, since we may replace $\p'_+$
with $\p'_-$, replace $h'$ by $-h'$, etc.

Since we are assuming $\Pb{\neg\ev J(r,R^*), \evZa\md\oconfr}>0$,
there is a path $\vartheta\subset\bal {R^*}$ such that
$\Pb{\gamma_*=\vartheta,\neg\ev J(r,R^*),\evZa\md\oconfr}>0$. Let
$\hat\Gamma$ be the collection of all such $\vartheta$, and fix some
$\vartheta\in\hat\Gamma$. Obviously, the length of $\vartheta$ is
$O(R^*)^2$. We start extending $\exter r\gamma$ starting at one of
the endpoints, say $x^r$, and consider the conditional probability
that each successive step follows $\vartheta$, given that the
previous steps have and given $\oconfr$. Each step is decided by the
sign of $h$ on a specific vertex $v$. When we condition on the
values of $h$ on the neighbors of $v$, the conditional
law of $h(v)$ is a
Gaussian with some constant positive variance. It follows
from~\eref{e.atv} that with high
probability (conditioned on the success of
the previous steps) the mean of this
Gaussian random variable is unlikely to be large. Thus, the
probability for either sign is bounded away from zero, which means
that each step is successful with probability bounded away from
zero. By~\eref{e.btv} it is unlikely that $h(v)$ will be very close
to zero. Proposition~\ref{p.hcont} therefore implies that if $R>R^*$
is very large, the probability for a successful one step extension
for $\gamma'$ is almost the same as for $\gamma$. Thus, we conclude
that for sufficiently large $R>R^*$
$$
(1-\delta)\,\Pb{\gamma'_*=\vartheta\md\oconfr'}\le
\Pb{\gamma_*=\vartheta\md \oconfr}\le (1+\delta)\,\Pb{\gamma'_*=\vartheta\md\oconfr'}
$$
holds for all $\vartheta\in\hat\Gamma$.
It is moreover clear that
$$
\Pb{\evZa,\neg\ev J(r,R^*)\md\gamma_*=\vartheta,\oconfr}
=
\Pb{\evZa',\neg\ev J'(r,R^*)\md\gamma'_*=\vartheta,\oconfr'},
$$
because under $\neg\ev J(r,R^*)$ the random walk $S$ cannot
get close to any place where $\exter r\gamma$ differs
from $\exter r\gamma'$ between the first visit
to $q^r$ and time $\tau_0$.
Thus,
$$
1-\delta\le
\frac
{
\Pb{\gamma_*=\vartheta,\evZa,\neg\ev J(r,R^*)\md \oconfr}
}{
\Pb{\gamma'_*=\vartheta,\evZa',\neg\ev J'(r,R^*)\md\oconfr'}
}
\le 1+\delta\,.
$$
The lemma follows (though perhaps $\delta$ needs to be readjusted).
\QED

The next lemma shows that given $\oconf R$ the events $\evZ$ have comparable
probabilities for different $\sig$.

\begin{lemma}\label{l.evZ}
As usual, assume~\iref{i.h}, \iref{i.D} and~\iref{i.S}.
There is a constant $c=c(\upperhco)\ge 1$ such that for all $R$ sufficiently large
and every $\sig,\sig'\in\{0,1,\dots,5\}$ we have
$$
\Pb{\evZ\md\oconf R}\le c\,\Pb{\evZs{\sig'}\md\oconf R}\,.
$$
\end{lemma}

\proof
The statement is clear when $R=100$, because in that case if it
is at all possible to extend $\oconf R$ in such a way that $\evZ$ holds,
then there is probability bounded away from zero
(by a function of $\upperhco$) that
$\evZs{\sig'}$ holds. (We may choose the continuations of
$\gamma$ and $S$ as we please, and as long as the continuations involve
a bounded number of steps, the probability for these continuations are
bounded away from zero,
as in the proof of Lemma~\ref{l.genseq}.)
When $R>100$, we may just condition on the corresponding extension of
$\oconf R$ up to radius $100$.
\QED

We now come to one of the main results in this section --- the existence of a limiting interface.

\begin{theorem}[Limit existence]\label{t.limit}
There is a (unique) probability measure $\mu_\infty$ on the space of
two-sided infinite simple $\TG^*$-paths $\dot\gamma$ which is the limit of the law
of $\gamma$ (unoriented) conditioned on $\evZa$ and $\oconf R$, in the following sense.
Assume~\iref{i.h}, \iref{i.D} and~\iref{i.S}.
For every finite set of $\TG^*$-edges $E_0$ and every $\delta>0$
there is an $R_0=R_0(\delta,E_0,\upperhco)$
such that if $R>R_0$, $\vv\notin\bal R\subset D$ and
$\Pb{\evZa\md\oconf R}>0$,
then
$$
\Bl|\Pb{E_0\subset\gamma\md\evZa,\oconf R}-\mu_\infty[E_0\subset\dot\gamma]\Br|<\delta\,.
$$
\end{theorem}

\proof
Clearly, it suffices to show that for every $r>0$ if $R$ is sufficiently large,
$\vv,\vv'\notin\bal R\subset D\cap D'$,
$\Pb{\evZa\md \oconf R}>0$ and $\Pb{\evZa'\md\oconf R'}>0$
then we may couple the conditioned laws of
$\gamma$ given $\oconf R$ and $\evZa$ and
$\gamma'$ given $\oconf R'$ and $\evZa'$
such that
\begin{equation}\label{e.match}
\Pb{\gamma\setminus\exter r\gamma = \gamma'\setminus\exter r\gamma' \md \oconf R,\evZa,\oconf R',\evZa'}>1-\delta\,.
\end{equation}
(Here, the equivalence is an equivalence of unoriented paths.)
Let $a_n$ be the constant $a$ given by Lemma~\ref{l.direct} when one takes $p=\delta_n :=2^{-n}\,\delta/8$.
We define a sequence of radii $r_0,r_1,\dots$ inductively, as follows.
Let $c_1$ be the constant $c$ given by Corollary~\ref{c.couple2}, and set $r_0:=r\vee10\vee c_1$.
Given $r_n$, let $\hat r_n$ be the $R$ promised by Lemma~\ref{l.genseq} when
we take $r_n$ for $r$, $\delta_n$ for $\delta$ and $a_n\,r_n$ for $R^*$.
Finally, set $r_{n+1}=4\,a_n\,\hat r_n$.
Let $c_2$ be the constant promised by Lemma~\ref{l.evZ}.
We assume, with no loss of generality, that $\delta<1/(4\,c_1)$.
Let $N\in\N$ be sufficiently large so that $(1-1/(72\,c_1\,c_2^2))^{N-1}<\delta/2$.
We will prove~\eref{e.match} on the assumption that $R>6\,r_N$.

The construction of the coupling is as follows.
First, we choose $\oconf {r_{N-1}}$ and $\oconf {r_{N-1}}'$ independently
according to their conditional distribution given
$\oconf R$, $\evZa$, $\oconf R'$ and $\evZa'$.
We proceed by reverse induction.
Suppose  that $n\in [1,N-1]\cap \N$ and
that $\oconf {r_n}$ and $\oconf {r_n}'$ have been determined.
If $\oconf {r_n}$ and $\oconf {r_n}'$ match inside $\bal {\hat r_n}$,
then we couple $\gamma$ and $\gamma'$ in such a way as
to maximize the probability that
$\gamma\setminus\exter {r_n}\gamma=\gamma'\setminus\exter{r_n}\gamma'$,
subject to maintaining their correct  conditional distributions given
the choices previously made.
If they do not match, then we couple
$\oconf{r_{n-1}}$ and $\oconf{r_{n-1}}'$ in such a way as to
maximize the probability that they match in $\bal {\hat r_{n-1}}$,
subject to their correct conditional distributions.
If $\oconf {r_0}$ and $\oconf {r_0}'$ have been determined,
but $\gamma$ and $\gamma'$ have not, then we couple
$\gamma$ and $\gamma'$ arbitrarily, subject to their correct
conditional distributions.

We claim that the coupling just described achieves the bound~\eref{e.match}.
Let $\ev M_n$ denote the event that
$\oconf {r_{n}}$ and $\oconf {r_{n}}'$ match inside
$\bal {\hat r_{n}}$, where $n\in\{1,2,\dots,N-1\}$,
and let $\ev M$ denote the union of these events $\ev M_n$.
It follows from the choice of $\hat r_n$ that
if $\ev M_n$ holds and $n$ is the largest
$n'\in\{1,2,\dots,N-1\}$ with that property,
then there is a coupling
of the appropriately conditioned laws of $\gamma$
and $\gamma'$ such that
\begin{multline*}
\Pb{\gamma\setminus\exter {r_n}\gamma\ne\gamma'\setminus\exter{r_n}\gamma'
\md \oconf {r_n},\oconf {r_n}',\evZa,\evZa'}
\\
\le \delta_n+
\Pb{\ev J(r_n,a_n\,r_n)\md \oconf {r_n},\evZa'}+
\Pb{\ev J'(r_n,a_n\,r_n)\md \oconf {r_n}',\evZa'},
\end{multline*}
and hence this also holds for our coupling.
By taking conditional expectation and summing over $n$, we get
\begin{align*}
&
\Pb{\ev M ,\gamma\setminus\exter r\gamma \ne \gamma'\setminus\exter r\gamma'
\md \oconf R,\evZa,\oconf R',\evZa'}
\\&\qquad
\le \delta/8+\sum_{n=1}^{N-1}
\Bl(
\Pb{ \ev J(r_n,a_n\,r_n)\md \oconf R,\evZa}
+
\Pb{ \ev J'(r_n,a_n\,r_n)\md \oconf R',\evZa'}
\Br)
\\&\qquad
\le
3\,\delta/8\qquad\text{(by the choice of $a_n$)}\,.
\end{align*}

Now fix some $n\in\{1,2,\dots,N-2\}$, and suppose that
none of the events $\ev M_{n'}$, $n'>n$ occurs.
Fix some arbitrary $\sig\in\{0,1,\dots,5\}$.
Conditional on $\oconf {r_{n+1}}$ and on $\evZa$,
by the choice of $c_2$, there
is probability at least $1/(6\,c_2)$ that $\evZ$ holds,
and the same is true for the system in $D'$.
If we additionally condition on $\evZ$ and $\evZp$,
then by the choice of $c_1$ there is a coupling
of the appropriate conditioned laws of
$\conf {r_{n+1}/4}$ and $\conf {r_{n+1}/4}'$ such that
$$
\Pb{
\conf {r_{n+1}/4}=\conf {r_{n+1}/4}'\md \oconf {r_{n+1}},
\oconf {r_{n+1}}',\evZ,\evZp}\ge 1/c_1\,.
$$
But note that if $\conf {r_{n+1}/4}=\conf {r_{n+1}/4}'$
and
$\neg \bl(\ev J(\hat r_n,r_{n+1}/4)\cup \ev J'(\hat r_n,r_{n+1}/4)\br)$
both hold (as well as $\evZ\cap\evZp$), then $\ev M_n$ holds as well.
We may then consider a coupling of the two systems which first decides
the two events $\evZ$ and $\evZp$ independently,
and if both hold (which happens with probability at least
$1/(6\,c_2)^2$), then with conditional probability at least $1/c_1$
we also have $\conf {r_{n+1}/4}=\conf {r_{n+1}/4}'$.
Thus, under this coupling
\begin{align*}
&
(6\,c_2)^2\,\Pb{\ev M_n
\md \oconf {r_{n+1}},\oconf{r_{n+1}}',\evZa,\evZa'}
\\&\qquad
\ge c_1^{-1} - \Pb{
\ev J(\hat r_n,r_{n+1}/4)\cup \ev J'(\hat r_n,r_{n+1}/4)
\md \oconf {r_{n+1}},\oconf{r_{n+1}}',\evZ,\evZp}
\\&\qquad
\ge c_1^{-1}-2\,\delta_n
\qquad\text{ (by the choices of $r_{n+1}$ and $a_n$)}
\\&\qquad
\ge c_1^{-1}/2
\qquad\text{ (by our assumption $4\,c_1\,\delta<1$).}
\end{align*}
This must hold for our coupling as well.
Consequently, induction gives
$$
\Pb{\neg\ev M}\le \Bigl(1-\frac 1{72\, c_2^2\, c_1}\Bigr)^{N-1}\,,
$$
which is less than $\delta/2$ by the choice of $N$.
Thus~\eref{e.match} follows, and the proof is complete.
\QED

\subsection{Boundary values of the interface}
\label{ss.boundaryvalues}

Consider the random path $\dot\gamma$ whose law is the
measure $\mu_\infty$ provided by Theorem~\ref{t.limit}.
We orient $\dot\gamma$ so that the edges $e_\sig^*$,
$\sig=0,1,\dots,5$ that are in $\dot\gamma$ are oriented
clockwise around the hexagon $\bigcup_{\sig=0}^5e_\sig^*$
centered at $0$.
Let $U_+$ denote the set of $\TG$-vertices
adjacent to $\dot\gamma$ on its right hand side,
and let $U_-$ denote the set of
vertices adjacent to $\dot\gamma$ on its left hand side.
Using the heights interface continuity (Proposition~\ref{p.hcont}),
it is clear that given $\dot\gamma$ we may define
the DGFF $\dot h$ on all of $\TG$ conditioned to be positive
on $U_+$ and negative on $U_-$, as a limit of an appropriately
conditioned DGFF on bounded domains.
Moreover, many properties of the DGFF on bounded domains easily
transfer to $\dot h$. In particular,~\eref{e.atv} applies,
to give $\Eb{|\dot h(0)|\md\dot\gamma}=O(1)$.
Set
\begin{equation}\label{e.hcodef}
\hco:=\Eb{\dot h(0)}.
\end{equation}
Clearly, $0<\hco<\infty$.

Recall that $\tau$ is the first time $t$ such that $S_t\in\p\DD(\gamma)$ and recall
the notation $\dist(\cdot,\cdot;\cdot)$ from~\eref{e.distdef}.
In this subsection we will show that in the limit as $\dist(\vv,\p D)\to\infty$ we have
$$
\EB{(\pm h(\zz)-\hco)
\,1_{\{\zz\in V_{\pm}(\gamma)\}}\md \gamma}\to0
$$
in probability, under the assumption that
\begin{enumerate}
\hitem{i.p}{($\p$)} $h_\p  (\p_+\cap V_{\p})\subset[-\lowerhco,\upperhco]$
and $h_\p(\p_-\cap V_{\p})\subset[-\upperhco,\lowerhco]$,
\end{enumerate}
where $\lowerhco=\lowerhco(\upperhco)>0$ is the constant given by Lemma~\ref{l.bdnarrows}.
The importance of the assumption~\iref{i.p} is that
by Remark~\ref{r.bdbarrier} it enables the application of the Barriers Theorem~\ref{t.barrier} to barriers
with endpoints on $\p D$, provided that the barriers in $Y_+$ do not have
endpoints in $\p_-$ and those in $Y_-$ do not have endpoints in $\p_+$.
This will allow us to prove:

\begin{theorem}\label{t.nohitbd}
Assume~\iref{i.h}, \iref{i.D}, \iref{i.S} and~\iref{i.p}.
There are positive constants $\expo2=\expo2(\upperhco)$
and $c=c(\upperhco)$ such that the following holds true.
Let $z$ be any point on $\p_+$. Then for every $r>1$,
$$
\Pb{\dist(z,\gamma;D)<r} \le c\,\bl(r/\dist(z,\p_-;D)\br)^\expo2.
$$
\end{theorem}

\proof
Let $\beta_1$ and $\beta_2$ be the two components of $\p_+\setminus\{z\}$.
Let $R=\dist(z,\p_-;D)$.
For each $\rho\in (0,R)$, let $A(\rho)$ denote the
connected component of $B(z,\rho)\cap D$ that has $z$ in its boundary,
and let $\alpha(\rho)$ denote the connected component
of $\p A(\rho)\setminus\p D$ that separates $z$ from $\p_-$ in $D$.
Using this construction,
the proof proceeds as in the proof~\eref{e.pen},
except that the barriers start from the outside and get closer to $z$,
and we appeal to Remark~\ref{r.bdbarrier} instead of the Barriers Theorem.
We leave it to the reader to verify that the proof carries over with no
other significant modifications.
\QED

Our next lemma shows that it is unlikely that $\zz$
is adjacent to $\gamma$ and is near $\p D$.

\begin{lemma}\label{l.farbd}
Assume~\iref{i.h}, \iref{i.D}, \iref{i.S} and~\iref{i.p}.
For every $p>0$ there is some $\delta=\delta(p,\upperhco)>0$
such that
$$
\Pb{0<\dist(\zz,\p D)<\delta\,\dist(\vv,\p D)}<p\,.
$$
\end{lemma}

\proof Let ${\tau_D}$ be the first $t$ such that $S_t\in\p D$, and
let $M=M(S,\gamma):=\{S_t: S_t\in \p D(\gamma), 0\le t<\tau_D\}$. We
will prove the stronger statement
\begin{equation}\label{e.M}
\Pb{\dist(M,\p D)<\delta\,\dist(\vv,\p D)}<p\,.
\end{equation}
(By convention $\dist(\emptyset,\p D)=\infty$.) Fix $r>0$ and set
$R=\dist(\vv,\p D)$. Conditioned on $\dist(M,\p D)<r$, we have
$\dist(S_{\tau_D},\gamma;D)<4\,r$ with probability bounded away from
zero, since the random walk started at any $v\in M$ such that
$\dist(v,\p D)<r$ has probability bounded away from $0$ to surround
the closest point to $v$ on $\p D$ (and therefore hit $\p D$) before
exiting the ball of radius $2\,r$ about that point. It therefore
suffices to prove
\begin{equation}\label{e.zbdgam}
\Pb{\dist(S_{{\tau_D}},\gamma;D)<\delta\,R}<p
\end{equation}
for $\delta=\delta(p,\upperhco)>0$.
Let $\ev A_+=\ev A_+(\delta)$ denote the event
$ \dist(S_{{\tau_D}},\p_-; D)>\delta^{1/2}\,R$,
and similarly define $\ev A_-$ with $\p_-$ replaced by $\p_+$.
By conditioning on $S_{{\tau_D}}$,
Theorem~\ref{t.nohitbd} shows that
$$
\Pb{\ev A_+,\dist(S_{\tau_D},\gamma;D)<\delta\,R}<p/3
$$
for an appropriate choice of $\delta$.
 A symmetric argument
applies on $\ev A_-$.
Consequently, it is enough to prove that
$\Pb{\neg (\ev A_+\cup\ev A_-)}<p/3$
for an appropriate choice of $\delta$.

\begin{figure}
\SetLabels
(.3*.8)$\vv$\\
(.2*.6)$D_*$\\
\T(.1*-0.01)$x_\p$\\
\L(1*.85)$y_\p$\\
\endSetLabels
\centerline{\epsfysize=2.3in%
\AffixLabels{%
\epsfbox{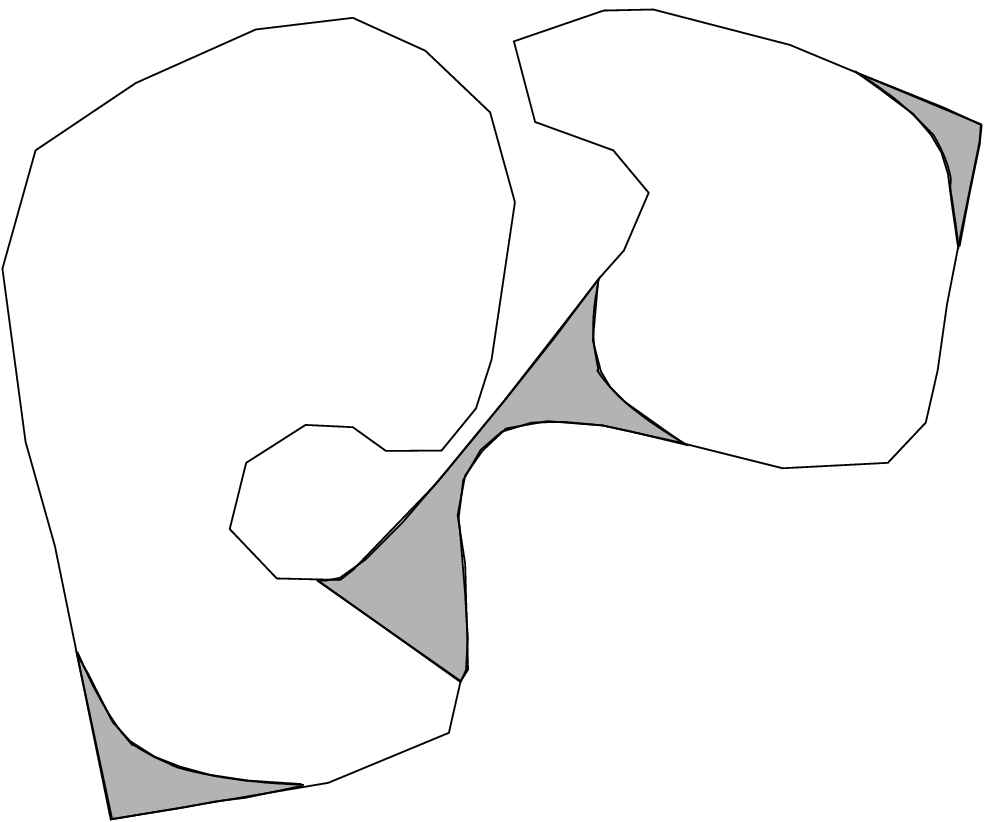}}
}
\begin{caption} {\label{f.L}The set $L_\rrho$ (shaded) and $D_*$.}
\end{caption} \end{figure}

Fix $\rrho:=\delta^{1/2}\,R$, let $L_\rrho$ denote the
set of points that lie on some path in $\closure D$
of diameter at most $\rrho$ connecting
$\p_+$ and $\p_-$,
and let $D_*$ be the connected component of
$\vv$ in $D\setminus L_\rrho$.
(See Figure~\ref{f.L}.)
We now prove that
\begin{equation}\label{e.twoballs}
\p D_*\setminus \p D \text{ is contained in the union of two balls of radius }2\,\rrho\,.
\end{equation}
Every path connecting $\p_+$ and $\p_-$ in
$\closure D$ must separate
$\vv$ from $x_\p$ or from $y_\p$ in $\closure D$
(because, by Jordan's theorem, it separates $x_\p$ from $y_\p$ in $\closure D$).
Let $\Gamma_1$ [respectively, $\Gamma_2$] denote the collection of paths in
$\closure D$ of
diameter at most $\rrho$ that
connect $\p_+$ and $\p_-$ and separate
$x_\p$ [respectively, $y_\p$] from $\vv$.
Then $\p D_*\setminus \p D$ is
contained in the union of the set of points
belonging to a path in $\Gamma_1$ and the set of
points belonging to a path in $\Gamma_2$.
Suppose that $\alpha$ and $\alpha'$ are two paths
in $\Gamma_1$, both of which intersect
$\p D_*$.
Let $\beta$ be a path connecting $x_\p$ with
$\alpha\cup\alpha'$ in $\closure D$, which is
disjoint from $\alpha\cup\alpha'$, except for its endpoint.
If $\beta\cap\alpha\ne\emptyset$, then
one can connect $x_\p$ to $\vv$ in $\beta\cup\alpha\cup D_*$,
and therefore $\alpha\cap\alpha'\ne\emptyset$
(since $\alpha'$ separates $\vv$ from $x_\p$ in $\closure D$).
Similar reasoning applies if $\beta\cap\alpha'\ne\emptyset$.
It follows that any two paths in $\Gamma_1$
that intersect $\p D_*$ must intersect each other,
and hence the collection of all such paths is covered by the ball of radius $2\,\rrho$
centered at any point on any such path.
Since a similar argument applies to $\Gamma_2$,~\eref{e.twoballs} follows.

By~\eref{e.twoballs} and
Lemma Hit Near~\ref{l.hitnear}
$$
\Pb{\exists t\le\tau_D: S_t\in L_\rrho}<p/3
$$
for all sufficiently small $\delta>0$.
Thus $\Pb{\neg (\ev A_+\cup\ev A_-)} = \Pb{S_{\tau_D}\in L_\rrho}<p/3$,
and the proof is complete.
\QED

Next, we show that $\zz$ is unlikely to be close to $\vv$ by proving the same
for $\gamma$.

\begin{lemma}\label{l.farvv}
Assume~\iref{i.h} and \iref{i.D}.
There are constants $c>0$, $\expo3>0$, both depending only on
$\upperhco$, such that for every $\delta>0$
$$
\Pb{\dist(\vv,\gamma)<\delta\,\dist(\vv,\p D)}<c\,\delta^{\expo3}\,.
$$
\end{lemma}

We expect that the left hand side is bounded by
$\delta^{1/2+o(1)}$, using the corresponding result~\cite{\RSsle}
for \SLEkk 4/.

\proof
Let $\Upsilon_n$ denote the circle of radius $2^{-n-1}\,\dist(\vv,\p D)$ about $\vv$.
As in the proof of~\eref{e.pen}, the Barriers Theorem~\ref{t.barrier}
implies that given that $\gamma$ intersects $\Upsilon_n$ (where $n\in\N$),
the conditional probability that $\gamma$ does not intersect $\Upsilon_{n+1}$
is bounded away from zero by a function of $\upperhco$, provided that
$2^{-n-1}\,\dist(\vv,\p D)>10$, say.
The lemma follows by induction.
\QED

\begin{proposition}\label{p.sides}
Assume~\iref{i.h}, \iref{i.D}, \iref{i.S} and \iref{i.p}.
For every $\eps>0$ there is an $R=R(\eps,\upperhco)$ such that
\begin{equation}\label{e.bdvalR}
\EB{\Eb{\bl(|h(\zz)|-\hco\br)\,1_{\{\zz\notin\p D\}}\md \gamma}^2}<\eps
\end{equation}
holds, provided that $\dist(\vv,\p D)>R$,
where $\hco$ is the constant given by~\eref{e.hcodef}.
\end{proposition}

The proof is based on the simple idea that given a single
instance of $\gamma$ we consider two independent
copies of $(h,S)$.

\proof
Set
$$
X:=\Eb{(|h(\zz)|-\hco)\,1_{\{\zz\notin\p D\}}\md \gamma}.
$$
To get a handle on $\Es{X^2}$,
let $(h',S')$ be independent from $(h,S)$ given $\gamma$
and have the same conditional law as that
of $(h,S)$ given $\gamma$. Thus,
$(h,S,\gamma)$ has the same law as $(h',S',\gamma)$.
Let $\tau':=\min\bl\{t:S'_t\in\p \DD(\gamma) \br\}$,
$y:=\Eb{|h(\zz)|-\hco\md \zz,\gamma}$ and $y':=\Eb{|h'(\zzp)|-\hco\md \zzp,\gamma}$.
Then
$ X^2 = \Eb{y\,y'\,1_{\{\zz,\zzp\notin\p D\}}\md\gamma}$ and hence
\begin{equation}
\label{e.xx}
\Eb{X^2} =
 \Eb{y\,y'\,1_{\{\zz,\zzp\notin\p D\}}}.
\end{equation}
Fix some $r_3\gg r_2\gg r_1\gg 0$,
 and assume that $\dist\bl(0,\{\vv,\p D\}\br)>r_3$.
Suppose that we condition on $\evZa$; that is, on $\zz=0$.
Then with high conditional probability $|\zzp|>2\,r_2$
and moreover $\dist\bl(0,\{S'_0,S'_1,\dots,\zzp\}\br)>2\,r_2$.
By the heights interface continuity (Proposition~\ref{p.hcont}),
given $\exter {r_1}\gamma$
and $\zzp$ and $\gamma\setminus\exter{r_1}\gamma\subset\bal{r_2}$, the
actual choice of $\gamma\setminus\exter{r_1}\gamma$ can change the value of $y'$ by
very little if $|\zzp|>2\,r_2$.
Thus, we conclude that $y'\,1_{\zzp\notin\p D}$ is nearly
independent from $y$ given $\evZa$, $\neg \ev J(r_1,r_2)$ and $\exter {r_1}\gamma$.
Since $y$ and $y'$ are bounded, in the limit as $r_3\to\infty$
\begin{multline}
\label{e.yyp}
 \Eb{y\,y'\,1_{\{\zzp\notin\p D\}}\md \evZa,\exter {r_1}\gamma,\neg\ev J}
\\
= o(1)+
\Eb{y'\,1_{\{\zzp\notin\p D\}}\md \evZa,\exter {r_1}\gamma,\neg\ev J}\,
 \Eb{y\md \evZa,\exter {r_1}\gamma,\neg\ev J},
\end{multline}
where $\ev J=\ev J(r_1,r_2)$.
By Lemma~\ref{l.direct}, $\Pb{\ev J\md\evZa}\to 0$ as $r_2/r_1\to\infty$.
Thus
$$
 \Eb{y\md \evZa,\exter {r_1}\gamma,\neg\ev J}-
 \Eb{y\md \evZa,\exter {r_1}\gamma}\to 0
$$
in probability as $r_2/r_1\to\infty$. A similar remark applies to
the other terms in~\eref{e.yyp}. Since $y$ and $y'$ are bounded,
taking conditional expectation given $\evZa$ in~\eref{e.yyp} gives
\begin{multline}
\label{e.yyy}
 \Eb{y\,y'\,1_{\{\zzp\notin\p D\}}\md \evZa}
\\
=
o(1)+
\EB{\Es{y'\,1_{\{\zzp\notin\p D\}}\md \evZa,\exter {r_1}\gamma}\,
 \Es{y\md \evZa,\exter {r_1}\gamma}\md \evZa}.
\end{multline}

By the Limit Existence Theorem~\ref{t.limit}, given $\evZa$ and $\exter{r_1}\gamma$,
near $0$ the path $\gamma_*$ is
close in distribution to $\dot\gamma$ when $r_1$ is large.
Consequently, Proposition~\ref{p.hcont} implies that
$\Eb{|h(0)|\md \evZa,\exter {r_1}\gamma}-\hco=o(1)$ as $r_1\to\infty$,
which gives
$$
\Eb{y\md \evZa,\exter {r_1}\gamma}=o(1)\,.
$$
Now~\eref{e.yyy} implies that
\begin{equation}
\label{e.yyz}
 \Eb{y\,y'\,1_{\{\zzp\notin\p D\}}\md \evZa} \to 0
\end{equation}
as $\dist\bl(0,\p D\cup\{\vv\}\br)\to\infty$.
Lemma~\ref{l.farbd} and Lemma~\ref{l.farvv} tell us that
for $r_0<\infty$ fixed
$\Pb{\zz\notin\p D,
\dist\bl(\zz,\p D\cup\{\vv\}\br)
<r_0}\to 0$ as $\dist(\vv,\p D)\to\infty$.
Since there is nothing special about the vertex at $0$,
except for our assumption that $\dist\bl(0,\p D\cup\{\vv\}\br)$ is large,
we conclude from~\eref{e.yyz} that
$$
 \lim\Eb{y\,y'\,1_{\{\zz,\zzp\notin\p D\}}\md \zz} = 0
$$
in probability, as $\dist(\vv,\p D)\to\infty$.
Now,~\eref{e.xx} implies that $\Eb{X^2}\to 0$, since $y\,y'\,1_{\{\zz,\zzp\notin\p D\}}$ is bounded.
This gives~\eref{e.bdvalR} and completes the proof.
\QED

Let $F$ be the function that is equal to $h_\p$ on $\p D$, $\hco$ on
$V_+(\gamma)$, $-\hco$ on $V_-(\gamma)$ and is discrete harmonic on
all other $\TG$ vertices in $D$. Since $\Eb{h(\vv)\md \gamma} =
\Eb{h(\zz)\md\gamma}$, Proposition~\ref{p.sides} gives
\begin{equation}
\label{e.hF}
\Eb{h(\vv)\md\gamma}-F(\vv)\to 0
\end{equation}
in probability as $\dist(\vv,\p D)\to \infty$.

\medskip

We now need to generalize the Proposition and~\eref{e.hF} to apply
when $\gamma$ is replaced by an appropriate initial segment
of $\gamma$.

Let $T$ be some stopping time for $\gamma$ started at $x_\p$
and let $\gamma^T$ denote $\gamma$ stopped at $T$.
(Note that the relevant filtration here, the one generated by
intitial segments of $\gamma$, only reveals the signs
of $h$ on vertices adjacent to these initial segments, but not
the actual values of $h$.)
Let $\zzT$ denote the vertex in $\p D(\gamma^T)$ first visited
by $S$, and
let $S^T$ denote the initial segment of $S$ up to its
first visit to $\zzT$. \label{zzTdef}

\begin{lemma}
\label{l.zzT}
Assume~\iref{i.h}, \iref{i.D}, \iref{i.S} and \iref{i.p}.
For every $p_0>0$ there is  some $s=s(p_0,\upperhco)>0$
such that
$$
\Pb{\dist(\zzT,\gamma\setminus\gamma^T)<s\, \dist(\vv,\p D), \zzT\notin\p D}<p_0\,.
$$
\end{lemma}

Note that we could rather easily prove the estimate with
$\dist(\zzT,\gamma\setminus\gamma^T;\DD(\gamma^T))$
instead of $\dist(\zzT,\gamma\setminus\gamma^T)$,
by the argument giving~\eref{e.pen},
but this is not sufficient for our purposes.
The idea of the proof of the lemma is to first show that $\zz=\zzT$ is usually
not too unlikely given $\gamma^T$ and $S$.
Then Lemma~\ref{l.direct} may be used in conjunction with the argument
giving~\eref{e.pen} to deduce the required result.
Note, that the event $\zz=\zzT$ is the event that $\gamma\setminus\gamma^T$
is not adjacent to any vertex visited by $S$ prior to $\zzT$.

\proof
We first show that for every $\eps>0$ there is a $p>0$ and an $R>0$,
both depending only on
$\eps$ and $\upperhco$ such that
\begin{equation}\label{e.zzTzz}
\PB{\Ps{\zzT=\zz\md S,\gamma^T}<p}<\eps
\end{equation}
holds provided that $\dist(\vv,\p D)>R$.

We choose $\delta=\delta(\eps,\upperhco)>0$ very small.
Set $r=\dist(\vv,\p D)$, and assume that $r> 100\,\delta^{-2}$.
Set $D_T:=D(\gamma^T)$.
Let $b_T$ be the point
on $\p D_T$ near the tip of $\gamma^T$ that is at equal distance
from $V_+(\gamma^T)$ and $V_-(\gamma^T)$ along $\p D_T$.
Let $\p_+^T$ and $\p_-^T$ denote the two connected components
of $\p D_T\setminus \{y_\p,b_T\}$ that have $y_\p$ and $b_T$ as their
endpoints, with $\p_+^T$ the one containing vertices in
$\p_+$.

Let $\ev A_1$ be the event
$\dist(\zzT,\p_+^T; D_T)\vee\dist(\zzT,\p_-^T;D_T)\ge 2\, \delta\,r$,
let $\ev A_2$ be the event
$\diam(S^T)< \delta ^{-1} \,r$,
let $\ev A_3$ be the event that
the diameter of the segment of $S^T$ after the first time at which
it is distance at most $\delta^2\,r$ from $\p D_T$ is less than $\delta\,r/2$,
and let $\ev A_4$ be the event $\dist(\vv,\gamma^T)\ge \delta^{1/2}\,r$.

Lemma~\ref{l.farvv} shows that if $\delta=\delta(\upperhco,\eps)$ is
sufficiently small, then $\Pb{\neg\ev A_4}<\eps/4$.
Lemma Hit Near~\ref{l.hitnear} implies that, by choosing $\delta$ sufficiently small,
one can ensure that $\Pb{\neg\ev A_j\md\gamma^T}<\eps/4$ for $j=2,3$.
We now prove the same for $j=1$.
Assume that $\ev A_4$ holds.
Let $L$ be the set of points in $D_T$ that lie on a path of diameter at
most $4\,\delta\,r$ in
$\closure D_T$ connecting $\p_+^T$ and $\p_-^T$,
and let $D_*$ be the connected component of $D_T\setminus L$
that contains $\vv$ (we know that $\vv\notin L$, since
$\delta$ is small and $\ev A_4$ holds).
By~\eref{e.twoballs} applied to $D_T$ in place of $D$,
$\p D_*\cap D_T$ may be covered by two balls
of radius $8\,\delta\,r$. Thus, Lemma Hit Near~\ref{l.hitnear}
shows that if $\delta=\delta(\eps)$ is chosen sufficiently
small $\Pb{S\text{ hits }L\text{ before }\p D_T,\,\ev A_4}<\eps/4$.
This implies $\Pb{\neg\ev A_1,\ev A_4}<\eps/4$.
Thus $\Pb{\neg\ev A}<\eps$, where $\ev A:=\ev A_1\cap\ev A_2\cap\ev A_3\cap A_4$.

We now complete the proof of inequality~\eref{e.zzTzz} by showing
that the event $\Pb{\zz=\zzT\md S,\gamma^T}<p$ is
contained in $\neg \ev A$ if $p=p(\upperhco,\delta)>0$
is chosen sufficiently small.
The latter is equivalent to showing that
the random variable
$\Pb{\zz=\zzT\md S,\gamma^T}$ is bounded away
from zero on $\ev A$ by a function of $\delta$ and $\upperhco$.

Suppose that $\ev A$
holds, and that $\zzT\in \p_+^T$.
Let $\bbeta_1$ and $\bbeta_2$ be the two connected components
of $\p_+^T\setminus \{\zzT\}$.
The construction of $\Upsilon_n$ in the proof~\eref{e.pen}
shows that there is a path $\Upsilon_0$ connecting $\bbeta_1$ and
$\bbeta_2$ in $B(\zzT,\delta\,r)\setminus B(\zzT,\delta\,r/2)$
that separates $\zzT$ from $\p_-^T$ in $D_T$
such that $\Upsilon_0$ is an $\bl(s,\diam(\Upsilon_0)\br)$-barrier
for $(D,\gamma_T)$ for every $s\in(0,\delta']$, where $\delta'\in(0,1)$ is a
universal constant.
(The assumption that $\ev A_1$ holds is used here.)
If $S^T$ never visits a vertex adjacent to $\Upsilon_0$,
then $\Upsilon_0$ separates $S^T$ from $\p_-^T$.
In this situation, if $\gamma\setminus\gamma^T$ does
not hit $\Upsilon_0$, then $S^T$ does not visit a vertex
adjacent to it, and therefore $\zzT=\zz$.
Thus, the Barriers Theorem~\ref{t.barrier}
(or Remark~\ref{r.bdbarrier}) applies to give
the needed lower bound on $\Pb{\zz=\zzT\md S,\gamma^T}$
when $S^T$ does not visit a vertex adjacent to $\Upsilon_0$.

Suppose now that $S^T$ does visit vertices adjacent to $\Upsilon_0$.
We can then construct a path $\Upsilon$ whose image is $\Upsilon_0$
as well as all the boundaries of hexagons visited by $S^T$
that are not separated from $\p_-^T$ by $\Upsilon_0$.
Since we are assuming that $\ev A_2$ holds,
$\diam(\Upsilon)\le 2\,\delta ^{-1} \,r$.
Since $\ev A_3$ holds,
$\dist(\Upsilon\setminus\Upsilon_0,\p D_T)>\delta^2\,r/2$
and therefore also $\diam(\Upsilon_0)>\delta^2\,r/2$.
Consequently, $\Upsilon$ is a $(\delta'\,\delta^3/2, 2\,\delta ^{-1} \,r)$-barrier.
Now Theorem~\ref{t.barrier} and Remark~\ref{r.bdbarrier} may be used again to give a similar
lower bound on $\Pb{\zz=\zzT\md S,\gamma^T}$.
As a similar argument applies when $\zz\in\p_-^T$, the proof
of~\eref{e.zzTzz} is now complete.
\smallskip

We now choose $\eps=p_0/9$ and take a $p>0$
and $R>0$ depending only on $\eps$ and $\upperhco$ and satisfying~\eref{e.zzTzz}.
Let $a$ be such that the estimate given in Lemma~\ref{l.direct}
holds with the $p$ there replaced by $p_0\,p/9$.
Let $\delta=\delta(p_0,\upperhco)>0$
be sufficiently small so that $\delta^{-1}>R\vee a$.
We assume that $r>10\,\delta^{-5}$.
For $z\in D$, let $\ev J_{z}$ denote the event that there are
more than two
disjoint arcs in $\gamma$ joining the two circles
$\p B(z,\delta^4\,r)$ and $\p B(z, \delta^5\,r)$,
and let $\ev J_{z}^T$ denote the event that there are
more than two such
arcs in $\gamma^T$.
By the choice of $a$ and $\delta$, the probability that
 $\dist\bl(\zz,\p D\cup\{\vv\}\br)>4\,\delta^4\,r$
and $\ev J_{\zz}$ holds is at most $p_0\,p/9$.
Consequently, the same bound applies for the probability that
$\dist\bl(\zzT,\p D\cup\{\vv\}\br)>4\,\delta^4\,r$,  $\zzT=\zz$
and $\ev J_{\zzT}^T$ holds.
Thus, as the events
$\dist\bl(\zzT,\p D\cup\{\vv\}\br)>4\,\delta^4\,r$ and
$\ev J_{\zzT}^T$ are $(\gamma^T,S)$ measurable,
$$
\begin{aligned}
p_0\,p/9
&
\ge
\Pb{
\zzT=\zz,
\dist(\zzT,\p D\cup\{\vv\})>4\,\delta^4\,r,
\ev J_{\zzT}^T}
\\&
=
\EB{\Ps{
\zzT=\zz,
\dist(\zzT,\p D\cup\{\vv\})>4\,\delta^4\,r
,
\ev J_{\zzT}^T
\md\gamma^T,S}}
\\&
=
\EB{\Ps{
\zzT=\zz\md \gamma^T,S}\,
1_{\{\dist(\zzT,\p D\cup\{\vv\})>4\,\delta^4\,r\}}\,
1_{\ev J_{\zzT}^T}}\,.
\end{aligned}
$$
By~\eref{e.zzTzz} and our choice of $\eps$, we therefore have
\begin{equation}\label{e.kr}
\Pb{
\dist(\zzT,\p D\cup\{\vv\})>4\,\delta^4\,r,
\ev J_{\zzT}^T}\le 2\,p_0/9\,.
\end{equation}

Let $\ev H$ denote the event
$(\gamma\setminus \gamma^T)\cap \p B(\zzT,\delta^5\,r)\ne\emptyset$.
Now condition on $\gamma^T$ and $S$ such that
$\dist(\zzT,\p D)>\delta\,r$ and $\neg \ev J_{\zzT}^T$ holds.
Suppose also that $\dist(\zzT,b_T)> \delta^3\,r$.
Then there are precisely two connected components of
$B(\zzT,\delta^4\,r)\cap D_T$ that intersect $\p B(\zzT,\delta^5\,r)$.
Let $w_1,w_2\in\p B(\zzT,\delta^5\,r)\cap D_T$ be points in
each of these two
connected components.
By constructing barriers as in the proof of~\eref{e.pen}
it is easy to see that if $\delta=\delta(p_0,\upperhco)$
is sufficiently small, then
$$
\Pb{\dist\bl(w_j,\gamma\setminus\gamma^T;D(\gamma^T)\br)\le \delta^4\,r)\md\gamma^T,S}
\le p_0/9
$$
for $j=1,2$.
Now note that if $\dist\bl(w_j,\gamma\setminus\gamma^T;D(\gamma^T)\br)>\delta^4\,r$
for $j=1,2$, then $\neg\ev H$ holds.
Consequently,
$$
\Pb{\ev H, \neg\ev J_{\zzT}^T, \dist(\zzT,b_T)> \delta^3\,r, \dist(\zzT,\p D)>\delta\,r}
\le 2\,p_0/9\,.
$$
We combine this with~\eref{e.kr}, and get
\begin{multline}\label{e.JzzT}
\Pb{\ev H,
\dist(\vv,\zzT)>4\,\delta^4\,r,
  \dist(\zzT,b_T)> \delta^3\,r, \dist(\zzT,\p D)>\delta\,r}
\\
\le 4\,p_0/9\,.
\end{multline}
Since $r=\dist(\vv,\p D)$,
provided that we take $\delta=\delta(p_0,\upperhco)>0$
sufficiently small, Lemma~\ref{l.farvv} gives
$\Pb{\dist(\vv,\zzT)\le\delta\,r}<p_0/9$,
 Lemma Hit Near~\ref{l.hitnear} gives
$ \Pb{ \dist(\vv,\zzT)> \delta\,r, \dist(\zzT,b_T)\le\delta^2\,r} <p_0/9 $
and~\eref{e.M} gives
$\Pb{\dist(\zzT,\p D)\le\delta\,r,\zzT\notin \p D}<p_0/9$.
These last three estimates may be combined with~\eref{e.JzzT}, to yield
$\Pb{\ev H,\zzT\notin\p D}\le 7\,p_0/9$,
which completes the proof.
\QED

We now prove the analog of~\eref{e.hF} with $\gamma^T$ replacing
$\gamma$. Let $F_T$ denote the function that is $+\hco$ on
$V_+(\gamma^T)$, $-\hco$ on $V_-(\gamma^T)$, equal to $h_\p$ on
$\TG$-vertices in $\p D$, and is discrete harmonic at all other
vertices in $\closure D$.

\begin{proposition}\label{p.stoppingbd}
Assume~\iref{i.h}, \iref{i.D}, \iref{i.S} and \iref{i.p}.
$$
\Eb{h(\vv)\md\gamma^T}-F_T(\vv)\to 0
$$
in probability as
$\dist(\vv,\p D)\to\infty$ while $\upperhco$ is held fixed.
\end{proposition}

\proof
Fix $\eps>0$, and set $r:=\dist(\vv,\p D)$.
We have by the heights interface continuity (Proposition~\ref{p.hcont})
$$
\Bl|\Eb{h(\zzT)\md \gamma,\zzT}-\Eb{h(\zzT)\md\gamma^T,\zzT}\Br|<\eps
$$
if $\dist(\zzT,\gamma\setminus\gamma_T)>R_0$, where $R_0=R_0(\eps,\upperhco)$.
Consequently, Lemma~\ref{l.zzT} with $p_0=\eps/\OC(1)>0$ gives
\begin{equation}\label{e.eee}
\EB{\Bl|\Eb{h(\zzT)\md \gamma,\zzT}-\Eb{h(\zzT)\md\gamma^T,\zzT}\Br|}<2\,\eps
\end{equation}
when $r>s^{-1}\,R_0$, and $s$ is as given by the lemma.
(Note that
$\Eb{h(\zzT)\md \gamma,\zzT}=\Eb{h(\zzT)\md\gamma^T,\zzT}$
when $\zzT\in\p D$.)

In the following, we will use a parameter $\delta>0$.
The notation $o(1)$ will be shorthand for any quantity $g$ satisfying
$\lim_{\delta\to 0}\lim_{r\to\infty} |g|=0$
while $\upperhco$ is fixed.
Let $Z$ be a maximal set of $\TG$-vertices in
$D\cap B(\vv,2\,\delta ^{-1} \,r)$ such that the distance between
any two such vertices is at least $\delta^2\,r$ and the distance
between any such vertex to $\p D$ is at least $\delta^3\,r$.
Then $|Z|=O(\delta^{-6})$.
By~\eref{e.hF} (with each $u\in Z$ in place of $\vv$) we therefore have
\begin{equation}\label{e.eachZ}
\Pb{\exists u\in Z: \Es{ h(u)-F(u)\md\gamma}^2> \eps} = o(1)\,.
\end{equation}

Let $t_0$ be the first $t\in\N$ such that
$\dist(S_t, \p D(\gamma^T))\le\delta\,r$.
By Lemma Hit Near~\ref{l.hitnear},
\begin{equation}\label{e.zznearS}
\Pb{\dist(\zzT,S_{t_0};D)\ge \delta^{1/2}\,r}=o(1)\,.
\end{equation}
By~\eref{e.zbdgam},
$$
\Pb{\dist(\zzT,\gamma;D)<\delta^{1/3}\,r,\zzT\in\p D}=o(1)\,,
$$
while Lemma~\ref{l.zzT} gives
$$
\Pb{\dist(\zzT,\gamma\setminus\gamma_T)<\delta^{1/3}\,r,\zzT\notin\p D}=o(1)\,.
$$
Thus,
$$
\Pb{\dist(\zzT,\gamma\setminus\gamma_T;D)<\delta^{1/3}\,r}=o(1)\,.
$$
This and~\eref{e.zznearS} imply
\begin{equation}\label{e.farga}
\Pb{\dist(S_{t_0},\gamma\setminus\gamma^T;D)<\delta^{1/3}\,r/2}=o(1)\,.
\end{equation}
Since $\dist(S_{t_0},\p D(\gamma^T))<\delta\,r$,
this and Lemma Hit Near~\ref{l.hitnear} imply that with probability $1-o(1)$
the $L^1$ norm of the  difference between
the discrete harmonic measure from $S_{t_0}$  on $\p D(\gamma^T)$
and the discrete harmonic measure from $S_{t_0}$  on $\p D(\gamma)$
is $o(1)$.
Because $\Eb{h(S_{t_0})\md\gamma,S_{t_0}}$ is the average of
$\Eb{h(z)\md\gamma}$, where $z$ is selected according to harmonic
measure on $\p D(\gamma)$ from $S_{t_0}$ and similarly for
$\gamma^T$, we conclude from the above and~\eref{e.eee} that
$$
\Eb{h(S_{t_0})\md \gamma^T,S_{t_0}}-
\Eb{h(S_{t_0})\md \gamma,S_{t_0}}=o(1)
$$
in probability.  By Lemma Hit Near~\ref{l.hitnear},
$\Pb{\dist(S_{t_0},\vv)>\delta^{-1}\,r}=o(1)$. On the event
$\dist(S_{t_0},\vv)\le \delta ^{-1} \,r$, fix some $z_0\in Z$ within
distance $\delta^2\,r$ from $S_{t_0}$ (if there is more than one
such $z_0$, let $z_0$ be chosen uniformly at random among these given
$(h,\gamma,S)$). By the Discrete Harnack Principle~\ref{l.harnack}
and~\eref{e.farga}, we have
$\Eb{h(S_{t_0})\md\gamma,S_{t_0}}-\Eb{h(z_0)\md\gamma,z_0}=o(1)$ in
probability, which in conjunction with~\eref{e.eachZ} yields $
\Eb{h(S_{t_0})\md\gamma,S_{t_0}}- F(z_0)=o(1)$. The Discrete Harnack
Principle now implies $F(z_0)-F(S_{t_0})=o(1)$ in probability,
and~\eref{e.farga} gives $F(S_{t_0})-F_T(S_{t_0})=o(1)$ in
probability. Consequently, $\Eb{h(S_{t_0})\md\gamma^T,S_{t_0}}-
F_T(S_{t_0})=o(1)$ in probability. Since $\Eb{h(\cdot)\md\gamma^T}$
and $F_T$ are discrete harmonic in $D(\gamma^T)$, the proposition
follows. \QED

We can now prove the height gap theorem:

\begin{theorem}\label{t.gap}
Assume~\iref{i.h}, \iref{i.D}, \iref{i.S} and \iref{i.p}.
As above, let $T$ denote a stopping time for $\gamma$,
let $\gamma^T:=\gamma[0,T]$, let
$D_T$ be the complement in $D$ of the
closed triangles meeting $\gamma[0,T)$.
Let $h_T$ denote the restriction of $h$
to $V\cap \p D_T$ and let $\vv$ be some
vertex in $D$.
Then
$$
\Eb{h(\vv)\md\gamma^T,\, h_T}-F_T(\vv)\to 0
$$
in probability as
$\dist(\vv,\p D)\to\infty$ while $\upperhco$ is held fixed,
where $F_T$ is as in Proposition~\ref{p.stoppingbd}.
\end{theorem}

\proof
Set
$$
X:=
\Eb{h(\vv)\md\gamma^T,\, h_T}-
\Eb{h(\vv)\md\gamma^T}.
$$
By Proposition~\ref{p.stoppingbd}, it suffices to show
that $X\to 0$ in probability as
$\dist(\vv,\p D)\to\infty$.
For $v\in V\cap \p D_T$, let $a_v$ denote the
conditional probability that simple random walk
started at $\vv$ first hits $\p D_T$ at $v$,
given $\gamma^T$.
Then
$$
X = \sum_{v\in V\cap \p D_T}
a_v \bigl(h(v)-\Es{h(v)\md \gamma^T}\bigr).
$$
Consequently,
$$
X^2 = \sum_{v,u} a_v\,a_u
\bigl(h(v)-\Es{h(v)\md \gamma^T}\bigr)
\bigl(h(u)-\Es{h(u)\md \gamma^T}\bigr).
$$
Now Corollary~\ref{c.corr} implies
that it suffices to show that for every $r>0$
$$
\sum_{v,u} a_v\,a_u\,1_{\{|v-u|<r\}}\,1_{\{v,u\notin V_\p\}}\to 0
$$
in probability.
This follows by Lemmas~\ref{l.hitnear} and~\ref{l.farvv}.
\QED

\section{Recognizing the driving term}\label{s.drive}

In this section we use a technique introduced in~\cite{\LSWlesl}
and used again in~\cite{\SchrammSheffieldHE} in order to show that
the driving term for the Loewner evolution given by the DGFF interface
with boundary values $-a$ and $b$ converges to
the driving term of \SLEab/ if $a,b\in[-\lowerhco,\upperhco]$.
The reader unfamiliar with this method is advised to
first learn the technique from~\cite[\S4]{\SchrammSheffieldHE}
or~\cite[\S3.3]{\LSWlesl}.
The account in~\cite{\SchrammSheffieldHE} is closer to the present
setup and somewhat simpler, but some parts of the argument there
are referred back to~\cite{\LSWlesl}.

The present argument is more involved than those of the
above mentioned papers, because we prove convergence to
an instance of \SLEr \rho_1,\rho_2/  rather than just
plain \SLEkk4/. The main added difficulty comes from the fact that
the drift term in \SLEr\rho_1,\rho_2/ becomes unbounded
as $W_t$ comes close to the force points.
These difficulties disappear if $a=b=\hco$, in which case
the convergence is to ordinary SLE and the argument giving the
convergence of the driving term to scaled Brownian motion
is easily established with minor adaptations of the established
method. We therefore forego dwelling on this simpler case, and
move on to the more general setting, assuming that the reader is
already familiar with the fundamentals of the method.

\subsection{About the definition of \SLEr \rho_1,\rho_2/}

Throughout this section, given a Loewner evolution defined by a
continuous $W_t$, we will let $x_t$ and $y_t$ be defined as in
Section \ref{ss.definingsle} by $x_t:=\sup\{g_t(x):x<0, x\notin
K_t\}$ and $y_t:=\inf\{g_t(x):x>0, x\notin K_t\}$, and we make use
of the definition of \SLEr \rho_1, \rho_2/ in Section~\ref{ss.definingsle}
by means of the SDE \eref{e.sdesystem}.  As we
mentioned in Section \ref{ss.definingsle}, some subtlety is
involved in extending the definition of \SLEr \rho_1, \rho_2/
 beyond times when
$W_t$ hits the force points, and in starting the process from the
natural initial values $x_0=W_0=y_0 = 0$. This is closely related
to the issues involved in defining the Bessel process, which we
presently recall.

The Bessel process $Z_t$ of dimension $\delta > 0$ and initial value $x\ne 0$ satisfies the SDE
\begin{equation}\label{e.BesselSDE} dZ_t = \frac{\delta-1}{2Z_t}\,dt + dB_t, \qquad Z_0 = x\,,
\end{equation} which we also write in integral form as \begin{equation}\label{e.integralBesselSDE}
Z_t = x + \int_0^t \frac{\delta-1}{2Z_s}\,ds + B_t-B_0\,, \end{equation}
up until the first time $t$ for which $Z_t = 0$.  When defining $Z_t$ for all times, this SDE
is awkward to work with directly since the drift blows up whenever $Z_t$ gets
close to zero (and some of the standard existence and uniqueness
theorems for SDE solutions, as given, e.g., in~\cite{\RevuzYor}, do not apply
in this situation).  However, for every $\delta > 0$, the square of
the Bessel process $Z_t^2$ turns out to satisfy an SDE whose drift remains bounded
and for which existence and uniqueness of solutions follow easily from standard theorems.
For this reason, many authors construct the Bessel process by first defining
the square of the Bessel process via an SDE that it satisfies and then taking its square root~\cite{\RevuzYor}.
(Recall also that when $\delta \leq 1$ the Bessel itself does not satisfy \eref{e.integralBesselSDE}
at all without a
principal value correction.  Even when $1 < \delta < 2$, which, as we will see
below, is the case that corresponds to
\SLEr \rho/ that hit the boundary and can be continued after hitting the boundary,
the solution to \eref{e.integralBesselSDE} is not unique
unless we restrict attention to non-negative solutions.)

The formal definition for \SLEr\rho/ with one force point (i.e., $\rho_1 = \rho$ and
$\rho_2 = 0$) was given in~\cite{\LSWrestriction}.  It was observed there
that in this case,  \eref{e.sdesystem} implies that the process $W_t-x_t$
satisfies the same SDE as the Bessel process of dimension
$\delta = 1 + \frac{2(\rho+2)}{\kappa}$ up until the first time $t$ for which $W_t=x_t$.
Thus, to define \SLEr\rho/, the paper~\cite{\LSWrestriction} starts with a
constant multiple of a Bessel
process $Z_t$ of the appropriate dimension and defines
the evolution of the force point $x_t$ by $x_t=x_0+\int_0^t  2\,ds/Z_s$ and
the driving term by $W_t=x_t+Z_t$.

Defining \SLEr \rho_1,\rho_2/ is a slightly more delicate matter
since neither $W_t-x_t$ nor $y_t-W_t$ is
 exactly a Bessel (although each one is
quite close to a Bessel when the other
force point is relatively far away).
Although this is not a very difficult issue,
 it seems that there does not yet exist in the literature an
adequate definition of \SLEr \rho_1,\rho_2/ that is valid
beyond the time that the driving term hits a force point.
Since we prove the convergence to \SLEr \rho_1,\rho_2/, we have to
define it.

The approach we adopt
is basically similar to the way in which the Bessel process (and hence
\SLEr \rho/) is usually defined:
we pass to a coordinate system in which the corresponding SDE becomes tractable.
We will describe the coordinate change we use in Section \ref{s.coord}.
Within this new coordinate system, we then prove the convergence of the
Loewner driving parameters of our discrete processes to those of the corresponding
\SLEr \rho_1,\rho_2/ in Sections \ref{ss.loewnerfordgff} and \ref{s.approximatediffusions}.
Section \ref{s.backtochordal} then describes the reverse coordinate transformation
and use it to give a formal definition of \SLEr \rho_1,\rho_2/, Definition \ref{d.slekr}.

We remark that there are many equivalent ways to define \SLEr \rho_1,\rho_2/
(for example, one can probably show directly that \eref{e.sdesystem} has a unique strong solution
for which $x_t \leq W_t \leq y_t$ for all $t$), but ours seems most efficient given that
the coordinate change also simplifies the proofs in
Sections \ref{ss.loewnerfordgff} and \ref{s.approximatediffusions}.

\subsection{A coordinate change}\label{s.coord}

In this subsection, we recall a different coordinate system for Loewner evolutions,
which is virtually identical to the setup used in~\cite[\S3]{\LSWiii}.
Suppose that $\gamma:[0,\infty)\to\closure\H$ is a continuous simple path that starts at
$\gamma(0)=0$, does not hit $\R\setminus\{0\}$, satisfies
$\lim_{t\to\infty}|\gamma(t)|=\infty$ and is parameterized by half-plane
capacity from $\infty$.
Let $g_t:\H\setminus \gamma[0,t]\to\H$ be the conformal map
satisfying the hydrodynamic normalization at $\infty$,
let $W_t=g_t(\gamma(t))$ be the corresponding Loewner driving term.
Loewner's theorem says that $g_t$ satisfies Loewner's chordal equation~\eref{e.chordal}.
Now we introduce a one parameter family of maps $G^*_t:\H\setminus\gamma[0,t]\to\H$
satisfying the normalization for $t>0$
$$
G^*_t(\infty)=\infty\,,\qquad G^*_t((0,\infty))=(1,\infty)\,,\qquad G^*_t((-\infty,0))=(-\infty,-1)\,.
$$
That is,
\begin{equation}\label{e.Gg}
G^*_t(z)=\frac {2\,g_t(z)-x_t-y_t}{y_t-x_t}\,,
\end{equation}
where $x_t$ and $y_t$ (as defined earlier) are the two images under $g_t$ of $0$
and $x_t<y_t$.
Set
$$
W^*_t:=G^*_t(\gamma(t))= \frac{2 \,W_t-x_t-y_t}{y_t-x_t}\,.
$$
By differentiating~\eref{e.Gg} and using~\eref{e.chordal} and~\eref{e.xy} it is immediate to verify
that $G^*_t$ satisfies
$$
\frac{dG^*_t(z)}{dt} =
\frac{8}{(y_t-x_t)^2}\,
\frac{1-G^*_t(z)^2}{(G^*_t(z)-W^*_t)(1-(W^*_t)^2)}\,.
$$
We now define a new time parameter
$$
s(t)=\log(y_t-x_t)=\log 2-\log {G^*_t}'(\infty)\,.
$$
It is easy to verify that $s$ is continuous and monotone increasing and
$s\bl((0,\infty)\br)=(-\infty,\infty)$.
Set $G_s=G^*_t$ and $\tilde W_s=W^*_t$ when $s=s(t)$.
Differentiation gives
\begin{equation}
\label{e.ds}
ds= \frac{8\,(y_t-x_t)^{-2}\,dt}{1-(W^*_t)^2}\,.
\end{equation}
Consequently, this change of time variable
allows us to write the ODE satisfied by $G$ as
\begin{equation}\label{e.dG}
\frac{dG_s(z)}{ds}=
\frac{1-G_s(z)^2}{G_s(z)-\tilde W_s}\,,
\end{equation}
where all the terms come from the new coordinate system.
Later, in subsection~\ref{s.backtochordal},
we explain how to go back to the standard chordal coordinate system.

\subsection{The Loewner evolution of the DGFF interface} \label{ss.loewnerfordgff}

In addition to our previous assumptions~\iref{i.h} and~\iref{i.D}
about the domain $D$ and the boundary conditions, we
now add the assumption that
\begin{enumerate}
\hitem{i.ab}{(ab)} there are constants $a,b$ such that
$h_\p=b$ on $\p_+$, $h_\p=-a$ on $\p_-$
and $\min\{a,b\}>-\lowerhco$, where $\lowerhco>0$
is given by Lemma~\ref{l.bdnarrows} with
$\upperhco:=\max\{|a|,|b|\}$.
\end{enumerate}
In this case, clearly~\iref{i.p} holds.
In the following, $a$ and $b$ will be considered as constants,
and the dependence of various constants on $a$ and $b$ will sometimes
be suppressed (for example, when using the $O(\cdot)$ notation).

Let
$\phi:D\to \H$ be a conformal map that corresponds $\p_+$ with the positive real ray.
Let $\gamma$ be the zero height interface of $h$ joining the endpoints
of $\p_+$, and let $\gamma^{\phi}$ denote the image of $\gamma$ under
$\phi$.
Now, $\gamma^{\phi}$ satisfies the assumptions in the
previous subsection. Consequently, we may parameterize $\gamma^\phi$
according to the time parameter $s=s(t)$ and consider
the conformal maps $G_s:\H\setminus\gamma^\phi(-\infty,s]\to \H$
as defined in \S\ref{s.coord}.
As above, we set $\tilde W_s=G_s(\gamma^\phi(s))$ and have the
differential equation~\eref{e.dG}.
Our goal now is to determine the limit of the law of $\tilde W$
as $\inr{\phi^{-1}(i)}(D)\to\infty$.
Set for $x\in[-1,1]$,
\begin{equation}\label{e.q}
q_1(x):=2\,(1-x^2),
\qquad
q_2(x):=-\frac{\hco+a}{2\hco}\,(x-1)-\frac{\hco+b}{2\hco}\,(x+1)\,.
\end{equation}
We extend the definitions of $q_1$ and $q_2$ to all
of $\R$ by taking each $q_j$ to be constant in each of
the two intervals $(-\infty,-1]$ and $[1,\infty)$.  Consider the SDE
\begin{equation}\label{e.SDE}
dY_s = q_2(Y_s)\,ds + q_1(Y_s)^{1/2}\,d B_s\,,
\end{equation}
where $B$ is a standard one-dimensional Brownian motion.
A weak solution is known to exist (see~\cite[\S 5.4.D]{\KaratsasShreve}).
We also recall that the weak solution is strong and pathwise
unique (see~\cite[\S IX, Theorems 1.7 \& 3.5]{\RevuzYor}).

\begin{theorem}\label{t.lim}
There is a time-stationary solution $Y:(-\infty,\infty)\to[-1,1]$
of~\eref{e.SDE}.
Moreover, for every finite $S>1$ and $\eps>0$
there is an $R_0=R_0(S,\eps)$ such that
if $R:= \inr{\phi^{-1}(i)}(D)>R_0$ and
the assumptions~\iref{i.h}, \iref{i.D} and~\iref{i.ab} hold,
then there is a coupling of $Y_s$ with $h$ such that
$$
\PB{\sup\bl\{|Y_s-\tilde W_s\br|:s\in[-S,S]\br\} >\eps}<\eps\,.
$$
\end{theorem}

The following proposition is key in the
proof of the theorem.  In essence, it states that $\tilde W_s$ satisfies a discrete
version of~\eref{e.SDE}.
Let $\ev F_{s}$ be the $\sigma$-field generated by
$(\tilde W_r:r\le s)$.
(Note that although the filtration defining $\tilde W_s$ is discrete, there
is no problem in considering $\ev F_s$ for arbitrary $s$,
though the behavior of $\tilde W_r$ for $r$ in some neighborhood of $s$
might be determined by $\ev F_s$.)

\begin{proposition}\label{p.diffusion}
Assume~\iref{i.h}, \iref{i.D} and~\iref{i.ab}.
Fix some $S>1$ large and some $\delta,\errprb>0$ small.
There is a constant $C>0$, depending only on $a,b$ and $S$,
and there is a
function $R_0=R_0(S,\delta,\errprb)$, depending only on $a,b,S,\delta$ and
$\errprb$, such that the following holds.
If $R:= \inr{\phi^{-1}(i)}(D)>R_0$
and $s_0,s_1$ are two stopping times for $\tilde W_s$ such that a.s.\
$-S\le s_0\le s_1\le S$, $\Delta s:=s_1-s_0\le \delta^2$ and
$\sup_{s\in[s_0,s_1]} |\tilde W_s-\tilde W_{s_0}|\le\delta$,
then the following two estimates hold with probability at least $1-\errprb$
\begin{align}
\label{e.dW}
\Bigl|
\Eb{
\Delta\tilde W
-
q_2(\tilde W_{s_0})\,\Delta s
\,
 \md \ev F_{s_0}}\Bigr|
&
\le C\,\delta^3,
\\
\label{e.dWW}
\Bigl|\Eb{
(\Delta \tilde W)^2-
q_1(\tilde W_{s_0})
\,
\Delta s
 \md \ev F_{s_0}}\Bigr|
&
\le C\, \delta^3,
\end{align}
where $\Delta\tilde W:=\tilde W_{s_1}-\tilde W_{s_0}$.
\end{proposition}

To prepare for the proof of the proposition, we need the following easy lemma.
The first two statements in this lemma should be rather obvious to anyone with
a solid background on conformal mappings.

\begin{lemma}\label{l.choosept}
Set $c_1=c_1(S)=100\,e^S$.
There are finite constants $c_2=c_2(S)>0$
and $R_0=R_0(S)>0$, depending
only on $S$, such that
if $R:= \inr{\phi^{-1}(i)}(D)>R_0$ and
if $z\in\H$ satisfies $5\,c_1\ge\Im z \ge c_1 \ge |\Re z|$,
then the following holds true:
\begin{enumerate}
\item\label{i.far} $\inr{\phi^{-1}(z)}(D)> c_2\,R$,
\item\label{i.a}
      there is a $\TG$ vertex $v\in D$ satisfying $|\phi(v)-z|<|z|/100$,
\item\label{i.bdd}
  $\Im G_s(z)\ge e^{-s}\, c_1/2$ for $s\in[-S,S]$
      (and in particular, $G_s(z)$ is well defined in that range),
   and
\item\label{i.bdd1}
  $|G_s(z)-2\,e^{-s}\,z|\le 2$ for $s\in[-S,S]$.
\end{enumerate}
\end{lemma}
\proof
Consider the conformal map $\psi(z)=(z-i)/(z+i)$ from $\H$ onto the unit disk $\U$
taking $i$ to $0$ and set $f=(\psi\circ\phi)^{-1}$.
The Schwarz lemma applied to
the map $z\mapsto f^{-1}\bl(f(0)+R\,z\br)$
restricted to $\U$ gives $1/|f'(0)|=|(f^{-1})'(f(0))|\le 1/R$.
Thus $|f'(0)|\ge R$.
For a fixed $c_1$ the set of possible $z$ is a compact subset of $\H$,
and its image under $\psi$ is a compact subset of $\U$.
Consequently, the Koebe distortion theorem (see, e.g.,~\cite[Theorem 1.3]{\PommeBDRY})
implies that $|f'(\psi(z))|\ge c_2'\,R$ for some $c_2'$ depending only on
$c_1$.
Now the Koebe $1/4$ theorem (see, e.g.,~\cite[Corollary 1.4]{\PommeBDRY}) gives
$\inr{f(\psi(z))}(D)> c_2\,R$ for some $c_2$ depending on $c_1$.
This takes care of Statement~\ref{i.far}.

Let $B$ be the open disk of radius $|z|/200$ about $z$.
Clearly, $B\subset\H$.
We conclude from $|f'(\psi(z))|\ge c_2'\,R$ that
$|(\phi^{-1})'(z)|> c_2''\,R$ for some $c_2''$ depending only on $c_1$.
Thus, the Koebe $1/4$ theorem  implies
$$
\inr {\phi^{-1}(z)} (\phi^{-1}(B))\ge c_2''\,\rad(B)\,R/4\,.
$$
Consequently, Statement~\ref{i.a} holds once $R_0>4/(c_2''\,\rad(B))$.
This takes care of~\ref{i.a}, because $1/\rad(B)$ is bounded
by a function of $c_1$.

It is easy to check that~\ref{i.bdd} follows from~\ref{i.bdd1}.
It remains to prove the latter.
Let $x_t$ and $y_t$ be as in \S\ref{s.coord}.
Note that $x_t<W_t<y_t$ for all $t> 0$
and $\lim_{t\searrow 0} x_t=W_0=0=\lim_{t\searrow 0} y_t$.
Therefore,~\eref{e.xy} implies $x_t<0<y_t$ for all $t>0$.
By~\eref{e.xy},
$$
\p_t(y_t-x_t)= \frac {2\,(y_t-x_t)}{(y_t-W_t)\,(W_t-x_t)}
\ge \frac 8{y_t-x_t}\,.
$$
Therefore,
$\p_t \bl((y_t-x_t)^2\br) \ge 16$, which gives
\begin{equation}
\label{e.ts}
  e^{2\,s(t)}\ge 16\,t\,.
\end{equation}

Observe by~\eref{e.chordal}
$$
\p_t \bl((\Im g_t(z))^2\br) \ge -4\,.
$$
Thus, $\bl(\Im g_t(z)\br)^2 \ge (\Im z)^2-4\,t\ge c_1^2-4\,t$.
Another appeal to~\eref{e.chordal} now gives
$$
\bl|g_t(z)-z\br|\le \frac{2\,t}{\sqrt{c_1^2-4\,t}}\,.
$$
By~\eref{e.Gg} and the definition of $s(t)$, the above gives
$$
\bl| G_s(z)-2\,e^{-s}\,z\br|
\le
\Bl|\frac{x_t+y_t}{y_t-x_t}\Br|+
\frac{2\,t\,e^{-s}}{\sqrt{c_1^2-4\,t}}\,.
$$
Now, the first summand on the right hand side is at most $1$, because $y_t>0>x_t$.
The second summand is also at most $1$ in the range
$s\in[-S,S]$, by~\eref{e.ts} and the choice of $c_1$.
This completes the proof of the lemma.
\QED

\proofof{Proposition~\ref{p.diffusion}}
With the notations of Lemma~\ref{l.choosept}, let $z_j:=2\, c_1\,i+c_1\,j/2$, $j=0,1$,
and let $v_j$ be a $\TG$-vertex satisfying Condition~\ref{i.a} of the lemma with $z_j$ in place of $z$.
Then $z_j':=\phi(v_j)$ satisfies in turn the assumptions required of $z$ in the lemma.
For $k=0,1$, note that there is a stopping time $T_k$
for $\gamma$ such that 
$\phi\circ\gamma(T_k)=\gamma^\phi(s_k)$.
Fix $j\in\{0,1\}$ and set
$$
X_k=X_k(j):= \Eb{h(v_j)\md \ev F_{s_k}}\,, \qquad k=0,1.
$$
Clearly,
\begin{equation}
\label{e.mart}
\Eb{X_1\md\ev F_{s_0}}=X_0\,.
\end{equation}
Recall the definition of the function $F_T$ from
Proposition~\ref{p.stoppingbd}.
Let $\ev A_k$ be the event $|X_k-F_{T_k}(v_j)|\ge\delta^5$.
By that proposition and the fact that $z_j'$ satisfies
Condition~\ref{i.far} of Lemma~\ref{l.choosept}, if $R$ is
chosen sufficiently large then $\Pb{\ev A_k}<\errprb\,\delta^5/4$.
Since $F_{T_k}(v_j)=O_{a,b}(1)$, and likewise
$X_k=O_{a,b}(1)$ by~\eref{e.atv}, we get
$$
\Bigl|
\Eb{X_{1}-F_{T_{1}}(v_j)\md
\ev F_{s_0}}
\Bigr|
\le \delta^5+ O_{a,b}(1)\,
\Pb{\ev A_{1} \md \ev F_{s_0}}.
$$
Let $\tilde{ \ev A}$ be the event that
$\Pb{\ev A_{1} \md \ev F_{s_0}}>\delta^5$.
Then $\Pb{\tilde{\ev A}}<\errprb/4$ (since we are assuming
$\Pb{\ev A_1}<\errprb\,\delta^5/4$) and we have
$$
\Bigl|
\Eb{X_{1}-F_{T_{1}}(v_j)\md
\ev F_{s_0}}
\Bigr|
\le
O_{a,b}(\delta^5)
\qquad \text{on }\neg \tilde{\ev A}\,.
$$
Thus, we have from~\eref{e.mart}
\begin{equation}
\label{e.mart2}
\Eb{F_{T_{1}}(v_j)\md \ev F_{s_0}}
- F_{T_0}(v_j) =
O_{a,b}(\delta^5)
\qquad \text{on }\neg (\tilde{\ev A}\cup {\ev A_0})\,.
\end{equation}
Let $H_k$ be the
bounded function that is harmonic (not discrete harmonic) in
$D\setminus\gamma[0,T_k]$, has boundary values $b$ on $\p_+$, $-a$
on $\p_-$, $+\hco$ on the right hand side of $\gamma[0,T_k]$, and
$-\hco$ on the left hand side of $\gamma[0,T_k]$. We claim that
$H_k(v_j)-F_{T_k}(v_j)$ is small if $\dist(v_j,\p
D\cup\gamma[0,T_k])$ is large. Indeed, this easily follows by
coupling the simple random walk on $\TG$ to stay with high
probability relatively close to a Brownian motion and
using~\eref{e.zbdgam} and Lemmas~\ref{l.farbd}
and~\ref{l.hitnear} to show
that with high probability the
boundary value sampled by the hitting point of the Brownian motion
is the same as that sampled by the hitting vertex of
the simple random walk. Now,
Lemma~\ref{l.farvv} guarantees that if $R$ is sufficiently large,
then with high probability $\dist(v_j,\p D\cup\gamma[0,T_k])$ is
large as well. Consequently, if $R$ is chosen sufficiently large we
have $\Pb{|H_k(v_j)-F_{T_k}(v_j)|\ge \delta^5}<\errprb\,\delta^5/4$. Let
$\ev B_k$ be the event $|H_k(v_j)-F_{T_k}(v_j)|\ge \delta^5$, let
$\tilde{\ev B}$ be the event $\Pb{\ev B_{1}\md\ev F_{s_0}}\ge
\delta^5$ and let
$\bar{\ev A}:=\ev A_0\cup\tilde{\ev A}\cup\ev B_0\cup\tilde{\ev B}$.
Note that $\Ps{\bar{\ev A}}<\errprb$.
The above proof of~\eref{e.mart2}
from~\eref{e.mart} now gives
\begin{equation}
\label{e.mart3}
\Eb{H_{1}(v_j)\md \ev F_{s_0}}
- H_{0}(v_j) =
O_{a,b}(\delta^5)
\qquad \text{on }\neg \bar{\ev A}\,.
\end{equation}

Now, the point is that $H_k(v_j)$ can easily be expressed analytically in terms
of $Z_k=Z_{k,j}:=G_{s_{k}}(z_j')=G_{s_{k}}\bl(\phi(v_j)\br)$
and $\tilde W_{s_{k}}$.
Indeed, conformal invariance implies that the harmonic measure of $\p_+$ in
$D\setminus\gamma[0,T_k]$ from $v_j$ is the same as the harmonic
measure of $[1,\infty)$ from $Z_k$, which is $1-\arg(Z_k-1)/\pi$, because
$G_s\circ \phi$ corresponds $\p_+$ with $[1,\infty)$.
Likewise, the harmonic measure of the right hand side of $\gamma[0,T_k]$ is
$\bl(\arg(Z_k-1)-\arg(Z_k-\WW_k)\br)/\pi$.
Similar expressions hold for the harmonic measure of the left hand side of
$\gamma[0,T_k]$ and of $\p_-$.
These give
\begin{equation}
\label{e.Hn}
{\pi}\,(H_k(v_j)-b)=
(\hco-b)\,\arg(Z_k-1)-2\,\hco\,\arg(Z_k-\WW_k)+(\hco-a)\, \arg(Z_k+1).
\end{equation}

Recall that $Z_k=G_{s_{k}}(z_j')$.
By~\eref{e.dG}, we have in the interval $s\in[s_0,s_{1}]$
\begin{equation}
\label{e.dGi}
G_s(z_j')=Z_0+
\int_{s_0}^{s}
\frac{1-G_r(z_j')^2}{G_r(z_j')-\tilde W_r}\,dr\,.
\end{equation}
Note that Conditions~\ref{i.bdd} and~\ref{i.bdd1} of Lemma~\ref{l.choosept} imply that the
integrand is $O_S(1)$.
Therefore $G_s(z_j')-Z_0=O_S(\Delta s)=O_S(\delta^2)$ for
$s\in[s_0,s_{1}]$. Moreover, we have
$|\tilde W_s-\WW_0|\le \delta$ in that range.
Thus, it follows from~\eref{e.dGi} and Condition~\ref{i.bdd} of the lemma
that
\begin{equation}
\label{e.dZ}
Z_{1}-Z_0 = \Delta s\,\Bigl(\frac {1-Z_0^2}{Z_0-\WW_0} +  O_S(\delta)\Bigr)
= \Delta s\,\frac {1-Z_0^2}{Z_0-\WW_0} +  O_S(\delta^3)\,.
\end{equation}
We will now write an expression for $H_{1}(v_j)-H_{0}(v_j)$ and then use~\eref{e.mart3}
to complete the proof.
Let us first look at the term $\arg(Z_k-\WW_k)$ on the right hand side
of~\eref{e.Hn}, and see how it changes from $k=0$ to $k=1$.
For this purpose, we expand $\log(Z-\tilde W)$ in Taylor series up to first order
in $Z-Z_0$ and up to second order in $\tilde W-\WW_0$ (since
$Z_{1}-Z_0=O_S(\delta^2)$ while
$\WW_{1}-\WW_{0}=O(\delta)$), as follows:
\begin{multline*}
\arg(Z_{1}-\WW_{1})-
\arg(Z_0-\WW_{0})
=
\Im\bigl(\log (Z_{1}-\WW_{1})-
\log (Z_0-\WW_{0})\br)\\
=
\Im\Bigl(
\frac{Z_{1}-Z_0}{Z_0-\WW_{0}}
-
\frac{\WW_{1}-\WW_0}{Z_0-\WW_{0}}
-
\frac{(\WW_{1}-\WW_0)^2}{2\,(Z_0-\WW_{0})^2}
\Bigr)+O_S(\delta^3)\,.
\end{multline*}
Similar (but simpler) expansions apply to the other arguments in~\eref{e.Hn}.
We use these expansions as well as~\eref{e.Hn} and~\eref{e.dZ} to write
\begin{align*}
&
\pi\,( H_{1}(v_j)- H_0(v_j))=
(\hco-b)\,
\Im\Bigl(
\frac{\Delta s\,(1-Z_0^2)}{(Z_0-\WW_0)(Z_0-1)}
\Bigr) +{}
\\
&\qquad
-2\,\hco\,
\Im\Bigl(
\frac{\Delta s\,(1-Z_0^2)-\Delta \tilde W (Z_0-\WW_0)-(\Delta\tilde W)^2/2}{(Z_0-\WW_{0})^2}
\Bigr)+{}\\
&\qquad
{}+
(\hco- a)\,
\Im\Bigl(
\frac{\Delta s\,(1-Z_0^2)}{(Z_0-\WW_0)(Z_0+1)}
\Bigr) + O_S(\delta^3)\,.
\end{align*}
With the abbreviations $\tilde x:= \Re(Z_0-\WW_0)$, $y:= \Im Z_0$,
the above simplifies to
\begin{align*}
&\frac{\pi}{\hco}\,( H_{1}(v_j)-H_0(v_j))=\\
&\quad
\frac {2\,y}{\tilde x^2+y^2}\,\Bigl(
\frac{\tilde x}{\tilde x^2 + y^2}\,\bl( \Delta s\,q_1(\WW_0)-(\Delta\tilde W)^2\br)
+
\Delta s\,q_2(\WW_0)
-\Delta\tilde W\Bigr) + O_S(\delta^3)\,.
\end{align*}
We know from~\eref{e.mart3} that on $\neg \bar{\ev A}$ the conditioned expectation given $\ev F_{s_0}$
of the left hand side is $O_{a,b}(\delta^5)$.
Since $(\tilde x^2+y^2)/y=O_S(1)$ (by Statement~\ref{i.bdd1} of Lemma~\ref{l.choosept}),
we have on $\neg\bar{\ev A}$
\begin{equation}\label{e.ee}
\EB{
\frac{\tilde x}{\tilde x^2 + y^2}\,\bigl(
\Delta s \, q_1(\WW_0)- (\Delta\tilde W)^2\bigr)
+
\Delta s\, q_2(\WW_0)
-\Delta\tilde W \md \ev F_{s_0}}
=O_{a,b,S}(\delta^3)\,.
\end{equation}

Now, this is valid for $z_j'$, with $j=0,1$.
The choice of $j$ only affects the left hand side in the term
$\tilde x/(\tilde x^2+y^2)$.
By the choice of the points $z_j'$ and by Statement~\ref{i.bdd1}
of Lemma~\ref{l.choosept}, the factor
${\tilde x}/(\tilde x^2 + y^2)=\Re\bl((Z_0-\WW_0)^{-1}\br)$ differs between the two $z_j'$ by an
amount that is bounded away from zero by a constant depending on $S$.
Subtracting the above relation~\eref{e.ee} for $z_0'$ from that of $z_1'$,
we therefore get~\eref{e.dWW} on $\neg\bar{\ev A}$.
When this is used in conjunction with~\eref{e.ee} again, one obtains~\eref{e.dW}
on $\neg\bar{\ev A}$.
This concludes the proof of the proposition.
\QED

\subsection{Approximate diffusions} \label{s.approximatediffusions}

In this subsection we embark on the general study of
random processes satisfying the conclusions of Proposition~\ref{p.diffusion}
and show that the proposition essentially characterizes the macroscopic
behavior of the process.
As one of the referees of this paper pointed out, one can try
to do this more \lq\lq traditionally\rq\rq\ by
proving tightness of the driving term and characterizing
the subsequential fine mesh limit using the appropriate martingale problem.
However, our approach is somewhat different (though not necessarily better).
Motivated by the proposition, we say that a continuous
random $W:[0,S]\to [-1,1]$ is a {\it $(C,\delta,\errprb)$-approximate
$(q_1,q_2)$-diffusion\/} if it satisfies the conclusion of the proposition; namely,
for every pair of stopping times $s_0$ and
$s_1$ such that a.s.\ $0\le s_0\le s_1\le S$, $\Delta s:=s_1-s_0\le \delta^2$
and $\sup_{s\in[s_0,s_1]} |W_s-W_{s_0}|\le \delta$ we have with probability at least
$1-\errprb$ that~\eref{e.dW} and~\eref{e.dWW} hold with $W$ in place of
$\tilde W$.

\begin{lemma}\label{l.Ys}
Suppose that $a,b\ge -\hco$ and that $Y_s$ satisfies~\eref{e.SDE},
where $q_1$ and $q_2$ given by~\eref{e.q}.
Suppose that $Y_0\in[-1,1]$ a.s.
Then there is a $C>0$ such that
$Y_s$ is a $(C,\delta,0)$-approximate $(q_1,q_2)$-diffusion
for every $\delta\in(0,1)$ and in every time interval $[0,S]$.
\end{lemma}
\proof
First, note that $q_1(x)=0$ and $x\,q_2(x)\le 0$ for $|x|\ge 1$.
This clearly implies that $\{Y_s:s\ge 0\}\subset[-1,1]$ a.s.
Now fix some $\delta>0$ and
two stopping times $s_0\le s_1$ satisfying the assumptions in the definition
of approximate diffusions.
Let $\ev F_{s_0}$ denote the $\sigma$-field generated by $(Y_s:s\le s_0)$.
Then
$$
Y_{s}-Y_{s_0}=\int_{s_0}^{s} q_2(Y_r)\,dr
+\int_{s_0}^{s} q_1(Y_r)^{1/2}\,dB_r\,.
$$
The second summand is a martingale, and therefore
$$
\Eb{Y_{s_1}-Y_{s_0}\md \ev F_{s_0}}=\EB{\int_{s_0}^{s_1} q_2(Y_s)\,ds\md \ev F_{s_0}}.
$$
Since $|Y_s-Y_{s_0}|\le\delta$ for $s\in[s_0,s_1]$ and $q_2$ is
a Lipschitz function, we conclude that
$$
\Eb{Y_{s_1}-Y_{s_0}\md \ev F_{s_0}}=\Eb{\Delta s\, q_2(Y_{s_0})\md \ev F_{s_0}}+O(\delta)\,\Eb{\Delta s\md\ev F_{s_0}}.
$$
Thus, $Y_s$ satisfies~\eref{e.dW}.

We now use \Ito/'s formula to calculate $(Y_{s_1}-Y_{s_0})^2$:
\begin{multline*}
(Y_{s_1}-Y_{s_0})^2
=
\int_{s_0}^{s_1} 2\,(Y_s-Y_{s_0})\,dY_s+
\langle Y\rangle_{s_1}-\langle Y\rangle_{s_0}
\\
=
\int_{s_0}^{s_1} 2\,(Y_s-Y_{s_0})\,q_2(Y_s)\,ds
+
\int_{s_0}^{s_1} 2\,(Y_s-Y_{s_0})\,q_1(Y_s)^{1/2}\,dB_s
+
\int_{s_0}^{s_1} q_1(Y_s)\,ds\,.
\end{multline*}
The left summand is $O(\delta^3)$ and the middle summand is
a martingale and therefore its expectation given $\ev F_{s_0}$
is zero. Thus
$$
\Eb{
(Y_{s_1}-Y_{s_0})^2\md\ev F_{s_0}}=
\EB{ \int_{s_0}^{s_1} q_1(Y_s)\,ds}+O(\delta^3)
= \Eb{\Delta s\,q_1(Y_{s_0})}+O(\delta^3)\,,
$$
because $q_1$ is Lipschitz.
This shows that $Y_s$ satisfies~\eref{e.dWW}, and completes the proof.
\QED

\begin{proposition}\label{p.diffusioncouple}
Fix $S>2$. Let $q_1,q_2:[-1,1]\to\R$ be defined as in~\eref{e.q},
where we assume that $a,b> -\hco$.
Suppose that $W^1$ is a $(C,\delta,\errprb)$-approximate
$(q_1,q_2)$-diffusion $W^1:[0,S]\to[-1,1]$
 and $W^2:[0,S]\to[-1,1]$ is a solution of~\eref{e.SDE}
with the same $q_1$ and $q_2$ and $W^1(0)=W^2(0)$ a.s.
Also assume $\errprb<\delta^5/S^2$.
Then there is a coupling of $W^1$ and $W^2$ such that
$\sup_{s\in[0,S-1]} |W^1_s-W^2_s|\to 0$
in probability as $\delta\to 0$ while $C$ is fixed.
Namely, for every $\eps>0$ there is a $\delta_0>0$,
depending only on $a,b,S,C$ and $\eps$ such that
$\sup_{s\in[0,S-1]} |W^1_s-W^2_s|\le \eps$ with
probability at least $1-\eps$ if $\delta<\delta_0$.
\end{proposition}

A useful tool in the proof of the proposition is the following lemma.

\begin{lemma}\label{l.dito}
Let $W$ be a $(C,\delta,\errprb)$-approximate $(q_1,q_2)$-diffusion,
$W:[0,S]\to[-1,1]$, and let $\tau_0$ and $\tau_1$
be two stopping times for $W$ satisfying
$0\le\tau_0\le\tau_1\le S$.
Assume that $C\,\delta<1/2$.
Let $f:[-1,1]\to\R$ be a function whose second derivative
is Lipschitz with Lipschitz constant $1$ and which satisfies
$\|f'\|_\infty,\|f''\|_\infty\le 1$.
Set
$$
Lf(x):= \frac 12\, q_1(x)\,f''(x)+q_2(x)\,f'(x)\,.
$$
Then there is a stopping time $\tau_1'$
satisfying $\tau_0\le\tau_1'\le\tau_1$ a.s.\ and
$\Pb{\tau_1'\ne\tau_1}\le\errprb$ such that
\begin{multline*}
\EB{f(W_{\tau'_1})-f(W_{\tau_0})
-\int_{\tau_0}^{\tau'_1} Lf(W_s)\,ds
\md
\ev F_{\tau_0}}
\\
=
O(C+1)\,\delta\,\Eb{\delta^2+\tau_1'-\tau_0\md\ev F_{\tau_0}}.
\end{multline*}
Moreover, in the above the function $f$ may be random, provided that it is
$\ev F_{\tau_0}$ measurable.
\end{lemma}

\proof
We inductively define the stopping times $s_j$ as follows.
Set $s_0:=\tau_0$, and
$s_{j+1}:=\min\{s\ge s_j:s=s_j+\delta^2
\text{ or }|W_s-W_{s_j}|=\delta\text{ or }s_j=\tau_1\}$.
If there is a $j\in\N$ such that $W$ does not satisfy~\eref{e.dW}
or~\eref{e.dWW} for the stopping times $(s_j,s_{j+1})$
in place of $(s_0,s_1)$, then let $n$ be the minimal such $j$.
(Note that the event that $W$ does not satisfy~\eref{e.dW} or~\eref{e.dWW}
for $(s_j,s_{j+1})$ is $\ev F_{s_j}$-measurable.)
Otherwise, let $n$ be the minimal
$j$ such that $s_j=\tau_1$.
Note that $(s_n,s_{n+1})$ is a pair of stopping times
and they do not satisfy both~\eref{e.dW} and~\eref{e.dWW} unless
$s_n=\tau_1$. Consequently,
\begin{equation}\label{e.gettoend}
\Pb{s_n\ne\tau_1}\le\errprb\,.
\end{equation}
Since $|W_{s_{j+1}}-W_{s_j}|\le\delta$ 
using a Taylor series for $f$ around $W_{s_j}$ we have
$$
f(W_{s_{j+1}})-f(W_{s_j})=f'(W_{s_j})\,\Delta_j W+(1/2)\,f''(W_{s_j})\,(\Delta_jW)^2+O(\delta^3)\,,
$$
where $\Delta_j W:=W_{s_{j+1}}-W_{s_j}$.
We may use~\eref{e.dW} and~\eref{e.dWW} to estimate the conditioned expectation of
$\Delta_j W$ and $(\Delta_jW)^2$ given $\ev F_{s_j}$ and get for $j<n$
$$
\Eb{f(W_{s_{j+1}})-f(W_{s_j})\md\ev F_{s_j}}
= Lf(W_{s_j})\,\Eb{s_{j+1}-s_j\md\ev F_{s_j}}+
O(1+C)\,\delta^3.
$$
By our assumptions about $f$ this may also be written as
$$
\EB{f(W_{s_{j+1}})-f(W_{s_j})-\int_{s_j}^{s_{j+1}} Lf(W_s) \,ds\md\ev F_{s_j}}= O(1+C)\,\delta^3.
$$
We sum this over $j$ from $0$ to $n-1$, then take expectations conditioned on $F_{\tau_0}$, to obtain
\begin{equation}\label{e.dfsn}
\EB{f(W_{s_{n}})-f(W_{\tau_0})-\int_{\tau_0}^{s_{n}} Lf(W_s) \,ds\md\ev F_{\tau_0}}
=
O(1+C)\,\delta^3\,\Eb{n\md\ev F_{\tau_0}}.
\end{equation}
Now fix some $j\in\N$.
On the event $j+1<n$, we have $(\Delta_j W)^2=\delta^2$ or $s_{j+1}-s_j=\delta^2$.
Therefore,
$$
\EB{\bl((\Delta_j W)^2+s_{j+1}-s_j\br)
\,1_{\{j<n\}}
\md\ev F_{s_j}}
\ge\delta^2\, \Pb{j+1<n\md\ev F_{s_j}}.
$$
By~\eref{e.dWW}, this gives
\begin{multline*}
\EB{\bl(1+q_1(W_{s_j})\br)\,(s_{j+1}-s_j)\,1_{\{j<n\}}\md\ev F_{s_j}} 
\\
\ge \delta^2\,
\Pb{j+1<n\md \ev F_{s_j}}-C\,\delta^3\,1_{\{j<n\}}\,.
\end{multline*}
We take expectation conditioned on $\ev F_{\tau_0}$ and use the fact that $q_1$ is bounded, to obtain
\begin{multline*}
O(1)\,\Eb{(s_{j+1}-s_j)\,1_{\{j<n\}}\md\ev F_{\tau_0}}
\\
 \ge \delta^2\,\Pb{j+1<n\md\ev F_{\tau_0}}
-C\,\delta^3\,\Pb{j<n\md \ev F_{\tau_0}}.
\end{multline*}
We sum this over all $j\in\N$, to get
$$
O(1)\,\Eb{s_n-\tau_0\md\ev F_{\tau_0}} \ge \delta^2\,(1-C\,\delta)\,\Eb{n-1\md\ev F_{\tau_0}}
-C\,\delta^3.
$$
By our assumption that $C\,\delta<1/2$, this implies
$$
O(1)\,\Eb{\delta^2+s_n-\tau_0\md\ev F_{\tau_0}}
\ge \delta^2\,\Eb{n\md\ev F_{\tau_0}}\,.
$$
When combined with~\eref{e.dfsn}, this gives
\begin{multline*}
\EB{f(W_{s_{n}})-f(W_{\tau_0})-\int_{\tau_0}^{s_{n}} Lf(W_s) \,ds\md\ev F_{\tau_0}}
\\
=
O(1+C)\,\delta\,\Eb{\delta^2+s_n-\tau_0 \md\ev F_{\tau_0}}.
\end{multline*}
By~\eref{e.gettoend}, this completes the proof with $\tau_1'=s_n$.
\QED

The next lemma bounds the expected time that $W_s$ spends close to $\pm1$.

\begin{lemma}\label{l.occupation}
Let $W$ be a $(C,\delta,\errprb)$-approximate $(q_1,q_2)$-diffusion
$W:[0,S]\to[-1,1]$,
where $q_1$ and $q_2$ are given by~\eref{e.q}, $C\,\delta<1/2$ and $b> -\hco$.
Suppose that $S\ge 1$. Given any $\eps>0$ there is some $x_0<1$,
$\delta'>0$ and $\errprb'>0$ all depending only on $\eps,a,b$ and $S$
such that if $\delta<\delta'$ and $\errprb<\errprb'$,
then
$$
\EB {\int_0^S 1_{\{W_s>x_0\}}\,ds}<\eps\,.
$$
A similar statement holds for the set of times such that
$W_s$ is near $-1$, provided that $a>-\hco$.
\end{lemma}
\proof
Set $\mu(A):=\EB{ \int_0^S 1_{\{W_s\in A\}}\,ds}$.
Note that $q_2(1)<0$. Fix some $y_0\in [0,1)$ such that
$q_2(x)\le q_2(1)/2$ throughout $[y_0,1]$ and set
$y_n:= 1-(1-y_0)\,2^{-n}$.
Let $f(x)$ be the twice continuously differentiable function that is
zero on $[-1,y_0]$, satisfies $f''(x)= \min\{x-y_0,y_1-x\}$ on $[y_0,y_1]$
and satisfies $f''(x)=0$ on $[y_1,1]$.
We apply Lemma~\ref{l.dito} with to $f$ with $\tau_0=0$ and
$\tau_1=S$.
Clearly $Lf(x)=0$ in $[-1,y_0]$.
On the interval $[y_0,y_1]$, we have
$f''(x)\le (y_1-y_0)/2$, $f'(x)\ge 0$ and $q_2(x)<0$.
Consequently, $Lf(x)\le q_1(x)(y_1-y_0)/4 \le (1-x)(y_1-y_0)\le (1-y_0)^2$
on $[y_0,y_1]$.
On the interval $[y_1,1]$, we have
$f''(x)=0$, $f'(x)= (y_1-y_0)^2/4$ and $q_2(x)\le q_2(1)/2<0$.
Consequently,
$Lf(x)\le -c\,(1-y_0)^2$, where $c>0$ depends only on $q_2(1)$.
Also note that
$$
|f(W_S)-f(W_0)|\le \sup_x f(x)-\inf_x f(x)=f(1)-0\le (1-y_0)^3.
$$
Therefore, Lemma~\ref{l.dito} gives
\begin{multline*}
(1-y_0)^2\,\mu([y_0,y_1))-c\,(1-y_0)^2\,\mu([y_1,1])
\\
\ge
-(1-y_0)^3+
O(C+1)\,\delta \,S+O(S)\,\errprb\,.
\end{multline*}
We assume that $\delta'$ and $\errprb'$ are sufficiently small so that
the right hand side is larger than $-2\,(1-y_0)^3$.
Then we get
$$
c\,\mu([y_1,1])\le 2\,(1-y_0)+
\mu([y_0,y_1))\,.
$$
This implies
$$
\mu([y_1,1])\le 2\,(1-y_0)+
\mu([y_0,1])/(1+c)\,.
$$
A similar inequality applies to $y_n$ and $y_{n+1}$. Induction therefore gives
$$
\mu[y_n,1]\le
\frac{\mu([y_0,1])}{(1+c)^{n}} +
\sum_{j=0}^{n-1} 2\,(1-y_j) <
\frac{\mu([y_0,1])}{(1+c)^{n}} +
 4 \,(1-y_0)\,,
$$
provided that $\delta'$
and $\errprb'$ are smaller than some functions of $n,S,C$ and $y_0$.
Consequently, we first choose $y_0$ such that
in addition to the requirements stated in the beginning of the proof,
$4\,(1-y_0)\le \eps/2$. Then we take $n$
sufficiently large so
that $(1-c)^{-n}\, S<\eps/2$. Then $\delta'$ and $\errprb'$ are determined.
This proves the first claim. The second one follows by symmetry.
\QED

The next lemma estimates the conditional expectation and  conditional
second moment of the time it takes $W_s$ to move a distance of $\delta$ beyond its
location at a stopping time.

\begin{lemma}\label{l.hoptime}
Let $W$ be a $(C,\delta,\errprb)$-approximate $(q_1,q_2)$-diffusion
$W:[0,S]\to[-1,1]$,
where $q_1$ and $q_2$ are given by~\eref{e.q} and $S>1$.
Let $x_0\in(0,1)$. There is a function $\delta_0>0$,
depending only on $x_0,C,a$ and $b$ such that the
following holds if $\delta<\delta_0\wedge (1/2)$ and
$\errprb<\delta^5/S^2$.
Let $\tau_0$ be a stopping time for $W$ and
let $\tau_1:=\inf\bl\{s\ge \tau_0: s=S\text{ or } |W_s-W_{\tau_0}|=\delta\br\}$.
Let $\ev A$ denote the event $\{\tau_0<S-1/2\}\cap\{|W_{\tau_0}|<x_0\}$.
Then
\begin{equation}
\label{e.timeexp}
\EB{\bl|\Es{\tau_1-\tau_0\md \ev F_{\tau_0}}-\delta^2/q_1(W_{\tau_0}) \br|\,1_{\ev A}}=
O_{x_0}(C+1)\,\delta^3 
\,.
\end{equation}
Moreover,
\begin{equation}\label{e.timemoment}
\EB{\Es{(\tau_1-\tau_0)^2\md \ev F_{\tau_0}}\,1_{\ev A}}
=O_{x_0}(\delta^4) 
\,.
\end{equation}
\end{lemma}
\proof
Let $f(x)=(x-W_{\tau_0})^2/4$, and let $L$ be as in Lemma~\ref{l.dito}.
Then $Lf(x)=q_1(W_{\tau_0})/4+O(\delta)$ for $x\in[W_{\tau_0}-\delta,W_{\tau_0}+\delta]$.
Therefore, Lemma~\ref{l.dito} gives
$$
\Eb{f(W_{\tau_1'})-q_1(W_{\tau_0})\,(\tau_1'-\tau_0)/4\md \ev F_{\tau_0}}
=O(C+1)\,\delta\,\Eb{\delta^2+\tau_1'-\tau_0\md \ev F_{\tau_0}}.
$$
That is,
\begin{equation}\label{e.time}
\Eb{
f(W_{\tau_1'})
\md \ev F_{\tau_0}}
-
\frac{q_1(W_{\tau_0})+O(C+1)\,\delta}4\,
\Eb{
\tau_1'-\tau_0
\md \ev F_{\tau_0}}
=O(C+1)\,\delta^3.
\end{equation}
By choosing $\delta_0$ sufficiently small,
we make sure that $O(C+1)\,\delta<q_1(W_{\tau_0})/2$
on $\ev A$.
Since $|f(W_{\tau_1'})|\le \delta^2/4$,
the above gives
$$
\Pb{\tau_1'=S\md\ev F_{\tau_0}}
\,
1_{\ev A}
=O(\delta^2)/q_1(W_{\tau_0})\,.
$$
Recall that $f(W_{\tau'_1})=\delta^2/4$ unless $\tau_1'=S$ or
$\tau_1'<\tau_1$.
Therefore, on $\ev A$,
$$
\Eb{
f(W_{\tau_1'})
\md \ev F_{\tau_0}}
=\delta^2/4+O_{x_0}(\delta^4)+
O(\delta^2)\,\Pb{\tau_1'<\tau_1\md\ev F_{\tau_0}}.
$$
We plug this and
$$
\Eb{
\tau_1'-\tau_1
\md \ev F_{\tau_0}}=O(S)\,\Pb{\tau_1'\ne\tau_1\md\ev F_{\tau_0}}
$$
 into~\eref{e.time}, simplify, and get
\begin{multline}\label{e.mixed}
\delta^2
-
\bl(q_1(W_{\tau_0})+O(C+1)\,\delta\br)\,
\Eb{
\tau_1-\tau_0
\md \ev F_{\tau_0}}
\\
=O_{x_0}(C+1)\,\delta^3
+
O(S)\,\Pb{\tau_1'<\tau_1\md\ev F_{\tau_0}} .
\end{multline}
on $\ev A$.
Now~\eref{e.timeexp} follows by dividing~\eref{e.mixed} by $q_1(W_{\tau_0})+O(C+1)\,\delta$, taking expectation
and recalling that $\Pb{\tau_1'\ne\tau_1}\le\errprb$.

Now define $t_n:=\tau_1\wedge \bl(\tau_0+16\,n\,\delta^2/q_1(W_{\tau_0})\br)$.
If $\delta_0$ is sufficiently small, then $Lf(x)>Lf(W_{\tau_0})/2=q_1(W_{\tau_0})/8$
throughout $[W_{\tau_0}-\delta,W_{\tau_0}+\delta]$ on the event $\ev A$.
Thus, we get by applying Lemma~\ref{l.dito} to the stopping times $t_n$ and $t_{n+1}$
$$
\bl(q_1(W_{\tau_0})-O(C+1)\,\delta\br)\,\Eb{t_{n+1}'-t_n\md\ev F_{t_n}}/8\le \delta^2/4+O(C+1)\,\delta^3,
$$
where $t_{n+1}'$ is the stopping time provided by the lemma.
Again, on $\ev A$ we may assume that $O(C+1)\,\delta<(q_1(W_{\tau_0})/2)\wedge (1/4)$.
Thus,
$$
q_1(W_{\tau_0})\,\Eb{t_{n+1}'-t_n\md\ev F_{t_n}}/16\le \delta^2/2\,,
$$
which implies
$$
\Pb{t'_{n+1}=t_n+16\,n\,\delta^2/q_1(W_{\tau_0})\md\ev F_{t_n}}\le 1/2\,.
$$
But if $t'_{n+1}\ne t_n+16\,n\,\delta^2/q_1(W_{\tau_0})$, then $t_{n+1}'=\tau_1$ or
$t'_{n+1}\ne t_{n+1}$. If $\ev B_n$ denotes the event that $t_{j}'=t_{j}$
for all $j=1,\dots,n$, then induction gives
$$
\Pb{t_n\ne \tau_1,\,\ev B_n\md\ev F_{\tau_0}} \,1_{\ev A} \le 2^{-n}.
$$
Lemma~\ref{l.dito} gives $\Pb{t_{n+1}'\ne t_{n+1}}\le\errprb$ and
therefore $\Pb{\neg\ev B_n}\le n\,\errprb$.
(In fact, it is not hard to get the better estimate
$\Pb{\neg\ev B_n}\le \errprb$.)
Consequently for $n\in\N$
$$
\Pb{\tau_1-\tau_0>16\,n\,\delta^2/q_1(W_{\tau_0})\md\ev F_{\tau_0}}\,1_{\ev A}\le 2^{-n}+\Pb{\neg\ev B_n\md\ev F_{\tau_0}}.
$$
 The above applies with
$\tilde n:=n\wedge\lceil -\log_2\errprb\rceil$
in place of $n$, and hence
\begin{multline*}
\PB{\tau_1-\tau_0>\frac{16\,n\,\delta^2}{q_1(W_{\tau_0})}\md\ev F_{\tau_0}}\,1_{\ev A}
\\
\le
\PB{\tau_1-\tau_0>\frac{16\,\tilde n\,\delta^2}{q_1(W_{\tau_0})}\md\ev F_{\tau_0}}\,1_{\ev A}\le 
2^{-\tilde n}+\Pb{\neg\ev B_{\tilde n}\md\ev F_{\tau_0}}.
\end{multline*}
We multiply both sides by $2\,(n+1)\,\bl(16\,\delta^2/{q_1(W_{\tau_0})}\br)^2$,
and sum over $n$ from $n=0$ to the least $m$ such that $16\,m\,\delta^2/{q_1(W_{\tau_0})}\ge S$.
The result on the left hand side bounds $\Eb{(\tau_1-\tau_0)^2\md \ev F_{\tau_0}}\,1_{\ev A}$.
Consequently, the required bound~\eref{e.timemoment} follows by taking expectations
and using our assumed upper bound for $\errprb$.
\QED

The following lemma shows that when we discretize the approximate diffusion
the resulting random walk has transition probabilities that can be well estimated
from $q_1$ and $q_2$ away from the boundary.

\begin{lemma}\label{l.transitions}
Fix some $x_0\in(0,1)$.
Let $W$ be a $(C,\delta,\errprb)$-approximate $(q_1,q_2)$-diffusion
$W:[0,S]\to[-1,1]$,
where $q_1$ and $q_2$ are given by~\eref{e.q},
$S\ge 1$ and $\errprb\le\delta^5/S$.
Set
$$
Z:=\bl\{k\,\delta:k\in\Z,\ |k\,\delta|< x_0\br\},
$$
$s_0:= \inf\{s\ge 0:W_s\in Z\text{ or }s=S\}$ and inductively
$$
s_{n+1}:= \inf\bl\{s\ge s_n: W_{s_n}\ne W_s\in Z\text{ or }s=S \br\}.
$$
Also set $X_n:=W_{s_n}$ and $Z^0:=Z\setminus\{ \min Z, \max Z\}$.
Let
$$
p^{\pm}_n:=\Pb{X_{n+1}=X_n\pm\delta\md \ev F_{s_n}}
$$
and
$$
r^{\pm}_n:=
\frac 12\pm\delta\,\frac{q_2(X_n)}{2\,q_1(X_n)}.
$$
There is a $\delta_0>0$, depending only on $C,x_0,a$ and $b$,
such that if $\delta<\delta_0$, then for all $n\in\N$
$$
\EB{\bl|p^{\pm}_n-r^{\pm}_n\br|\, 1_{\{s_n<S-1/2\}}\,1_{\{X_n\in Z^0\}}}
\le O_{x_0}(C+1)\,\delta^2.
$$
\end{lemma}
\proof
We now use a different test function:
$$
f(x):=\alpha\,(x-X_n)^2+\beta\,(x-X_n)\,,
$$
where $\alpha:= -q_2(X_n)\,\beta/q_1(X_n)$ and
$\beta:=|q_1(X_n)/(6\,q_2(X_n))|\wedge (1/3)$ with
$\beta=1/2$ if $q_2(X_n)=0$.
The choice of $\alpha$ and $\beta$ above is tailored to give
$Lf(X_n)=0$ and $|4\,\alpha|+|\beta|\le 1$.
The latter implies that $\|f''\|_\infty,\|f'\|_\infty\le 1$
when $f$ is restricted to the interval $[-1,1]$.
We now apply Lemma~\ref{l.dito} again with this $f$
and stopping times $s_n$ and $s_{n+1}$.
Note that $Lf=O(\delta)$ in the interval $[X_n-\delta,X_n+\delta]$.
Hence $Lf(W_s)=O(\delta)$ for $s\in[s_n,s_{n+1}]$.
Lemma~\ref{l.dito} gives
$$
\Eb{f(W_{s'_{n+1}})-f(W_{s_n})\md\ev F_{s_n}} =
O(C+1)\,\delta\,\Eb{\delta^2+s_{n+1}'-s_n\md\ev F_{s_n}},
$$
on the event $X_n\in Z^0$, where $s'_{n+1}$ is the stopping time produced by the lemma.
The above may be written
\begin{multline*}
(\beta\,\delta+\alpha\,\delta^2)\,p_n^+ + (-\beta\,\delta+\alpha\,\delta^2)\,p_n^-
= -\Eb{ 1_{\{ W_{s'_{n+1}}\notin Z\}} f(W_{s'_{n+1}})\md\ev F_{s_n}}+{}
\\
{}+
O(C+1)\,\delta\,\Eb{\delta^2+s_{n+1}'-s_n\md\ev F_{s_n}}.
\end{multline*}
Set $p_n^o:= \Pb{  W_{s'_{n+1}}\notin \{X_n-\delta,X_n+\delta\} \md\ev F_{s_n}}$.
Then $1\ge p_n^++p_n^-\ge 1-p_n^o$.
Hence, the above gives
$$
 (2\,p_n^+-1)\,\beta
 +\alpha\,\delta
 \\
 = O(p_n^o)+ O(C+1)\,\delta^2+
O(C+1)\,\Eb{s_{n+1}-s_n\md\ev F_{s_n}},
$$
which, by the definitions of $\alpha$ and $r_n^+$ may be rewritten,
\begin{equation}
\label{e.moof}
2\,\beta\, (p_n^+-r_n^+)
= O(p_n^o)+ O(C+1)\,\delta^2+
O(C+1)\,\Eb{s_{n+1}-s_n\md\ev F_{s_n}}.
\end{equation}
By~\eref{e.timeexp} and our assumption $\errprb\le\delta^5/S$, we have
\begin{equation}
\label{e.goof}
\Eb{(s_{n+1}-s_n)\,1_{X_n\in Z^0}\,1_{s_n<S-1/2}} = O_{x_0}(\delta^2)\,,
\end{equation}
provided that $\delta_0$ is sufficiently small.
Since $W_{s'_{n+1}}\notin Z$ only when $s'_{n+1}\ne s_{n+1}$
or $s'_{n+1}=s_{n+1}=S$,
on the event $ \{X_n\in Z^0\}\cap\{s_n<S-1/2\}$,
\begin{equation}
\label{e.split}
\begin{aligned}
 \Eb{p_n^o \md \ev F_{s_n}}
&\le
\Pb{ s'_{n+1}\ne s_{n+1}\md\ev F_{s_n}}+
\Pb{ s'_{n+1}\ge s_n+1/2 \md\ev F_{s_n}}
\\&
\le
\Pb{ s'_{n+1}\ne s_{n+1}\md\ev F_{s_n}}+
2\,
\Eb{s_{n+1}-s_n\md \ev F_{s_n}}.
\end{aligned}
\end{equation}
Note that $\beta^{-1}=O_{x_0}(1)$
and $ \Pb{ s'_{n+1}\ne s_{n+1}}\le\errprb$.
Hence, we now obtain
the result for $p_n^+$ by taking expectation
on the event
 $\{X_n\in Z^0\}\cap\{s_n<S-1/2\}$
 in~\eref{e.moof} and
using~\eref{e.goof} and~\eref{e.split}.
A symmetric argument applies to  $p_n^-$ and $r_n^-$, and the
proof is complete.
\QED

Next, we show that the time parameterization of $W$ can be
well approximated by a function of the discretized walk trajectory.

\begin{lemma}\label{l.whatsthetime}
Assume the setting and notation of Lemma~\ref{l.transitions}
in addition to $\delta<\delta_0$.
For $n\in\N$ let $t_n$ denote the time spent
up to time $s_n$ in segments
$[s_j,s_{j+1}]$ such that $X_j=W_{s_j}\in Z^0$; that is,
$$
t_n:=
\sum_{j=0}^{n-1}
1_{\{X_j\in Z^0\}}\, (s_{j+1}-s_j)\,.
$$
Also let
$$
\sigma_n:=
\sum_{j=0}^{n-1}
1_{\{X_j\in Z^0\}}\,\delta^2/q_1(X_j)\,.
$$
Let $N_0:=\min\{n\in \N:s_n\ge S-1/2\}$.
Then for all $n\in\N$
\begin{equation}
\label{e.whats}
\Eb{\max\{|\sigma_j-t_j|:j=1,2,\dots,n\wedge N_0\}}
\le
O_{x_0}(C+1)\,\bl(\delta^2\,n^{1/2}+\delta^3\,n\br)
\,.
\end{equation}
\end{lemma}
\proof
Let
$v_j:= (s_{j+1}-s_j)\,1_{\{X_j\in Z^0\}}\,1_{\{j<N_0\}}$,
$u_j:=\Eb{v_j\md \ev F_{s_j}}$ and
$w_j:= \delta^2\,q_1(X_j)^{-1}\,1_{\{X_j\in Z^0\}}\,1_{\{j<N_0\}}$.
Now, $M_n:=\sum_{j=0}^{n-1}(v_j-u_j)$ is clearly a martingale.
Consequently, Doob's maximal inequality
for $L^2$ martingales~\cite[II.1.6]{\RevuzYor}
gives
$$
\Eb{\max\{|M_j|:j=1,\dots, n\}}^2\le O(1)\,\Eb{M_n^2}
=O(1)\,\sum_{j=0}^{n-1}\Eb{(v_j-u_j)^2}.
$$
Since $u_j=\Eb{v_j\md\ev F_{s_j}}$, we have $\Eb{(v_j-u_j)^2}\le \Eb{v_j^2}$.
By Lemma~\ref{l.hoptime}, $\Eb{v_j^2}$ is bounded by the right hand
side of~\eref{e.timemoment}.
Now, the right hand side of~\eref{e.timeexp} bounds
$\Eb{|u_j-w_j|}$.
The result follows by our assumption
$\errprb\le \delta^5/S^2$, since for every $m\le n\wedge N_0$
\begin{multline*}
|\sigma_m-t_m|
=
\Bl|\sum_{j=0}^{m-1} (v_j-w_j)\Br|
\le|M_m|+
\sum_{j=0}^{m-1}
|u_j-w_j|
\\
\le
\max\bl\{|M_j|:j=1,\dots, n\br\}+\sum_{j=0}^{n-1}
|u_j-w_j|.
\QED
\end{multline*}

\proofof{Proposition~\eref{p.diffusioncouple}}
By Lemma~\ref{l.Ys}, $W^2$ is a $(C',\delta,0)$-approximate
$(q_1,q_2)$-diffusion for some fixed constant $C'>0$ and
every $\delta>0$. We may assume, with no loss of generality,
that $C\ge C'$.
Let $\eps>0$.
Let $x_0'\in (1-\eps,1)$ satisfy
Lemma~\ref{l.occupation} with this given
$\eps$, and assume that $\delta$ is sufficiently small so that that lemma is valid.
Take $x_0=(1+x_0')/2$.
Let $Z$ and $Z^0$ be as in Lemma~\ref{l.transitions}, let
$s_j^k$ be the corresponding stopping times introduced there for
$W^k$
and let $p_{k,j}^\pm$ denote the random transition probabilities
for $W^k$ defined there.
 Also abbreviate $X^k_j:=W^k_{s_j}$.
Let $\ev F^k_s$ denote the filtration of $W^k$, $k=1,2$.
Let $Y^k_j=1$ if $X^k_{j+1}-X^k_j=\delta$,
$Y^k_j=-1$ if $X^k_{j+1}-X^k_j=-\delta$
and $Y^k_j=0$ if $|X^k_{j+1}-X^k_j|\ne\delta$.
Then $\Pb{Y^k_j=\pm1\md\ev F^k_{s^k_j}}=p_{k,j}^\pm$
and $\Pb{Y^k_j=0\md\ev F^k_{s^k_j}}=1-p_{k,j}^+-p_{k,j}^-$
if $X^k_j\in Z^0$.

For the coupling of $W^1$ and $W^2$ we use an i.i.d.\ sequence
$U_j$ of uniform random variables in $[0,1]$.
The coupling proceeds as follows.
Up to their corresponding stopping times $s_0^k$, $k=1,2$, let them
run independently. Inductively, we suppose that the coupling has been
constructed up to their corresponding stopping times
$s_j^k$, $k=1,2$.
For each $k=1,2$ we take $Y^k_{j}= 1$ if $U_j\le p^+_{k,j}$,
$Y^k_{j}=-1$ if $U_j\ge 1-p^-_{k,j}$,
and $Y^k_{j}=0$ if $U_j\in (p^+_{k,j},1-p^-_{k,j})$.
(In other words, we try and match up $X^1_{j+1}-X^1_j$
with $X^2_{j+1}-X^2_j$ as much as possible.)
These choices respect the correct conditional distributions for these
variables.
Now we sample the restriction of $W^1$ to $[s_j^1,s_{j+1}^1]$ and
the restriction of $W^2$ to $[s_j^2,s_{j+1}^2]$ independently
from their corresponding conditional distribution given
$(\ev F^1_{s_j},Y^1_j)$ and $(\ev F^2_{s_j},Y^2_j)$, respectively.
This completes the description of the coupling.

Let $N:=\min\{n: s^1_n\vee s^2_n\ge S-1/2\}$.
Let $\ev A_j$ be the event $\{X_j^1,X_j^2\in Z^0\}$ and set
$$
Q_j:=1_{\ev A_j}\,|Y_j^1- Y_j^2|,\qquad
\hat Q_n:=\sum_{j=0}^{n} Q_j\,1_{\{j<N\}}.
$$
Note that on $\neg \ev A_j$ we have $|X^1_{j+1}-X^2_{j+1}|\le |X^1_j-X^2_j|$
unless $s^1_{j+1}\vee s^2_{j+1}=S$.
Moreover, $|X^1_0-X^2_0|\le\delta$ and when $Y^k_j=0$ we have $s^k_{j+1}=S$.
Consequently,
\begin{equation}
\label{e.Xdiff}
|X_n^1-X_n^2|\,1_{\{n<N\}}\le \delta+ \delta\,\hat Q_{n-1}\,.
\end{equation}
We now proceed to estimate
$ \theta_j:=\Eb{Q_j\,1_{\{j<N\}}} $.  Clearly,
$$
\Pb{Q_j\ne 0\md\ev F^1_{s^1_j}\vee \ev F^2_{s^2_j}}
\le \bl(|p_{1,j}^+-p_{2,j}^+|+
|p_{1,j}^--p_{2,j}^-|\br)\,1_{\ev A_j}\,.
$$
Let $r^\pm_{k,j}$ be the $r^\pm_j$ in Lemma~\ref{l.transitions}
corresponding to the process $W^k$. Then
$$
|p_{1,j}^\pm-p_{2,j}^\pm|
\le
|p_{1,j}^\pm-r_{1,j}^\pm|+ |r_{1,j}^\pm-r_{2,j}^\pm| + |r_{2,j}^\pm-p_{2,j}^\pm|.
$$
By that lemma,
$$
\Eb{
|p_{k,j}^\pm-r_{k,j}^\pm|\,1_{\ev A_j}\,1_{\{j<N\}}
} \le
O_{x_0}(C+1)\,\delta^2.
$$
Using the expression given for $r_{k,j}^\pm$, we deduce that
$$
|r_{1,j}^\pm-r_{2,j}^\pm|\le O_{x_0}(\delta)\,|X^1_j-X^2_j|\,.
$$
Thus, we get
$$
\theta_n=
\Eb{Q_n\,1_{\{n<N\}}}
\le
O_{x_0}(C+1)\,\delta^2 + O_{x_0}(\delta)\,\Eb{|X^1_n-X^2_n|\,1_{\{n<N\}}}.
$$
In conjunction with~\eref{e.Xdiff}, this gives
$$
\theta_n \le O_{x_0}(C+1)\,\delta^2\,\Bl(1+\sum_{j=0}^{n-1}\theta_j\Br).
$$
Induction therefore implies
$$
\theta_n\le O_{x_0}(C+1)\,\delta^2\,\bl(1+O_{x_0}(C+1)\,\delta^2\br)^n.
$$
Taking note of~\eref{e.Xdiff}, we infer
\begin{equation}
\label{e.distbd}
\EB{\max_{j<n\wedge N} |X^1_j-X^2_j|}
\le
\delta+
O_{x_0}(C+1)\,\delta^3\,n\,\bl(1+O_{x_0}(C+1)\,\delta^2\br)^n.
\end{equation}
Now let $m_0:=\lceil 4\,S\,\delta^{-2}\,\max\{q_1(x):|x|\le x_0\}\rceil$.
Observe that at least one from every two consecutive $j\in\N$  satisfies $X_j^k\in Z^0$ or $s_j^k=S$.
Consequently, $m_0<N$ implies
$\sigma_{m_0}^k\ge 2\,S$ for $k=1,2$, where $\sigma^k_j$ denotes the
$\sigma_j$ from Lemma~\ref{l.whatsthetime} corresponding to $W^k$.
Note that in that lemma $t_j\le s_j\le S$.
Therefore,
taking $n=m_0$ in~\eref{e.whats} implies
$$
\Pb{N>m_0}\le
O_{x_0}(C+1)\,S^{-1}\,\bl(\delta^2\,m_0^{1/2}+\delta^3\,m_0\br)
=
O_{x_0}(C+1)\, \delta\,.
$$
Set $X^*:= \max_{j<N } |X_j^1-X_j^2|$.
The above and~\eref{e.distbd} with $n=m_0$ imply
\begin{equation}
\label{e.Xs}
\Pb{X^*>\delta^{1/2}\text{ or }N>m_0}
\le O_{x_0,C,S}(\delta^{1/2}).
\end{equation}
Now for each $s\le S$ let $J(s):=\min\{j\in \N: s_j^1\ge s\}$.
Then
\begin{equation}
\label{e.Wdiff}
\begin{aligned}
&
 \sup_{s\in[0,S-1]}
|W_s^1-W_s^2|
\le
 \sup_{s\in[0,S-1]}
|W_s^1-X_{J(s)}^1|+{}
\\
&\qquad
{}+
 \sup_{s\in[0,S-1]}
|X_{J(s)}^1-X_{J(s)}^2|+
 \sup_{s\in[0,S-1]}
|W_s^2-X_{J(s)}^2|.
\end{aligned}
\end{equation}
First, it is clear that
$$
\sup_{s\in[0,S-1]}|W_s^1-X_{J(s)}^1|\le \eps +\delta\,.
$$
(The left hand side is usually at most $\delta$ but can be
as large as $1-x_0+\delta$ if, for example,
$X_{J(s)-1}=\max Z=x_0$.)
We leave aside, for now, the estimation of
the second summand in~\eref{e.Wdiff} and
consider the last.
Set
$$
t^*:= \sup_{s\in[0,S-1]} |s- s_{J(s)}^2|.
$$
Since $X_{J(s)}^2=W_{s_{J(s)}^2}^2$,
we have
\begin{equation}
\label{e.eqco}
\sup_{s\in[0,S-1]}|W_s^2-X_{J(s)}^2|\le
\sup_{s\in[0,S-1]}\, \sup_{t\in[0,t^*]}
|W_s^2-W_{s+t}^2|.
\end{equation}
Observe from~\eref{e.whats} that for $k=1,2$,
$$
\max\bl\{|\sigma^k_j-t^k_j|:j\le N\wedge m_0\br\}\to 0
,
$$
in probability as $\delta\to 0$,
where $t^k_j$ is the $t_j$ of Lemma~\ref{l.whatsthetime} corresponding to
$W^k$. Now the choice of $x_0$ (via Lemma~\ref{l.occupation}) implies that
$$
\PB{\max_{j\in\N } |s^k_0+t^k_j-s_j^k| >\sqrt\eps} <\sqrt{\eps}\,.
$$
Lemmas~\ref{l.occupation} and~\ref{l.hoptime} imply that
$\Eb{s_0^k}\le \eps \vee O_{x_0,C}(\delta^2)$.
Consequently, we have
\begin{equation}
\label{e.sigs}
\max\{|\sigma^k_j-s^k_j|:j\le N\wedge m_0\} \to 0
\end{equation}
in probability as $\eps,\delta\to 0$.
Now
\begin{multline*}
|\sigma^1_n-\sigma^2_n|
\le
\sum_{j=0}^{n-1}
\Bl|
1_{\{X^1_j\in Z^0\}}\,\delta^2/q_1(X^1_j)
-
1_{\{X^2_j\in Z^0\}}\,\delta^2/q_1(X^2_j)\Br|
\\
\le
\sum_{j=0}^{n-1}
\Bl|
\,\delta^2/q_1(X^1_j)
-
\delta^2/q_1(X^2_j)\Br|
+
\sum_{k=1}^2 \sum_{j=0}^{n-1} 1_{\{X^k_j\notin Z^0\}}
\,
\delta^2/q_1(X^{k}_j)
\,.
\end{multline*}
The right hand side is monotone non-decreasing in $n$.
When $n\le N\wedge m_0$ the first sum
is at most $O_{x_0}(\delta^2)\,m_0\,X^*$.
It is easy to see that for $n= N\wedge m_0$ the iterated
sum on the right tends to $0$ in probability:
this follows from the proof of~\eref{e.sigs},
because if we replace $x_0$ by $x_0+\delta$, the terms appearing
in this iterated sum are included in $\sigma_n$.
We now get from~\eref{e.Xs}
$$
\max_{j\le N}|\sigma^1_j-\sigma^2_j| \to 0
\qquad \text{in probability as $\eps,\delta\to 0$.}
$$
 By~\eref{e.sigs}, this also gives
$$
\max_{j\le N}|s^1_j-s^2_j| \to 0
\qquad \text{in probability as $\eps,\delta\to 0$.}
$$
In particular, \eref{e.Xs} and~\eref{e.sigs} imply
$\sup_{j\le N-1} s^k_{j+1}-s^k_j\to 0$
in probability, because $\sigma^k_{j+1}-\sigma^k_j\le O_{x_0}(\delta^2)$.
One consequence is that $\Pb{J(S-1)\ge N}\to 0$.
Furthermore, we now have
$$
t^* \le \sup_{s\in[0,S-1]} |s- s_{J(s)}^1|
+
\sup_{s\in[0,S-1]} |s^2_{J(s)}-s^1_{J(s)}|\to 0
$$
in probability.
Now~\eref{e.eqco} implies that
$ \sup_{s\in[0,S-1]}|W_s^2-X_{J(s)}^2|\to 0$ in probability,
because the right hand side in~\eref{e.eqco}
is smaller than $(t^*)^{1/3}$ with probability going to
$1$ as $t^*\to 0$, since $W^2$ is a solution of~\eref{e.SDE}.
This takes care of the last summand on the right hand side of~\eref{e.Wdiff}.

The middle summand on the right hand side of~\eref{e.Wdiff}
also tends to $0$ in probability because, as we have seen,
$\Pb{J(S-1)< N}\to 1$ and $X^*\to 0$ in probability.
This completes the proof.
\QED

\proofof{Theorem~\ref{t.lim}}
Let $x_1,x_2\in[-1,1]$ be two arbitrary points, and let
$Y^1_s,Y^2_s:[0,\infty)\to[-1,1]$ be two
independent solutions of~\eref{e.SDE} (with respect
to two independent Brownian motions)
which start at $x_1$ and $x_2$, respectively.
We claim that $s_{=}:=\min\{s:Y^1_s=Y^2_s\}<\infty$ a.s.
The argument is quite standard.
Suppose without loss of generality that $x_2>x_1$.
By Lemma~\ref{l.occupation}, it is unlikely
that $Y^2_s$ stays very  close to $1$ for a long
time and unlikely that $Y^1_s$ stays very close to
$-1$ for a long time. It is therefore easy to
conclude from Lemmas~\ref{l.hoptime} and~\ref{l.transitions}
that there are constants $s_0,c_0>0$  (which do not
depend on $x_1$ or $x_2$) such that
$\Pb{Y^2_{s_0}<0}>c_0$ and $\Pb{Y^1_{s_0}>0}>c_0$.
This implies $\Pb{s_{=}<s_0}>c_0^2$.
By strong uniqueness of solutions of~\eref{e.SDE}, it follows
that the solutions are Markov and have stationary transition
probabilities. Consequently, we get by induction
$\Pb{s_{=}>n\,s_0}< (1-c_0^2)^n$ for all $n\in\N$,
which proves that $s_=<\infty$ a.s.

We now argue that $\Pb{s_=<\infty}=1$ and the uniqueness in
law of solutions of~\eref{e.SDE} implies that for every
Borel subset $A\subset [-1,1]$ the limit
$$
\mu(A):= \lim_{r\to\infty}\Pb{Y^1_r\in A}
 $$
 exists.
We may couple a solution started at some time $s_0<r$ such that
$r-s_0$ is a large constant to be independent from $Y^1_s$ until the
first time in $[s_0,\infty)$ in which they meet and to agree with
$Y^1_s$ afterwards. Because these solutions are likely to meet prior
to time $r$, it follows that $\Pb{Y^1_r\in A}$ is close to the
probability that the solution started at time $s_0$ is in $A$ at
time $r$, proving the existence of $\mu$. (In our setting, $\mu$ may
be explicitly described. Its density with respect to Lebesgue
measure is proportional to
$(1+x)^{(b-\hco)/(2\hco)}\,(1-x)^{(a-\hco)/(2\hco)}$.)
  Since solutions of~\eref{e.SDE} are Markov,
 a solution $\tilde Y:[0,\infty)\to[-1,1]$
 of~\eref{e.SDE} such that the distribution of $\tilde Y_0$ is
 given by $\mu$ is time-stationary.
 To get a time-stationary solution $Y:(-\infty,\infty)\to[-1,1]$,
 we may take the weak limit of time-translations of $\tilde Y$.

Now let $S'$ be much larger than $S$.
By Proposition~\ref{p.diffusion} with $S'$ instead of
$S$ and Proposition~\ref{p.diffusioncouple} translated to start
at time $-S'$ and an appropriate choice of the $S$ appearing there,
we may couple $\tilde W_s$ so that with probability close to $1$
it stays close to a solution $W^2_s$ of~\eref{e.SDE} starting at
$\tilde W_{-S'}$ throughout $[-S',S']$.
We may at the same time couple $W^2_s$ so that with high probability
it agrees with $Y_s$ inside the interval $[-S,S]$.
Then with high probability $\tilde W_s$ stays close to $Y_s$
in $[-S,S]$, which concludes the proof of the theorem.
\QED

\begin{remark}
At this point, it may be worthwhile to point out which properties of
the functions $q_1$ and $q_2$ played a part in the proof.
The only properties that are essential for the above proof
are that $q_1$ and $q_2$ are both Lipschitz continuous
in $[-1,1]$, that $q_1>0$ in $(-1,1)$ and
$q_1=0$ on $\{-1,1\}$, and that $q_2(1)<0<q_2(-1)$.
\end{remark}

\subsection{Back to chordal} \label{s.backtochordal}

In Subsection~\ref{s.coord}, we described the transition from the chordal Loewner system
to the setup with the points $\pm 1$ fixed. We now describe the reverse transformation.
We start with some continuous $Y:(-\infty,\infty)\to[-1,1]$.
Set
\begin{equation}
\label{e.ws}
w_s:=\frac{e^s}2\,Y_s+ \frac 12\int_{-\infty}^s e^{u}\, Y_u\,du\,.
\end{equation}
Also define
$$
t^*(s):=\frac 18\int_{-\infty}^s
e^{2u}\,\bl(1-Y_u^2\br)\,du\,,
\quad s^*(t):=\sup\{s\in(-\infty,\infty):t^*(s)\le t\}\,.
$$
Now set $\hat Y_t:=w_{s^*(t)}$ for $t>0$, $\hat Y_0:=0$, and observe that $\hat Y$ is continuous
provided that there is no nontrivial time interval in which $Y\in\{\pm 1\}$.

\begin{lemma}\label{l.fullcircle}
If $Y_s=\tilde W_s$ is defined from $W_t$ as in Subsection~\ref{s.coord},
then $\hat Y_t=W_t$.
\end{lemma}
\proof
By the definition of $s(t)$ in Subsection~\ref{s.coord}, we have
$e^s=y_t-x_t$. Consequently,~\eref{e.ds} implies
$\p_t \bl(t^*(s(t))\br)=1$. Since $s(0)=-\infty$ and
$t^*(-\infty)=0$, it follows that $t^*(s(t))=t$
for all
$t\ge 0$
and $s^*(t^*(s))=s$ for all $s\in(-\infty,\infty)$.
Next,~\eref{e.xy} and the definition of $W^*$ give
$$
\p_t(y_t+x_t)=
\frac {8\,W^*_t}{(y_t-x_t)\,\bl(1-(W^*_t)^2\br)}\,.
$$
Since $y_0=0=x_0$, and $\tilde W_{s(t)}=W^*_t$ this implies
$$
y_t+x_t=
\int_0^t
\frac {8\,\tilde W_{s(r)}}{(y_r-x_r)\,\bl(1-\tilde W_{s(r)}^2\br)}\,dr
=
\int_{-\infty}^{s(t)}e^u\,\tilde W_u\,du
\,,
$$
where the second equality follows by a change of variable.
Now $\hat Y_t=w_{s(t)}=W_t$ follows from the definition
of $W^*_t$. The proof of the lemma is therefore complete.
\QED

We now discuss the behavior of solutions of~\eref{e.SDE}
in the chordal coordinate system, but generalize to the case $\kappa\ne 4$.

\begin{lemma}\label{l.Ychordal}
Let $\tilde a,\tilde b>0$, let $Y:(-\infty,\infty)\to[-1,1]$ be a solution
of
$$
dY_s =
(-\tilde a\,(Y_s-1)-\tilde b\,(Y_s+1))\, ds +
\sqrt{\kappa\,(1-Y_s^2)/2}\,dB_s\,,
$$
and let $\hat Y_t$ be the corresponding process, as described
following~\eref{e.ws}.
Then on any time interval which avoids
$\bl\{t:Y_{s^*(t)}\notin\{\pm 1\}\br\}$,
the process $\hat Y_t$ satisfies the SDE of the driving term for
\SLEkr\kappa;2\,(\tilde a-1),2\,(\tilde b-1)/:
$$
d\hat Y_t=
\frac{2\,(\tilde a-1)}{\hat Y_t-x_t}
+
\frac{2\,(\tilde b-1)}{\hat Y_t-y_t}+ \sqrt\kappa \,d\hat B_t
$$
for some Brownian motion $\hat B_t$, where $x_t$ and $y_t$
satisfy~\eref{e.xy} with $\hat Y_t$ in place of $W_t$.
\end{lemma}
\proof
Let $\beta_s=\int_{-\infty}^s e^u\,Y_u\,du$,
$y^*_s:= (\beta_s+e^s)/2$ and $x^*_s:=(\beta_s-e^s)/2$.
\Ito/'s formula and the definition of $t^*$ give
\begin{multline*}
dw_s=
e^s\,Y_s\,ds+\frac {e^s}2\,\sqrt{\kappa\,(1-Y_s^2)/2 }\,dB_s
+ \frac {e^s}2 \,
(-\tilde a\,(Y_s-1)-\tilde b\,(Y_s+1))\, ds
\\=
\sqrt{\kappa\,{t^*}'(s) }\,dB_s
+ 4\,{e^{-s}} \,
\Bl(\frac{\tilde a-1}{Y_s+1}+\frac{\tilde b-1}{Y_s-1}\Br)\,
{t^*}'(s)\, ds
\\=
\sqrt{\kappa\,{t^*}'(s) }\,dB_s
+ 2\,
\Bl(\frac{\tilde a-1}{w_s-x^*_s}+\frac{\tilde b-1}{w_s-y^*_s}\Br)\,
{t^*}'(s)\, ds \,.
\end{multline*}
Now set
$$
\hat B_t = \int_{-\infty}^{s^*(t)}\sqrt{{t^*}'(u)}\,dB_u\,.
$$
Then $\hat B_t$ is clearly a continuous martingale.
Since also $\langle \hat B\rangle_t=\int_{-\infty}^{s^*(t)} {t^*}'(u)\,du=t$,
we find that $\hat B$ is a Brownian motion with respect to $t$.
The above formula for $dw$ gives
$$
d\hat Y_t=\sqrt\kappa\,d\hat B_t
+ 2\,
\Bl(\frac{\tilde a-1}{\hat Y_t-x_t}+\frac{\tilde b-1}{\hat Y_t-y_t}\Br)\,dt\,,
$$
where $x_t:= x^*_{s^*(t)}$ and $y_t:= y^*_{s^*(t)}$.
Now
$$
\p_t x_t = \frac {\p_s x_s^*}{{t^*}'(s)}=\frac 2{x_t-\hat Y_t}\,,
$$
and similarly for $y_t$. This concludes the proof.
\QED

As mentioned at the beginning of this section (\S\ref{s.drive}), existence
and uniqueness of solutions to the usual SDE defining
\SLEkr\kappa;\rho_1,\rho_2/ have not been proved beyond times
when the driving term $W_t$ meets the force points. We now offer the following.

\begin{definition}\label{d.slekr}
If $\tilde a,\tilde b>0$, then
 the Loewner equation driven by the process $\hat Y$
of Lemma~\ref{l.Ychordal} is called
\SLEkr\kappa;2\,(\tilde a-1),2\,(\tilde b-1)/
(starting from $(0,0_-,0_+)$).
\end{definition}

\subsection{Loewner driving term convergence}

In this section, we complete the proof of Theorem~\ref{absmall}.
The theorem will follow quite easily from Theorem~\ref{t.lim}.

\proofof{Theorem \ref{absmall}}
Fix $T,\eps,\eps_0,\eps'>0$.
Let $Y_s$ and $\tilde W_s$ be coupled as in Theorem~\ref{t.lim},
but with $\eps'$ in place of $\eps$.
Let $\hat Y_t$, $t^*(s)$ and $s^*(t)$
be defined from $Y_s$ as in the beginning of \S\ref{s.backtochordal}.
Since the interior of the set
of times for which $Y_s\in\{-1,1\}$ is empty,
it follows that there is some positive $t_0>0$ such that
with probability at least $1-\eps_0$ we have
$t^*(1)-t^*(0)>t_0$. Because $Y_s$ is stationary,
it follows that
$$
\Pb{t^*(s+1)-t^*(s)>e^{2s} t_0}\ge 1-\eps_0\,.
$$
In particular, there is some $S_0>0$ such
that $\Pb{t^*(S_0)>T+1}\ge 1-\eps_0$.
Fix an $S_0$ satisfying this and additionally
$e^{- S_0}<\eps_0/2$.

Let
$$
\tilde t(s):=\frac 18 \int_{-\infty}^s e^{2u}\,(1-\tilde W_u^2)\,du\,,
$$
which is the equivalent of $t^*(s)$ with $\tilde W$ replacing $Y$.
It is clear that if $\eps'=\eps'(S_0,\eps_0)$ is sufficiently small
and
\begin{equation}\label{e.goodcouple}
\sup\bl\{|\tilde W_s-Y_s|: s\in[-S_0,S_0]\br\}<\eps',
\end{equation}
then for every $s\in[-\infty,S_0]$
the right hand side in~\eref{e.ws} differs from the corresponding
quantity where $\tilde W$ replaces $Y$ by at most $\eps_0$.
Lemma~\ref{l.fullcircle} then gives
$$
\sup\bl\{|\hat W_{\tilde t(s)}-W_{t^*(s)}|:s\le [0,S_0]\br\}<\eps_0\,,
$$
where $W_t$ is the chordal driving term for \SLEab/ and
$\hat W$ is the chordal driving term for $\phi_D\circ \gamma$.
(Here, we also use the fact that $s^*(t^*(s))=s$.)
If we assume~\eref{e.goodcouple} with $\eps'$ sufficiently small,
we also get $\sup\bl\{|t^*(s)-\tilde t(s)|:s\le S_0\br\}<\eps_0$.
Let $T_0$ be the obvious upper bound for
$\tilde t$ and $t^*$ in $(-\infty,S_0]$;
that is, $T_0:=e^{2S_0}/16$.
Also set
$$
M(\eps_0):=
\sup\bl\{|W_{t_0}-W_{t_1}|:0\le t_0\le t_1\le t_0+\eps_0\le T_0+\eps_0\br\}\,.
$$
Since $W_t$ is a.s.\ continuous, $M(\eps_0)\to0$ in probability as $\eps_0\to 0$.
Now the triangle inequality
$|\hat W_{\tilde t(s)}-W_{\tilde t(s)}|\le
|\hat W_{\tilde t(s)}-W_{t^*(s)}|+
|W_{t^*(s)}-W_{\tilde t(s)}|$
shows that when~\eref{e.goodcouple} holds
we have
$ \sup\bl\{|\hat W_{\tilde t(s)}-W_{\tilde t(s)}|:s\le S_0\br\}<\eps_0+M(\eps_0) $.
Hence,
$ \sup\bl\{|\hat W_{t}-W_{t}|:t\in[0,\tilde t(S_0)]\br\}<\eps_0+ M(\eps_0)$.
But we have seen that $t^*(S_0)$ is very likely to be larger than $T+1$
and that $|t^*(S_0)-\tilde t(S_0)|<\eps_0$ when~\eref{e.goodcouple} holds.
By Theorem~\ref{t.lim}, when $r_D$ is large~\eref{e.goodcouple}
holds with high probability.
This concludes the proof.
\QED

\subsection{Carath\'eodory convergence}\label{CKC}

For $K\subset\closure\H$
let $N_\eps(K)$ denote the $\eps$-neighborhood
of $K$ in $\overline\H$, and let $\cH_\eps(K)$
denote the unbounded connected component
of $\closure{\H}\setminus N_\eps ( K)$.
Set
$$
\dCKC(K,K'):=\inf\{\eps>0: K\cap \cH_\eps(K')=\emptyset=K'\cap\cH_\eps(K)\}\,.
$$
It is easy to see that $\dCKC$ is a metric
on the collection of compact connected $K\subset\closure\H$
such that $\H\setminus K$ is connected.
This metric is related to the
Carath\'eodory Kernel convergence topology, which is of
central importance in the theory of conformal mappings.

Let $K_t$ and $K'_t$ denote the evolving hulls
corresponding to two Loewner evolutions generated by
continuous driving terms $W_t$ and $W'_t$, respectively
(as defined in~\S\ref{ss.definingsle}).
Such evolving hulls are also sometimes called {\bf  Loewner chains}.
We set for $T\ge 0$,
$$
\dCKC^T(K,K'):=\sup_{t\in[0,T]} \dCKC(K_t,K'_t)\,.
$$
Following is a simple lemma relating uniform convergence
of driving terms to $\dCKC$-convergence of the
corresponding Loewner chains.

\begin{lemma}\label{l.WCKC}
The Loewner transform $W\mapsto K$ is a continuous map
from the space of continuous paths $W$ with the topology of
uniform convergence to the space of Loewner chains with
$\dCKC^T$-convergence.

In other words, for every $\eps>0,T>0$ and
$W:[0,T]\to\R$ continuous, there is some $\delta=\delta(\eps,T,W)>0$
such that if $\tilde W:[0,T]\to\R$
is continuous and satisfies $\sup_{t\in[0,T]}|W_t-\tilde W_t|<\delta$,
then the corresponding Loewner chains satisfy
$\dCKC^T(K,\tilde K)<\eps$.
\end{lemma}

This lemma is similar in spirit to~\cite[Proposition 4.47]{\LawlerConformalBook}.
As is well known, $\tilde K_T\to K_T$ in the Hausdorff metric
does not follow from $\tilde W\to W$ uniformly in $[0,T]$.

\proof
Fix $T>0$. Suppose that $W^n\to W$ uniformly in $[0,T]$.
Let $K^n$ denote the Loewner chain corresponding to $W^n$
and let $g^{(n)}_t:\H\setminus K^n_t\to\H$ denote the
corresponding Loewner evolution.
Fix some $t_0\in[0,T]$.
Since $\diam(K^n_{t_0})$ is clearly bounded
by a function of $T$ and $\|W^n\|_\infty$
\cite[Lemma 4.13]{\LawlerConformalBook},
the closure of $\{K^n_{t_0}:n\in\N_+\}$
is compact with respect to the Hausdorff metric
on nonempty compact subsets of $\closure\H$.
Consider some integer sequence $n_j\to\infty$
for which the Hausdorff limit $K':=\lim_{j\to\infty} K^{n_j}_{t_0}$
exists.
If $z\in\closure \H\setminus K_{t_0}$,
then there is a neighborhood  $U$ of $z$ in $\closure\H$
such that $U\cap K^n_{t_0}=\emptyset$ for all
sufficiently large $n$, by the continuity of
solutions of ODE's in the vector field specifying the
ODE.
It follows that $z\notin K'$, and hence $K'\subset K_{t_0}$.

With the intention of reaching a contradiction, suppose that
there is some point $z\in(\p K_{t_0})\setminus (K'\cup \R)$.
Let $z'$ be a point in $\H\setminus K_{t_0}$ satisfying
$|z-z'|< \dist(z',K'\cup\R)/3$.
The above argument shows that
$\lim_{j\to\infty}g_{t_0}^{(n_j)}(z')= g_{t_0}(z')\in\H$.
On the other hand, for all sufficiently large $j$ we
have $\dist(z',K^{n_j}_{t_0})> 2\,|z-z'|$.
Now the Koebe distortion theorem (e.g.~\cite[Cor.~1.4]{\PommeBDRY})
applied to the restriction of $g^{(n_j)}_{t_0}$ to the
disk of radius $2\,|z-z'|$ about $z'$
(once with $z'$ and again with $z$) shows that
$\Im g^{(n_j)}_{t_0}(z)\ge (3/16) \,\Im g^{(n_j)}_{t_0}(z')$.
However, since $z\in K_{t_0}$, for every
$\eps>0$ there is some first $t_1\in[0,t_0)$
such that $\Im g_{t_1}(z) \le\eps$.
The convergence argument above shows that for
all arbitrarily large $n$
$\Im g_{t_1}^{(n)}(z)\le 2\,\eps$. Since $\Im g_t^{(n)}(z)$
decreases monotonically in $t$, it follows that
$\Im g^{(n)}_{t_0}(z)\le 2\,\eps$ for all
sufficiently large $n$.
This contradicts our previous conclusion
$\Im g^{(n_j)}_{t_0}(z)\ge (3/16) \,\Im g^{(n_j)}_{t_0}(z') \to (3/16)\,\Im g_{t_0}(z')>0$,
and proves $K'\supset (\p K_{t_0})\setminus\R$.

Let $\tilde K$ be the union of $K'$ and the bounded connected components of $\closure\H\setminus K'$.
The above implies $\tilde K\supset K_{t_0}\setminus\R$.
Now note that $K_{t_0}\setminus\R$ is dense in $K_{t_0}$.
(This follows from the easy direction 2 $\Rightarrow$ 1
in~\cite[Theorem 2.6]{\LSWi} and from the
fact that $\H\cap K_t\setminus K_{t'}$ is nonempty
when $t>t'$.)
Consequently, $\tilde K= K_{t_0}$,
which implies that $\dCKC(K^n_t,K_t)\to 0$
for every fixed $t\in [0,T]$.

Observe that the above proof also gives $\lim_{s\to t}\dCKC(K_s,K_t)=0$
for $t\in [0,T]$ and $s$ tending to $t$
in $[0,T]$. Thus, $K_t$ is continuous in
$t$ with respect to $\dCKC$.
Since
$\dCKC(L,\tilde L)\le \dCKC(L',\tilde L)\vee \dCKC (L'',\tilde L)$
when $L'\subset L\subset L''$, the $\dCKC^T$ convergence
easily follows from the pointwise convergence,
from continuity of $K_t$ and from monotonicity
of $K^n_t$ in $t$.
\QED

\subsection{Improving the convergence topology}\label{ss.improve}

In this subsection we complete the proof of Theorem~\ref{preciseSLEconvergence}.
There are examples showing that the convergence of the Loewner driving term
does not imply the uniform convergence of the paths parameterized by capacity.
(See~\cite[\S3.4]{\LSWlesl}.) Therefore, we will need to apply other considerations.
Before embarking on the proof we note that when $a,b\ge \hco$
the trace of \SLEab/ is a simple path
that does not hit $\R$, except at its starting point.
Indeed, note first that the force points are moving monotonically
away from one another. By comparison with a Bessel process,
for example, it is easy to see that the trace does not hit the real line
at any time $t>0$.
It also does not hit itself, since $t\mapsto g_s(\gamma(t+s))$
has law that is mutually absolutely continuous with the path of
\SLEkk4/.

\begin{lemma}\label{l.Hconv}
Let $T>0$ and let $W_n:[0,T]\to\R$ be
a sequence of continuous functions converging uniformly to
a function $W:[0,T]\to\R$.
Suppose that each $W_n$ is the driving term
of a Loewner evolution of a path $\gamma_n:[0,T]\to\closure{\H}$
and $W$ is the driving term of a Loewner evolution of a simple path 
$\gamma:[0,T]\to\closure\H$
satisfying $\gamma(0,T]\cap\R=\emptyset$.
Then 
$\lim_{n\to\infty}\sup_{t\in[0,T]}\dhaus(\gamma_n[0,t],\gamma[0,t])= 0$,
where $\dhaus$ denotes the
the Hausdorff metric.
\end{lemma}
\proof
First note that $\diam \gamma_n[0,T]$ is bounded, because
$\|W_n\|_\infty$ is bounded.
Fix some $t\in[0,T]$, and let $\Gamma$ denote a subsequential Hausdorff limit of
$\gamma_n[0,t]$. It suffices to prove that $\Gamma=\gamma[0,t]$.
By Lemma~\ref{l.WCKC}, we know that for every $\eps>0$
we have $\Gamma\cap\cH_\eps(\gamma[0,t])=\emptyset$.
Since $\gamma$ is a simple path satisfying $\gamma(0,T]\cap\R=\emptyset$,
it follows that
$\bigcup_{\eps>0} \cH_\eps(\gamma[0,t])= \H\setminus \gamma[0,t]$,
which implies that $\Gamma\subset \R\cup\gamma[0,t]$.
Fix some $z_1\in\R\setminus \gamma[0,t]=\R\setminus\{\gamma(0)\}$.
By the continuity in $W$ and $z$ of the solutions of Loewner's equation~\eref{e.chordal},
it follows that there is a neighborhood
$V$ of $z_1$ such that $V\cap \gamma_n[0,T]=\emptyset$ for all sufficiently
large $n$. This implies $z_1\notin\Gamma$, and hence
$\Gamma\subset\gamma[0,t]$.

Now let $t'\in[0,t]$.
By Lemma~\ref{l.WCKC} again, for every $\eps>0$ and
every $n$ sufficiently large $\gamma(t')\notin\cH_\eps(\gamma_n[0,t])$,
which means that every path connecting $\gamma(t')$ to $\infty $ in $\cH$
must come within distance $\eps$ from $\gamma_n[0,t]$.
Thus, every such path must intersect $\Gamma$. 
Since $\Gamma\subset\gamma[0,t]$ is closed, this implies $\gamma(t')\in\Gamma$.
Therefore, $\Gamma=\gamma[0,t]$; that is,
$\lim_{n\to\infty}\dhaus(\gamma_n[0,t],\gamma[0,t])=0$.
Since $\gamma_n[0,t]$ and $\gamma[0,t]$ are monotone increasing in $t$
and $\gamma[0,t]$ is continuous in $t$ with respect to $\dhaus$,
it easily follows that 
$\lim_{n\to\infty}\sup_{t\in[0,T]}\dhaus(\gamma_n[0,t],\gamma[0,t])= 0$.
\QED

Here is an outline of the main ideas going into the proof
of Theorem~\ref{preciseSLEconvergence}.
Let $\gamma^\phi$ be the path $\phi\circ\gamma$
parametrized by half-plane capacity.
The main step in the proof is to show that if we fix $T>0$
we have $\sup_{t\in[0,T]}|\gamma^\phi(t)-\gSLE(t)|\to 0$ in probability.
By Theorem~\ref{absmall} and Lemma~\ref{l.Hconv}, we get 
$\sup_{t\in[0,T]}\dhaus(\gamma^\phi[0,t],\gSLE[0,t])\to 0$
in probability (since $\gSLE$ is a simple path).
 We only need to rule out the possibility that
$\gamma^\phi$ has significant (and fast) backtracking
along $\gSLE$. This is ruled out by invoking Lemma~\ref{l.direct}
and observing that the $\phi$-image of the place where
simple random walk (starting from a vertex near $\phi^{-1}(i)$)
hits $\p\DD(\gamma)$ can be close to any fixed segment of $\gSLE(0,T]$.

\proofof{Theorem~\ref{preciseSLEconvergence}}
Let $\gamma^\phi$ be the path $\phi\circ\gamma$
parametrized by half-plane capacity.
Let $\delta$, $T>0$ and $\rr= \inr{\phi^{-1}(i)}(D)$.
Let $W$ denote the Loewner driving term of $\gSLE$.
Since $\gSLE$ is a.s.\ a simple path,
Lemma~\ref{l.Hconv} implies that for every $\eps>0$
there is some $\eps'=\eps'(\eps,\gSLE)>0$ such that
if $\tilde W$ is the driving term of a continuous
path $\tilde\gamma$ and
$\sup_{t\in[0,T]} |\tilde W_t-W_t|<\eps'$,
then $\sup_{t\in[0,T]} \dhaus(\gSLE[0,t],\tilde\gamma[0,t])<\eps$.
Moreover, it is not hard to see that $\eps'$ can be chosen as a measurable
function of $W$.
Hence,
Theorem~\ref{absmall} implies that for every $\eps_0>0$
if $\rr$ is larger than some
function of $\eps_0,\delta,T,a$ and $b$, then there is a coupling
of $h$ and \SLEab/ such that
$\rho:= \sup_{t\in [0,T]} \dhaus\bl(\gamma^\phi[0,t],\gSLE[0,t]\br)<\eps_0$
with probability at least $1-\delta$.
Without yet specifying $\eps_0$,
we assume that indeed $\gamma^\phi[0,T]$ and $\gSLE$ are so coupled.
Let $\ev A_0$ be the event that $\rho<\eps_0$.
Then $\Pb{\ev A_0}\ge 1-\delta$.

Let $0<t_0<t_1<t_2<t_3\le T$.
We will show that under this coupling, if $\rr$ is large,
then with high probability $\gamma^\phi(t)$
is close to $\gSLE[t_0,T]$ for every $t\in [t_3,T]$.
This will then imply that $\sup_{t\in[0,T]}|\gSLE(t)-\gamma^\phi(t)|$ is
small.

Since $\gSLE$ is a.s.\ a simple path disjoint from $\{i\}$,
there is a constant $s_0=s_0(t_0,t_1,t_2,t_3,T,\delta)>0$ such that
$\Pb{\ev A_1}\ge 1-\delta$, where
$\ev A_1$ is the event that
\begin{enumerate}
\item \label{i.harm}
the harmonic measure of $\gSLE[t_1,t_2]$ from $i$
with respect to the domain
 $\H\setminus \gSLE[0,\infty)$
is at least $s_0$,
\item \label{i.distR} $\dist\bl(\R,\gSLE[t_0,T]\br)>s_0$,
\item $\dist\bl(\gSLE[0,t_j],\gSLE[t_{j+1},T]\br)>s_0$ for $j=0,1,2$,
\item \label{i.diamA} $\diam\gSLE[0,T]<1/s_0$, and
\item  \label{i.ifar} $i\in\cH_{s_0}(\gSLE[0,T])$.
\end{enumerate}
(It is tedious, but straightforward, to check that $\ev A_1$
is measurable.)

Consider a simple random walk $S$ independent from $h$
starting at a $\TG$-vertex closest to $\phi^{-1}(i)$.
Let $\tau_T$ be the first time $t$ when
$S(t)\in \p \DD(\phi^{-1}\circ \gamma^\phi[0,T])$,
and as in \S\ref{zzTdef}, let $\zzT:=S(\tau_T)$.
We claim that for every $\eps_1>0$ if
$\rr$ is sufficiently large and $\eps_0$
is sufficiently small, then
\begin{equation}\label{e.bingo}
\PB{\dist\bl(\phi(\zzT),\gSLE[t_1,t_2]\br)<\eps_1\md \gSLE,\gamma^\phi[0,T]}
> s_0/2 \qquad\text{on }\ev A_0\cap \ev A_1\,.
\end{equation}

To prove~\eref{e.bingo} first observe
that conditional on $\gSLE$ such that $\ev A_1$ holds
 a two-dimensional Brownian motion $\hat S$ started at $i$
has probability at least $s_0$ to first hit
$\gSLE\cup\R$ in $\gSLE[t_1,t_2]$.
Moreover, if this happens, the Brownian motion is
likely to stay within a compact subset $L\subset\H$
before hitting $\gSLE[t_1,t_2]$
and not to come arbitrarily close to $\gSLE$ far from
its hitting point.
On compact subsets of $\H$, the map $\rr^{-1}\phi^{-1}$
distorts distances by a bounded factor,
by the Koebe distortion theorem~\cite[Thm.~1.3 \& Cor.~1.4]{\PommeBDRY}.
Since $\phi^{-1}$ takes a Brownian motion to a
monotonically time-changed Brownian motion,
by taking $\rr$ large
we may couple $(\rr)^{-1}\,S$ and $(\rr)^{-1}\,\phi^{-1}\circ\hat S$
to stay arbitrarily close (until $\phi^{-1}\circ\hat S$ hits
$\p D$) with high
probability, up to a time change.
Assuming that $\rho$
is arbitrarily small and taking~\ref{i.ifar}
 into account,
we find that on $\ev A_1$ and given
$(\gSLE,\gamma^\phi[0,T])$,
with conditional probability 
at least $(5/6)\,s_0$
the random walk gets to a vertex $v$
where $\dist\bl(\phi(v),\gSLE[t_1,t_2]\br)$
is arbitrarily small before time $\tau_T$.
Now, on the event $\ev A_0$, $\gamma^\phi[0,T]$ has to be close by,
and so we find from
Lemma~\ref{l.hitnear} that on $\ev A_0\cap \ev A_1$
and given $(\gSLE,\gamma^\phi[0,T])$, with
conditional probability at least
$(2/3)\,s_0$ we have
$\dist\bl(\phi(z_T),\gSLE[t_1,t_2]\br)$ arbitrarily small.
This proves~\eref{e.bingo}.

We take $\gamma^\phi[T,\infty)$ independent from $\gSLE$ given
$\gamma^\phi[0,T]$ in the coupling of $\gSLE$ with $\gamma$.
As in \S\ref{ss.enterRW}, let $\tau$ be the hitting time
of $S$ on $\p\DD(\gamma)$.
In~\eref{e.zzTzz} we choose $\eps=\delta$, and get a corresponding
$p>0$. Let $\ev A_2$ denote the event
$$\Pb{\zzT=\zz\md S,\gamma^\phi[0,T]}\ge p\,.
$$
Then~\eref{e.zzTzz} reads $\Pb{\ev A_2}\ge 1-\delta$.
In conjunction with~\eref{e.bingo} and $\Pb{\ev A_j}\ge 1-\delta$, $j=0,1$,
this implies that
\begin{equation}\label{e.separatingpoint}
\PB{
\PB{\dist\bl(\phi(\zz),\gSLE[t_1,t_2]\br)<\eps_1\md \gSLE,\gamma^\phi[0,T]}
> p\,s_0/2}\ge 1-3\,\delta \,.
\end{equation}

Conditional on $\zz$, on the event
$\ev A_1\cap\{\dist\bl(\phi(\zz),\gSLE[t_1,t_2]\br)<\eps_1\}$
we invoke Lemma~\ref{l.direct} with $\zz$ translated to
$0$ where the $p$ in that lemma is chosen as $p\,s_0\,\delta/2$
and the $R$ is taken to be $s_2\,\rr$,
where $s_2=s_2(s_0)>0$ is a small constant depending
only on $s_0$.
Note that the assumption necessary for the lemma that $\dist(\zz,\{\vv\}\cup\p D)> 4\,R$
 holds by~\ref{i.distR}, \ref{i.diamA} and~\ref{i.ifar} in the
definition of $\ev A_1$ and the fact that the
distance distortion of $\rr^{-1}\,\phi^{-1}$ is bounded on compact subsets of $\H$.
Let $\tilde a$ denote
the $a$ provided by the lemma, which is a function of
$s_0,a,b$ and $\delta$.
Set $r=R/(\tilde a+1)$.
Then the lemma together with~\eref{e.separatingpoint}
imply that with probability $1-O(\delta)$
there is within distance $\eps_1$ from $\gSLE[t_1,t_2]$
the $\phi$-image of a vertex $v\in \p \DD(\gamma)$
such that $\phi^{-1}\circ\gamma^\phi[0,T]$ has precisely two disjoint crossings
of the annulus $\{z: r\le |z-v|\le R\}$.
If this happens, let 
$z_1$ be a point in $\gSLE[t_1,t_2]$ closest to such a $\phi(v)$.

Let $r'$ denote the lower bound we get
on $\bl\{|\phi(v)-\phi(z)|:|v-z|\ge r\br\}$
which follows from the bounded distortion
of $\rr\,\phi$. 
We may also assume that $|\phi(v)-\phi (z)|\le s_0/8$
when $|v-z|\le R$. Now take $\eps_1=r'/4$ and
let $s_3=s_3(\delta,r')\in(0,r'/8)$ be so small that with probability
at least $1-\delta$ for every ball of radius $r'/4$ centered
at a point $z_0\in \gSLE[t_0,T]$ the distance outside of the ball $B(z_0,r'/4)$
between the two
connected components of $\R\cup\gSLE[0,T]\setminus\{z_0\}$
is at least $s_3$.
Note that when this is the case, every path connecting these
two components outside of $B(z_0,r'/4)$ has to intersect
$\cH_{s_3/3}(\gSLE[0,T])$.
Consequently, if additionally $\rho<s_3/3$,
then $\gamma^\phi[0,T]$ cannot contain an arc whose
endpoints are within distance $s_3/4$ of these two components
unless the arc visits the ball $B(z_0,r'/4)$.
If $\rho$ is sufficiently small, then
there is some $t_3'\le t_3$ such that
$\bl|\gSLE(t_3)-\gamma^\phi(t_3')\br|< s_3/4$.
Now choose $z_0:=z_1$, for the previous paragraph.
The path $\gamma^\phi[0,t_3']$ must pass through
the ball $B(z_0,r'/4)$, and therefore
$\phi^{-1}\circ \gamma^\phi[0,t_3']$
contains two disjoint crossings of the annulus $\{z: r\le |z-v|\le R\}$.
If we assume that $\phi^{-1}\circ\gamma^\phi[0,T]$ has no more than two disjoint
crossings of this annulus (which happens with probability
at least $1-O(\delta)$),
it follows that for $t\in[t_3,T]$ the point
$\gamma^\phi(t)$ is closer to $\gSLE[t_0,T]$ than to
$\gSLE[0,t_0]\cup\R$.
Since in the above $\delta$, $t_0$ and $t_3$ are arbitrary
subject to the constraint $0<t_0<t_3$ and $\delta>0$,
the claimed uniform convergence in $[0,T]$ follows.

To prove convergence in law with respect to the uniform $d_*$
metric, it suffices to show that for 
every radius $r_1>0$ there is some $r_2>r_1$ such that
$\gamma^\phi$ is unlikely to return to $B(0,r_1)$
after its first exit from $B(0,r_2)$.
For this proof, we will use the conformal
invariance of extremal length
(see~\cite{\Ahlfors}).

Fix some $r_1>0$.
The extremal length of the collection of arcs
in the half-annulus $A:=\{z\in \cH: r_1<|z|<r_2\}$ which connect
$\R_+$ with $\R_-$ tends to zero as $r_2\to\infty$.
Let $A':=\phi^{-1}(A)$,
$\p_j:=\phi^{-1}\bl(\{z\in\cH:|z|=r_j\})$, $j=1,2$,
$L:=\dist (\p_1,\p_2;A')$
and $L':=\dist (\p_+,\p_-;A')$.
By conformal invariance of extremal length,
it follows that  $L'/L\to 0$ as $r_2\to\infty$,
uniformly in $D$.
(Otherwise, the metric which is equal
to the Euclidean metric in
the ball of radius $3\,L$ centered
on a point in an arc of length at most $2\,L$
from $\p_1$ to $\p_2$ in $A'$
and is zero outside this ball
contradicts the extremal length going to zero.)

Let $\beta\subset A'$ be an arc of diameter
at most $2\,L'$ connecting $\p_+$ and $\p_-$.
Let $L_1:=\dist(\beta,\p_1;A')$ and
$L_2:=\dist(\beta,\p_2;A')$.
Since $L'/L\to0$ as $r_2\to\infty$,
we have $L_1\ge L/3$ or $L_2\ge L/3$
if $r_2$ is sufficiently large large.
Suppose first that $L_2\ge L/3$.
Let $s_0$ denote the first time such that
$|\phi\circ\gamma(s_0)|=r_2$. Then there are two
connected components $\beta_1$ and $\beta_2$ of
$\beta\setminus \gamma[0,s_0]$ such that
$\gamma[0,s_0]\cup \beta_1\cup\beta_2$
separates $\phi^{-1}\bl(B(0,r_1)\br)$
from $y_\p$ in $D$.
Now the proof of Theorem~\ref{t.nohitbd} shows that
$$
\Pb{\gamma\cap\beta_1\ne\emptyset\md \gamma[0,s_0]} <\delta
$$
if $L_2/L'$ is sufficiently large.
(Hence when $L_2\ge L/3$ and $r_2$ is sufficiently large.)
A similar estimate holds with $\beta_2$.
Thus, with probability at most $O(\delta)$,
$\gamma^\phi$ contains two
disjoint crossings of the annulus $r_1\le |z|\le r_2$.
On the other hand, if $L_2<L/3$ and $L_1\ge L/3$,
then we may apply the same argument to the reversal of
$\gamma$ (or else slightly modify the
way the analog of Theorem~\ref{t.nohitbd} is proved)
to reach the same conclusion. This completes the proof.
\QED

\section{Other lattices}\label{s.finalremarks}

In this section we describe the modifications necessary to adapt the proof
of Theorems~\ref{absmall} and~\ref{preciseSLEconvergence} to the more general
framework of Theorem~\ref{t.otherlattices}.

Before we go into the actual proof, a few words need to be said
about the properties of the weighted random walk on $\GG$ and its convergence
to Brownian motion. Fix some vertex $v_0$ in $\GG$, and let
$V_0$ denote its orbit under the group generated by the
two translations $T_1$ and $T_2$ preserving $\GG$.
If the walk starts at $v_0$, then a new Markov chain is obtained
by looking at the sequence of vertices in $V_0$ which the walk visits.
A simple path reversal argument shows that for this new Markov chain the
transition probability from $v$ to $u$ is the same as the transition
probability from $u$ to $v$, for every pair of vertices $v,u\in V_0$.
Also observe that the $\R^2$-length of a single step has an exponential tail.
This is enough to show that the Markov chain on $V_0$, rescaled
appropriately in time and space, converges to a linear image of
Brownian motion (and it is not hard to verify that the linear transformation
is non-singular). Moreover, the few properties of the simple random walk on
$\TG$ which we have used in the course of the paper are easily verified
for this Markov chain on $V_0$ and easily translated to the
weighted walk on $\GG$.

\proofof{Theorem~\ref{t.otherlattices}}
Very few changes are needed to adapt the proof.
Let $\GG$ denote the original lattice consisting of
only edges of positive weight, and let $\bar\GG$ denote
the triangulation of $\GG$, as described in~\S\ref{precisestatementsection}.
Let $\bar\GG^*$ denote the planar dual of $\bar \GG$.

The statement and proof of Lemma~\ref{l.ibd} requires some changes, because
in the more general setup it is not true that every vertex
adjacent to an interface on the right has a $\GG$-neighbor
on the left of the interface (and similarly in the other direction).
Thus, in the revised version of the lemma,
the assumption that each vertex in $V_+$ neighbors with a vertex in $V_-\cup V_\p$
and every vertex in $V_-$ neighbors with a vertex in $V_+\cup V_\p$ needs
to be replaced by the assumption that for some constant $m$, depending on the
lattice, for every vertex $v\in V_+$  the $\GG$-graph-distance from
$v$ to $V_-\cup V_\p$ is at most $m$, and symmetrically for vertices in $V_-$.
This change requires a few extra lines in the proof of~\eref{e.hbexp}.
Let $M_j$ be the maximum of $\Eb{e^{h(v)}\md \ev K}$
for vertices in $V_+$ at $\GG$-distance at most $j$ from
$V_-\cup V_\p$. Every vertex at $\GG$-distance $j>0$ from
$V_-\cup V_\p$ has a $\GG$-neighbor at $\GG$-distance $j-1$
from $V_-\cup V_\p$. Therefore, the proof of~\eref{e.Mbd2} now gives
$$
M_j \le O(1) \, M_m^{c}\,M_{j-1}^{1-c}+O(1)\,,
$$
where $c<1$ is some constant depending only on the
lattice and its edge weights.
We can certainly drop the trailing additive $O(1)$.
Induction on $j$ now gives
$$
M_j \le O(1)^{q_j}\,M_m^{c\;q_j}\,M_0^{(1-c)^j},
$$
where $q_j=1+(1-c)+\cdots+(1-c)^{j-1}=\bl(1-(1-c)^j\br)/c$.
When $j=m$ this reads
$$
M_m^{(1-c)^m} \le O(1)^{q_m}\,M_0^{(1-c)^m},
$$
Clearly, $M_0\le e^{\upperhco}$, and the
bound $M_m=O_{m,\upperhco}(1)$ follows.
A corresponding bound clearly also holds for
$\Eb{e^{-h(v)}\md\ev K}$ when $v\in V_-$.
The remainder of the proof of
the analog of Lemma~\ref{l.ibd} proceeds without
difficulty.

The proof of Lemma~\ref{l.ibd2} needs to be similarly
adapted, but essentially the same argument works.

The next point which requires adaptation is
the definition of $\evZ$ in \S\ref{ss.enterRW}.
Let $\Sigma$ denote the set of pairs $(v,e^*)$,
where $v$ is a vertex in $\GG$ and
$e^*$ is an edge in $\bar\GG^*$ that is dual to one of the edges incident
with $v$ in $\bar\GG$.
If $\sigma=(v,e^*)\in\Sigma$, let $\ev Z^\sig$ denote the
event that the first vertex adjacent to $\gamma$ that
$S$ hits is $v$ and moreover $e^*\in\gamma$.
Let $\Sigma'$ be a collection of elements of $\Sigma$,
one from each orbit under the group generated by the
translations $T_1$ and $T_2$ preserving $\GG$.
Let $\ev Z_0:=\bigcup_{\sig\in\Sigma'}\ev Z^\sig$.
The proof then proceeds essentially unchanged, with
$\ev Z^\sig$ in place of $\evZ$ and with the modified
definition for $\ev Z_0$.
\QED

\addcontentsline{toc}{section}{Bibliography}

\bibliographystyle{halpha}
\bibliography{mr,prep,notmr,DGFFcontours}

\end{document}